\documentclass[leqno]{report}
\usepackage{amsmath,amsthm,amsfonts,amssymb,setspace,graphicx,
diagrams,float, floatflt,wrapfig, verbatim}
\newcommand{\frontmatter}{\clearpage \pagenumbering{roman}}
\newcommand{\mainmatter}{\clearpage \pagenumbering{arabic}}

\usepackage{graphicx}
\usepackage{amssymb}
\usepackage{relsize}
\usepackage{amsmath}
\usepackage[mathscr]{eucal}
\newtheorem{thm}{\hspace{1cm}Theorem}[subsection]

\newtheorem{cor}[thm]{\hspace{1cm}Corollary}
\newtheorem{lem}[thm]{\hspace{1cm}Lemma}
\newtheorem{sublem}[thm]{\hspace{1cm}Sublemma}
\newtheorem{prop}[thm]{\hspace{1cm}Proposition}

\theoremstyle{definition}
\newtheorem{defn}[thm]{\hspace{1cm}Definition}
\theoremstyle{remark}
\newtheorem{rem}[thm]{\hspace{1cm}Remark}
\numberwithin{equation}{subsection}

\newenvironment{myproof}[1]{\emph{#1}.  }{\qed\vspace{.3cm}}

 \newcommand{\s}{\mathcal{S}}

 \newcommand{\Real}{\mathbb{R}}
 \newcommand{\Natural}{\mathbb{N}}
 
 \newcommand{\Complex}{\mathbb{C}}
 \newcommand{\Rational}{\mathbb{Q}}
 
 \newcommand{\Int}{\mathbb{Z}}

 \newcommand{\Gh}{\mathfrak{h}}

 \newcommand{\GL}[2]{\mathrm{GL}_{#1}(#2)}
 
 \newcommand{\SL}[2]{\mathrm{SL}_{#1}(#2)}

  \newcommand{\Mat}[2]{\mathrm{M}(#1,#2)}
 \newcommand{\SO}[2]{\mathrm{SO}_{#1}(#2)}
  
 \newcommand{\SOR}[1]{\mathrm{SO}(#1)}
 \newcommand{\SOM}[2]{\mathrm{SO}(#1,#2)}

 \newcommand{\romtwo}{{\mathrm{I \hspace{-0.2ex}I}}}

 \newcommand{\intd}{\,\mathrm{d}}

 \newcommand{\half}{\frac{1}{2}}

  \newcommand{\vectx}{\mathbf{x}}
 \newcommand{\vectj}{\mathbf{j}}
 \newcommand{\vecti}{\mathbf{i}}

 \newcommand{\vectu}{\mathbf{u}}
 \newcommand{\vectv}{\mathbf{v}}
  \newcommand{\vectw}{\mathbf{w}}

 \newcommand{\conj}{\mathbf{c}}

 \newcommand{\vol}{\mathrm{vol}}

 \newcommand{\diag}{\mathrm{diag}}

 \newcommand{\Id}{\mathrm{Id}}

 \newcommand{\Go}{\mathfrak{o}}
 \newcommand{\Mats}{\mathrm{Mat}}

 \newcommand{\inv}{^{-1}}
 \newcommand{\scrF}{\mathscr{F}}
  \newcommand{\scrG}{\mathscr{G}}
 
 \newcommand{\scrB}{\mathscr{B}}
  
 \newcommand{\scrS}{\mathscr{S}}

 \newcommand{\scrM}{\mathscr{M}}
 
 \newcommand{\scrC}{\mathscr{C}}

 \newcommand{\Bflambda}{\mathbf{\Lambda}}
 \newcommand{\res}{\mathrm{res}}

 \newcommand{\Spos}{\mathrm{\Spos}}

 \newcommand{\Rept}{\mathrm{Re}}
 \newcommand{\Impt}{\mathrm{Im}}
 \newcommand{\bbH}{\mathbb{H}}
 \newcommand{\bbS}{\mathbb{S}}

 \newcommand{\coh}{\mathscr{C}_{\mathbf{H}}}
 \newcommand{\intrr}{\mathrm{Int}}
 \newcommand{\finishdisp}{\vspace{0.3cm} \]}
 \newcommand{\startdisp}{\vspace{0.3cm}\[ }
 \newcommand{\starteqn}{\vspace{0.3cm}\begin{equation}}
 \newcommand{\finisheqn}{\vspace{0.3cm}\end{equation}}
 \newcommand{\starteqnarraynn}{\vspace*{0.1cm}\begin{eqnarray*}}
 \newcommand{\finisheqnarraynn}{\vspace*{2cm}\end{eqnarray*}}
  \newcommand{\starteqnarray}{\vspace*{0.1cm}\begin{eqnarray}}
 \newcommand{\finisheqnarray}{\vspace*{2cm}\end{eqnarray}}
 
 \newcommand{\Fix}{\mathrm{Fix}}

 \newcommand{\ord}{\mathrm{ord}}

 \newcommand{\forallindisp}{\quad\text{for all}\;}

 \newcommand{\bmth}{\boldmath}
 \newcommand{\ubmth}{\unboldmath}

  \newcommand{\red}{\mathrm{red}}
  \newcommand{\MtwoN}{\mathrm{M}_2^{\mathrm{N}}}
  \newcommand{\rt}{\mathrm{rt}}
  \newcommand{\sq}{\mathrm{sq}}
  \newcommand{\scrK}{\mathscr{K}}
  \newcommand{\Iso}{\mathrm{Iso}}
  \newcommand{\vectt}{\mathbf{t}}
  \newcommand{\bigsubseteq}{\mbox{\Large\,$\subseteq$\,}}
  
  \newcommand{\Graph}{\mathrm{Graph}}
  
\begin{document}
\pagestyle{headings}
\frontmatter
\begin{abstract}
\begin{center}
A fundamental domain of Ford type for some subgroups of orthogonal
groups\\
April    2006
\end{center}
We initiate a study of the spectral theory of the locally symmetric space
$X=\Gamma\backslash G/K$, where $G=\SO{3}{\Complex}$,
$\Gamma=\SO{3}{\Int[\vecti]}$, $K=\SOR{3}$.
We write down explicit equations defining a fundamental
domain for the action of $\Gamma$ on $G/K$.  The fundamental
domain is well-adapted for studying the theory of $\Gamma$-invariant
functions on $G/K$.  We write down equations defining a fundamental domain
for the subgroup $\Gamma_{\Int}=\SOM{2}{1}_{\Int}$ of
$\Gamma$ acting on the symmetric space $G_{\Real}/K_{\Real}$,
where $G_{\Real}$ is the split real form $\SOM{2}{1}$ of $G$
and $K_{\Real}$ is its maximal compact subgroup $\SOR{2}$.
We formulate a simple geometric relation between the fundamental
domains of $\Gamma$ and $\Gamma_{\Int}$ so described.
We
then use the previous results
compute the covolumes of of the lattices $\Gamma$ and $\Gamma_{\Int}$
in $G$ and $G_{\Real}$.
\end{abstract}
\begin{titlepage}
\begin{singlespace}
\vspace*{50pt}
\begin{center} \Huge \bfseries
A fundamental domain of Ford type for some subgroups of orthogonal
groups
\end{center}
\vspace{2.25in}
\begin{center}
A Preprint
\end{center}
\vspace{2.25in}
\begin{center}
by \\ Eliot Brenner
\end{center}
\begin{center}
Center for Advanced Studies in Mathematics\\
at Ben Gurion University
\end{center}
\begin{center}
May 2006
\end{center}
\end{singlespace}
\end{titlepage}
\clearpage
\thispagestyle{empty}
\begin{center}
\large\vspace*{1in}
ACKNOWLEDGMENTS
\end{center}
\vspace{0.25in}
The Institute for the Advanced Study of Mathematics
at Ben-Gurion University provided support
and a pleasant and stimulating working environment
during the writing of this memoir.  Additionally,
the author  thanks Mr. Tony Petrello for financial assistance.
\tableofcontents
\mainmatter
\chapter{Fundamental Domain}
\section{Representation of $\SO{3}{\Int[\vecti]}$ as a lattice
in $\SL{2}{\Complex}$}\label{sec:lattice}
We begin by establishing some basic notational conventions.

Let $n$ be a positive integer and $\Go$ a ring.  We will
use $\mathrm{Mat}_n(\Go)$ to denote the set of all\linebreak
$n$-by-$n$ square matrices with coefficients in $\Go$.
We reserve use the Greek letters $\alpha$, and so on,
for the elements of $\mathrm{Mat}_n(\Go)$,
and the roman letters $a,b,c,d$ and so on,
for the entries of the matrices.
We will denote scalar mutliplication
on $\mathrm{Mat}_n(\Go)$ by simple juxtaposition.
Thus, if $\Go=\Int[\vecti]$,
$\ell\in\Int[\vecti]$ and $\alpha\in\mathrm{Mat}_2(\Int[\vecti])$,
then
\startdisp
\alpha=\begin{pmatrix}a&b\\c&d\end{pmatrix}\;\text{implies}\;
\ell\alpha=\begin{pmatrix}\ell a&\ell b\\\ell c&\ell d\end{pmatrix}.
\finishdisp
The letters $p,q,r,s$ will be reserved
to denote a quadruple of elements of $\Go$
such that $ps-rq=1$.
In what follows, we normally have $\Go=\Int[\vecti]$,
whenever $\alpha$ is written
with entries $p$ through $s$.  Therefore,
\startdisp\alpha=
\begin{pmatrix}p&q\\
r&s\end{pmatrix}\in \SL{2}{\Int[\vecti]},
\finishdisp
unless stated otherwise.

\subsection{Realization of $\SO{3}{\Int[\vecti]}$ as a group
of fractional linear transformations}\label{subsec:flts}
We will denote a conjugation action of a group on a space $V$
by $\conj_V$, when the context makes clear
what this action is.
For example, if $H$ is a linear Lie group and $\Gh$ the Lie
algebra of $H$, then we have
\startdisp
\conj_{\Gh}(h)X=hXh\inv,\forallindisp h\in H,\;X\in\Gh.
\finishdisp
Note that the morphism $\conj_{\Gh}(h)$ is the image under the Lie
functor of the
usual conjugation $\conj_H(h)$ on the group level.
Using $\mathrm{SL}(V)$ to denote the group of unimodular transformations
of a vector space $V$, it is easy to see that
\starteqn\label{eqn:conjalwaysmorphism}
\conj_{\Gh}: H\rightarrow SL(\Gh)\;\text{is a Lie group morphism}.
\finisheqn

Henceforth, whenever $H$ is a group acting on a Lie algebra
$\Gh$ by conjugation, we will omit the subscript $\Gh$.
Thus, we define
\startdisp
\conj:=\conj_{\Gh},
\finishdisp
when we are in the situation of \eqref{eqn:conjalwaysmorphism}.

Throughout (except in a few sections where it is prominently
noted otherwise), we will use the notation $G=\SO{3}{\Complex}$, $\Gamma=\SO{3}
{\Int[\vecti]}$.
We use $B$ to denote the half-trace form on $\mathfrak{sl}_2(\Complex)$, the Lie
algebra of traceless $2$-by-$2$ matrices.  That is,
\startdisp
B(X,Y)=\half\mathrm{Tr}(XY).
\finishdisp
We use the notation $\beta'=\{X'_1,X'_2,Y'\}$ for the ``standard"
basis of $\mathfrak{sl}_2(\Complex)$, where
\starteqn\label{eqn:betaprime}
\begin{gathered}
X_1'=\begin{pmatrix}0&1\\
0&0
\end{pmatrix},\qquad X_2'=\begin{pmatrix}0&0\\
1&0\end{pmatrix},\\
\text{and}\quad
Y'=\begin{pmatrix}1&0\\
0&-1
\end{pmatrix}
\end{gathered}
\finisheqn

The following properties of $B$ are verified either immediately
from the definition or by straightforward calculations.
\begin{itemize}
\item[\textbf{B1}] $B$ is nondegenerate.
\item[\textbf{B2}]  Setting
\starteqn\label{eqn:beta}
\begin{gathered}
X_1=X_1'+X_2',\qquad X_2=\vecti(X_1'-X_2'),\\
\text{and}\quad Y=Y',
\end{gathered}
\finisheqn
we obtain an orthonormal basis $\beta=\{X_1,X_2,Y\}$, with
respect to the bilinear form $B$.
\item[\textbf{B3}] $B$ is invariant under the conjugation action of $\SL{2}{\Complex}$, meaning
that
\startdisp
B(X,Y)=B(\conj(g)Z,\conj(g)W),\forallindisp Z,W\in\mathfrak{sl}_2(\Complex),\,
g\in\SL{2}{\Complex}.
\finishdisp
\end{itemize}
By \textbf{B3}, $\conj$ is a morphism of $\SL{2}{\Complex}$ into $G$.
The content of part (a) of Proposition \ref{prop:conjsurjection} below
is that the morphism
$\conj$ just described is an epimorphism.

As a consequence of \textbf{B1} and \textbf{B2}, we have that
\starteqn\label{eqn:halftracequadratic}
B(x^1_1 X_1+x^1_2 X_2+y^1 Y,x^2_1 X_1+x^2_2 X_2+y^2Y)=
x^1_1x^2_1+x^1_2x^2_2+y^1y^2,\;
x^i_j, y\in\Complex.
\finisheqn
For any bilinear form $B$ on a vector space $V$,
we use $\mathrm{O}(B)$ to denote the group of linear transformations of $V$
preserving $B$, and we use $\mathrm{SO}(B)$ to denote the unimodular
subgroup of $\mathrm{O}(B)$.
If $B$ is as in \eqref{eqn:halftracequadratic}, then the isomorphism,
\starteqn\label{eqn:concreteorthequiv}
\mathrm{SO}(B)\cong G,
\finisheqn
induced by the identification of the vector space $\mathfrak{sl}_2(\Complex)$
with $\Complex\langle X_1,X_2,Y\rangle$, puts a system of coordinates
on $G$.  Part (b) of Proposition \ref{prop:conjsurjection},
below, will describe the epimorphism $\conj:\SL{2}{\Complex}\rightarrow G$
in terms of these coordinates.
\begin{prop} \label{prop:conjsurjection} With $G$, $\conj$
as above, we have\begin{itemize}
\item[(a)] The map $\conj$ induces an isomorphism
\startdisp
\SL{2}{\Complex}/\{\pm I\}\stackrel{\cong}{\longrightarrow} G
\finishdisp
of Lie groups.
\item[(b)]  Relative to the standard coordinates
on $\SL{2}{\Complex}$ and the coordinates on $G$ induced from
the orthonormal basis $\beta$ of $\mathfrak{sl}_2(\Complex)$, as defined
in \eqref{eqn:beta}, the
epimorphism $\conj: \SL{2}{\Complex}\rightarrow G$
has the following coordinate expression.
\starteqn\label{eqn:imagematrixconj}
\conj\left(\begin{pmatrix}a&b\\c&d\end{pmatrix}\right)=
\begin{pmatrix}\frac{a^2-c^2+d^2-b^2}{2}&
\frac{\vecti(a^2-c^2+b^2-d^2)}{2}&cd-ab\vspace*{0,2cm}\\
\frac{\vecti(b^2+d^2-a^2-c^2)}{2}&\frac{a^2+c^2+b^2+d^2}{2}&\vecti(ab+cd)
\vspace*{0.2cm}\\
-ac+bd&\vecti(ac+bd)&ad+bc\end{pmatrix}.
\finisheqn
\end{itemize}
\end{prop}
Before giving the proof of the Proposition, we establish some further
notational conventions regarding conjugation mappings.  Whenever
a matrix group $H$ has a conjugation action $\conj_V$ on a \textit{finite dimensional
vector space} $V$ over a field $F$,
each basis $\beta$ of $V$ naturally induces a morphism
\starteqn\label{eqn:conjwithresptobasis}
\conj_{V,\beta}:H\rightarrow\GL{N}{F},\;\text{where}\; N=\dim V.
\finisheqn
Let $\beta,\beta'$ be two bases of $V$.  Write $\alpha^{\beta\mapsto\beta'}$
for the change-of-basis matrix from $\beta$ to $\beta'$.  That is,
if $\beta$, $\beta'$ are written as $N$-entry row-vectors, then
\starteqn\label{eqn:changeofbasis}
\beta\alpha^{\beta\mapsto\beta'}=\beta'.
\finisheqn
Then elementary linear algebra tells us that
\starteqn\label{eqn:conjchangeofbasis}
\begin{aligned}
\conj_{V,\beta}&=&&\conj_{\GL{N}{F}}\left(\left(\alpha^{\beta\mapsto\beta'}\right)\inv\right)\conj_{V,\beta'}\\
&=&&\conj_{\GL{N}{F}}\left(\alpha^{\beta'\mapsto\beta}\right)\conj_{V,\beta'}.
\end{aligned}
\finisheqn
Assuming that $c_V$ is injective, and writing $c_V\inv$ for the
left-inverse of $c_V$, we calculate from \eqref{eqn:conjchangeofbasis}
that
\starteqn\label{eqn:changeofbasisconj}
\conj_{V,\beta}\conj_{V,\beta'}\inv\in\mathrm{Aut}
(\GL{N}{F})\;\text{is given by}\;
\conj_{\GL{N}{F}}\left(\alpha^{\beta\mapsto\beta'}\right).
\finisheqn
In keeping with the practice established after \eqref{eqn:conjalwaysmorphism},
we will omit the subscript $\Gh$ when $H$ is a Lie group
acting on its Lie algebra by conjugation.  Thus, for any basis $\beta$
of $\Gh$,
\startdisp
\conj_{\beta}:=\conj_{\Gh,\beta}.
\finishdisp
Further, except in situations such as the proofs of Proposition
\ref{prop:conjsurjection} and \ref{prop:gammaziso} below, where
we are dealing with $\conj_{\beta}$ for different bases
$\beta$ of $\Gh$ at the same time, we will fix a single basis $\beta$
for $\Gh$ and blur the distinction between $\conj$ and $\conj_{\beta}$.
For example, in this chapter, whenever $H=\SL{2}{\Complex}$ and $V=\mathrm{Lie}
(H)$, we will write $\conj$ to denote both the ``abstract" morphism $\conj$
of $H$ into $\mathrm{Aut}(V)$ and the linear morphism $\conj_{\beta}$
of $H$ into $\GL{3}{\Complex}$, where $\beta$ is the orthonormal
basis for $\mathrm{Lie}(H)$ defined in \eqref{eqn:beta}.  Whenever
the linear morphism into $\GL{3}{\Complex}$ is induced
by a basis $\beta'\neq\beta$, the notation $\conj_{\beta'}$ will
be used.

\begin{myproof}{Proof of Proposition \ref{prop:conjsurjection}}
For (a), it is clear that the kernel of $\conj$ is the center
of $\SL{2}{\Complex}$, which equals $\{\pm I\}$.  The surjectivity
of $\conj$ follows from a comparison of the dimensions of $\SL{2}{\Complex}$
and $G$ (both have complex dimension $3$),
and from the well-known fact that $G$ is connected.

Let $h\in\SL{2}{\Complex}$, so that
\startdisp
h=\begin{pmatrix}a&b\\c&d\end{pmatrix},\quad\text{with}\;ad-bc=1.
\finishdisp
Recall the basis
$\beta'=\{X_1', X_2', Y'\}$ for $\mathrm{Lie}(\SL{2}{\Complex})$.
A routine calculation using \eqref{eqn:betaprime} shows that we have
\starteqn\label{eqn:conjholdbasis}
\conj_{\beta'}(h)=\begin{pmatrix}a^2&-b^2&-2ab\\
-c^2&d^2&2cd\\
-ac&bd&ad+bc\end{pmatrix}.
\finisheqn
Let $\alpha^{\beta\mapsto\beta'}$ be the change-of-basis
matrix, satisfying relation \eqref{eqn:changeofbasis}.
As a consequence of the the relations \eqref{eqn:beta}, we obtain
\starteqn
\label{eqn:basistransformationequals}
(\alpha^{\beta\mapsto\beta'})\inv=\begin{pmatrix}1&\vecti&0\\
1&-\vecti&0\\
0&0&1\end{pmatrix}.
\finisheqn
Using \eqref{eqn:conjchangeofbasis},
\eqref{eqn:basistransformationequals}, and \eqref{eqn:conjholdbasis}
we obtain
\startdisp
\conj(h):=\conj_{\beta}(h)=\begin{pmatrix}\frac{a^2-c^2+d^2-b^2}{2}&
\frac{\vecti(a^2-c^2+b^2-d^2)}{2}&cd-ab\vspace*{0.2cm}\\
\frac{\vecti(b^2+d^2-a^2-c^2)}{2}&\frac{a^2+c^2+b^2+d^2}{2}&\vecti(ab+cd)
\vspace*{0.2cm}\\
-ac+bd&\vecti(ac+bd)&ad+bc\end{pmatrix}.
\finishdisp
This gives part (b) of the proposition.
\end{myproof}

Whenever $G_1\subseteq\GL{N_1}{F}$ and $G_2\subseteq\GL{N_2}{F}$
are linear groups and $\varphi$ is a morphism expressed
by a formula involving only polynomial functions of the entries,
$\varphi$ may be extended to a rational map
\startdisp
\tilde{\varphi}: \Mat_{N_1}(F)\rightarrow\Mat_{N_2}(F).
\finishdisp
For example, since all the entries in \eqref{eqn:imagematrixconj}
are polynomials in the entries of the matrix in $\SL{2}{\Complex}$,
we may extend $\conj$ to a map
\starteqn\label{eqn:rationalextension}
\tilde{\conj}: \mathrm{Mat}_2(\Complex)\rightarrow\mathrm{Mat}_3(\Complex).
\finisheqn
We note some properties of $\tilde{\conj}$ that will be useful
later.  Let $\ell$ denote a complex number and $\alpha$ an arbitrary
2-by-2 matrix.
\starteqn\label{eqn:multiplehasintegral}
\text{If $\tilde{\conj}(\alpha)$ has integral entries and $\ell^2\in\Int[\vecti]$
then $\tilde{\conj}(\ell\alpha)$ has integral entries.}
\finisheqn
We deduce \eqref{eqn:multiplehasintegral} by observing
that the entries of $\tilde{\conj}(\alpha)$ are
homogeneous polynomials of degree $2$ in the entries of $\alpha$.
As an immediate corollary of \eqref{eqn:multiplehasintegral},
one has
\starteqn\label{eqn:multiplehasintegralunits}
\text{If $\ell^2\in\Int[\vecti]^*$, then $\tilde{\conj}(\alpha)$
has integral entries if and only if $\tilde{\conj}(\ell\alpha)$
has integral entries.
}
\finisheqn

We now wish to describe the inverse image $\conj\inv(\Gamma)$
as a subset of $\SL{2}{\Complex}/\{\pm I\}$ with respect
to the standard coordinates of $\SL{2}{\Complex}$.  According to Proposition
\ref{prop:conjsurjection}, this amounts to describing the quadruples
\starteqn\label{eqn:quadruplecond}
(a,b,c,d)\in \Complex^4,\;\text{with $ad-bc=1$, and the entries
of the right-side of \eqref{eqn:imagematrixconj} integers.}
\finisheqn
Describing the quadruples meeting conditions
\eqref{eqn:quadruplecond} will
be the subject of the remainder of this section and the
next, culminating in the proof of Proposition \ref{prop:inversematrixexplicitdescription}.

\vspace*{0.3cm}\noindent\boldmath\textbf{Conventions regarding multiplicative
structure of $\Int[\vecti]$}.  \unboldmath
Before stating the proposition, we establish
certain conventions we will use when dealing with the multiplicative
properties of the Euclidean ring  $\Int[\vecti]$.  First,
it is well-known that $\Int[\vecti]$ is a Euclidean, hence principal, ring.
That $\Int[\vecti]$ is principal
means that all ideals $\mathscr{I}$ of $\Int[\vecti]$ are generated by a single
element $m\in\Int[\vecti]$, so that every $\mathscr{I}$ is of the form $(m)$.
However, there is an unavoidable ambiguity in the choice of generators
caused by the presence in $\Int[\vecti]$
of four units, $\vecti^{j}$, for $j\in\{0,\ldots, 3\}$, in $\Int[\vecti]$.  We
will adopt the following convention to sidestep the ambiguity caused
by the group of units.
\begin{defn}  We refer to the following subset of $\Complex^{\times}$
as the \textbf{standard subset}
\starteqn\label{eqn:standardset}
\{z\in\Complex^{\times}\;|\; \Rept(z)> 0,\, \Impt(z)\geq 0\}.
\finisheqn
That is, the standard subset of $\Complex^{\times}$
is the union of the interior of the first quadrant
and the positive real axis.  An element of $\Int[\vecti]$
in the standard subset will be referred to as a \textbf{standard Gaussian
integer}, or more simply as a \textbf{standard integer} when the context
is clear.
\end{defn}
Because of the units in $\Int[\vecti]$, each nonzero
ideal $\mathscr{I}$ of $\Int[\vecti]$ has precisely one generator which
is a standard integer.   Henceforth, we refer to generator of $\mathscr{I}$
which is a standard integer as the \textbf{standard
generator} of $\mathscr{I}$.
Unless otherwise stated, whenever we write $\mathscr{I}=(m)$,
to indicate the ideal $\mathscr{I}$ generated by an $m\in\Int[\vecti]$, it will
be understood that $m$ is standard.  Conversely, whenever
we write an ideal $\mathscr{I}$ in the form $(m)$, it
will be understood that $m$ is the standard generator
of $\mathscr{I}$.
Thus, for example, since $(1-\vecti)=\vecti^3(1+\vecti)$ with
$1+\vecti$ standard,
we write $\mathscr{I}=:(1-\vecti)\Int[\vecti]$, defined
as the ideal of Gaussian integers
 divisible
by $1-\vecti$, in the form $\mathscr{I}=(1+\vecti)$.

Similar comments apply to Gaussian primes, factorization, and greatest
common divisor in $\Int[\vecti]$.  By a ``prime in $\Int[\vecti]$", we will
always mean a \textit{standard prime}.  By ``prime factorization"
in $\Int[\vecti]$ we will always mean \textit{factorization
into a product of standard primes}, multiplied by the appropriate unit factor.
Note that the convention regarding standard primes uniquely
determines the unit factor in a prime factorization.  For example,
since
\startdisp
2=\vecti^{3}(1+\vecti)^2
\finishdisp
and $(1+\vecti)^3$ is standard, the above expression is the standard
factorization of the Gaussian integer $2$, and
$\vecti^{3}$ is uniquely determined as the \textit{standard unit factor}
in the prime factorization of $2\in\Int[\vecti]$.

\begin{defn} \label{defn:GCD} For $x,\,y\in\Int[\vecti]$, consider $(x,y)\Int[\vecti]$,
the ideal generated by $x$ and $y$.  Define
\startdisp
\text{\boldmath$\mathrm{GCD}(x,y)$\unboldmath to be the unique}\;
m\in\Int[\vecti]\;\text{such that $m$ is standard and}\;
(m)=(x,y)\Int[\vecti].
\finishdisp
\end{defn}
We will use the following three basic properties of $\mathrm{GCD}(x,y)$.
In each case, the proof is the same as for the corresponding
properties of the greatest common divisor for the rational integers.
Let $x,\,y,\,z\in\Int[\vecti]$.  Then we have
\begin{itemize}
\item[\textbf{GCD1}]\hspace*{1cm}$\mathrm{GCD}(x,y)| x$,
\item[\textbf{GCD2}]\hspace*{1cm}$
\mathrm{GCD}\left(\frac{x}{\mathrm{GCD}(x,y)},
\frac{y}{\mathrm{GCD}(x,y)}\right)=1$,
\item[\textbf{GCD3}]\hspace*{1cm}There exist $z,w\in\Int[\vecti]$
such that $xz-yw=\mathrm{GCD}(x,y)$.
\end{itemize}

By convention, unless stated otherwise, the ``trivial ideal"
$\Int[\vecti]$ will be understood to belong to the set of ideals
of $\Int[\vecti]$.  The standard generator of
the trivial ideal $\Int[\vecti]$ is, of course, $1$.

To facilitate the statement of Proposition
\ref{prop:inversematrixexplicitdescription}, we
estblish the following conventions.  First, we use $\omega_8$
to denote the unique primitive eighth root of unity in the standard set
of $\Complex^{\times}$.  Observe that
\starteqn\label{eqn:omegaeightdefn}
\omega_8=\frac{\sqrt{2}}{2}(1+\vecti),\quad\text{and}\quad
\omega_8^2=\vecti.
\finisheqn

\vspace*{0.3cm}\noindent\boldmath\textbf{The $\SL{2}{\Int[\vecti]}$-space
$\mathrm{M}_2^{\mathrm{N}}$}.\unboldmath
\begin{defn}  For $N\in\Int[\vecti]$, \boldmath
$\mathrm{M}_2^{\mathrm{N}}\;$\unboldmath
will the subset of $\mathrm{Mat}_2(\Int[\vecti])$ consisting
of the elements with determinant $N$.  Since the group $\SL{2}{\Int[\vecti]}$ acts
on $\MtwoN$ by multiplication on the left, $\MtwoN$ is a
$\SL{2}{\Int[\vecti]}$-space.
\end{defn}

It is not difficult to see that the
action of $\SL{2}{\Int[\vecti]}$ on $\MtwoN$ fails to be
transitive, so $\MtwoN$ is not a $\SL{2}{\Int[\vecti]}$-homogeneous
space.  The purpose of the subsequent definitions and results
is to give a description of the orbit structure of the
$\SL{2}{\Int[\vecti]}$-space $\MtwoN$.

Let
\starteqn\label{eqn:omegamdefn}
\Omega_{y}:=\text{a fixed set of representatives of
$\Int[\vecti]/(y)$},\;\text{for all}\;y\in\Int[\vecti]
\finisheqn
It is clear that, for each
$y\in\Int[\vecti]$, a number of possible $\Omega_{y}$ exist.
For the general result,
Proposition \ref{prop:heckedecomp}, below, the choice
of $\Omega_y$ does not matter, and we leave it unspecified.
However, in the specific applications
of Proposition \ref{prop:heckedecomp}, where $y$ is always
of the form $y=(1+\vecti)^n$ for $n$ a positive integer,
it will be essential
to give an $\Omega_y$ explicitly, which we now do.

So let $n\in\Natural$, $n\geq 1$.
In the definition of $\Omega_{(1+\vecti)^n}$,
we use the ``ceiling" notation, defined as follows:
\startdisp
\lceil q\rceil\;=\text{smallest integer $\geq \frac{x}{y}$.}
\finishdisp
Now set
\starteqn\label{eqn:omegaspecific}
\Omega_{(1+\vecti)^n}=\left\{r+s\vecti\;\text{with $r,\,s\in\Int$,
$0\leq r<2^{\lceil\frac{n}{2}\rceil}$, $0\leq s< 2^{n-
\lceil\frac{n}{2} \rceil}$}\right\}.
\finisheqn
The definition is justified by Lemma \ref{lem:standardresiduerep},
below.
In order not to interrupt the flow of paper, we delay
the proof of Lemma \ref{lem:standardresiduerep} until
the end of \S\ref{subsec:flts}.

\begin{lem}\label{lem:standardresiduerep}
For $n\geq 1$ an integer, let $\Omega_{(1+\vecti)^n}$
be defined as \eqref{eqn:omegaspecific}.  Then
\startdisp
\Omega_{(1+\vecti)^n}\;\text{is a complete set of representatives of
$\Int[\vecti]/\hspace*{-.5mm}\left((1+\vecti)^n\right)$}\;\text{
for all $n$}.
\finishdisp
\end{lem}

\begin{defn}\label{defn:alphamat}  Let $N\in\Int[\vecti]$ be fixed, and
for each $y\in\Int[\vecti]$ let $\Omega_y$ be as in
\eqref{eqn:omegamdefn}.  Define the matrix $\alpha^\mathrm{N}(m,x)
\in \MtwoN$
as follows,
\starteqn\label{eqn:MofNmxdefn}
\alpha^\mathrm{N}(m,x)=\begin{pmatrix}m&x\\
0&\frac{N}{m}\end{pmatrix},\;\text{for}\;m\in\Int[\vecti],\,
m|N,\, x\in\Omega_{\frac{N}{m}}.
\finisheqn
\end{defn}
It is trivial to verify that $\alpha^\mathrm{N}(m,x)$, as given by
\eqref{eqn:MofNmxdefn}, indeed has determinant $N$, \textit{i.e.}
$\alpha^\mathrm{N}(m,x)\in \MtwoN$.  The point of Definition
\ref{defn:alphamat} is given by the following proposition.

\begin{prop}\label{prop:heckedecomp}  For $N\in\Int[\vecti]-\{0\}$,
let $\MtwoN$ be the $\SL{2}{\Int[\vecti]}$-space of matrices with entries
in $\Int[\vecti]$ and determinant $N$.
Define the matrices $\alpha^\mathrm{N}(m,x)$ as in \eqref{eqn:MofNmxdefn}.
Then
\starteqn\label{eqn:heckedecomp}
\mathrm \MtwoN=
\bigcup_{\left\{\stackrel{m\in\Int[\vecti]| \; m|N,}
{\frac{N}{m}\;\text{standard}}
\right\}}
\hspace*{-1.07cm}\cdot\hspace*{1cm}
\bigcup_{x\in\Omega_{\frac{N}{m}}}
\hspace*{-.55cm}\cdot\hspace*{0.55cm}
\SL{2}{\Int[\vecti]}\alpha^\mathrm{N}(m,x),
\finisheqn
and \eqref{eqn:heckedecomp} gives the decomposition of
the $\SL{2}{\Int[\vecti]}$-space $\MtwoN$
into distinct $\SL{2}{\Int[\vecti]}$-orbits.
\end{prop}
We delay the proof of Proposition \ref{prop:heckedecomp} until
the end of \S\ref{subsec:flts}, and here restrict ourselves to some
concerning the significance of Proposition \ref{prop:heckedecomp}.
First, a statement equivalent to Proposition \ref{prop:heckedecomp}
is that an arbitrary $\alpha\in M_2^{\rm N}$
has a uniquely determined product decomposition of the form
\starteqn\label{eqn:heckeexplicit}\alpha=
\begin{pmatrix}a&b\\
c&d\end{pmatrix}=\begin{pmatrix}p&q\\
r&s\end{pmatrix}\begin{pmatrix}m&x\\0&\frac{N}{m}\end{pmatrix},\;
\text{with}\, m\in\Go,\; m|N,\; \frac{N}{m}\;\text{standard},\; x\in\Omega_
{\frac{N}{m}},\, pr-qs=1.
\finisheqn
The uniqueness is derived from Proposition \ref{prop:heckedecomp}
as follows.  The
second matrix in the product of \eqref{eqn:heckeexplicit}
is uniquely determined by the matrix decomposition because
of the disjointness of the union in \eqref{eqn:heckedecomp}.
The first matrix in the product appearing in \eqref{eqn:heckeexplicit}
is therefore also uniquely determined.

The second remark is that Proposition \ref{prop:heckedecomp}
may be thought of as the Gaussian-integer
version of the decomposition of elements
of $\mathrm{Mat}_2(\Int)$ of fixed determinant $N$,
sometimes known as the Hecke decomposition.  Occasionally
we refer to \eqref{eqn:heckeexplicit} as the \textit{Gaussian}
Hecke decomposition, to distinguish it from this \textit{classical}
Hecke decomposition in the context of the rational integers.
Readers familiar with the the proof
of the classical Hecke decomposition may skip the proof
of Proposition \eqref{prop:heckedecomp}, since
the proof is the same as that of the classical decomposition
except for some care that has to be taken because of the presence
of additional units in $\Int[\vecti]$.
\newcounter{ranrom}\setcounter{ranrom}{7}
For the classical
Hecke decomposition, see page 110, \S\Roman{ranrom}.4, of \cite{langmodforms},
which is the source of our notation for the Gaussian version.

\vspace*{0.3cm}\noindent\boldmath\textbf{Statement
of the Main Result of \S\ref{sec:lattice}}.  \unboldmath
Let $\Xi$ be an arbitrary subset of $\SL{2}{\Int[\vecti]}$.
Suppose, at first, that $\Xi$ is actually a subgroup
of $\SL{2}{\Int[\vecti]}$.
Since $\SL{2}{\Int[\vecti]}\alpha^\mathrm{N}(m,x)$ is an $\SL{2}{\Int[\vecti]}$-space, it is
also a $\Xi$-space.  For general, $\Xi$, however, the action of $\Xi$
on $\SL{2}{\Int[\vecti]}\alpha^{\rm N}(m,x)$
fails to be transitive, \textit{i.e.},
$\SL{2}{\Int[\vecti]}\alpha^\mathrm{N}(m,x)$ is
not a $\Xi$-homogeneous space.  We will now describe
the orbit structure of $\SL{2}{\Int[\vecti]}\alpha^\mathrm{N}(m,x)$ for a specific
subgroup $\Xi$.  In order to make the description of the
subgroup and some related subsets of $\SL{2}{\Int[\vecti]}$
easier, we introduce the epimorphism
\startdisp
\mathrm{red}_{1+\vecti}: \SL{2}{\Int[\vecti]}\rightarrow\SL{2}{\Int[\vecti]
/(1+\vecti)}
\finishdisp
by inducing from the reduction map
\startdisp
\mathrm{red}_{1+\vecti}: \Int[\vecti]\rightarrow\Int[\vecti]/(1+\vecti)
\finishdisp
That is, we ``extend" $\mathrm{red}_{1+\vecti}$ from elements
to matrices
by setting
\starteqn\label{eqn:matrixredmapdefn}
\mathrm{red}_{1+\vecti}\left(\begin{pmatrix}p&q\\ r&s\end{pmatrix}\right)=
\begin{pmatrix}\mathrm{red}_{1+\vecti}p&\mathrm{red}_{1+\vecti}q\\
\mathrm{red}_{1+\vecti}r&\mathrm{red}_{1+\vecti}s\end{pmatrix}.
\finisheqn
Since $\Omega_{1+\vecti}=\{0,1\}$, we may identify
$\Int[\vecti]/(1+\vecti)$ with $\{0,1\}$.  Similarly
to the convention with $p,q,r,s\in\Int[\vecti]$,
we use $(\overline{p},\overline{q},\overline{r},\overline{s})$
to denote a quadruple of elements of $\Int[\vecti]/(1+\vecti)$
such that
\startdisp
\overline{p}\hspace*{0.5mm}\overline{s}-\overline{r}\hspace*{0.5mm}\overline{q}=1.
\finishdisp
Here are two elements of $\SL{2}{\Int[\vecti]/(1+\vecti)}$ of particular
interest.
\starteqn\label{eqn:xi12elements}
\overline{I}:=\begin{pmatrix}1&0\\0&1\end{pmatrix},\;
\overline{S}:=\begin{pmatrix}0&1\\1&0
\end{pmatrix}\in\SL{2}{\Int[\vecti]/(1+\vecti)}.
\finisheqn
The notation in \eqref{eqn:xi12elements}
is chosen to remind the reader that $\overline{I}=\red_{1+\vecti}(I)$
and $\overline{S}=\red_{1+\vecti}(S)$, where $I,\, S$ are the standard
generators of $\SL{2}{\Int}$, as in \setcounter{ranrom}{6}
\S\Roman{ranrom}.1
of \cite{jol05}.
Since $\overline{S}^2=\overline{I}$,
it is easy to see that $\{\overline{I},\overline{S}\}$ is a subgroup
of $\SL{2}{\Int[\vecti]/(1+\vecti)}$.
Now define
\starteqn\label{eqn:xi12defn}
\Xi_{12}=\red_{1+\vecti}\inv(\{\overline{I},\, \overline{S}\}).
\finisheqn
Since $\red_{1+\vecti}$ is a morphism, $\Xi_{12}$ is a subgroup
of $\SL{2}{\Int[\vecti]}$.

Also, using the epimorphism $\red_{1+\vecti}$
we define the following subsets of $\SL{2}{\Int[\vecti]}$:
\starteqn
\label{eqn:residuematrices}
\begin{aligned}
\Xi_1&=&\red_{1+\vecti}^{-1}\left(\left\{\begin{pmatrix}0&1\\1&1\end{pmatrix},\,
\begin{pmatrix}1&1\\0&1\end{pmatrix}\right\}\right),\\
\Xi_2&=&\red_{1+\vecti}^{-1}\left(\left\{\begin{pmatrix}1&1\\1&0\end{pmatrix},\,
\begin{pmatrix}1&0\\1&1\end{pmatrix}\right\}\right).
\end{aligned}
\finisheqn
(The subscripts on the $\Xi$ of \eqref{eqn:xi12defn}
and \eqref{eqn:residuematrices}
are chosen in order to
remind the reader of the column in which zeros appear
in the matrices of $\red_{1+\vecti}(\Xi)$.)
Since $\SL{2}{\Int[\vecti]/(1+\vecti)}$
consists of the elements $\overline{I},\overline{S}$ and the four elements
appearing on the right-hand side of \eqref{eqn:residuematrices},
and $\red_{1+\vecti}$
is an epimorphism,
\starteqn\label{eqn:sl2xidecomp}
\SL{2}{\Int[\vecti]}=\Xi_1\bigcup\hspace*{-.325cm}\cdot\hspace*{.325cm}\Xi_2
\bigcup\hspace*{-.325cm}\cdot\hspace*{.325cm}\Xi_{12}.
\finisheqn
Unlike $\Xi_{12}$, the subsets $\Xi_1$ and $\Xi_2$ of $\SL{2}{\Int[\vecti]}$
are not subgroups.  All three subsets $\Xi$ in \eqref{eqn:xi12defn}
and \eqref{eqn:residuematrices}
though have the following property.
\starteqn\label{eqn:xisetdescription}
\begin{gathered}
\text{For fixed}\;
\begin{pmatrix}\overline{p}&\overline{q}\end{pmatrix},
\begin{pmatrix}\overline{r}&\overline{s}\end{pmatrix}\in
\left\{
\begin{array}{l}
\begin{pmatrix}1&1\end{pmatrix},\vspace*{0.15cm}\\
\begin{pmatrix}1&0\end{pmatrix},\vspace*{0.15cm}\\
\begin{pmatrix}0&1\end{pmatrix}
\end{array}
\right\}\subset(\SL{2}{\Int[\vecti]/(1+\vecti)})^2,\\
\Xi=\red_{1+\vecti}\inv\left(
\left\{\begin{pmatrix}\overline{p}&\overline{q}\\
\overline{r}&\overline{s}\end{pmatrix},\,
\begin{pmatrix}\overline{r}&\overline{s}\\
\overline{p}&\overline{q}
\end{pmatrix}\right\}\right).
\end{gathered}
\finisheqn
For example, we obtain $\Xi_{12}$ by taking
\startdisp
\begin{pmatrix}\overline{p}&\overline{q}\end{pmatrix}=
\begin{pmatrix}1&0\end{pmatrix}\;\text{and}
\begin{pmatrix}\overline{r}&\overline{s}\end{pmatrix}
=\begin{pmatrix}0&1\end{pmatrix}
\finishdisp
in \eqref{eqn:xisetdescription}.

The reason for introducing the subsets $\Xi$
of \eqref{eqn:residuematrices}
is that they allow us, in Sublemma \ref{sublem:finerdecomp} below
to describe precisely the orbit structure of the $\Xi_{12}$-space $\SL{2}{\Int[\vecti]}
\alpha^\mathrm{N}(m,x)$.
\begin{sublem}\label{sublem:finerdecomp}
Using the notation of
\eqref{eqn:MofNmxdefn} and \eqref{eqn:residuematrices}, we have
\starteqn\label{eqn:finerdecomp}
\SL{2}{\Int[\vecti]}\alpha^\mathrm{N}(m,x)=\bigcup_{\Xi\,=\,\Xi_{1},\,\Xi_{2},\Xi_{12}}
\hspace*{-1cm}\cdot\hspace*{0.7cm}\Xi \alpha^\mathrm{N}(m,x).
\finisheqn
Each of the three sets in the union \eqref{eqn:finerdecomp} is closed
under the action, by left-multiplication, of
$\Xi_{12}$ on $\SL{2}{\Int[\vecti]}\alpha^\mathrm{N}(m,x)$ and equals precisely
one $\Xi_{12}$-orbit in the space
$\SL{2}{\Int[\vecti]}\alpha^\mathrm{N}(m,x)$.
\end{sublem}
\begin{myproof}{Proof}  It is clear from \eqref{eqn:sl2xidecomp}
that $\SL{2}{\Int[\vecti]}\alpha^\mathrm{N}(m,x)$ equals the disjoint union three
sets as on the right-hand side of \eqref{eqn:finerdecomp}.  What
remains is to show that each of the three sets in the union
is indeed a $\Xi_{12}$-orbit.  It clearly suffices to show
that each $\Xi$-set is closed under the action of $\Xi_{12}$
and that the action of $\Xi_{12}$ on $\Xi$ is transitive.

Since the action of $\Xi_{12}$ on $\SL{2}{\Int[\vecti]}\alpha^{\rm N}(m,x)$
is by left-multiplication, it will suffice to prove, under
the action by left-multiplication of $\Xi_{12}$ on $\SL{2}{\Int[\vecti]}$,
each $\Xi$-subset is precisely one orbit.  Note that by \eqref{eqn:xi12defn},
the reduction map $\red_{1+\vecti}$ is a group epimorphism
of $\Xi_{12}$ onto $\{\overline{I},\overline{S}\}$.  The group
$\{\overline{I},\overline{S}\}$ has an induced action on
$\SL{2}{\Int[\vecti]/(1+\vecti)}$.  Further, the reduction map
\startdisp
\red_{1+\vecti}:\SL{2}{\Int[\vecti]}\rightarrow
\SL{2}{\Int[\vecti]/(1+\vecti)}
\finishdisp
respects the actions of $\Xi_{12}$ and its image $\{\overline{I},\overline{S}\}$.
Therefore, it will suffice to prove that
\starteqn\label{eqn:finerdecompinter}
\red_{1+\vecti}\Xi\;\text{is a $\{\overline{I},\overline{S}\}$ orbit,
for each}\; \Xi\in\{\Xi_{1},\,\Xi_2,\,\Xi_{12}\}.
\finisheqn
In order to prove \eqref{eqn:finerdecomp}, we compute the product
$\gamma_2\gamma_1\inv$ for $\gamma_1,\gamma_2\in\red_{1+\vecti}\Xi$.  Recall
the description of $\Xi$ given in \eqref{eqn:xisetdescription}.
It follows from \eqref{eqn:xisetdescription} that
\startdisp
\gamma_2\gamma_1\inv=\begin{cases}\overline{I}&\text{if}\;
\gamma_1=\gamma_2\\
\overline{S}&\text{if}\gamma_1\neq\gamma_2.
\end{cases}
\finishdisp
Therefore, for each $\Xi$-subset and fixed $\gamma\in\Xi$, we have
\starteqn\label{eqn:finerdecomp2}
\{\gamma_2\gamma\inv\;|\; \gamma_2\in\red_{1+\vecti}\Xi\}=
\{\overline{I},\overline{S}\}=\red_{1+\vecti}\Xi_{12}
\finisheqn
We have
\startdisp
\red_{1+\vecti}\Xi=\{(\gamma_2\gamma\inv)\gamma\;|\;\gamma_2\in\red_{1+\vecti}\Xi\}=
\{\overline{I},\overline{S}\}\gamma,\;\text{for each $\gamma\in\red_{1+\vecti}\Xi$},
\finishdisp
where the second equality follows from \eqref{eqn:finerdecomp2}.
We have verified \eqref{eqn:finerdecompinter}.  By the comments
preceding \eqref{eqn:finerdecompinter}, this completes
the proof of Sublemma \ref{sublem:finerdecomp}.
\end{myproof}
\begin{prop}\label{prop:inversematrixexplicitdescription}  Let $\conj$
be the morphism from $\SL{2}{\Complex}$ onto $G$
as in \eqref{eqn:imagematrixconj}.  Let $\Gamma=\SO{3}{\Int[\vecti]}$
be the group of integral points of $G$ in the coordinatization
of $G$ induced by the isomorphism \eqref{eqn:concreteorthequiv}.
Let the subsets $\Xi_1,\, \Xi_2,\,\Xi_{12}$ of $\SL{2}{\Int[\vecti]}$ be as
defined in \eqref{eqn:xi12defn} and \eqref{eqn:residuematrices}.
Let the matrices $\alpha^{\mathrm{N}}(m,x)$ be as in \eqref{eqn:MofNmxdefn}.
Let $\omega_8\in\Complex$ be as in
\eqref{eqn:omegaeightdefn}.
Then we have
\starteqn\label{eqn:inversematrixexplicitdescription}
\conj\inv(\Gamma)=\bigcup_{\delta,=0,1}
\hspace*{-0.5cm}\cdot\hspace*{.5cm}
\left(\frac{1}{\omega_8^{\delta}}
\Xi_{12}\alpha^{\vecti^{\delta}}\hspace{-0.5mm}(\vecti^{\delta},0)\bigcup\hspace*{-0.32cm}
\cdot\hspace*{0.32cm}
\left(\bigcup_{\epsilon=0,1}\hspace*{-0.45cm}\cdot\hspace*{.45cm}
\frac{1}{\omega_8^{\delta}(1+\vecti)}
\Xi_{2}\alpha^{2\vecti^{1+\delta}}\hspace*{-0.7mm}
(\vecti^{1+\delta},\vecti^{\epsilon})
\right)\right).
\finisheqn
\end{prop}
The proof of Proposition \eqref{prop:inversematrixexplicitdescription}
will be completed at the end of \S\ref{subsec:inverseimagedescription}.
\vspace*{.3cm}

\noindent\textbf{Remarks}\begin{itemize}
\item[(a)]  We use $\Int[\omega_8]$ to denote the ring
generated over $\Int$ by $\omega_8$.  By \eqref{eqn:omegaeightdefn}
we have $\Int[\vecti]\subset\Int[\omega_8]$ and $\Int[\omega_8]=
\Int[\omega_8,\vecti]$.  It follows from Proposition
 \ref{prop:inversematrixexplicitdescription} that
$\conj\inv(\Gamma)\subseteq \SL{2}{\Complex}$ is in fact a subset of
$\SL{2}{\Rational(\omega)}$.  More precisely, of the two
parts of the right-hand side of \eqref{eqn:inversematrixexplicitdescription},
we have
\starteqn\label{eqn:firstsubsetofinvim} \frac{1}{\omega_8^{\delta}}
\Xi_{12}\alpha^{\vecti^{\delta}}
(\vecti^{\delta},0)\subseteq\SL{2}{\Int[\vecti,\omega_8]}
\quad\text{for}\;\delta\in\{0,\,1\}, \finisheqn while
\starteqn\label{eqn:secondsubsetofinvim}
\left(\bigcup_{\epsilon=0,1}\hspace*{-0.45cm}\cdot\hspace*{.45cm}
\frac{1}{\omega_8^{\delta}(1+\vecti)}
\Xi_{2}\alpha^{2\vecti^{1+\delta}}(\vecti^{1+\delta},\vecti^{\epsilon})
\right)\subseteq\mathrm{SL}_2\left
(\Int\left[\vecti,\omega_8,\frac{1}{1+\vecti}\right]\right)
\quad\text{for}\;\delta\in\{0,\,1\} \finisheqn
\item[(b)]  One can easily verify that the set on the left-hand
side of \eqref{eqn:firstsubsetofinvim}
is closed under multiplication, while
the set on the left-hand side of
\eqref{eqn:secondsubsetofinvim} is not.  More precisely,
through a rather lengthy calculation,
not included here, one verifies that
\starteqn\label{eqn:twosubsetproductalt}
\text{for $(x,y)$ a pair of elements of the form of
\eqref{eqn:secondsubsetofinvim},
$xy$ is}\begin{cases}\text{of form \eqref{eqn:secondsubsetofinvim}}\\
\quad\quad\text{or}\\
\text{of form \eqref{eqn:firstsubsetofinvim}.}
\end{cases}
\finisheqn
with each possibility in \eqref{eqn:twosubsetproductalt} being
realized for an suitable pair $(x,y)$.
These calculations amount to a brute-force
verification of the fact that the right-hand
side of \eqref{eqn:inversematrixexplicitdescription} is closed under multiplication.
But, because $\Gamma$ is a group and $\conj$ a morphism, this fact
also follows from
Proposition \ref{prop:inversematrixexplicitdescription},
after \S\ref{subsec:inverseimagedescription}.
\end{itemize}

The explicit representation
of $\conj\inv(\Gamma)$ in \ref{prop:inversematrixexplicitdescription}
allows us to read off certain group-theoretic facts relating
$\conj\inv(\Gamma)$ to $\SL{2}{\Int[\vecti]}$.  In Lemma \ref{lem:indexlemma}
below we use the notation
\startdisp
[G:H]\;\text{is the index of $H$ in $G$, for any group $G$ with subgroup $H$}.
\finishdisp
\begin{lem}\label{lem:indexlemma}  Let $\conj\inv(\Gamma)$ be the subgroup
of $\SL{2}{\Complex}$ described above,
given explicitly in matrix form
in \eqref{eqn:inversematrixexplicitdescription}.  All the
other notation is also as in
Proposition \ref{prop:inversematrixexplicitdescription}.
\begin{itemize}\item[(a)]  We have
\startdisp
\conj\inv(\Gamma)\cap\SL{2}{\Int[\vecti]}=\Xi_{12}.
\finishdisp
\item[(b)]  We have
\starteqn\label{eqn:index6}
[\conj\inv(\Gamma):\Xi_{12}]=6,\quad [\SL{2}{\Int[\vecti]}:\Xi_{12}]=3.
\finisheqn
Explicitly, the six right cosets of $\Xi_{12}$ in $\conj\inv(\Gamma)$
are the two cosets obtained by letting $\delta$ range
over $\{0,1\}$ in
\startdisp
\frac{1}{\omega_8^{\delta}}
\Xi_{12}\alpha^{\vecti^{\delta}}\hspace{-0.5mm}(\vecti^{\delta},0)
\finishdisp
and the four cosets obtained by
letting $\delta, \epsilon$ range over $\{0,1\}$ independently in
\startdisp
\frac{1}{\omega_8^{\delta}(1+\vecti)}
\Xi_{12}\begin{pmatrix}1&1\\0&1\end{pmatrix}\alpha^{2\vecti^{1+\delta}}\hspace*{-0.7mm}
(\vecti^{1+\delta},\vecti^{\epsilon}).
\finishdisp
\end{itemize}
\end{lem}
\begin{myproof} {Proof} Letting $\delta,\, \epsilon$
range over $\{0,1\}$ independently, we see that the expression
for $\conj\inv(\Gamma)$
on the right-hand side \eqref{eqn:inversematrixexplicitdescription}
does indeed represent $\conj\inv(\Gamma)$ as the union
of six disjoint sets.  Now,
Sublemma \ref{sublem:finerdecomp} says
that each of the six subsets in the union is closed under
left-multiplication by $\Xi_{12}$.
Therefore, Sublemma \ref{sublem:finerdecomp} implies that
that each of the six sets in the union is actually a
right $\Xi_{12}$-coset.  These observations prove part (b),
except for the statement that
\startdisp
[\SL{2}{\Int[\vecti]}:\Xi_{12}]=3,
\finishdisp
which follows immediately from \eqref{eqn:xi12defn}.
\startdisp
\frac{1}{\omega_8^{\delta}(1+\vecti)}
\Xi_{12}\Xi_{2}\alpha^{2\vecti^{1+\delta}}\hspace*{-0.7mm}
(\vecti^{1+\delta},\vecti^{\epsilon}).
\finishdisp
We now prove part (a).  The coset
obtained by taking $\delta=0$ in
\startdisp
\frac{1}{\omega_8^{\delta}}
\Xi_{12}\alpha^{\vecti^{\delta}}\hspace{-0.5mm}(\vecti^{\delta},0)
\finishdisp
is exactly $\Xi_{12}$.
Now, since $\Xi_{12}\subset\SL{2}{\Int[\vecti]}$ by definition
of $\Xi_{12}$,
any of the six $\Xi_{12}$ subsets mentioned in Part (b)
is either contained within $\SL{2}{\Int[\vecti]}$,
or else is disjoint from $\SL{2}{\Int[\vecti]}$.  Therefore,
to complete the proof of Part (a), it suffices to
observe that
\startdisp
\frac{1}{\omega_8}
\Xi_{12}\alpha^{\vecti}\hspace{-0.5mm}(\vecti,0)
\finishdisp
is not contained in $\SL{2}{\Int[\vecti]}$, nor is
\startdisp
\frac{1}{\omega_8^{\delta}(1+\vecti)}
\Xi_{12}\Xi_{2}\alpha^{2\vecti^{1+\delta}}\hspace*{-0.7mm}
(\vecti^{1+\delta},\vecti^{\epsilon}).
\finishdisp
when $\epsilon,\delta$ take values in $\{0,1\}$.
\end{myproof}

\vspace*{0.3cm}\noindent\boldmath\textbf{Proofs
of Lemma \ref{lem:standardresiduerep} and Proposition \ref{prop:heckedecomp}}.
\ubmth We close \S\ref{subsec:flts} by giving two proofs
that were deferred in the course of the main exposition.
\vspace*{0.3cm}

\begin{myproof} {Proof of Lemma \ref{lem:standardresiduerep}} Since
\startdisp
2^{\lceil\frac{n}{2}\rceil}\equiv 0\mod (1+\vecti)^n,
\finishdisp
it is easy to see that every residue class modulo $(1+\vecti)^n$
has at least one representative of the form
\starteqn\label{eqn:standardresiduerepproof}
r+s\vecti\;\text{with $0\leq r,\,s<2^{\lceil\frac{n}{2}\rceil}$}.
\finisheqn
If $n$ is even, then $\lceil\frac{n}{2}\rceil=\frac{n}{2}$.  Therefore,
the set of elements
of the form \eqref{eqn:standardresiduerepproof}
is precisely $\Omega_{(1+\vecti)^n}$.  If $n$ is odd, then
$n-\lceil\frac{n}{2}\rceil=\frac{n}{2}-\frac{1}{2}$, so that
\startdisp
2^{n-\lceil\frac{n}{2}\rceil}(1+\vecti)\equiv 0\mod (1+\vecti)^n.
\finishdisp
If the representative of the form
\eqref{eqn:standardresiduerepproof} isnot in
$\Omega_{(1+\vecti)^n}$, then
subtract $2^{n-\lceil\frac{n}{2}\rceil}(1+\vecti)$
from the initial representative, and add back $2^{\lceil\frac{n}{2}\rceil}$
to the result, if necessary.  In this manner, we produce a representative
of the given residue class in $\Omega_{(1+\vecti)^n}$,
and we are done with the case $n$ even.  So,
we have now shown that every residue class modulo $(1+\vecti)^n$
has at least one representative in $\Omega_{(1+\vecti)^n}$.
It is now routine to very that there is exactly one
such representative in $\Omega_{(1+\vecti)^n}$.  One method is to
calculate the number of residue classes modulo $(1+\vecti)^n$
(\textit{i.e.}, $2^n$)
and compare that number with $\#\Omega_{(1+\vecti)^n}$,
which can be calculated directly from \eqref{eqn:omegaspecific}.
\end{myproof}

\begin{myproof}{Proof of Proposition \ref{prop:heckedecomp}}
In the remarks following the statement of the Proposition, we gave
a more explicit, but equivalent, form of Proposition \ref{eqn:heckedecomp}.
The more explicit form is the statement
that each element $\alpha\in\MtwoN$
has a unique expression of the form \eqref{eqn:heckeexplicit}.
It is the explicit form of Proposition \ref{prop:heckedecomp}
given in \eqref{eqn:heckeexplicit} that we will now prove,
\starteqn\label{eqn:heckeesplicitmultliplied}
\begin{pmatrix}a&b\\
c&d\end{pmatrix}=\begin{pmatrix}p&q\\
r&s\end{pmatrix}\begin{pmatrix}m&x\\0&\frac{N}{m}\end{pmatrix},\;
\text{with}\, m\in\Go,\; m|N,\; \frac{N}{m}\;\text{standard},\; x\in\Omega_
{\frac{N}{m}},\, pr-qs=1.
\finisheqn
We now give an explicit procedure for finding $m,x$ in which
each step is uniquely determined.
Then we fine $q,s$ so that \eqref{eqn:heckeesplicitmultliplied}
is verified.  The uniqueness in the proposition
follows from the fact that the procedure
uniquely determines $m,x$, hence the first matrix in the explicit
decomposition \eqref{eqn:heckeesplicitmultliplied},
hence also the second matrix in the
explicit decomposition \eqref{eqn:heckeesplicitmultliplied}.

Choose $\epsilon\in\{0,1,2,3\}$ so that
\startdisp
\frac{N}{\vecti^{\epsilon}\mathrm{GCD}(a,c)}\quad\text{is standard}.
\finishdisp
Set
\starteqn\label{eqn:mdefn}
m=\vecti^{\epsilon}\mathrm{GCD}(a,c).
\finisheqn
With this definition of $m$, we have $\frac{N}{m}$ standard.

Next, choose $p,r\in\Int[\vecti]$ satisfying
\starteqn\label{eqn:prdefn}
\begin{pmatrix}a\\ c\end{pmatrix}=m\begin{pmatrix}p\\ r\end{pmatrix}.
\finisheqn
By \eqref{eqn:mdefn} and \eqref{eqn:prdefn},
it is clear that $p,r\in\Int[\vecti]$ satisfying
\eqref{eqn:prdefn} exist and are uniquely determined.
From the determinant condition $ad-bc=N$ and \eqref{eqn:prdefn}
we have
\starteqn\label{eqn:detNconsequence}
pd-rb=\frac{N}{m}
\finisheqn
By property \textbf{GCD 2} and \eqref{eqn:prdefn}, we have
\starteqn\label{eqn:prrelprime}
\mathrm{GCD}(p,r)=1.
\finisheqn
By Property \textbf{GCD 3} there exist $q_0,s_0\in\Int[\vecti]$
with
\starteqn\label{eqn:unimodularityfirstmatrix}
ps_0-rq_0=1.
\finisheqn
By \eqref{eqn:omegamdefn}, we can define an element $x\in\Int[\vecti]$
by the conditions
\starteqn\label{eqn:xdefn}
x\in\Omega_{\frac{N}{m}}\;\text{such that}\; \red_{\frac{N}{m}}(x)=
\red_{\frac{N}{m}}(-qd+sb).
\finisheqn
At this point we have given the procedure for finding $m,\,x$
with each step uniquely determined, so it is clear that
$m,\,x$ are uniquely determined.  Now it remains to find $q,\,s$
so that \eqref{eqn:heckeesplicitmultliplied} is verified.

By \eqref{eqn:xdefn} there is a unique $\ell\in\Int[\vecti]$ such that
\starteqn\label{eqn:xdefnconsequence}
x=q_0d+s_0b+\ell\frac{N}{m}.
\finisheqn
Now define $q,\,s$ by the condition
\starteqn\label{eqn:sqdefn}
\begin{pmatrix}q\\s
\end{pmatrix}=
\begin{pmatrix}-p\ell+q_0\\-r\ell+s_0
\end{pmatrix}.
\finisheqn
Note first
that according to the definitions of $m,\,x$ and $p$ through $s$,
\begin{equation}\label{eqn:heckecomputation1}
\begin{aligned}
\begin{pmatrix}p&q\\
r&s\end{pmatrix}\begin{pmatrix}m&x\\0&\frac{N}{m}\end{pmatrix}&=&&
\begin{pmatrix}p&q_0\\
r&s_0\end{pmatrix}
\begin{pmatrix}1&-\ell\\
0&1\\
\end{pmatrix}
\begin{pmatrix}1&\ell\\
0&1
\end{pmatrix}
\begin{pmatrix}m&-q_0d+s_0b\\
0&\frac{N}{m}
\end{pmatrix}\\&=&&
\begin{pmatrix}p&q_0\\
r&s_0\end{pmatrix}
\begin{pmatrix}m&-q_0d+s_0b\\
0&\frac{N}{m}
\end{pmatrix}.
\end{aligned}
\end{equation}
Multiplying the matrices on the right side of \eqref{eqn:heckecomputation1},
then applying \eqref{eqn:unimodularityfirstmatrix}, \eqref{eqn:prdefn},
and \eqref{eqn:detNconsequence} we obtain
\starteqn\label{eqn:heckecomputation2}
\begin{pmatrix}pm&-pq_0d+ps_0b+q_0\frac{N}{m}\vspace*{0.2cm}\\
rm&-rq_0d+rs_0b+s_0\frac{N}{m}
\end{pmatrix}=\begin{pmatrix}a&q_0(-pd+rb+\frac{N}{m})+b\vspace*{0.2cm}\\
c&s_0(pd+rb+\frac{N}{m})+d\end{pmatrix}=\begin{pmatrix}a&b\\
c&d
\end{pmatrix}.
\finisheqn
Putting \eqref{eqn:heckecomputation1} and \eqref{eqn:heckecomputation2}
together, we have completed the verification of \eqref{eqn:heckeexplicit}.
By the above comments, this completes the proof of Proposition
\ref{prop:heckedecomp}.
\end{myproof}

\subsection{Proof of Proposition \ref{prop:inversematrixexplicitdescription}.}
\label{subsec:inverseimagedescription}
The present section is devoted to developing the machinery
used in proving  Proposition \ref{prop:inversematrixexplicitdescription}
and then completing the proof.
Since this machinery will not be used again in this chapter,
the reader may wish to skip this section on a first reading.
For easier reference in the course of these lemmas,
it will be useful to identify and label several properties
that may belong to the quadruple
\starteqn\label{eqn:quadruple}
(a,\,b,\,c,\,d)\in\Complex^4.
\finisheqn
In order to state these properties more easily, we make
the following definition.
\begin{defn}  A \textbf{permutation of the quadruple
\bmth$(a,\,b,\,c,\,d)$\ubmth} is a quadruple with the same
entries as $(a,\,b,\,c,\,d)$, but possibly in a different
order.
\end{defn}
It will be our convention to use $(x,\,y,\,z,\,w)$
for a permutation of $(a,\,b,\,c,\,d)$, and for the remainder
of the section, it will always be assumed, unless stated
otherwise, that the letters $x,y,z,w$ denote the elements
of an arbitrary permutation of $(a,\,b,\,c,\,d)$.
In order to illustrate the reason
for formalizing the notion of the permutation
of a quadruple, note the difference between speaking
of, on the one hand, \textit{the pair of elements} $(x,y)$ of the permutation $(x,y,z,w)$
of the quadruple $(a,\,b,\,c,\,d)$, and, on the other hand,
\textit{a pair of elements} $(x,y)$ of the quadruple $(a,\,b,\,c,\,d)$.
  Assuming that the entries $a,b,c,d$ are distinct complex numbers we may,
in the first case, infer that $x,y,z,w$ are also distinct complex numbers,
whereas in the second case, we cannot infer this because, for example,
it may be that $x=y=a$.

\vspace*{0.3cm}\noindent\boldmath\textbf{List of commonly
used properties of the quadruple $(a,b,c,d)\in\Complex^4$
of \eqref{eqn:quadruple}}.
\ubmth  In the list of properties,
we will always assume that $(x,\,y,\,z,\,w)$ is an arbitrary
permutation of \ref{eqn:quadruple}.
\begin{itemize}
\item[(i)] $ad-bc=1$.
\item[(ii)] $2xy\in \Int[\vecti]$.
\item[(iii)]  $x^2+y^2+z^2+w^2,\;x^2-y^2+z^2-w^2\in \left((1+\vecti)^2\right).$
\item[(iv)]  At least three of the elements
 of \eqref{eqn:quadruple} are nonzero.
\item[(v)]  Exactly two of the elements of \eqref{eqn:quadruple}
are nonzero.
\item[(vi)] $x^2\in\Rational(\vecti).$
\item[(vii)]  If $x\neq 0$ then
$x=\omega_8^{\pm j}\quad\text{for $j\in\{0,1,2,3,4\}$ depending only
on \eqref{eqn:quadruple}.}$
\end{itemize}

In certain situations, it will be convenient to refer to the following
weaker version of \linebreak Property (i):
\begin{itemize}
\item[(i$'$)] $ad-bc\neq 0$.
\end{itemize}
The reason for introducing Property (i$'$) can be glimpsed
from the following simple observation.
\begin{sublem}\label{sublem:multprops}
Let $(a,b,c,d)\in\Complex^4$ be a quadruple, and let $z\in\Complex$
such that $z^2\in\Complex-\{0\}$.
\begin{itemize}
\item[(a)]  The quadruple
\startdisp
z(a,b,c,d):=(za,zb,zc,zd)\in\Complex^4
\finishdisp
has Property (iv), respectively Property (v), if and only if
$(a,b,c,d)$ has Property (iv), respectively Property (v).
\item[(b)]  Assume in addition that $z\in\Int[\vecti]-\{0\}$.
Then, if $(a,b,c,d)$
has Property (i$'$), resp., Property (ii), resp. Property (iii),
then the quadruple $z(a,b,c,d)$
has Property (i$'$), resp., Property (ii), resp. Property (iii).
\end{itemize}
\end{sublem}

Sublemma \ref{sublem:multprops} is verified immediately
from the definitions of Properties (i$'$), (ii) and (iii).

\begin{sublem}\label{sublem:sqrational}  Let
$(a,\, b,\, c,\, d)\in\Complex^4$ be a quadruple as in
\eqref{eqn:quadruple}.  Then we have the following implications.
\begin{itemize}
\item[(a)] Properties (ii) and (iv) imply Property (vi).
\item[(b)] Properties (iii) and (v) imply Property (vi).
\item[(c)] Properties (i$'$), (iii), (v) together imply Property (vii).
\end{itemize}
\end{sublem}
\begin{myproof} {Proof} For (a), by Property (iv) we may assume
without loss of generality that $x,z,w\neq 0$.
We may express $x^2$ in the following form.
\starteqn\label{eqn:sqrationala}
x^2=\frac{(2xy)(2xz)}{(2zw)}.
\finisheqn
By the choice of $z,w$, the denominator of the right hand side of \eqref{eqn:sqrationala}
is nonzero, so \eqref{eqn:sqrationala} makes sense.  By property (ii),
each of the factors in the numerator and denominator of \eqref{eqn:sqrationala}
are in $\Int[\vecti]$, so \eqref{eqn:sqrationala} lies in $\Rational(\vecti)$.

For (b), if $x=0$, then obviously $x^2\in\Rational(\vecti)$.
By Property (v) we may assume that $x,y\neq 0$, while $z=w=0$.
Write
\starteqn\label{eqn:sqrationalb}
2x^2=(x^2+y^2)+(x^2-y^2).
\finisheqn
Since $z,w=0$, Property (iii) implies that the two terms in parentheses
on the right side of \eqref{eqn:sqrationalb} belong to $\left((1+\vecti)^2\right)$.
So \eqref{eqn:sqrationalb} implies that
we have $2x^2\in \left((1+\vecti)^2\right)$.  Therefore,
\starteqn\label{eqn:sublempartb}
\text{Properties (iii) and (v) together imply}\; x^2\in\Int[\vecti]
\subseteq \Rational(\vecti).
\finisheqn
By \eqref{eqn:sublempartb}, Properties (iii) and (v) together certainly
imply Property (vi).

For (c), by Property (v), we may assume that $x,y\neq 0$
and $z=w=0$.  By (i$'$), we must have that $\{x,y\}=\{a,d\}$ or $\{b,c\}$.
In either case, substitution into (i$'$) yields
\starteqn\label{eqn:unitxy}
xy=\pm 1,\quad\text{so that}\; (x^2)(y^2)=1.
\finisheqn
By Properties (iii), (v), and \eqref{eqn:sublempartb},
we have that $x^2,\, y^2\in\Int[\vecti]$.  Therefore, \eqref{eqn:unitxy}
implies that
$x^2$ is a unit in $\Int[\vecti]$.  But the units
of $\Int[\vecti]$ are $\vecti^{\pm \epsilon}$, $\epsilon\in\{0,1,2\}$,
so that by \eqref{eqn:omegaeightdefn},
$x=\omega_8^{\pm\delta}$ for $\delta\in\{0,1,2,3,4\}$.  This
completes the proof of (c).
\end{myproof}

\begin{sublem} \label{sublem:props45} Let $(a,b,c,d)\in\Complex^4$ as
in \eqref{eqn:quadruple}.  Then we have the implicaton
\startdisp
\text{Property (i$'$) implies either property (iv) or (v).}
\finishdisp
\end{sublem}
\begin{myproof}{Proof}
It is clear that Property (i) implies that at least
2 of the entries of $(a,\,b,\,c,\,d)$
are nonzero, so \eqref{eqn:quadruple} has Property (iv) or Property (v).
\end{myproof}

\begin{sublem}\label{sublem:props34}
Let $(a,\,b,\,c,\,d)$ as
in \eqref{eqn:quadruple} be a quadruple
of complex numbers satisfying the conditions of \eqref{eqn:quadruplecond}.  Then
\begin{itemize}
\item[(a)]  We have that \eqref{eqn:quadruple} satisfies both properties (ii)
and (iii).
\item[(b)]  We have that \eqref{eqn:quadruple} satisfies either property
(iv) or (v).
\end{itemize}
\end{sublem}
\begin{myproof}{Proof}  The conditions of \eqref{eqn:quadruplecond} consist
of Property (i) above and the statement that the entries of
\eqref{eqn:imagematrixconj} belong to $\Int[\vecti]$.  Property (i)
obviously implies Property (i$'$).  Therefore, we can apply Sublemma
\ref{sublem:props45} to obtain part (b).

For part (a), we explain why the statement that the entries of
\eqref{eqn:imagematrixconj} belong to $\Int[\vecti]$ implies both Properties
(ii) and (iii).  As a consequence of Property (i), we have
\starteqn\label{eqn:detintegral}
ad-bc\in\Int[\vecti].
\finisheqn
Assume as above that $(x,y,z,w)$ is an arbitrary permutation of \eqref{eqn:quadruple}.
By examining the entries in the third
row and column of \eqref{eqn:imagematrixconj}
and by using \eqref{eqn:detintegral}, we see that
\starteqn\label{eqn:squareslemmaa1}
(xy+zw),\, (xy-zw)\in \Int[\vecti].
\finisheqn
A routine calculation shows that
\starteqn\label{eqn:squareslemmaa2}
2xy\in\Int\left[\vecti,\,(xy+zw), \,(xy-zw)\right].
\finisheqn
Together \eqref{eqn:squareslemmaa1} and \eqref{eqn:squareslemmaa2} imply
that $2xy\in\Int[\vecti]$, which is Property (ii).  For Property (iii),
we see by direct inspection of the upper left $2$-by-$2$
block of entries of \eqref{eqn:imagematrixconj} that
\startdisp
\frac{x^2+y^2+z^2+w^2}{2},\;\frac{(x^2-y^2+z^2-w^2)}{2}\,\in\,\Int[\vecti].
\finishdisp
We have
\startdisp
x^2+y^2+z^2+w^2, x^2-y^2+z^2-w^2\in 2\Int[\vecti]=\left((1+\vecti)^2\right).
\finishdisp
So \eqref{eqn:quadruple} satisfies Property (iii).
\end{myproof}

\begin{lem}  \label{lem:squaresrational} Let $(a,\,b,\,c,\,d)$ as
in \eqref{eqn:quadruple} be a quadruple
of complex numbers satisfying the conditions of \eqref{eqn:quadruplecond}.
\begin{itemize}
\item[(a)] We have \eqref{eqn:quadruple} has property (vi).
\item[(b)] In particular if
\eqref{eqn:quadruple} has property (v), then \eqref{eqn:quadruple}
has property (vii).
\end{itemize}
\end{lem}
\begin{myproof}{Proof}
By Sublemma \ref{sublem:props34}, part (a),
\eqref{eqn:quadruple} satisfies both properties (ii)
and (iii).  By Sublemma \ref{sublem:props34}, part (b),
\eqref{eqn:quadruple} satisfies either (iv) or (v).
If \eqref{eqn:quadruple} satisfies (iv), then we apply
part (a) of Sublemma \ref{sublem:sqrational}, and if \eqref{eqn:quadruple}
satisfies (v), then we apply part (b) of Sublemma \ref{sublem:sqrational},
to conclude that \eqref{eqn:quadruple} satisfies Property (vi).
This completes the proof of (a).

Part (b) follows from Sublemma \ref{sublem:props34}, part (a) and
Sublemma \ref{sublem:sqrational},
part (c).
\end{myproof}


\vspace*{0.3cm}\noindent\boldmath\textbf{The map $\ord_{\nu}$
and its properties}.
\ubmth
Let $\Go$ be a PID with field of fractions $k$.
For any prime $\nu\in\Go$, the $\nu$-adic
valuation on $k$, written $\ord_{\nu}$ has its usual meaning.  Namely,
\startdisp
\ord_{\nu}: k\rightarrow \Int\cup\infty,
\finishdisp
such that for all $x\in k^*$,
\startdisp
x=(\nu)^{\ord_{\nu}(x)}\frac{p_x}{q_x},\quad\text{with $p_x,q_x\in\Go$ and
$x\not|\ p_xq_x$},
\finishdisp
and by convention,
\startdisp
\ord_{\nu}(0)=\infty.
\finishdisp
Thus, $\ord_{\nu}(x)$ is positive (resp. negative, resp. zero)
when $\nu$ divides the numerator (resp. denominator, resp. neither numerator
nor denominator) of $x\in k^*$ written ``in lowest terms".
We recall the following elementary properties of  $\ord_{\nu}$.
In each property $x,y$ are arbitrary elements of $\Go$.
For any two elements $x,y\in\Go$, we use the notation $x\sim y$
to indicate that $x$ and $y$ differ by multiplication by a unit in $\Go$.
\begin{itemize}
\item[\textbf{ORD 1}] $\ord_{\nu}(xy)=\ord_{\nu}(x)+\ord_{\nu}(y)$.
\item[\textbf{ORD 2}] $\ord_{\nu}(x+y)\geq\min\left(\ord_{\nu}(x),
\ord_{\nu}(y)\right).$
\item[\mbox{\textbf{ORD 3}}] $\ord_{\nu}(x)>\ord_{\nu}(y)\;
\text{implies}\;\ord_{\nu}(x+y)=\ord_{\nu}(y).$
\item[\textbf{ORD 4}]  We have \startdisp
y\sim\prod_{\nu}\nu^{\ord_{\nu}(y)},
\;\textit{i.e.},\;
y=\vecti^{\epsilon(y)}\prod_{\nu}\nu^{\ord_{\nu}(y)}, \; \text{for some}\,
u\in\Go^*,\finishdisp
where the product is taken over all standard primes $\nu$ and
is actually finite, because \linebreak $\ord_{\nu}(y)=0$ for almost all $\nu$,
and $\epsilon(y)\in\{0,1,2\}$.
\item[\textbf{ORD 5}]  Let $n$ be a positive integer, $x\in\Go$,
$u\in\Go^*$.
Then we have
\startdisp
\min\left(\ord_{\nu}(x),\ord_{\nu}(x+u\nu^n)\right)\leq n.
\finishdisp
\end{itemize}
Properties \textbf{ORD 1}, \textbf{ORD 2}, and are
covered in all standard treatments of the subject, see
\textit{e.g.} \cite{hindrysilverman}, p. 170, and \textbf{ORD 4}
is easy to verify from the definitions,
so we omit the proofs of these three properties.

\begin{myproof}{Proof of \textbf{ORD 3} and \textbf{ORD 5}}
Property \textbf{ORD 3} is clear whenever either $x$ or $y$ is zero,
so we may assume that $x,y\in k^*$.
By the definition of $\ord_{\nu}$ we have
\starteqn\label{eqn:valpropelemproof1}
x=\nu^{\ord_{\nu}(x)}\frac{p_x}{q_x},\;y=\nu^{\ord_{\nu}(y)}\frac{p_y}{q_y},
\quad\text{with}\; \nu\not|\ p_xq_xp_yq_y.
\finisheqn
One sees that
\starteqn\label{eqn:valpropelemproof2}
x+y=\nu^{\ord_{\nu}(y)}\frac{p_yq_x+\nu^{\ord_{\nu}(x)-\ord_{\nu}(y)}p_xq_y}
{q_xq_y}
\finisheqn
By \eqref{eqn:valpropelemproof1}, $\nu \not|\ q_xq_y$.  We claim that also
\starteqn\label{eqn:valpropfalse}
\nu\not|\ \left(p_yq_x+\nu^{\ord_{\nu}(x)-\ord_{\nu}(y)}p_xq_y\right).
\finisheqn
For suppose that \eqref{eqn:valpropfalse} is false.
Then because $\ord_{\nu}(x)>\ord_{\nu}(y)$,
we have
\startdisp
\nu | \nu^{\ord_{\nu}(n)-\ord_{\nu}(y)}p_xq_y.
\finishdisp
We obtain $\nu | p_yq_x$.  But $\nu | p_yq_x$ contradicts
\eqref{eqn:valpropelemproof1}.
This contradiction proves \eqref{eqn:valpropfalse}.  By
\eqref{eqn:valpropfalse},
\eqref{eqn:valpropelemproof2} can be written as
\startdisp
x+y=\nu^{\ord_{\nu}(y)}\frac{p_{x+y}}{q_{x+y}}\quad\text{with $\nu\not|\
p_{x+y}q_{x+y}$}.
\finishdisp
Here,
\startdisp
p_{x+y}=p_yq_x+\nu^{\ord_{\nu}(x)-\ord_{\nu}(y)}p_xq_y,\quad\text{and}\;
q_{x+y}=q_xq_y.
\finishdisp
This completes the proof of Property \textbf{ORD 3}.

Before proving \textbf{ORD 5}, note that \textbf{ORD 4} trivially
implies that
\starteqn\label{eqn:unitshaveorderzero}
\ord_{\nu}(u)=0\forallindisp\;\text{$\nu$ prime, $u\in\Go^*$.}
\finisheqn

In proving Property \textbf{ORD 5}, we may assume without loss
of generality that
\startdisp
\ord_{\nu}(x)>n.
\finishdisp
By \textbf{ORD 1} and \eqref{eqn:unitshaveorderzero},
we have $\ord_{\nu}(u\nu^n)=n$.  We may
therefore apply \textbf{ORD 3} with $y=u\nu^n$, and we
deduce that
\startdisp
\ord_{\nu}(x+u\nu^n)=\ord_{\nu}(u\nu^n)=n.
\finishdisp
Thus, we have verified Property \textbf{ORD 5}.
\end{myproof}
\begin{rem}  We must assume strict inequality
in the hypothesis of \textbf{ORD 3}, as is easily seen
from the example of $\Go=\Int[\vecti]$, $\nu=(1+\vecti)$,
$x=y\in\Int[\vecti]-\{0\}$.
\end{rem}

\vspace*{0.3cm}\noindent\boldmath\textbf{Application
of $\ord_{\nu}$ to the quadruple \eqref{eqn:quadruple}}.
\ubmth
For Sublemma \ref{sublem:primefactors}, it will be more convenient to use
the following alternate forms of Properties (ii) and (iii).
\begin{itemize}
\item[(ii$'$)] For each pair $(x,y)$ of elements of
\eqref{eqn:quadruple} we have
\startdisp
4x^2y^2\in \Int[\vecti].
\finishdisp
\item[(iii$'$)] For each pair $(x,y)$ of elements of \eqref{eqn:quadruple} we have
\startdisp
x^2+y^2\in \Int[\vecti].
\finishdisp
\end{itemize}
It is easy to see that each of Properties (ii) and (iii)
implies its ``primed" form.

We shall use the above notation and the \textbf{ORD}
properties primarily in the case when $\Go=\Int[\vecti]$,
so that $k=\Rational(\vecti)$.
\begin{sublem}\label{sublem:primefactors}  Assume for all that
follows that \eqref{eqn:quadruple}
satisfies Property (vi).  Let $(x,y,z,w)$ be an arbitrary
permutation of \eqref{eqn:quadruple}.
\begin{itemize}\item[(a)]
Assume that \eqref{eqn:quadruple} has Property (ii$'$).
Let $\nu\in\Int[\vecti]$ be a (Gaussian) prime.   Then we have
\starteqn\label{eqn:sublemprimefactorspta}
\ord_{\nu}(x^2)+\ord_{\nu}(y^2)\geq -\ord_{\nu}(4).
\finisheqn
\item[(b)]  Assume that \eqref{eqn:quadruple} has Property (iii$'$).  Then
\starteqn\label{eqn:sublemprimefactorsptb}
\ord_{\nu}(x^2)<0\;\text{implies}\; \ord_{\nu}(x^2)=\ord_{\nu}(y^2).
\finisheqn
\item[(c)]  Assume that \eqref{eqn:quadruple} has both Properties
(ii$'$) and (iii$'$).  Then we have
\starteqn\label{eqn:sublemmaprimefactorsptcnot2}
\ord_{\nu}(x^2)\geq 0  \quad\text{for $\nu\neq 1+\vecti$},
\finisheqn
and
\starteqn\label{eqn:sublemmaprimefactorsptc2}
\ord_{1+\vecti}(x^2)\geq -2.
\finisheqn
\end{itemize}
\end{sublem}
\begin{myproof}{Proof}  Property (vi) guarantees that all the quantities
mentioned in the lemma belong to $\Rational(\vecti)$, that is,
these quantities are in the domain of $\ord_{\nu}$.

Part (a) is an application of the fact that $\ord_{\nu}$ is positive
on $\Int[\vecti]$, combined with \textbf{ORD 1},
applied to the statement of Property (ii$'$).

For (b), suppose otherwise.  Then \eqref{eqn:quadruple}
satisfies property (iii$'$), and we have both \linebreak
\mbox{$\ord_{\nu}(x^2)<0$},
and $\ord_{\nu}(y^2)\neq\ord_{\nu}(x^2)$.
Since $\ord_{\nu}(y^2)<\ord_{\nu}(x^2)$
we may assume without loss of generality that
$\ord_{\nu}(y^2)>\ord_{\nu}(x^2)$.  By \textbf{ORD 3},
we have
\starteqn\label{eqn:falseassumptiononval}
\ord_{\nu}(x^2+y^2)=\ord_{\nu}(x^2),\;\text{so that}\;
\ord_{\nu}(x^2+y^2)<0
\finisheqn
Property (iii$'$) implies that $\ord_{\nu}(x^2+y^2)\geq0$.  But
\eqref{eqn:falseassumptiononval} says that contrary.  Therefore,
our assumptions imply a contradiction.  This proves part (b)
of the Sublemma.

We now prove (c).  First, since $4\sim (1+\vecti)^4$, we compute that
\starteqn\label{eqn:foursvaluation}
\ord_{\nu}(4)=0,\,\text{for $\nu\neq 1+\vecti$},\;\text{whereas}
\;\ord_{1+\vecti}(4)=4.
\finisheqn
Clearly, to prove (c), we may assume that $\ord_{\nu}(x^2)<0$, because
if $\ord_{\nu}(x^2)\geq 0$, (c) is trivial.  Therefore,
(b), gives $\ord_{\nu}(x^2)=\ord_{\nu}(y^2)$.  Substituting this
equality into \eqref{eqn:sublemprimefactorspta}, and solving for $\ord_{\nu}(x^2)$
we obtain
\starteqn\label{eqn:ordnuxsquaredestimate}
\ord_{\nu}(x^2)\geq \frac{\ord_{\nu}(4)}{2}.
\finisheqn
Applying \eqref{eqn:foursvaluation} to \eqref{eqn:ordnuxsquaredestimate} gives
\eqref{eqn:sublemmaprimefactorsptcnot2} and \eqref{eqn:sublemmaprimefactorsptc2}.
\end{myproof}

\begin{sublem} \label{sublem:parity} Suppose that \eqref{eqn:quadruple}
satisfies Property (ii).  Let $(x,y,z,w)$ be an arbitrary permutation
of \eqref{eqn:quadruple}.
\begin{itemize}
\item[(a)]
Then
\starteqn\label{eqn:ordersquareeven}
\ord_{\nu}(4x^2y^2)\in 2\Int,\forallindisp \nu\;\text{prime.}
\finisheqn
\item[(b)]  Assume that \eqref{eqn:quadruple} satisfies
Property (vi).  Then
\starteqn\label{eqn:orderssameparity}
\ord_{\nu}(x^2)\equiv\ord_{\nu}(y^2)\!\!\!\!\mod 2,\,\,
\text{for all}\,\, \nu\;\text{prime.}
\finisheqn
\item[(c)]  Assume that \eqref{eqn:quadruple} satisfies
Property (i).  Then we have
\startdisp
\min(\ord_{\nu}(ad),\ord_{\nu}(bc))\leq 0\forallindisp\nu\;\text{prime.}
\finishdisp
\end{itemize}
\end{sublem}
\begin{myproof}{Proof}  For (a), Property (ii) implies that
$2xy$ is in the domain of $\ord_{\nu}$.  Since $4x^2y^2=(2xy)^2$,
we have from \textbf{ORD 1} that
\startdisp
\ord_{\nu}(4x^2y^2)=2\ord_{\nu}(2xy)\in 2\Int.
\finishdisp

For (b), Property (vi) implies that each entry $x^2$ lies in the
domain of $\ord_{\nu}$.  By \textbf{ORD 1}
and \eqref{eqn:foursvaluation}, we have
\starteqn\label{eqn:sublemparityptbpf}
\ord_{\nu}(4x^2y^2)=\ord_{\nu}(4)+\ord_{\nu}(x^2)+\ord_{\nu}(y^2)
\equiv \ord_{\nu}(x^2)+\ord_{\nu}(y^2)\mod 2.
\finisheqn
By \eqref{eqn:ordersquareeven}, the left side of
\eqref{eqn:sublemparityptbpf} is even.  Therefore,
\eqref{eqn:sublemparityptbpf} implies that
\startdisp
\ord_{\nu}(x^2)+\ord_{\nu}(y^2)\equiv 0\mod 2,
\finishdisp
which gives \eqref{eqn:orderssameparity}.

For (c), Property (ii) implies that $2ad,\, 2bc\in\Int[\vecti]$,
from which it certainly follows that
\startdisp
ad,\,bc\in \Rational(\vecti).
\finishdisp
Thus, $ad,\, bc$ are in the domain of $\ord_{\nu}$'s.
By \textbf{ORD 2}, we have
\starteqn\label{eqn:paritysublemc1}
\min(\ord_{\nu}(ad),\ord_{\nu}(bc))\leq \ord_{\nu}(ad-bc).
\finisheqn
But by Property (i),
\starteqn\label{eqn:paritysublemc2}
\ord_{\nu}(ad-bc)=\ord_{\nu}(1)=0.
\finisheqn
Combining \eqref{eqn:paritysublemc1} and \eqref{eqn:paritysublemc2},
we obtain part (c).
\end{myproof}

\vspace*{0.3cm}\noindent\boldmath\textbf{The three
possibilities for the order at $1+\vecti$
of the entries of \eqref{eqn:quadruple} satisfying \eqref{eqn:quadruplecond}}.
\ubmth
With $(a,\,b,\,c,\,d)$ a quadruple as in \eqref{eqn:quadruple},
we define the $\textbf{P}i$, $i\in\{0,1,2\}$,
possibilities that \eqref{eqn:quadruple} may
satisfy.  In the definition of the $\textbf{P}i$ Possibilities,
it is assumed as usual that $(x,\,y,\,z,\,w)$ is an arbitrary
permutation of \eqref{eqn:quadruple}.  Also, each $\textbf{P}i$
Possibility includes the implicit condition that all
of the quantities mentioned are rational,
\textit{i.e.} lie in the domain of $\ord_{1+\vecti}$.
\begin{itemize}
\item[$\textbf{P}$0]  $\ord_{1+\vecti}(x^2)\in 2\Int_{\geq 0}$.  Further,
\starteqn\label{eqn:p1}
\ord_{1+\vecti}(a^2d^2)\ord_{1+\vecti}(b^2c^2)=0.
\finisheqn
\item[$\textbf{P}$1]  $\ord_{1+\vecti}(x^2)=-1$.
\item[$\textbf{P}$2]  $\ord_{1+\vecti}(x^2)=-2$.
\end{itemize}
The notation reminds us that for $(a,b,c,d)$ satisfying
$\textbf{P}i$ for $i=1,2$ we have $\ord_{1+\vecti}(x^2)=i$
and $\ord_{1+\vecti}(x^2)\geq i$ when $i=0$, with equality
satisfied for at least one permutation $(x,y,z,w)$.
It is obvious that \eqref{eqn:quadruple} satisfies
at most one of the $\textbf{P}i$ possibilities.  Thus, the point of Lemma
\ref{lem:squaresorders}, part (b), is that, under the hypothesis
of \eqref{eqn:quadruplecond}, \eqref{eqn:quadruple}
satisfies at least one of the $\textbf{P}i$ possibilities.
\begin{lem} \label{lem:squaresorders} Let $(a,\,b,\,c,\,d)$ as
in \eqref{eqn:quadruple} be a quadruple
of complex numbers satisfying the conditions of \eqref{eqn:quadruplecond}.
As usual, $(x,y,z,w)$ denotes an arbitrary permutation of $(a,b,c,d)$
\begin{itemize}
\item[(a)]   We have
\starteqn\label{eqn:oddprimesordersposeven}
\ord_{\nu}(x^2)\in 2\Int_{\geq 0},\forallindisp \nu\;\text{prime},\;\nu\neq 1+\vecti.
\finisheqn
\item[(b)]  Exactly one of the possibilities (P0) through (P2) hold.
\end{itemize}
\end{lem}
\begin{myproof}{Proof}  By Lemma \ref{lem:squaresrational}, part (b),
\eqref{eqn:quadruple} satisfies Property (iv).  By Sublemma
\ref{sublem:props34}, \eqref{eqn:quadruple} has Properties (ii)
and (iii), so \eqref{eqn:quadruple} also has Properties (ii$'$)
and (iii$'$). Further, by \eqref{eqn:quadruplecond}, \eqref{eqn:quadruple}
has property (i).  Therefore, all the parts of Sublemmas
\ref{sublem:primefactors} \ref{sublem:parity} apply.

With $x,\,\nu$ as in part (a), we have by Sublemma \ref{sublem:primefactors},
part (c),
\starteqn\label{eqn:ordnuxsqpos}
\ord_{\nu}(x^2)\in\Int_{\geq 0}.
\finisheqn
To see that actually $\ord_{\nu}(x^2)\in 2\Int_{\geq 0}$, assume
otherwise.  By \eqref{eqn:ordnuxsqpos}, we may assume that
\starteqn
\label{eqn:parityimpossibility}
\ord_{\nu}(x^2)\in 1+2\Int_{\geq 0}.
\finisheqn
By Sublemma \ref{sublem:parity}, part (b), we then have,
\startdisp
\ord_{\nu}(y^2), \ord_{\nu}(z^2),\ord_{\nu}(w^2)\in 1+\Int_{\geq 0}.
\finishdisp
Therefore, it follows from \textbf{ORD 1} that
\startdisp
\ord_{\nu}(x^2y^2)\in 2+\Int_{\geq 0}\;\text{for an arbitrary
permutation of \eqref{eqn:quadruple}.}
\finishdisp
In particular, we deduce that
\starteqn\label{eqn:adbcimpossibleorders}
\ord_{\nu}(a^2d^2),\,\ord_{\nu}(b^2c^2)>0.
\finisheqn
Using Sublemma \ref{sublem:parity},
part (c), together with \textbf{ORD 1}, we deduce that
\startdisp
\min\left(\ord_{\nu}(a^2d^2),\ord_{\nu}(b^2c^2)\right)\leq 0.
\finishdisp
Together with \eqref{eqn:adbcimpossibleorders}, we have reached
a contradiction.  So \eqref{eqn:adbcimpossibleorders}
is false, as is assumption \eqref{eqn:parityimpossibility},
since \eqref{eqn:parityimpossibility}
implies \eqref{eqn:adbcimpossibleorders}.
This concludes the proof of (a).

By \eqref{eqn:sublemmaprimefactorsptc2} in
Sublemma \ref{sublem:primefactors}, part (b), we have
\startdisp
\ord_{1+\vecti}(x^2)\geq -2,
\finishdisp
so that $\ord_{1+\vecti}(x^2)=-1$ or $-2$, or $\ord_{1+\vecti}(x^2)\geq0$.
If
\starteqn\label{eqn:someordernegative}
\ord_{1+\vecti}(x^2)=-1\;\text{or}\;-2,
\finisheqn
then by Sublemma \ref{sublem:primefactors}, part (b), we have
\startdisp
\ord_{1+\vecti}(y^2)=\ord_{1+\vecti}(x^2)\forallindisp\,y\;\text{entry
of \eqref{eqn:quadruple}}.
\finishdisp
Since the permutation $(x,\,y,\,z,\,w)$,
is arbitrary, the two possibilities in \eqref{eqn:someordernegative}
correspond to Possibilities $\textbf{P}1$
and $\textbf{P}2$.  Aside from
\eqref{eqn:someordernegative}, the only possibility is that
\startdisp
\ord_{1+\vecti}(x^2)\geq 0.
\finishdisp
In this case, it follows from Sublemma \ref{sublem:parity} that
the conditions of $\textbf{P}0$ are satisfied.  First, by part (c)
of Sublemma \ref{sublem:parity}, we must have \eqref{eqn:p1}.
In order to satisfy \eqref{eqn:p1}, we must have
$\ord_{1+\vecti}(x^2)=0$ for some permutation
$(x,\,y,\,z,\,w)$ of \eqref{eqn:quadruple}.  Then
$\ord_{1+\vecti}(x^2)\in 2\Int_{\geq 0}$, and Sublemma \ref{sublem:parity},
part (b) implies that $\ord_{1+\vecti}(x^2)\in 2\Int_{\geq 0}$
for all elements $x$ of \eqref{eqn:quadruple}.
\end{myproof}

Assessing our progress towards the proof of
Proposition \ref{prop:inversematrixexplicitdescription},
we have up until this point, given necessary, but not sufficient
conditions, for a quadruple \eqref{eqn:quadruple} to be the entries
of a matrix in $\conj\inv(\Gamma)$, equivalently
to satisfy \eqref{eqn:quadruplecond}.  The necessary but
not sufficient conditions are that \eqref{eqn:quadruple} satisfies
exactly of the $\textbf{P}i$ possibilities,
$i\in\{0,1,2\}$, which are given in Lemma \ref{lem:squaresorders},
part (b).  We now proceed by analyzing what the necessary
and sufficient conditions \eqref{eqn:quadruplecond} tell us
about a quadruple
\eqref{eqn:quadruple} satisfying one of the $\textbf{P}i$ possibilities.
Once again, we start with some easy and general observations.

\vspace*{0.3cm}\noindent\boldmath\textbf{The units in the factorization
of $x^2$ and $xy$}.
\ubmth In what follows, we make use of the square map $\sq$
sending $z$ to $z^2$, on various domains and on various
restricted parts of these domains.  For Sublemma \ref{sublem:unitssquareroot}
below we only need the fact that
\starteqn\label{eqn:squaremapcstar}
\mathrm{sq}|_{\Complex^{\times}}: \Complex^{\times}\rightarrow\Complex^{\times}\;
\text{is an epimorphism with kernel $\{\vecti^{\epsilon}\;|\; \epsilon
\in\{0,2\}\}$}.
\finisheqn

\begin{sublem}\label{sublem:unitssquareroot}
Let $x\in\Complex$ such that $x^2\in\Rational(\vecti)$.  Assume that
\starteqn\label{eqn:supposeorderseven}
\ord_{\nu}(x^2)\in 2\Int,\forallindisp\nu\;\text{prime}.
\finisheqn
Let $\omega_8$ be as in \eqref{eqn:omegaeightdefn}.
Let $\epsilon(x^2)\in\{0,1,2,3\}$ be as defined in \textbf{ORD 4}.
Define $\delta(x)$ by
\starteqn\label{eqn:deltaxdefn}
\delta(x)=\red_{2}(\epsilon(x^2)).
\finisheqn
so that $\delta(x)=0$ or $1$ depending on whether $\epsilon(x^2)$ is
even or odd.
Then,we have
\startdisp
\omega_{8}^{\delta(x)}x\in\Rational(\vecti).
\finishdisp
Further,
\starteqn\label{eqn:omegaxdescription}
\omega_{8}^{\delta(x)}x\sim\prod_{\nu}\nu^{\ord_{\nu}(x^2)/2}
\finisheqn
\end{sublem}
\begin{myproof}{Proof}
Using \textbf{ORD 4}, \eqref{eqn:squaremapcstar},
and \eqref{eqn:omegaeightdefn} we verify that
\starteqn\label{eqn:squarerootslemmainter}
x=\pm\omega_8^{\epsilon(x^2)}\prod_{\nu}\nu^{\ord_{\nu}(x^2)/2}.
\finisheqn
Here, the hypothesis \eqref{eqn:supposeorderseven}
guarantees that the numbers $\ord_{\nu}(x^2)/2$
belong to $\Int$.
 With $\delta(x)$ defined as in \eqref{eqn:deltaxdefn},
we clearly have
\startdisp
(\epsilon(x^2)+\delta(x))/2\in \Int.
\finishdisp
Multiplying both sides of \eqref{eqn:squarerootslemmainter}
by $\omega_8^{\delta(x)}$ and applying \eqref{eqn:omegaeightdefn},
we obtain
\startdisp
\omega_8^{\delta(x)}x=\pm\vecti^{(\epsilon(x^2)+\delta(x))/2}
\prod_{\nu}\nu^{\ord_{\nu}(x^2)/2}.
\finishdisp
The above expression for $\omega_8^{\delta(x)}x$ shows that it is
in $\Rational(\vecti)$ and differs from $\prod_{\nu}\nu^{\ord_{\nu}(x^2)/2}$
only by multiplication by the unit $\pm\vecti^{(\epsilon(x^2)+\delta(x))/2}$.
This completes the proof of Sublemma \ref{sublem:sqrational}.
\end{myproof}

\begin{sublem} \label{sublem:sqrationalindep}  Suppose \eqref{eqn:quadruple}
satisfies Properties (ii) and (iv).  Let $(x,y,z,w)$ be a permutation
of \eqref{eqn:quadruple}.  Set $\delta=\delta(x)$,
where $\delta(x)$ is as defined in \eqref{eqn:deltaxdefn}.
Then
\starteqn
\label{eqn:sqrationalindepsublem}
\omega_8^{\delta}y\in\Rational(\vecti).
\finisheqn
Further,
\starteqn
\omega_{8}^{\delta}y\sim\prod_{\nu}\nu^{\ord_{\nu}(y^2)/2}.
\finisheqn
\end{sublem}
\begin{myproof}{Proof}  Because \eqref{eqn:quadruple}
satisfies Property (iv),
$x,y$ satisfy
the hypotheses of Sublemma \ref{sublem:sqrational}.  Therefore,
$\delta(x), \delta(y)$ exist, as in the conclusion of Sublemma \ref{sublem:sqrational}.
By the choice of $\delta(x), \delta(y)$
\startdisp
\left(\omega_8^{\delta(x)}x\right)\left(\omega_8^{\delta(y)}y\right)=
\omega_8^{\delta(x)+\delta(y)}xy\in\Rational(\vecti).
\finishdisp
On the other hand, Property (ii) implies that $xy\in\Rational(\vecti)$.
We deduce that
\starteqn\label{eqn:omega8xy}
\omega_8^{\delta(x)+\delta(y)}\in\Rational(\vecti).
\finisheqn
Since $\delta(x), \delta(y)\in\{0,1\}$,
$\delta(x)+\delta(y)$ is $0$, $1$, or $2$.  Since
$\omega_8\notin\Rational(\vecti)$, $\delta(x)+\delta(y)=0$ or $2$.
We deduce that $\delta(x)=\delta(y)$.  Therefore, we may take
$\delta=\delta(x)=\delta(y)$.  The Lemma follows
by substituting $\delta$ for $\delta(y)$ in the conditions defining
$\delta(y)$.
\end{myproof}

Since a quadruple satisfying \eqref{eqn:quadruplecond} has Property
(iv), Sublemma \ref{sublem:sqrational} says that there
exists for each permutation $(x,\,y,\,z,\,w)$
of \eqref{eqn:quadruple} a $\delta(x)\in\{0,1\}$,
such that $\vecti^{\delta(x)}
x\in\Rational(\vecti)$.
Sublemma \ref{sublem:sqrationalindep} is significant
addition because it says that $\delta(x)$ is actually independent
of the permutation.  Sublemma \ref{sublem:sqrationalindep}
implies that for every quadruple
\eqref{eqn:quadruple} satisfying \eqref{eqn:quadruplecond},
we may define
\starteqn\label{eqn:indepdeltadefn}
\delta=\delta\big((a,\,b,\,c,\,d)\big)\in\{0,1\}\;\text{such that}
\; \delta(a,\,b,\,c,\,d)\in(\Rational(\vecti))^4.
\finisheqn

\vspace*{0.3cm}\noindent\boldmath\textbf{The square map $\sq$
and the group of invertible residues modulo $(1+\vecti)^n$}.
\ubmth
Before proceeding with the  analysis of the set of quadruples
satisfying \eqref{eqn:quadruplecond},
we must establish certain facts pertaining to the square
map applied to the group of invertible residues
modulo $(1+\vecti)^n$.  Note first that for any $\nu$
prime and $\mu$ in $\Int[\vecti]$ such that $\nu|\mu$,
the function $\ord_{\nu}$ is well defined
on $\Int[\vecti]/(\mu)$.  In particular, then, we may define
\starteqn\label{eqn:Ugroupdefn}
U_{(1+\vecti)^n}:=\{x\in\Int[\vecti]/\left((1+\vecti)^n\right)\;|\;
\ord_{1+\vecti}(x)=0\}.
\finisheqn
The set $U_{(1+\vecti)^n}$ is closed under multiplication,
and $U_{(1+\vecti)^n}$ is precisely the group of invertible elements of $\Int[\vecti]\left/
(1+\vecti)^n\right.$.
The endomorphism $\mathrm{sq}$ of $\Int[\vecti]\left/((1+\vecti)^n\right)$
defined by
\startdisp
\mathrm{sq}(x)=x^2
\finishdisp
restricts to an endomorphism of $U_{(1+\vecti)^n}$.
In order to describe the image of the restriction
$\mathrm{sq}|_{U_{(1+\vecti)^n}}$, we make the following
observation.  Note that, since by Lemma \ref{lem:standardresiduerep},
the residue classes in
$\Int[\vecti]\left/((1+\vecti)^n\right)$ are in bijection with the elements
of $\Omega_{(1+\vecti)^n}$ we may identify $U_{(1+\vecti)^n}$
with its image in $\Omega_{(1+\vecti)^n}$.  We denote
this image by $\tilde{U}_{(1+\vecti)^n}$, and identify
$\tilde{U}_{(1+\vecti)^n}$ with
${U}_{(1+\vecti)^n}$.  The identification
of $\tilde{U}_{(1+\vecti)^n}$ with ${U}_{(1+\vecti)^n}$ makes
$\tilde{U}_{(1+\vecti)^n}$
into a group.
We note, for future use that, as a simple computation shows,
\starteqn\label{eqn:tildeUsize}
\#\tilde{U}_{(1+\vecti)^{n}}=\#{U}_{(1+\vecti)^{n}}=2^{n-1}.
\finisheqn
Likewise, the identification induces an endomorphism
$\mathrm{sq}|_{\tilde{U}_{(1+\vecti)^n}}$
of the group $\tilde{U}_{(1+\vecti)^n}$.
We usually refer to the induced endomorphism as `$\mathrm{sq}$', and as
`$\mathrm{sq}|_{\tilde{U}_{(1+\vecti)^n}}$\hspace{-0.07cm}'\hspace{0.07cm}
only when there is any danger of confusion.
\begin{defn}  We refer to the image $\mathrm{sq}(U_{(1+\vecti)^n})$
as the set of \textbf{quadratic residues in \bmth $U_{(1+\vecti)^n}$
\ubmth}.  Likewise, we refer to the image
$\mathrm{sq}(\tilde{U}_{(1+\vecti)^n})$
as the set of \textbf{quadratic residues in \bmth $
\tilde{U}_{(1+\vecti)^n}$
\ubmth}.
\end{defn}
\begin{sublem} \label{sublem:imptsqdescr}
Let
\starteqn\label{eqn:imptsqdescrstatement}
x=r+s\vecti\in\tilde{U}_{(1+\vecti)^n},\;\text{where}\; r,s\in\Int,
\finisheqn
be a quadratic residue.  Then $s\in 2\Int$.
\end{sublem}
\begin{myproof}{Proof}  The sublemma results from a straightforward
calculation.  If $n=1$, the result is trivial,
so we may assume that $n>1$.  The assumption that $x$ is a quadratic
residue in $\tilde{U}_{(1+\vecti)^n}$ implies that
\starteqn\label{eqn:imptsqdescrpf}
x\equiv y^2\mod (1+\vecti)^n,\quad\text{for some $y\in\tilde{U}_{(1+\vecti)^n}$}
\finisheqn
Supposing that
\startdisp
y=p+q\vecti\;\text{with}\; p,q\in\Int[\vecti],
\finishdisp
we have by \eqref{eqn:imptsqdescrstatement} and \eqref{eqn:imptsqdescrpf}
that
\startdisp
r+s\vecti\equiv p^2-q^2+2pq\vecti\mod(1+\vecti)^n.
\finishdisp
Thus $s$ differs from $2pq$ by the imaginary part of some
multiple $(1+\vecti)^n$.  Both $2pq$ and the imaginary part
of any multiple of $(1+\vecti)^n$ belong to $2\Int$,
so $s$ belongs to $2\Int$.
\end{myproof}

For Sublemma \ref{sublem:quadresidues} and the subsequent
discussion, we note that for $x,\,y\in\Int[\vecti]$ such that
$x|y$, the reduction map
$\red_{x}:\Int[\vecti]\rightarrow\Int[\vecti]/(x)$ naturally
induces a map,
\startdisp
\red_x:\Int[\vecti]/(y)\rightarrow\Int[\vecti]/(x).
\finishdisp
Note that the map $\red_x$ maps the set invertible elements
of $\Int[\vecti]/(y)$ onto the set of invertible
elements of $\Int[\vecti]/(x)$, and the restriction of $\red_x$
to the invertible elements is an epimorphism.
Applying this to the situation at hand, we see that for $n>1$,
\startdisp
\red_{(1+\vecti)^{n-1}}(U_{(1+\vecti)^n})=U_{(1+\vecti)^{n-1}}.
\finishdisp
Via the identification of $U_{(1+\vecti)^n}$
with $\tilde{U}_{(1+\vecti)^n} \subset\Omega_{(1+\vecti)^n}$,
for each $n>1$, we obtain a naturally induced epimorphism
\startdisp
\red_{(1+\vecti)^{n-1}}: \tilde{U}_{(1+\vecti)^n}\rightarrow
\tilde{U}_{(1+\vecti)^{n-1}}
\finishdisp
Since $\mathrm{sq}$ commutes with $\red_{(1+\vecti)^{n-1}}$
the reduction morphism restricts to a morphism of
$\ker\left(\mathrm{sq}_{(1+\vecti)^n}\right)$
into $\ker\left(\mathrm{sq}_{(1+\vecti)^{(n-1)}}\right)$.
For our analysis, we need more precise information
on the image of this morphism for certain low values
of $n$, which is subject of Sublemmas \ref{sublem:quadresidues}
and \ref{sublem:quadraticresn234}.

\begin{sublem}  \label{sublem:quadresidues} Let $n$ be
an integer greater than $1$.  Let
$\tilde{U}_{(1+\vecti)^n}$, $\mathrm{sq}|_{\tilde{U}_{(1+\vecti)^n}}$,
and\linebreak $\red_{(1+\vecti)^{n-1}}$ be as above.
Set
\startdisp
m=\begin{cases}
n-2&\text{if $n>3$}\\
n-1&\text{if $n=1$ or $n=2$}
\end{cases}
\finishdisp
\begin{itemize}
\item[(a)]
We have
\starteqn\label{eqn:nlargekernelsqmap}
\begin{aligned}
&&&\ker\left(\mathrm{sq}|_{\tilde{U}_{(1+\vecti)^n}}\right)=\\&&&
\hspace*{1cm}
\left\{x\in\Omega_{(1+\vecti)^n}\;|\;
x\equiv \ell(1+\vecti)^{m} \pm 1\mod(1+\vecti)^n,\;
\ell\in\Omega_{(1+\vecti)^{n-m}}\right\}.
\end{aligned}
\finisheqn
\item[(b)]  We have
\starteqn\label{eqn:nlargekernelsqmapred}
\begin{aligned}
&&&\res_{(1+\vecti)^{n-1}}\left(
\ker\left(\mathrm{sq}|_{\tilde{U}_{(1+\vecti)^n}}\right)\right)=\\&&&
\hspace*{1cm}
\left\{x\in\Omega_{(1+\vecti)^{n-1}}\;|\;
x\equiv \ell(1+\vecti)^{m} \pm 1\mod(1+\vecti)^{n-1},\;
\ell\in\Omega_{(1+\vecti)^{n-m-1}}\right\}.
\end{aligned}
\finisheqn
\end{itemize}
\end{sublem}
\begin{myproof}{Proof}
For (a), note that $z\in
\ker\left(\mathrm{sq}|_{\tilde{U}_{(1+\vecti)^n}}\right)$
if and only if
\startdisp
z^2-1=(z+1)(z-1)\equiv 0\mod (1+\vecti)^n.
\finishdisp
By taking $\ord_{(1+\vecti)^n}$ of both sides and applying
\textbf{ORD 1}, we obtain
\startdisp
z\in\ker\left(\mathrm{sq}|_{\tilde{U}_{(1+\vecti)^n}}\right)\;
\text{if and only if}\;
\ord_{1+\vecti}(z-1)
+\ord_{1+\vecti}(z+1)\geq n.
\finishdisp
Trivially, then
\starteqn\label{eqn:quadresidue1}
z\in\ker\left(\mathrm{sq}|_{\tilde{U}_{(1+\vecti)^n}}\right)\;
\text{if and only if}\;
(\min+\max)(\ord_{1+\vecti}(z-1)
,\ord_{1+\vecti}(z+1))\geq n.
\finisheqn
We now apply Property \textbf{ORD 5} with $x=z-1$,
$\nu=(1+\vecti)$, $n=2$ and $u=-\vecti$, to deduce that
\starteqn\label{eqn:quadresidueclaim}
\min\left(\ord_{1+\vecti}(z-1),\ord_{1+\vecti}(z+1)\right)\leq 2.
\finisheqn
Applying \eqref{eqn:quadresidue1} and \eqref{eqn:quadresidueclaim}, we
obtain
\starteqn\label{eqn:quadresiduengeq4}
z\in\ker\left(\mathrm{sq}|_{\tilde{U}_{(1+\vecti)^n}}\right)\;\text{
implies}\;
\max\left(\ord_{1+\vecti}(z-1),\ord_{1+\vecti}(z+1)\right)\geq n-2.
\finisheqn
We next claim that
\starteqn\label{eqn:quadresidueleq3}
\text{Under the assumption that $n<4$, we have strict inequality
in \eqref{eqn:quadresiduengeq4}.}
\finisheqn
In order to prove \eqref{eqn:quadresidueleq3}, assume
that
\startdisp
n<4\;\text{and}\;\max\left(\ord_{1+\vecti}(z-1),\ord_{1+\vecti}(z+1)\right)=n-2.
\finishdisp
The condition that $n<4$ is equivalent to the condition $n-2<2$.
Therefore, we may apply Property \textbf{ORD 3} with
$x=\pm 2$, $y=z\mp 1$, to obtain
\startdisp
\min\left(\ord_{1+\vecti}(z-1),\ord_{1+\vecti}(z+1)\right)=n-2.
\finishdisp
Therefore,
\startdisp
(\max+\min)\left(\ord_{1+\vecti}(z-1),\ord_{1+\vecti}(z+1)\right)=
2n-4<n,
\finishdisp
where the inequality on the right results from the assumption $n<4$.
By \eqref{eqn:quadresidue1}, we deduce that
$z\notin\ker\left(\mathrm{sq}|_{\tilde{U}_{(1+\vecti)^n}}\right)$
This completes the proof of \eqref{eqn:quadresidueleq3}.

Let $m$ be as in the statement of the sublemma.  Then,
combining \eqref{eqn:quadresidueleq3} with \eqref{eqn:quadresiduengeq4},
we obtain
\starteqn\label{eqn:quadresiduewithm}
z\in\ker\left(\mathrm{sq}|_{\tilde{U}_{(1+\vecti)^n}}\right)\;\text{
implies}\;
\max\left(\ord_{1+\vecti}(z-1),\ord_{1+\vecti}(z+1)\right)\geq m.
\finisheqn
In order to prove the converse of \eqref{eqn:quadresiduewithm}, assume
that $\max\left(\ord_{1+\vecti}(z-1),\ord_{1+\vecti}(z+1)\right)\geq m$.
Note that $m\geq 2$, \textit{i.e.}
\startdisp
\min\left(\max\left(\ord_{1+\vecti}(z-1),\ord_{1+\vecti}(z+1)\right),
2\right)=2.
\finishdisp
Therefore, applying \textbf{ORD 2} with $x=z\pm 1$ and $y=\mp 2$, we
obtain
\startdisp
\min\left(\ord_{1+\vecti}(z-1),\ord_{1+\vecti}(z+1)\right)\geq 2.
\finishdisp
We calculate that
\startdisp
(\min+\max)\left(\ord_{1+\vecti}(z-1),\ord_{1+\vecti}(z+1)\right)\geq m+2\geq n.
\finishdisp
By \eqref{eqn:quadresidue1}, we have
$z\in\ker\left(\mathrm{sq}|_{\tilde{U}_{(1+\vecti)^n}}\right)$,
and this completes the proof of the converse of \eqref{eqn:quadresiduewithm}.
So we have
\starteqn\label{eqn:quadresiduewithmiff}
z\in\ker\left(\mathrm{sq}|_{\tilde{U}_{(1+\vecti)^n}}\right)\;\text{
if and only if}\;
\max\left(\ord_{1+\vecti}(z-1),\ord_{1+\vecti}(z+1)\right)\geq m.
\finisheqn

From \eqref{eqn:quadresiduewithmiff} and the definition
of $\ord_{(1+\vecti)}$ we deduce that
\startdisp
\ker\left(\mathrm{sq}|_{\tilde{U}_{(1+\vecti)^n}}\right)=
\left\{x\in\Omega_{(1+\vecti)^n}\;|\;
x\equiv y(1+\vecti)^{m} \pm 1\mod (1+\vecti)^n,\;
y\in\Int[\vecti]\right\}.
\finishdisp
Using Lemma \ref{lem:standardresiduerep}, it is easy to verify
that for $y\in\Int$, there is a unique $\ell\in\Omega_{(1+\vecti)^{n-m}}$
such that
\startdisp
y(1+\vecti)^{m}\equiv\ell(1+\vecti)^m\mod(1+\vecti)^n.
\finishdisp
Therefore, we have deduced the description of
$\ker\left(\mathrm{sq}|_{\tilde{U}_{(1+\vecti)^n}}\right)$
given in \eqref{eqn:nlargekernelsqmap}.

Part (b) is obtained by applying the map $\red_{(1+\vecti)^{n-1}}$
to the elements on the right-hand side of \eqref{eqn:nlargekernelsqmap}.
Using Lemma \ref{lem:standardresiduerep},
it is easily verified that for each $y\in\Omega_{(1+\vecti)^{m-n}}$
there is a unique $\ell\in\Omega_{(1+\vecti)^{m-n-1}}$
such that
\startdisp
y(1+\vecti)^{m}\equiv\ell(1+\vecti)^m\mod(1+\vecti)^{n-1}.
\finishdisp
We thereby obtain the description of
$\red_{(1+\vecti)^{n-1}}\ker\left
(\mathrm{sq}|_{\tilde{U}_{(1+\vecti)^n}}\right)$ given in
\eqref{eqn:nlargekernelsqmapred}.
\end{myproof}

Given the general description of the quadratic residues
in Sublemma \ref{sublem:quadresidues}, we
now explicitly list, in Sublemma \ref{sublem:quadraticresn234},
the quadratic residues in $\tilde{U}_{(1+\vecti)^n}$
for the first few integer values of $n$.  This
is the information we actually need to carry out our analysis.

\begin{sublem}\label{sublem:quadraticresn234}  For small
values of $n$ in Sublemma \ref{sublem:quadresidues},
we calculate that
\begin{itemize}
\item[(a)]
$\ker\left(\mathrm{sq}|_{\tilde{U}_{(1+\vecti)^2}}\right)=
\{1,\vecti\}$,
$\res_{1+\vecti}\left(
\ker\left(\mathrm{sq}|_{\tilde{U}_{(1+\vecti)^2}}\right)\right)=\{1\}$.
\item[(b)]
$\ker\left(\mathrm{sq}|_{\tilde{U}_{(1+\vecti)^3}}\right)=
\{1,3\}$,
$\res_{(1+\vecti)^{2}}\left(
\ker\left(\mathrm{sq}|_{\tilde{U}_{(1+\vecti)^3}}\right)\right)=\{1\}$.
\item[(c)]
$\ker\left(\mathrm{sq}|_{\tilde{U}_{(1+\vecti)^4}}\right)=
\{1,3,1+2\vecti, 3+2\vecti\}$,
$\res_{(1+\vecti)^{3}}\left(
\ker\left(\mathrm{sq}|_{\tilde{U}_{(1+\vecti)^4}}\right)\right)=\{1,3\}$.
\end{itemize}
\end{sublem}
\begin{myproof}{Proof}
Parts (a)--(c) are verified by substituting the values $n=2,3,4$
into \eqref{eqn:nlargekernelsqmap} and \eqref{eqn:nlargekernelsqmapred}.
To complete the calculations, one uses the definition
of $\Omega_{(1+\vecti)^m}$ for $m\geq 1$
given in \eqref{eqn:omegamdefn}.  The details are routine, so
we omit them.
\end{myproof}

Sublemma \ref{sublem:quadraticresn234} gives us the information
appearing in the four leftmost columns of Table \ref{table:sqresiduesinfo}.
In the paragraphs immediately following,
we will explain the meaning of the two rightmost columns.

\starteqn\label{table:sqresiduesinfo}
\begin{array}{|c|p{1.8cm}|p{1.75cm}|p{1.8cm}|p{1.8cm}|p{1.85cm}|p{1.55cm}|}\hline
{\scriptstyle n} & $\hspace*{-0.1cm}{\scriptstyle \ker\left(\mathrm{sq}|_{\tilde{U}_{(1+\vecti)^n}}
\right)}$
 & $\hspace*{-0.15cm}\begin{array}{c}{\scriptstyle\#}\\
 \hspace*{-0.15cm}{\scriptstyle
\ker\left(\mathrm{sq}|_{\tilde{U}_{(1+\vecti)^n}}\right)} \end{array}$&
\hspace*{-0.3cm}$\begin{array}{l}
{\scriptstyle \red_{(1+\vecti)^{n-1}}}\\ {\scriptstyle\ker\left(\mathrm{sq}|_
{\tilde{U}_{(1+\vecti)^n}}\right)}
\end{array}$&
\hspace*{-0.3cm}$\begin{array}{l}
{\scriptstyle \#\red_{(1+\vecti)^{n-1}}}\\ {\scriptstyle\ker\left(\mathrm{sq}|_
{\tilde{U}_{(1+\vecti)^n}}\right)}
\end{array}$&\hspace*{-0.15cm}${\scriptstyle
\mathrm{rt}_n\left(\sq\left(\tilde
{U}_{(1+\vecti)^n}\right)
\right)}$ &\hspace*{-.3cm}$
\begin{array}{l}
{\scriptstyle\red_{(1+\vecti)^{n-1}}\mathrm{rt}_n}\\
{\scriptstyle \left(\sq\left(\tilde
{U}_{(1+\vecti)^n}\right)\right)}
\end{array}$
\\ \hline \hline
2& $1,\,\vecti$ & $2$& $1$& 1& 1&1\\ \hline
3& $1,\,3$ & $2$& $1$& 1& 1,\,$\vecti$ &$1,\,\vecti$ \\ \hline
4&
\hspace*{-.2cm}$\begin{array}{p{1cm}}
$1,\, 3,$\\
$1+2\vecti,$\\ $3+2\vecti$
\end{array}$& $4$& $1,3$ & $2$&1,$\vecti$ &$1,\vecti$\\ \hline
\end{array}
\finisheqn

\vspace*{0.3cm}\noindent\boldmath\textbf{Sections of the square map}.
\ubmth
Let $\varphi$ be an epimorphism of a group $G$ onto a group $G'$.
Let $\psi: G'\rightarrow G$ be a right-inverse to $\varphi$, \textit{i.e.}
a injective set map of $G'$ to $G$ satisfying
\startdisp
\varphi\psi=I_{G'}.
\finishdisp
Then $\psi(G')$ is a subset of $G$
containing exactly one element of each fiber of $\varphi$.
We will refer to $\psi(G')$ as a \textbf{section
of \bmth$\varphi$ in $G$\ubmth}.
Note that since $\psi$ is injective and $\phi$ surjective,
\startdisp
\#\psi\varphi(G)=\#\psi(G')=\# G'.
\finishdisp
Therefore, we have a map of sets
\startdisp
\psi\varphi: G\rightarrow G,\quad\text{with image of size}\;\# G'.
\finishdisp
By a standard isomorphism of elementary group theory,
\startdisp
G'=G/\ker(\varphi),\;\text{so that}\;\#G'=\frac{\#G}{\#\ker(\varphi)}.
\finishdisp
We conclude that
\starteqn\label{eqn:sizeofvarphig}
\psi\varphi(G)\;\text{is a subset of $G$ of size}\;\frac{\#G}{\#\ker(\varphi)},
\finisheqn
and also
\starteqn\label{eqn:gpsiphiinvinker}
g\left(\psi\varphi(g)\right)\inv\in\ker(\phi).
\finisheqn

In the applications, we will be taking $G=\tilde{U}_{(1+\vecti)^n}$,
$\varphi=\mathrm{sq}|_{\tilde{U}_{(1+\vecti)^n}}$,
$G'=\mathrm{sq}\left(\tilde{U}_{(1+\vecti)^n}\right)$.
\begin{defn}  A fixed right-inverse of $\mathrm{sq}|_{\tilde{U}_{(1+\vecti)^n}}$
will be denoted
\bmth$\mathrm{rt}_n$\ubmth.  Thus,
\startdisp
\mathrm{rt}_n: \mathrm{sq}\left(\tilde{U}_{(1+\vecti)^n}\right)
\rightarrow \tilde{U}_{(1+\vecti)^n},
\finishdisp
such that
\startdisp
\sq\circ\rt_n=\Id_{\tilde{U}_{(1+\vecti)^n}},
\finishdisp
and $\rt_n\left(\sq\left(\tilde{U}_{(1+\vecti)^n}\right)\right)$
is a section of the $\sq$ map in $\tilde{U}_{(1+\vecti)^n}$.
\end{defn}
By specifying $\mathrm{rt}_n$, we are in effect fixing a branch of the
square root function on $\sq\left(\tilde{U}_{(1+\vecti)^n}\right)$,
hence the notation ``rt'' for \textit{root}.  Because
of the identification of $U_{(1+\vecti)^n}$ with
$\tilde{U}_{(1+\vecti)^n}$, $\mathrm{rt}_n$ may equally
well be thought of as a map from $\sq\left({U}_{(1+\vecti)^n}\right)$,
namely, a fixed branch of the square-root function on
$\sq\left({U}_{(1+\vecti)^n}\right)$.

Applying \eqref{eqn:sizeofvarphig} to the situation at hand, and
using \eqref{eqn:tildeUsize}, we see that
\starteqn\label{eqn:sizertsq}
\rt_n\sq(\tilde{U}_{(1+\vecti)^n})\;\text{is a subset of $\tilde{U}_{
(1+\vecti)^n}$ of size $\frac{2^{n-1}}{\#\ker\left(
\mathrm{sq}|_{\tilde{U}_{(1+\vecti)^n}}\right)}$.}
\finisheqn
Applying \eqref{eqn:gpsiphiinvinker}, we have
\starteqn\label{eqn:xrtxsqinvinker}
x\left(\rt_n\sq(x)\right)\inv\in
\ker\left(sq|_{\tilde{U}_{(1+\vecti)^n}}\right),\forallindisp\;
x\in \tilde{U}_{(1+\vecti)^n}
\finisheqn

\begin{sublem}\label{sublem:sqmapsection}  Let $n$ be a positive
integer, and let
$\tilde{U}_{(1+\vecti)^n}$,
$\mathrm{sq}$, and $\mathrm{rt}_n$ be as defined above.
Then (a)--(c) below give the
the sizes of $\rt_n\sq(\tilde{U}_{(1+\vecti)^n})$, following
\eqref{eqn:sizertsq}, as well as the defining equations for
one valid choice of section $\rt_n$ of
$\sq|_{\tilde{U}_{(1+\vecti)^n}}$.
\begin{itemize}
\item[(a)] $\#\rt_2\sq\left(\tilde{U}_{(1+\vecti)^2}\right)=1$;\quad $\rt_2(1)=1$.
\item[(b)] $\#\rt_3\sq\left(\tilde{U}_{(1+\vecti)^3}\right)=2$;\quad $\rt_3(1)=1$,
$\rt_3(3)=\vecti$
\item[(c)] $\#\rt_4\sq\left(\tilde{U}_{(1+\vecti)^4}\right)=2$;\quad $\rt_4(1)=1$,
$\rt_4(3)=\vecti$.
\end{itemize}
\end{sublem}
\begin{myproof}{Proof}  The sizes of $\rt_n\sq\left(\tilde{U}_{(1+\vecti)^n}
\right)$ given in (a)--(c) are calculated
using \eqref{eqn:tildeUsize} and \eqref{eqn:sizertsq}.

We now indicate why the equations in (a)--(c) define a valid
choice of section $\mathrm{\rt}_n$ for $\sq|_{\tilde{U}_{(1+\vecti)^n}}$,
$n=2,3,4$.  Since
$\#\rt_2\sq\left(\tilde{U}_{(1+\vecti)^2}\right)=1$,
any map from $\{1\}$ to $\tilde{U}_{(1+\vecti)^2}$
defines a section of $\sq|_{\tilde{U}_{(1+\vecti)^n}}$.  When
$n=3,4$, $\#\rt_2\sq\left(\tilde{U}_{(1+\vecti)^2}\right)=2$.
Therefore, if $(x,y)$ is any pair of elements of $\tilde{U}_{(1+\vecti)^n}$
such that
\startdisp
x\in\ker\left(\sq|_{\tilde{U}_{(1+\vecti)^n}}\right),\;
\text{and}\;y\notin\ker\left(\sq|_{\tilde{U}_{(1+\vecti)^n}}\right)
\finishdisp
then $\#\{\sq(x^2),\sq(y^2)\}=2$. Therefore, we have
\starteqn\label{eqn:sectioncard2cond}
\{\sq(x^2),\sq(x^2)\}=\sq(U_{(1+\vecti)^n}).
\finisheqn
From \eqref{eqn:sectioncard2cond} it is easily verified that
the following formulae define a valid section of the square map.
\starteqn\label{eqn:sectioncard2cond2}
\mathrm{rt}(1)=x;\;\mathrm{rt}(y^2)=y.
\finisheqn
For $n=3, 4$, it is easily verified that
$\vecti\in\tilde{U}_{(1+\vecti)^n}-
\ker(\sq_{\tilde{U}_{(1+\vecti)^n}})$.  Substituting
$(x,y)=(1,\vecti)$ into \eqref{eqn:sectioncard2cond2}, we obtain
the formulae defining $rt_n$ for $n=3,4$ in Parts (b) and (c).
\end{myproof}

Henceforth, it will be assumed that $rt_{(1+\vecti)^n}$ is defined
according to the formulas of Sublemma \eqref{sublem:sqmapsection}.
Furthermore, Sublemma \eqref{sublem:sqmapsection} provides
the information in the rightmost
two columns of Table \eqref{table:sqresiduesinfo}.

\vspace*{0.3cm}\noindent\boldmath\textbf{Properties concerned with the residues
of the entries of an integral quadruple}.
\ubmth
For the following list of \textbf{R} properties, let
\starteqn\label{eqn:quadrupleres}
(a',b',c',d')\in (\Int[\vecti])^4.
\finisheqn
As a matter of convenience, we introduce the piece of
notation, for $x,y,q\in\Go$,
\startdisp
x\cong y\mod q\;\text{means}\; q|x-uy,\;\text{for some $u\in\Go^*$}.
\finishdisp

The properties
\textbf{R}$\mathbf{i}(n)$, $i\in\{0,1,2\}$,
applying to \eqref{eqn:quadrupleres}, are denoted $R$ for
``residue".  Throughout the definition of the
properties, $(x',\,y',\,z',\,w')$ will be used to denote an
arbitrary permutation of the quadruple \eqref{eqn:quadrupleres},
in a way parallel to that in which $(x,\,y,\,z,\,w)$ denoted
an arbitrary permutation of the quadruple \eqref{eqn:quadruple}.
The \textbf{R}$\mathbf{i}(n)$ properties are indexed
by a variable $n$ that takes on integer values greater than
or equal to $2$.
\begin{itemize}
\item[\textbf{R}$\mathbf{0}(n)$]  $x'\equiv 1\mod (1+\vecti)$
\item[\textbf{R}$\mathbf{1}(n)$]
$a'd'-b'c'\cong (1+\vecti)^{n-2}\mod (1+\vecti)^{n}.$
\item[\textbf{R}$\mathbf{2}(n)$] ${x'}^2+{y'}^2\equiv{z'}^2+{w'}^2\mod
(1+\vecti)^n$.
\end{itemize}
Of the three properties, \textbf{R}$\mathbf{0}(n)$ is exceptional in that it
that actually does not depend on the parameter $n$,
since it just says that no element of \eqref{eqn:quadrupleres}is divisible
by $(1+\vecti)$.  We use the parameter $n$ in referring
to \textbf{R}$\mathbf{0}(n)$ as a matter of convenience, in order to
make it easier to refer to the \textbf{R}$\mathbf{i}(n)$ properties
collectively.

We will sometimes have occasion to refer to the following weaker version
of \textbf{R}$\mathbf{1}(n)$.
\begin{itemize}
\item[\textbf{R}$\mathbf{1'}(n)$]
$a'd'-b'c'\cong (1+\vecti)^{n-2}\mod (1+\vecti)^{n-1}.$
\end{itemize}
It is easy to see that \textbf{R}$\mathbf{1}(n)$ implies
\textbf{R}$\mathbf{1'}(n)$, by reducing \textbf{R}$\mathbf{1}(n)$
modulo $(1+\vecti)^{n-1}$.
\begin{sublem} \label{sublem:primedquadruples}  Let $(a',\,
b',\,c',\,d')\in(\Int[\vecti])^4$ be a quadruple as in \eqref{eqn:quadrupleres}.
Let the properties \textbf{R}$\mathbf{i}(n)$, $i\in\{0,1,2\}$
be defined as above.  Then we have
\begin{itemize}
\item[(a)]  Suppose that we have
\starteqn\label{eqn:sublemprimedquadrupleshypa}
{x'}^2\equiv {y'}^2\mod(1+\vecti)^n,
\finisheqn
for \textit{some} permutation $(x',\,y',\,z',\,w')$ of \eqref{eqn:quadrupleres}.
Then \eqref{eqn:quadrupleres} satisfies property \textbf{R}$\mathbf{2}(n)$
if and only if we also have
\startdisp
{z'}^2\equiv{w'}^2\mod(1+\vecti)^n.
\finishdisp
\item[(b)]  If
\starteqn\label{eqn:sublemprimedquadrupleshypb}
\#\rt_n\left(\sq(U_{(1+\vecti)^n})\right)\leq 2,
\;\text{\textbf{R}$\mathbf{2}(n)$ is satisfied},
\finisheqn
then for \textit{some} permutation
$(x',\,y',\,z',\,w')$ of \eqref{eqn:quadrupleres}, we have
\starteqn\label{eqn:sublemprimedquadrupleb}
x'^2\equiv y'^2\mod(1+\vecti)^n,\quad\text{and}\quad z'^2\equiv w'^2
\mod(1+\vecti)^n.
\finisheqn
\end{itemize}
\end{sublem}

\begin{myproof}{Proof}  For (a), using the condition
\eqref{eqn:sublemprimedquadrupleshypa}, we have
\startdisp
x'^2+z'^2\equiv y'^2+w'^2\mod(1+\vecti)^n\;\;\text{if and only if}\;
z'^2\equiv w'^2\mod(1+\vecti)^n.
\finishdisp

For (b), note that \eqref{eqn:sublemprimedquadrupleshypb}
implies that, for some permutation of \eqref{eqn:quadrupleres},
we have
\startdisp
x'^2\equiv y'^2\mod(1+\vecti)^n.
\finishdisp
Therefore, part (a) applies.  Since we are assuming
\textbf{R}$\mathbf{2}(n)$ is satisfied, part (a) implies that
\startdisp
z'^2\equiv w'^2\mod(1+\vecti)^n.
\finishdisp
This completes the proof of \eqref{eqn:sublemprimedquadrupleb}.
\end{myproof}

\begin{sublem}\label{sublem:Rcases}\begin{itemize}
Let $(a',\,
b',\,c',\,d')\in(\Int[\vecti])^4$ be a quadruple as in \eqref{eqn:quadrupleres}.
Let the properties \textbf{R}$\mathbf{i}(n)$, $i\in\{0,1,2\}$
be defined as above.
\item[(a)]
Properties \textbf{R}$\mathbf{1}(2)$ and \textbf{R}$\mathbf{2}(2)$ are
simultaneously satisfied
if and only if
\starteqn\label{eqn:primedquadrupleokn2}
(a',\,b',\,c',\,d')\equiv (1,0,0,1),\;\text{or}\;(0,1,1,0)\mod(1+\vecti).
\finisheqn
\item[(b)]  All three properties
\textbf{R}$\mathbf{0}(3)$ and \textbf{R}$\mathbf{1}(3)$
and \textbf{R}$\mathbf{2}(3)$ \textit{cannot} be satisfied simultaneously.
\end{itemize}
\end{sublem}
\begin{myproof}{Proof}
For part (a), begin by substituting $n=2$ into the definition
of Property \textbf{R}$\mathbf{1}(n)$.
\begin{itemize}
\item[\textbf{R}$\mathbf{1}(2)$] $a'd'-b'c'\cong 1\mod (1+\vecti)^2$,
or equivalently,
\starteqn\label{eqn:R1casen2}
a'd'-b'c'\equiv 1\mod (1+\vecti).
\finisheqn
\item[\textbf{R}$\mathbf{2}(2)$] ${x'}^2+{y'}^2\equiv
{z'}^2+{w'}^2\mod (1+\vecti)^2$.
\end{itemize}
As tabulated in \eqref{table:sqresiduesinfo},
$\mathrm{sq}\left(\tilde{U}_{(1+\vecti)^2}\right)$ consists of $1$ element,
namely $1$.  Note that
\startdisp
\Int[\vecti]/(1+\vecti)^2= U_{(1+\vecti)^2}\cup\{\red_{(1+\vecti)^2}(0),
\red_{(1+\vecti)^2}(1)\}.
\finishdisp
We therefore have
\starteqn\label{eqn:Rcasessublem1}
\red_{(1+\vecti)^2}\{x'^2,y'^2,z'^2,w'^2\}
\subseteq\sq\left(\tilde{U}_{(1+\vecti)^2}\cup \{0,1+\vecti\}\right)=
\{1,0\}.
\finisheqn
By \eqref{eqn:Rcasessublem1}, there is a permutation
of \eqref{eqn:quadrupleres} such that $x'^2=y'^2$.
Using Sublemma \ref{sublem:primedquadruples},
we have that
\starteqn\label{eqn:Rcasessublem2}
\text{\textbf{R}$\mathbf{2}(2)$
is satisfied if and only if $z'^2\equiv w'^2\mod (1+\vecti)^n$.}
\finisheqn
Fix a permutation such that
such that
\starteqn\label{eqn:equalityofsqprimes}
x'^2=y'^2\;\text{and}\;z'^2=w'^2.
\finisheqn  From Table \eqref{table:sqresiduesinfo},
we see that for any entry $v'$ of \eqref{eqn:quadrupleres}, we have
\startdisp
\red_{(1+\vecti)^2}(v')=\begin{cases} 1&\text{if and only if}\;
q'\equiv1\mod (1+\vecti).\\
 0&\text{if and only if}\; q'\equiv 0\mod (1+\vecti).
\end{cases}
\finishdisp
Therefore, for any pair of entries $(u',v')$ of \eqref{eqn:quadrupleres},
\starteqn\label{eqn:squaresequiviffelemsequiv}
\red_{(1+\vecti)^2}(u'^2)=\red_{(1+\vecti)^2}(v'^2)\;\text{if and only if}\;
\red_{(1+\vecti)^2}(u')=\red_{(1+\vecti)^2}(v').
\finisheqn
Further,
\starteqn\label{eqn:conditionforprodbeing1}
u'v'\equiv 1\mod (1+\vecti)\;\text{if and only if}\; (u',v')\equiv (1,1)
\mod (1+\vecti).
\finisheqn
Applying \eqref{eqn:squaresequiviffelemsequiv} to the pairs
$(u',v')=(x',y')$, $(u',v')=(z',w')$,
we deduce from \eqref{eqn:equalityofsqprimes} that
\starteqn\label{eqn:equalpairs}
x'\equiv y'\mod (1+\vecti)^2\;\text{and}\; z'\equiv w'\mod (1+\vecti)^2.
\finisheqn
Applying \eqref{eqn:conditionforprodbeing1}, we deduce that
\eqref{eqn:R1casen2}, and hence property
\textbf{R}$\mathbf{1}(2)$ is equivalent to
\starteqn\label{eqn:casen2alternative}
\text{Exactly one of $(a',\,d'),\; (b',\,c')\equiv (1,1)\mod(1+\vecti)$}.
\finisheqn
By \eqref{eqn:casen2alternative}, there can be a permutation of
$(a',\,b',\,c',\,d')$
satisfying \eqref{eqn:equalpairs} only if $(a',\,b',\,c',\,d')$
satisfies the condition of \eqref{eqn:primedquadrupleokn2}.
Therefore, we see that Properties \textbf{R}$\mathbf{1}(2)$
and \textbf{R}$\mathbf{2}(2)$ imply that
\eqref{eqn:quadrupleres} satisfies the condition
of \eqref{eqn:primedquadrupleokn2}.

The converse, namely, that if \eqref{eqn:quadrupleres} satisfies
the condition of \eqref{eqn:primedquadrupleokn2},
then it satisfies Properties \textbf{R}$\mathbf{1}(2)$
and \textbf{R}$\mathbf{2}(2)$, is verified by direct calculation,
using the specific form of \textbf{R}$\mathbf{1}(2)$
and \textbf{R}$\mathbf{2}(2)$ given at the beginning of the proof.
This completes the proof of part (a).

In order to prove (b), we first verify a claim.
As usual, we let $\rt_n$ be fixed
by the equations of Table \eqref{table:sqresiduesinfo}.
Then
\starteqn\label{eqn:Rcasessublemclaimptb}
\begin{array}{l}
\text{For}\; x,y\in \Int[\vecti],\;
\text{such that}\; \red_{(1+\vecti)^3}(x),
\,\red_{(1+\vecti)^3}(y)\in \tilde{U}_{(1+\vecti)^3},\;\text{we have}\;\vspace*{0.2cm}\\
\hspace*{0.4cm}
\red_{(1+\vecti)^2}(xy)=
\begin{cases}1&\text{if $\rt_3\left(\red_{(1+\vecti)^3}(x^2)
\right)=\rt_3\left(\red_{(1+\vecti)^3}(y^2)\right).$}\\
\vecti&\text{if $\rt_3\left(\red_{(1+\vecti)^3}(x^2)
\right)\neq\rt_3\left(\red_{(1+\vecti)^3}(y^2)\right).$}
\end{cases}
\end{array}
\finisheqn
In order to verify \eqref{eqn:Rcasessublemclaimptb}, note
that by \eqref{eqn:xrtxsqinvinker}, we have
\starteqn\label{eqn:Rcasessublem3}
\red_{(1+\vecti)^3}(x)\left(\rt_3\left(\sq(\red_{(1+\vecti)^3}(x))
\right)\right)\inv\in\ker\left(\sq|_{U_{(1+\vecti)^3}}\right),
\forallindisp\; x\in\Int[\vecti].
\finisheqn
By Table \eqref{table:sqresiduesinfo}, though we have
\startdisp
\red_{(1+\vecti)^2}\left(\ker\left(\sq|_{U_{(1+\vecti)^3}}\right)\right)
=\{1\}.
\finishdisp
By \eqref{eqn:Rcasessublem3}, we can therefore replace
$x$ and $y$ in \eqref{eqn:Rcasessublemclaimptb} with
$\rt_3\left(\red_{(1+\vecti)^3}(x^2)\right)$ and
$\rt_3\left(\red_{(1+\vecti)^3}(y^2)\right)$.  Also,
from Table \eqref{table:sqresiduesinfo}, we see
that the two possible values for $\rt_3\left(\red_{(1+\vecti)^3}(x^2)\right)$ and
and $\rt_3\left(\red_{(1+\vecti)^3}(y^2)\right)$
are $1$ and $\vecti$.  A trivial calculation now completes
the proof of \eqref{eqn:Rcasessublemclaimptb}.

We now complete the proof of (b) from \eqref{eqn:Rcasessublemclaimptb}.
Assume first that all three properties
\textbf{R}$\mathbf{0}(3)$ and \textbf{R}$\mathbf{1}(3)$
and \textbf{R}$\mathbf{2}(3)$ are satisfied simultaneously.
Since, as we have noted in Table \eqref{table:sqresiduesinfo},
\startdisp
\#\rt_3\left(\sq(U_{(1+\vecti)^3})\right)=\#
\left(\sq(U_{(1+\vecti)^3})\right)= 2,
\finishdisp
and since we are assuming that
\textbf{R}$\mathbf{2}(3)$ is satisfied,
Sublemma \ref{sublem:Rcases}, part (b) applies.  So,
for some permutation of \eqref{eqn:quadrupleres},
we have
\startdisp
x'^2\equiv y'^2\mod(1+\vecti)^n,\quad\text{and}\;z'^2\equiv w'^2
\mod(1+\vecti)^n.
\finishdisp
That is,
\starteqn\label{eqn:squaresofpairssame}
\rt_3\left(\red_{(1+\vecti)^3}(x'^2)
\right)=\rt_3\left(\red_{(1+\vecti)^3}(y'^2)\right),\;\text{and}\;
\rt_3\left(\red_{(1+\vecti)^3}(z'^2)
\right)=\rt_3\left(\red_{(1+\vecti)^3}(w'^2)\right).
\finisheqn
Using \eqref{eqn:squaresofpairssame} and \eqref{eqn:Rcasessublemclaimptb}
we deduce that
\startdisp
x'y'-z'w'\equiv x'z'-y'w'\equiv x'w'-y'z'\equiv 0\mod (1+\vecti)^2.
\finishdisp
In particular, since $(x',\,y',\,z',\,w')$ is a permutation
of \eqref{eqn:quadrupleres}, we have
\starteqn\label{eqn:determinantzero}
a'd'-b'c'\equiv 0\mod(1+\vecti)^2.
\finisheqn
Since \eqref{eqn:determinantzero} contradicts property
\textbf{R}$\mathbf{2'}(3)$.  Therefore, our assumption
is inconsistent and all three properties
\textbf{R}$\mathbf{0}(3)$ and \textbf{R}$\mathbf{1}(3)$
and \textbf{R}$\mathbf{2}(3)$ \textit{cannot} be
satisfied simultaneously.
\end{myproof}

\vspace*{0.3cm}\noindent\boldmath\textbf{A quadruple of
integers canonically associated with a quadruple \eqref{eqn:quadruple}
satisfying \eqref{eqn:quadruplecond}}.
\ubmth
Let $(a,b,c,d)$ as in \eqref{eqn:quadruple} be
a quadruple of complex numbers satisfying the conditions
of \eqref{eqn:quadruplecond}.  Recall that,
by Lemma \ref{lem:squaresorders}, part (b),
there is a unique $i\in\{0,1,2\}$ such that Possibility
\textbf{P}$i$ holds.


\begin{sublem}\label{sublem:multiplehasprops}
Let $(a,b,c,d)\in\Complex^4$ as in \eqref{eqn:quadruple} satisfy
conditions \eqref{eqn:quadruplecond}.  Let $i$ be the unique
element of $\{0,1,2\}$ such that Possibility $\textbf{P}i$ holds.
\begin{itemize}
\item[(a)]  For any permutation of (a,b,c,d),
we have
\starteqn\label{eqn:xprimedsquaredorders}
\ord_{\nu}\left(\left(\sqrt{1+\vecti}^ix\right)^2\right)
\in2\Int_{\geq 0},\forallindisp\;\text{Gaussian primes}\;\nu.
\finisheqn
\item[(b)]  The quadruple
\starteqn\label{eqn:multquad}
\sqrt{1+\vecti}^i(a,b,c,d)
\finisheqn
satisfies Properties
(i$'$), (ii), and (iii), and either (iv) or (v).  In particular,
the quadruple of \eqref{eqn:multquad} satisfies
(iv), resp. (v), if and only if $(a,b,c,d)$ satisfies
(iv), resp. (v).
\end{itemize}
\end{sublem}
\begin{myproof}{Proof}  For (a), let $\nu$ be a Gaussian
prime and use \textbf{ORD 1} to calculate
\startdisp
\ord_{\nu}\left(\left(\sqrt{1+\vecti}^ix\right)^2\right)=
\ord_{\nu}\left((1+\vecti)^i\right)+\ord_{\nu}(x^2).
\finishdisp
Using the definition of Property
\textbf{P}$i$, we obtain
\starteqn\label{eqn:ordersofxprimeindifferentcases}
\ord_{\nu}=\left(\left(\sqrt{1+\vecti}^ix\right)^2\right)=
\begin{cases}i+\ord_{1+\vecti}(x^2) &\text{if $\nu=1+\vecti$.}\\
\ord_{\nu}(x^2)&\text{if $\nu\neq 1+\vecti$.}
\end{cases}
\finisheqn
In order to complete the proof of (a), use
Lemma \ref{lem:squaresorders}, part (a), to obtain
$\ord_{\nu}(x^2)\in\Int_{\geq 0}$ for $\nu\neq 1+\vecti$.
That leaves only the case of $\nu=1+\vecti$.  By the \textbf{P}$i$
possibilities, we have
\startdisp\label{eqn:ordersxsqcases12}
\ord_{1+\vecti}(x^2)=\begin{cases}-i&\text{if}\;\textbf{P}1\;\text{or}
\;\textbf{P}2\;\text{is
satisfied.}\\
\text{an element of $2\Int$}&\text{if \textbf{P}2 is satisfied.}\end{cases}
\finishdisp
Therefore, in the first case of \eqref{eqn:ordersofxprimeindifferentcases}
we obtain
\startdisp
\ord_{1+\vecti}(\sqrt{1+\vecti}^ix)^2=\begin{cases}0&\text{if}\;\textbf{P}1\;\text{or}
\;\textbf{P}2\;\text{is
satisfied.}\\
\ord_{1+\vecti}(x)\in 2\Int&\text{if \textbf{P}2 is satisfied.}
\end{cases}
\finishdisp
This completes the proof of (a).

In order to prove the statements
of (b), first note that, \eqref{eqn:quadruple}
satisfies Property (i$'$),  The reason is we are assuming
it satisfies the conditions of \eqref{eqn:quadruplecond},
which include Property (i), and Property (i)
implies Property (i$'$).
Also, by Sublemma \ref{sublem:props34}, \eqref{eqn:quadruple} satisfies
both Properties (ii) and (iii).  Therefore,
we can apply Sublemma \ref{sublem:multiplehasprops}
with $z=\sqrt{1+\vecti}^i$ to conclude
that the multiple $\sqrt{1+\vecti}^i(a,b,c,d)$
has Properties (i$'$), (ii) and (iii).
Since $\sqrt{1+\vecti}^i(a,b,c,d)$ has Property (i$'$),
Sublemma \ref{sublem:props45} implies that $\sqrt{1+\vecti}^i(a,b,c,d)$
has either Property (iv) or Property (v).
\end{myproof}

In several of the lemmas that follow, we will make repeated
use of the same set of hypotheses.  In order to save
the trouble of reiterating them, we now define the following
hypotheses, \textbf{HYP}$(i,\delta)$,
which a quadruple $(a,b,c,d)\in\Complex^4$ may satisfy.  In the label \textbf{HYP}$(i,\delta)$
collectively given to these hypotheses, $i$ is a variable
taking values in the set $\{0,1,2\}$, and $\delta$
is a variable taking values in the set in the set $\{0,1\}$.  As
usual, $(x,y,z,w)$ denotes an arbitrary permutation of $(a,b,c,d)$.
\begin{itemize}
\item[\textbf{HYP}$(i,\delta)$]  First, $x^2\in\Rational(\vecti)$.
Next, with $\nu$ a Gaussian prime, we have
\starteqn\label{eqn:hypfirstpart}
\ord_{\nu}(x^2)=\begin{cases}\text{element of $2\Int_{\geq 0}$}&\text{for
$\nu\neq1+\vecti$}.\\
\text{-i}&\text{for $\nu=1+\vecti$ if $i=1$ or $2$}.\\
\text{element of $2\Int_{\geq 0}$}&\text{for $\nu=1+\vecti$ if $i=0$}.
\end{cases}
\finisheqn
Finally, we have
\starteqn\label{eqn:sqintegralindepa1}
(i,\delta)\;\text{minimal in $\Int\times\Int_{\geq 0}$ such
that}\;\omega_8^{\delta}\sqrt{1+\vecti}^ix\in\Int[\vecti].
\finisheqn
\end{itemize}
Before giving Lemma \ref{lem:sqintegralindep},
we give the following definitions to which the definition
of \textbf{HYP}$(i,\delta)$
naturally leads.
\begin{defn} \label{defn:primedquadruple}
Let $(a,b,c,d)\in\Complex^4$.  Suppose that
$(a,b,c,d)$ satisfies \textbf{HYP}$(i,\delta)$ for some $(i,\delta)
\in\{0,1,2\}\times\{0,1\}$
and set\bmth
\starteqn\label{eqn:ideltadefn}
\begin{aligned}&
\big(i,\delta\big)\left((a,b,c,d)\right)=\\
&\text{\hspace*{0.5cm}the uniquely determined pair \ubmth
$(i,\delta)$ such that $(a,b,c,d)$ satisfies \textbf{HYP}$(i,\delta)$}.
\end{aligned}
\finisheqn
\ubmth
Extend this notation from $\Complex^4$ to $\mathrm{Mat}_2(\Complex)$,
by setting
\bmth
\starteqn\label{eqn:ideltadefnmatrices}
\big(i,\delta\big)(\alpha)=\\
\text{\ubmth
$\big(i,\delta\big)\big((a,b,c,d)\big)$, where the latter exists and $\alpha=
\begin{pmatrix}a&b\\c&d\end{pmatrix}$}.
\finisheqn

Define the quadruple $(a',\,b',\,c',\,d')\in(\Int[\vecti])^4$ by
\starteqn\label{eqn:primedquadrupledefn}
(a',\,b',\,c',\,d'):=\omega_8^{\delta}\sqrt{1+\vecti}^i
(a,b,c,d).
\finisheqn
In the situation of \eqref{eqn:ideltadefnmatrices}, set
\starteqn\label{eqn:primematrixdefn}
\alpha'=\omega_8^{\delta}\sqrt{1+\vecti}^i\alpha.
\finisheqn
\end{defn}

Note that whenever $(a,b,c,d)$ satisfies \textbf{HYP}$(i,\delta)$,
and $(a',\,b',\,c',\,d')$ is as defined in \eqref{eqn:primedquadrupledefn},
we have
\starteqn\label{eqn:sqintegralindepa2}
x'
\sim(1+\vecti)^{(\ord_{1+\vecti}(x^2)+i)/2}\prod_{\nu\neq (1+\vecti)}
\nu^{\ord_{\nu}(x^2)/2},
\finisheqn
where the product extends over Gaussian primes $\nu$
other than $1+\vecti$, and all the exponents
appearing on the
right-hand side of \eqref{eqn:sqintegralindepa2} are in
$\Int_{\geq 0}$.

\begin{lem} \label{lem:sqintegralindep}
Let $(a,b,c,d)\in\Complex^4$ as in \eqref{eqn:quadruple} be
a quadruple of complex numbers satisfying the conditions
of \eqref{eqn:quadruplecond}.  Thus there exists uniquely determined pair
of integers,
\startdisp
(i,\delta)=\big(i,\delta\big)\big((a,b,c,d)\big)\in\{0,1,2\}\times \{0,1\}
\finishdisp
such that $(a,b,c,d)$ satisfies \textbf{HYP}$(i,\delta)$.
\end{lem}
\begin{myproof}{Proof}
Applying Sublemma \ref{sublem:multiplehasprops}, part (b),
we see that
\starteqn\label{eqn:multhaspropsapplication}
\text{The quadruple $(a,b,c,d)$ either has Property (iv),
or it has Property (v)
}
\finisheqn
Suppose first that $(a,b,c,d)$ has Property (iv).  Then
Lemma \ref{lem:squaresorders} implies that there
is a unique $i\in\{0,1,2\}$ such that $(a,b,c,d)$ satisfies
Possitibility \textbf{P}$i$.  By Sublemma \ref{sublem:multprops},
with $z=\sqrt{1+\vecti}^i$,
$\sqrt{1+\vecti}^i(a,b,c,d)$ has Property (iv).
Then Sublemma \ref{sublem:sqrationalindep} implies that,
for
\starteqn\label{eqn:deltatheright}
\delta=\delta\left(\sqrt{1+\vecti}^i(a,b,c,d)\right)\in\{0,1\}
\finisheqn
defined as in \eqref{eqn:deltaxdefn},
we have
\startdisp
\omega_8^{\delta}\sqrt{1+\vecti}^iy\in\Rational(\vecti),
\finishdisp
with
\starteqn\label{eqn:sqintegralindeplem1}
\omega_{8}^{\delta}x\sqrt{1+\vecti}^i
\sim\prod_{\nu}\nu^{\ord_{\nu}\left(\left(\sqrt{1+\vecti}^i
x\right)^2\right)/2}.
\finisheqn
In \eqref{eqn:sqintegralindeplem1}, the product
extends over the Gaussian primes.
By \eqref{eqn:xprimedsquaredorders}, each of the exponents
in the product on the right-hand side of \eqref{eqn:sqintegralindeplem1}
belongs to $\Int_+$.  So the
right-hand side of \eqref{eqn:xprimedsquaredorders}
belongs to $\Int[\vecti]$, which establishes \eqref{eqn:sqintegralindepa1}.
We have therefore shown that with $i$ chosen so that
$(a,b,c,d)$ satisfies Possibility \textbf{P}$i$ and
$\delta$ chosen as in \eqref{eqn:deltatheright},
$(a,b,c,d)$ satisfies \textbf{HYP}$(i,\delta)$.
This completes the proof of the sublemma
in the case when Property (iv) is satisfied.

By \eqref{eqn:multhaspropsapplication} we may henceforth
assume that $(a,b,c,d)$ has Property (v).
 By Lemma \ref{lem:squaresrational}, Part (b),
$(a,b,c,d)$ has Property (vii).
Let $j\in\Int$ be as in the statement of Property (vii), and set
\starteqn\label{eqn:deltajred}
\delta=\red_2(j),
\finisheqn
so that $\delta=0$ or $1$, depending on whether $j$ is even or odd.
Then Property (vii) says that either $x=0$, in which
case \eqref{eqn:sqintegralindepa1} and \eqref{eqn:sqintegralindepa2}
are reduced to trivialities, or
\startdisp
\omega_8^{\delta}x=\vecti^{\frac{\delta+j}{2}}.
\finishdisp
By the choice of $\delta$ in \eqref{eqn:deltajred},
$\delta+j\in2\Int$, so that
$\omega_8^{\delta}x$ is a unit of the Gaussian
integers.  Therefore, \eqref{eqn:sqintegralindepa1} is verified.
So, with $\delta$ chosen as in \eqref{eqn:deltajred},
$(a,b,c,d)$ satisfies \textbf{HYP}$(0,\delta)$.
This completes the proof of Lemma \ref{lem:sqintegralindep}.
\end{myproof}

\vspace*{0.3cm}\noindent\boldmath\textbf{The `R' properties
applied to $(a',\,b',\,c',\,d')$}.
\ubmth
We note for later use a few simple consequences of
the definition of \eqref{eqn:primedquadrupledefn}.
These consequences all spring from the relation
the quadruples \eqref{eqn:quadruple} and \eqref{eqn:primedquadrupledefn}.
\starteqn\label{eqn:scalarmult}
(a',\,b',\,c'\,d')\;\text{is a scalar multiple of}\;
(a,b,c,d),\;\text{by the complex number};\omega_8^{\delta}\sqrt{1+\vecti}^i.
\finisheqn
Throughout the following we will use $m$ to denote a nonegative
integer.  Further, in the following,
\startdisp
f_m(x,y,z,w)\;\text{denotes a polynomial in the variables $x,y,z,w$,
homogeneous of degree $m$.}
\finishdisp
From \eqref{eqn:scalarmult}, we deduce that for arbitrary $f_m$,
\starteqn\label{eqn:scalarmulthomogpoly}
f_m(a',\,b',\,c',\,d')=\left(\omega_8^{\delta}\sqrt{1+\vecti}^i\right)^m
f_m(a,b,c,d)
\finisheqn
From \eqref{eqn:scalarmulthomogpoly}, we deduce that for $\ell\in\Int$,
\starteqn\label{eqn:scalarmultpolyeven}
\text{if $f_{2\ell}(a,b,c,d)\equiv 0\mod (1+\vecti)^{n}$, then
$f_{2\ell}(a,b,c,d)\equiv 0\mod(1+\vecti)^{n+\ell i}$
}
\finisheqn
The claim \eqref{eqn:scalarmultpolyeven} is an immediate consequence
of \eqref{eqn:scalarmulthomogpoly} and the calculation
\startdisp
\left(\omega_8^{\delta}\sqrt{1+\vecti}\right)^{2\ell}=
\vecti^{\delta\ell}(1+\vecti)^{\ell}\in\Int[\vecti].
\finishdisp
We will apply \eqref{eqn:scalarmulthomogpoly} in particular
to the case when $m=2$ and $f_2$ is the ``determinant polynomial
$f_2(x,y,z,w)=xw-yz$, which is homogeneous of degree $2$.
\starteqn\label{eqn:determinantprimed}
a'd'-b'c'=\vecti^{\delta}(1+\vecti)^i ad-bc.
\finisheqn
Applying \eqref{eqn:determinantprimed} to \eqref{eqn:primematrixdefn}
we find that
\startdisp
\det(\alpha')=\vecti^{\delta}(1+\vecti)^i\det(\alpha).
\finishdisp
Recalling the set $M^{\rm N}_2$ consisting of 2-by-2 integer matrices
of determinant $N$, introduced in \S\ref{subsec:flts}, we are led to
the following important re-interpretation of the condition
$\textbf{HYP}(i,\delta)$ for even $i$.  Let $j\in\Int_{\geq 0}$, $\delta=0$
or $1$.  Then,
\starteqn\label{eqn:hypreinterpretationeven}
\{\alpha\in\SL{2}{\Complex}\;|\; \text{$\alpha$ satisfies
$\textbf{HYP}(2j,\delta)$}\}=\frac{1}{\omega_8^{\delta}(1+\vecti)^{j}}
M^{\vecti^{\delta}(1+\vecti)^j}_2.
\finisheqn

We will apply \eqref{eqn:scalarmultpolyeven} to the case when $\ell=1$
$f_2(x,y,z,w)=x^2+y^2-z^2-w^2$, a homogeneous
polynomial of degree $2\ell=2$.
Referring to the definition of Property \textbf{R}$\mathbf{2}(n)$
above, we see that in this case \eqref{eqn:scalarmultpolyeven} implies that
\starteqn\label{eqn:R2satbyprimedquad}\text{
If $(a,b,c,d)$ satisfies Property \textbf{R}$\mathbf{2}(n)$,
then $(a',\,b',\,c',\,d')$ satisfies \textbf{R}$\mathbf{2}(n+i)$.
}\finisheqn

\begin{lem}\label{lem:multiplesatsprops}  Let $(a,b,c,d)\in\Complex^4$,
as in \eqref{eqn:quadruple}.   Suppose that $(a,b,c,d)$ satisfies
\eqref{eqn:quadruplecond}.
\begin{itemize}\item[(a)]  Suppose that $(a,b,c,d)$ satisfies
\textbf{HYP}$(n,\delta)$ with
\startdisp
\big(n,\delta\big)\big((a,b,c,d)\big)\in\{1,2\}\times\{0,1\}.
\finishdisp
Then
\starteqn\label{eqn:multiplisatsprops1}
(a',\,b',\,c',\,d')\;\text{satisfies}\;\text{\textbf{R}$\mathbf{i}(n+2)$}\;
\text{for $i\in\{0,1,2\}$}.
\finisheqn
\item[(b)]  Suppose that $(a,b,c,d)$ satisfies \textbf{HYP}$(0,\delta)$
for $\delta\in\{0,1\}$.  Then
\starteqn\label{eqn:r12sat}
(a',\,b',\,c',\,d')\;\text{satisfies}\;\text{\textbf{R}$\mathbf{1}(2)$
and \textbf{R}$\mathbf{2}(2)$}.
\finisheqn
\end{itemize}
\end{lem}
\begin{myproof}{Proof}
Throughout, we may assume that $(a,b,c,d)$ satisfies \textbf{HYP}$(n,\delta)$.
We will first assume that
\startdisp
\big(n,\delta\big)\big((a,b,c,d)\big)\in\{1,2\}\times\{0,1\},
\finishdisp
as in Part (a), and establish \eqref{eqn:multiplisatsprops1}
in the case $i=0$.
By using \eqref{eqn:hypfirstpart} and \eqref{eqn:sqintegralindepa2}
we deduce that $\ord_{1+\vecti}(x')=0$.
This proves that $x'\equiv 1\mod (1+\vecti)$, and therefore
we have verified \eqref{eqn:multiplisatsprops1} in the case when $i=0$.

We now allow
\startdisp
\big(n,\delta\big)\big((a,b,c,d)\big)\;\text{to take any
value in}\;\{0,1,2\}\times\{0,1\},
\finishdisp
We will prove that
\startdisp
(a',\,b',\,c',\,d')\;\text{satisfies}\;\text{\textbf{R}$\mathbf{i}(n+2)$}\;
\text{for $i\in\{1,2\}$},
\finishdisp
and this will complete the proof of the lemma.
Since $(a,b,c,d)$ satisfies the conditions of \eqref{eqn:quadruplecond},
it is clear that $(a,b,c,d)$ has Property (i).  From
Property (i) and \eqref{eqn:determinantprimed}, we have
\startdisp
a'd'-b'c'=\vecti^{\delta}(1+\vecti)^n
\finishdisp
Reducing modulo $(1+\vecti)^{n+2}$, we deduce that
\startdisp
a'd'-b'c'\cong (1+\vecti)^n\mod(1+\vecti)^{n+2}.
\finishdisp
Therefore, $(a',\,b',\,c',\,d')$ satisfies \textbf{R}$\mathbf{1}(n+2)$.
For Property \textbf{R}$\mathbf{2}(n+2)$, note that
by Sublemma \ref{sublem:props34} $(a,b,c,d)$ satisfies
Property (iii).  By applying the map $\red_{(1+\vecti)^2}$
we see that Property (iii) implies Property \textbf{R}$\mathbf{2}(2)$.
Thus $(a,b,c,d)$ satisfies Property \textbf{R}$\mathbf{2}(2)$.
We now use \eqref{eqn:R2satbyprimedquad} to deduce
that $(a',\,b',\,c',\,d')$ satisfies Property \textbf{R}$\mathbf{2}(n+2)$.
This completes the proof of the lemma.
\end{myproof}

\vspace*{0.3cm}\noindent\boldmath\textbf{Completion of the proof
of Proposition \ref{prop:inversematrixexplicitdescription}}.
\ubmth
Applying Part (a) of Lemma \ref{lem:multiplesatsprops},
we will indicate in \eqref{eqn:idivision} below how
the two parts of the union \eqref{eqn:inversematrixexplicitdescription}
can be characterized by the datum $i$ in
$\big(i,\delta\big)\big((a,b,c,d)\big)$.  Let
\starteqn\label{eqn:alphaquadentries}
\alpha=\begin{pmatrix}a&b\\c&d\end{pmatrix}\in\conj\inv(\Gamma).
\finisheqn
Recall that by Lemma \ref{lem:sqintegralindep}, each $(a,b,c,d)$
satisfying the conditions of \eqref{eqn:quadruplecond}
satisfies $\textbf{HYP}(i,\delta)$ for a uniquely determined
pair
\startdisp
\big(i,\delta\big)\big((a,b,c,d)\big)\in\{0,1,2\}\times \{0,1\}.
\finishdisp
For convenience, define a function
\startdisp
(i,\delta)(\cdot):\conj\inv(\Gamma)\rightarrow
\{0,1,2\}\times\{0,1\},
\finishdisp
by setting
\startdisp
\big(i,\delta\big)(\alpha)=\big(i,\delta\big)\big((a,b,c,d)\big),
\finishdisp
where $\alpha$ as in
\eqref{eqn:alphaquadentries} belongs to $\conj\inv(\Gamma)$.
By combining previously established results, we deduce
that for any $\alpha\in\conj\inv(\Gamma)$, $(i,\delta)(\alpha)
\notin\{1\}\times\{0,1\}$, so that
\starteqn\label{eqn:hypothesis1imposs}
\alpha\in\conj\inv(\Gamma)\;\text{implies that}\;
\big(i,\delta\big)(\alpha)\in\{0,2\}\times\{0,1\}.
\finisheqn
The reason for \eqref{eqn:hypothesis1imposs} is that, supposing
that $\alpha$ satisfies $\textbf{HYP}(i,\delta)$, Lemma
\ref{lem:multiplesatsprops}, Part (a) implies that
\starteqn\label{eqn:multiplisatsprops1case1}
(a',\,b',\,c',\,d')\;\text{satisfies}\;\text{\textbf{R}$\mathbf{i}(3)$}\;
\text{for $i\in\{0,1,2\}$}.
\finisheqn
Yet, Sublemma \ref{sublem:Rcases} says that
\eqref{eqn:multiplisatsprops1case1} cannot occur.  This contradiction
proves \eqref{eqn:hypothesis1imposs}.
On the hand hand, we have the case when
$\big(i,\delta\big)(\alpha)=(0,\delta)$ for $\delta\in\{0,1\}$.
According to \eqref{eqn:sqintegralindepa1},
we then have
\startdisp
\omega_8^{\delta}(a,b,c,d)\in(\Int[\vecti])^4.
\finishdisp
Since $\alpha\in\SL{2}{\Complex}$, by assumption, we therefore have
in the case $i=0$, that
$\alpha\in\SL{2}{\Int[\omega_8]}$.  On the other hand, we have
the case when $\big(i,\delta\big)(\alpha)=(2,\delta)$.  In that
case, \eqref{eqn:sqintegralindepa1} implies that
\startdisp
\omega_8^{\delta}\sqrt{1+\vecti}^2(a,b,c,d)\in\Int[\vecti],\;\text{but}\;
\omega_8^{\delta}(a,b,c,d)\notin(\Int[\vecti])^4.
\finishdisp
Therefore, in the case that $i=2$ we have $\alpha\notin\SL{2}{\Int[\omega_8]}$.
Summing up the above discussion we have so far demonstrated that
\starteqn\label{eqn:idivision}
\begin{aligned}
\conj\inv(\Gamma)\cap\SL{2}{\Int[\omega_8]}&=&
\big(i,\delta\big)\inv(0\times\{0,1\}),\\
\conj\inv(\Gamma) - \conj\inv(\Gamma)\cap\SL{2}{\Int[\omega_8]}
&=&\big(i,\delta\big)\inv(2\times\{0,1\}).
\end{aligned}
\finisheqn
The equalities of \eqref{eqn:idivision} will be used in conjunction
with \eqref{eqn:hypreinterpretationeven} to prove Part (a)
of Lemma \ref{lem:xi12actions}, below.

Next, in Lemma \ref{lem:subsetconstr} below
we prove a partial converse to Part (b)
of Lemma \ref{lem:multiplesatsprops}.  In Lemma \ref{lem:subsetconstr}
and its proof, we make use of the rational extension
$\tilde{\conj}$ of $\conj$
defined in \eqref{eqn:rationalextension}.
The reader should also recall the properties of $\tilde{\conj}$ stated in
\eqref{eqn:multiplehasintegral} and \eqref{eqn:multiplehasintegralunits}.

\begin{lem}\label{lem:subsetconstr}  Suppose that
$(a,b,c,d)\in\Complex^4$, with
\startdisp
\alpha=\begin{pmatrix}a&b\\c&d\end{pmatrix}\in\mathrm{Mat}_2(\Complex).
\finishdisp
Assume that $(a,b,c,d)$ satisfies $\textbf{HYP}(0,\delta)$
with $\delta\in\{0,1\}$.
\begin{itemize}\item[(a)]  Assuming also
that $(a'\,b'\,c',\,d')$ satisfies \textbf{R}$\mathbf{2}(2)$, we have
$\tilde{\conj}(\alpha)\in\mathrm{Mat}_3(\Int[\vecti])$.
\item[(b)]  We have
\starteqn\label{eqn:determinantcondition}
\alpha\in\conj\inv(\Gamma)\;\text{if and only if}\;
(a',\,b',\,c',\,d')\;\text{satisfies}\;\textbf{R}\mathbf{2}(2)\;
\text{and}\;
ad-bc=1.
\finisheqn
\end{itemize}
\end{lem}
\begin{myproof}{Proof}  For Part (a).
assume that $(a,b,c,d)$ satisfies $\textbf{HYP}(0,\delta)$
with $\delta\in\{0,1\}$, and also
that $(a'\,b'\,c',\,d')$ satisfies \textbf{R}$\mathbf{2}(2)$.
Define
\startdisp
\alpha'\in\mathrm{Mat}_2(\Int[\vecti])\;\text{by}\;\,
\alpha'=\begin{pmatrix}a'&b'\\c'&d'\end{pmatrix}.
\finishdisp
By Definition \ref{defn:primedquadruple}, we have $x'=\omega_8^{\delta}x$.
Therefore, $\alpha'=\omega_8^{\delta}\alpha$.
By \eqref{eqn:omegaeightdefn}, we may apply \eqref{eqn:multiplehasintegralunits}
with $\ell=\omega_8^{\delta}$.  We obtain
\starteqn\label{eqn:constructionsublemequiv}
\tilde{\conj}(\alpha)\in\mathrm{Mat}_3(\Int[\vecti])\;\text{if and only if}\;
\tilde{\conj}(\alpha')\in\mathrm{Mat}_3(\Int[\vecti]).
\finisheqn
Therefore, have reduced the proof of Part (a)
to showing that the entries of $\tilde{\conj}(\alpha')$ are Gaussian
integers.
Since the entries in the third row and column of \eqref{eqn:imagematrixconj},
the relation defining $\tilde{\conj}$,
belong to $\Int[a,b,c,d]$, and since $\alpha'\in\mathrm{Mat}_2(\Int[\vecti])$,
the entries in the third row and colun of $\tilde{\conj}(\alpha')$
belong to $\Int[\vecti]$.
It remains to prove that the elements in the upper-left
2-by-2 block of $\tilde{\conj}(\alpha')$ belong to $\Int[\vecti]$
Examining \eqref{eqn:imagematrixconj} again, we see
that in order to prove that the remaining
entries of $\tilde{\conj}(\alpha')$ belong to $\Int[\vecti]$,
it suffices to prove that
\starteqn\label{eqn:subsetconstrlem3}
x'^2+y'^2\pm (z'^2+w'^2)\in \left((1+\vecti)^2\right),
\finisheqn
In \eqref{eqn:subsetconstrlem3},  $(x',\,y',\,z',\,w')$ denotes
an arbitrary permutation
of $(a',\,b',\,c',\,d')$, as usual.
But \eqref{eqn:subsetconstrlem3} follows immediate from
the hypothesis that $(a',\,b',\,c',\,d')$ has Property
\textbf{R}$\mathbf{2}(2)$.  So we may conclude that
$\tilde{\conj}(\alpha')\in\mathrm{Mat}_3(\Int[\vecti])$,
and by \eqref{eqn:constructionsublemequiv}, this completes
the proof of Part (a).

For Part (b), assume at first that
$(a',\,b',\,c',\,d')$ satisfies \textbf{R}$\mathbf{2}(2)$
and
$ad-bc=1$.
Recall that the conditions
of \eqref{eqn:quadruplecond}, which are equivalent to
$\alpha\in\conj\inv(\Gamma)$, are that the right
side of \eqref{eqn:imagematrixconj} has integer entries,
together with the determinant condition
$ad-bc=1$.  By the definition of $\tilde{\conj}(\alpha)$,
it is clear that the right
side of \eqref{eqn:imagematrixconj} has integer entries
if and only if $\tilde{\conj}(\alpha)\in\mathrm{Mat}_3(\Int[\vecti])$.
Therefore, by Part (a), the integer
condition on the entries is satisfied.  The determinant
condition is satisfied by the assumption.
Thus $\alpha\in\conj\inv(\Gamma)$.  Conversely, if $\alpha\in\conj\inv(\Gamma)$,
then $ad-bc=1$, because $\alpha\in\SL{2}{\Complex}$.
Further, by Lemma \eqref{lem:multiplesatsprops}, and the
assumption that $(a,b,c,d)$ satisfies $\textbf{HYP}(0,\delta)$, we have that
$(a',\,b',\,c',\,d')$ satisfies \textbf{R}$\mathbf{2}(2)$.
This completes the proof of Part (b).
\end{myproof}

Lemma \ref{lem:subsetconstr} is used in two slightly different
ways below.  In Corollary \ref{cor:subsetconstr} the Lemma
is used to explain how the parameter
$\delta\in\{0,1\}$ of $\big(i,\delta\big)(\alpha)$ shows up in the structure of
\eqref{eqn:inversematrixexplicitdescription},
in much
the same way that \eqref{eqn:idivision} explains how the parameter
$i\in\{0,2\}$ shows up in \eqref{eqn:inversematrixexplicitdescription}.
In Corollary \ref{cor:xi12} Lemma \ref{lem:subsetconstr} is used
to construct a ``sufficiently large" subgroup of $\conj\inv(\Gamma)$.
By ``sufficiently large" we mean that the subgroup constructed
in Corollary \ref{cor:xi12},
is large enough that, in conjunction with the tools developed
earlier in this section, we can determine all the elements of
$\conj\inv(\Gamma)$.

\begin{cor} \label{cor:subsetconstr} We have
\starteqn\label{eqn:deltarelation}
\{\alpha\in\conj\inv(\Gamma)\;|\; \big(i,\delta\big)
(\alpha)=(0,1)\}=\frac{1}{\omega_8}\begin{pmatrix}\vecti&0\\0&1\end{pmatrix}
\{\alpha\in\conj\inv(\Gamma)\;|\; \big(i,\delta\big)
(\alpha)=(0,0)\},
\finisheqn
so that
\starteqn\label{eqn:translatesubset}
\begin{aligned}&&&
\{\alpha\in\conj\inv(\Gamma)\;|\; \big(i,\delta\big)
(\alpha)\in\{0\}\times\{0,1\}\}=\\
&&&\hspace*{0.5cm} \{\alpha\in\conj\inv(\Gamma)\;|\; \big(i,\delta\big)
(\alpha)=(0,0)\}\bigcup\hspace*{-0.3cm}\cdot\hspace*{0.3cm}
\frac{1}{\omega_8}\begin{pmatrix}\vecti&0\\0&1\end{pmatrix}
\{\alpha\in\conj\inv(\Gamma)\;|\; \big(i,\delta\big)
(\alpha)=(0,0)\}.
\end{aligned}
\finisheqn
\end{cor}
\begin{myproof}{Proof}
Because $\big(i,\delta\big)
(\alpha)$ takes values in $\{0,2\}\times\{0,1\}$,
it is clear that \eqref{eqn:translatesubset}
follows from \eqref{eqn:deltarelation}.
We now prove \eqref{eqn:deltarelation}.
Suppose that $\alpha_0\in\SL{2}{\Complex}$,
and make the definition
\starteqn\label{eqn:alphaonedefn}
\alpha_1:=\frac{1}{\omega_8}\begin{pmatrix}\vecti&0\\0&1\end{pmatrix}\alpha_0.
\finisheqn
We will adopt the following notational conventions.  The subscripts
on roman letters representing entries will match the subscript
on the Greek letter representing the matrix to which the entries
belong.  That is, for $j=0,1$ the quadruple of entries of $\alpha_j$,
$\alpha_j$ as above, will be denoted $(a_j,b_j,c_j,d_j)$.  The
same practice will be observed with Roman letters permutations, meaning
that $(x_j,y_j,z_j,w_j)$ represents an arbitrary permutation
of $(a_j,b_j,c_j,d_j)$, and so on.
Whenever entries of $(x_j,y_j,z_j,w_j)$ for both $j=0$ and $j=1$
appear in the same equation, it will be assumed that only
the case of corresponding permutations are being considered,
so that if, for example $x_0$ is the image of $a_0$
under a permutation of the quadruple indexed by $0$, then $x_1$ is the image of $a_1$
under the permutation of the quadruple indexed by $1$.  Similar
comments apply to the primed forms of the matrices $\alpha'_j$
and the primed quadruples $(x_j',\,y_j',\,z_j',\,w_j')$, $j=0,1$
representing their

According to \eqref{eqn:alphaonedefn},
the relationship between the two quadruples is given by
\starteqn\label{eqn:subsetconstrlem4}
(a_1,b_1,c_1,d_1)=(a_0\omega_8,b_0\omega_8,
c_1\omega_8^{-1},d_1\omega_8^{-1}).
\finisheqn
It is clear that $\alpha_0$ satisfies $\textbf{HYP}(0,0)$
if and only if $\alpha_1$ satisfies $\textbf{HYP}(0,1)$.
Therefore, in order to complete the proof of \eqref{eqn:deltarelation},
it will suffice to show that
\startdisp
(a_0,b_0,c_0,d_0)\;\text{satisfies the hypotheses of Lemma
\ref{lem:subsetconstr}(b)
if and only if}\;(a_1,b_1,c_1,d_1)\;\text{does so.}
\finishdisp
With reference to \eqref{eqn:determinantcondition}
since
\startdisp
\det\left(\frac{1}{\omega_8}\begin{pmatrix}\vecti&0\\
0&1\end{pmatrix}\right)=1,
\finishdisp
it is easily verified, using \eqref{eqn:alphaonedefn},
that the determinant of $\alpha_j$
($=a_jd_j-b_jc_j$) is the same
for $j=0$ and $j=1$.  The hypotheses concerning $\textbf{HYP}(0,j)$
concerning has already been
dealt with above.
Assuming, without loss of generality that
$(a_j,b_j,c_j,d_j)$ satisfies
$\textbf{HYP}(0,j)$, we see that
our task is reduced to proving the following equivalence.
\starteqn\label{eqn:subsetconstrlem5}\text{
$(a_0',\,b_0'\,c_0'\,d_0')$ satisfies \textbf{R}$\mathbf{2}(2)$
if and only if $(a_1'\,b_1',\,c_1',\,d_1')$ satisfies
\textbf{R}$\mathbf{2}(2)$.
}\finisheqn
We have, by definition,
\startdisp
(a_0',\,b_0'\,c_0'\,d_0')=(a_0,b_0,c_0,d_0),
\finishdisp
and
\startdisp
(a_1',\,b_1'\,c_1'\,d_1')=\omega_8(a_1,b_1,c_1,d_1)=
(\vecti a_0,\vecti b_0, c_0, d_0).
\finishdisp
From these equalities, it follows immediately that
\starteqn\label{eqn:squaresequalmod2}
(x_0')^2\equiv (x_1')^2\mod 1+\vecti.
\finisheqn
From the definition of \textbf{R}$\mathbf{2}(2)$, we see
that \eqref{eqn:subsetconstrlem5} follows
from \eqref{eqn:squaresequalmod2}.  By the comments
preceding \eqref{eqn:subsetconstrlem5}, this
completes the proof of \eqref{eqn:deltarelation}.
\end{myproof}

\begin{cor}\label{cor:xi12}
We have the containment
\starteqn\label{eqn:containment}
\Xi_{12}\subseteq \conj\inv(\Gamma)\cap\SL{2}{\Int[\omega_8]}
\finisheqn
\end{cor}
\begin{myproof}{Proof}  We show that elements $\Xi_{12}$
satisfies the conditions of Lemma \ref{lem:subsetconstr}, Part (b),
for belonging to $\conj\inv(\Gamma)$.
It is clear that any element of $\SL{2}{\Int[\vecti]}$,
of which $\Xi_{12}$ is a subset,
satisfies hypothesis $\textbf{HYP}(0,0)$.
Since $\alpha\in\SL{2}{\Int[\vecti]}$, the determinant
condition $ad-bc=1$ is clearly satisfied.
That leaves the condition that $(a',\,b'\,c'\,d')$
satisfies \textbf{R}$\mathbf{2}(2)$.  Since
for $\alpha\in\Xi_{12}$, the quadruple of entries $(a,b,c,d)$
satisfies $\textbf{HYP}(0,0)$, we have
\startdisp
(a',\,b'\,c'\,d')=(a,b,c,d).
\finishdisp
Further, the $\Xi_{12}$ is simply defined as the subgroup
of elements of $\SL{2}{\Int[\vecti]}$ satisfying
\startdisp
(a,b,c,d)\equiv (0,1,0,1)\mod (1+\vecti)^2,
\finishdisp
from which we verify \textbf{R}$\mathbf{2}(2)$ directly.
\end{myproof}
\begin{lem} \label{lem:xi12actions} Let $N\in\Int[\vecti]$.
Recall the action of \hspace{0.03cm} $\Xi_{12}$ on $M_2^{\mathrm{N}}$
by left multiplication.
\begin{itemize}
\item[(a)]
We have the containments
\starteqn\begin{aligned}
\label{eqn:containmentseries}
\Xi_{12}&&&\subseteq&&& \conj\inv(\Gamma)\cap\SL{2}{\Int[\omega_8]}&\subseteq&&
\hspace*{-.3cm}\bigcup_{\delta= 0,1}\hspace*{-.45cm}\cdot
\hspace*{0.45cm}\frac{1}{\omega_8^{\delta}}M_2^{\vecti^{\delta}},\\&&&&&&
\conj\inv(\Gamma)-\conj\inv(\Gamma)\cap\SL{2}{\Int[\omega_8]}&\subseteq&&
\bigcup_{\delta= 0,1}\hspace*{-.45cm}\cdot
\hspace*{0.45cm}\frac{1}{\omega_8^{\delta}(1+\vecti)}M_2^{2\vecti^{1+\delta}}.
\end{aligned}
\finisheqn
\item[(b)]
As a consequence of \eqref{eqn:containmentseries}
the action of $\Xi_{12}$ on
\startdisp
\bigcup_{\delta= 0,1}\hspace*{-.45cm}\cdot
\hspace*{0.45cm}\frac{1}{\omega_8^{\delta}}M_2^{\vecti^{\delta}}
\finishdisp
preserves
$\conj\inv(\Gamma)\cap\SL{2}{\Int[\omega_8]}$,
and the action of $\Xi_{12}$ on
\startdisp
\bigcup_{\delta= 0,1}\hspace*{-.45cm}\cdot
\hspace*{0.45cm}\frac{1}{\omega_8^{\delta}(1+\vecti)}M_2^{2\vecti^{1+\delta}}
\finishdisp
preserves
$\conj\inv(\Gamma)-\conj\inv(\Gamma)\cap\SL{2}{\Int[\omega_8]}$.
\end{itemize}
\end{lem}
\begin{myproof}{Proof}   The containment of $\Xi_{12}$ in
$\conj\inv(\Gamma)\cap\SL{2}{\Int[\omega_8]}$ has
already been established in Corollary \ref{cor:subsetconstr}
above.  The remaining containments in \eqref{eqn:containmentseries}
follow immediately from
\eqref{eqn:idivision} and \eqref{eqn:hypreinterpretationeven}.
This completes the proof of Part (a).

Part (b) is readily derived from Part (a) using the facts
that $\conj\inv(\Gamma)$
and $\conj\inv(\Gamma)\cap\SL{2}{\Int[\omega_8]}$ are groups.
\end{myproof}

As will be seen below, the main point of Lemma \ref{lem:twopartsdescribed}, Part (a)
is that, by the use of certain facts established in \S\ref{subsec:flts},
the first containment in Part (a) of Lemma \ref{lem:subsetconstr}
can actually be replaced by an equality.

\begin{lem}\label{lem:twopartsdescribed}  For $N,m\in\Int[\vecti]$, such that $N/m$
is a standard integer and $x\in\Omega_{\frac{N}{m}}$, let
$\alpha^{\rm N}(m,x)$ be
as in \eqref{eqn:MofNmxdefn}.  We have
\starteqn\label{eqn:twopartsdescribed}
\begin{aligned}
\big(i,\delta\big)\inv(0\times\{0,1\})&=&&
\conj\inv(\Gamma)\cap\SL{2}{\Int[\omega_8]}&=&\bigcup_{\delta,=0,1}
\hspace*{-0.5cm}\cdot\hspace*{.5cm}
\frac{1}{\omega_8^{\delta}}
\Xi_{12}\alpha^{\vecti^{\delta}}(\vecti^{\delta},0),\\
\big(i,\delta\big)\inv(2\times\{0,1\})&=&&
\conj\inv(\Gamma)-\conj\inv(\Gamma)\cap\SL{2}{\Int[\omega_8]}&=&
\bigcup_{\delta,\epsilon=0,1}\hspace*{-0.55cm}\cdot\hspace*{.55cm}
\frac{1}{\omega_8^{\delta}(1+\vecti)}
\Xi_{2}\alpha^{2\vecti^{1+\delta}}\hspace*{-0.7mm}
(\vecti^{1+\delta},\vecti^{\epsilon}).
\end{aligned}
\finisheqn
\end{lem}
\begin{myproof} {Proof} We first note that the left equality
in each line of \eqref{eqn:twopartsdescribed} is a restatement of
the corresponding line in \eqref{eqn:idivision}.  So it remains
to prove the right equality in each line of \eqref{eqn:twopartsdescribed}.
The right equality in the first line  is derived
from the equality
\starteqn\label{eqn:xi12equality}
\conj\inv(\Gamma)\cap\SL{2}{\Int[\vecti]}=\Xi_{12}
\finisheqn
and previous results.  The derivation of the right equality in the
first line from \eqref{eqn:xi12equality}
proceeds by substituting the description of \eqref{eqn:idivision}
for the left side of Part (a).  Specializing
\eqref{eqn:hypreinterpretationeven} to the case $(j,\delta)=(0,0)$,
one has
\startdisp
\{\alpha\in\conj\inv(\Gamma)\;|\;\big(i,\delta\big)(\alpha)=(0,0)\}=
\conj\inv(\Gamma)\cap\SL{2}{\Int[\vecti]}.
\finishdisp
Making the substitution given by this equality, one has
\starteqn
\label{eqn:someparts}
\conj\inv(\Gamma)\cap\SL{2}{\Int[\omega_8]}=
\bigcup_{\delta,=0,1}
\hspace*{-0.5cm}\cdot\hspace*{.5cm}
\frac{1}{\omega_8^{\delta}}\begin{pmatrix}\vecti^{\delta}&0\\
0&1\end{pmatrix}
\left(\conj\inv(\Gamma)\cap\SL{2}{\Int[\vecti]}\right).
\finisheqn
Substituting \eqref{eqn:xi12equality} into \eqref{eqn:someparts},
one deduces that
\startdisp\label{eqn:twopartslemma1}
\conj\inv(\Gamma)\cap\SL{2}{\Int[\omega_8]}=
\bigcup_{\delta,=0,1}
\hspace*{-0.5cm}\cdot\hspace*{.5cm}
\frac{1}{\omega_8^{\delta}}\begin{pmatrix}\vecti^{\delta}&0\\
0&1\end{pmatrix}
\Xi_{12}.
\finishdisp
But it is easily verified that
\startdisp
\begin{pmatrix}\vecti^{\delta}&0\\
0&1\end{pmatrix}
\Xi_{12}=\Xi_{12}\alpha^{\vecti^{\delta}}(\vecti^{\delta},0),\;
\text{for $\delta=0,1$,}
\finishdisp
so we have completed the derivation of Part (a) from \eqref{eqn:xi12equality}.

We now verify \eqref{eqn:xi12equality}.  Part (b) of Lemma
\ref{lem:subsetconstr} and Sublemma \ref{sublem:finerdecomp} imply
that \starteqn\label{eqn:thisunion}
\conj\inv(\Gamma)\cap\SL{2}{\Int[\vecti]}=\bigcup
\hspace*{-.3cm}\cdot\hspace*{0.3cm}\Xi,\quad\text{$\Xi$ ranging over
a subset of $\{\Xi_{12}, \Xi_{1}, \Xi_{2}\}$}. \finisheqn By Lemma
\ref{lem:subsetconstr}, we already know that $\Xi_{12}$ is included
in the union on the right-hand side of \eqref{eqn:thisunion}.
Therefore, it will suffice to show that for a particular
$\alpha\in\Xi_1$, resp., $\Xi_2$, the
$\alpha\notin\conj\inv(\Gamma)$.
We choose
\starteqn\label{eqn:couplaelements}\alpha_1=
\begin{pmatrix}1&1\\0&1\end{pmatrix}\in\Xi_1,\quad\text{and}\;\alpha_2=
\begin{pmatrix}1&0\\1&1\end{pmatrix}\in\Xi_2.
\finisheqn
The elements $\alpha_1$, $\alpha_2$ satisfy
hypothesis \textbf{HYP}$(0,0)$, so that
\startdisp
(a'_i,\,b'_i,\,c'_i,\,d'_i)=(a_i,b_i,c_i,d_i),\;\text{for}\; i=1,2.
\finishdisp
Thus
\startdisp
\red_{(1+\vecti)^2}
(a'_i,\,b'_i,\,c'_i,\,d'_i)\equiv\red_{(1+\vecti)^2}(a_i,b_i,c_i,d_i)\equiv
\begin{cases}(1,1,0,1),\;\text{for $\alpha_1$}\\ \text{or}\\
(1,0,1,1),\;\text{for $\alpha_2$}
\end{cases}\;\mod(1+\vecti)^2.
\finishdisp
Therefore,
\startdisp
a_i'^2+d_i'^2\equiv 0\mod(1+\vecti)^2,\quad\text{while}\; b_i'^2+c_i'^2\equiv 1
\mod(1+\vecti)^2,\;\text{for $i=1,2$}.
\finishdisp
Therefore, $(a'_i,\,b'_i,\,c'_i,\,d'_i)$ fails to satisfy
Property \textbf{R2}$(2)$.  By Lemma \ref{lem:multiplesatsprops},
$\alpha_i\notin\conj\inv(\Gamma)$, for $i=1,2$.
By the above comments, the union
in \eqref{eqn:thisunion} ranges only over $\{\Xi_{12}\}$.
This completes the proof of \eqref{eqn:xi12equality} and therefore
of the first line of \eqref{eqn:twopartsdescribed}.

We have by
Lemma \ref{lem:subsetconstr} that
$\conj\inv(\Gamma)-\conj\inv(\Gamma)\cap\SL{2}{\Int[\omega_8]}$
is an $\Xi_{12}$-space, and the purpose of the second line of
\eqref{eqn:twopartsdescribed} is to give a precise
description of this $\Xi_{12}$-space.
Using the equality on the left in the second line of \eqref{eqn:twopartsdescribed}
we can restate the result of the second line as follows,
\startdisp\begin{aligned}
\text{The $\Xi_{12}$-space}&&&
\big(i,\delta\big)\inv(2,0)&&&\bigcup
\hspace*{-.3cm}\cdot\hspace*{.3cm}&\big(i,\delta\big)\inv(2,1)
&\textbf{(A)}\\
\text{equals}&&&
\bigcup_{\epsilon=0,1}\hspace*{-0.45cm}\cdot\hspace*{.45cm}
\frac{1}{1+\vecti}
\Xi_{2}\alpha^{2\vecti}\hspace*{-0.7mm}
(\vecti,\vecti^{\epsilon})
&&&\bigcup
\hspace*{-.3cm}\cdot\hspace*{.3cm}&
\bigcup_{\epsilon=0,1}
\hspace*{-0.45cm}\cdot\hspace*{.45cm}
\frac{1}{\omega_8(1+\vecti)}
\Xi_{2}\alpha^{-2}\hspace*{-0.7mm}
(-1,\vecti^{\epsilon})
.&\textbf{(B)}
\end{aligned}
\finishdisp
In order to verify the second line of \eqref{eqn:twopartsdescribed},
we must show that line \textbf{(A)}
equals line \textbf{(B)}.
Each of lines \textbf{(A)} and \textbf{(B)}
is expressed above as a disjoint union of two sets, making four sets
in all.
We claim that
\starteqn\label{eqn:bothsetsxi12space}\text{
Each of the two sets in the union in line \textbf{(A)},
resp. \textbf{(B)}
is a $\Xi_{12}$-space}.
\finisheqn
In order to verify \eqref{eqn:bothsetsxi12space},
for line \textbf{(A)}, note that
\starteqn\label{eqn:deltapreservation}\text{for $\alpha\in\conj\inv(\Gamma)$,
$\gamma\in\Xi_{12}$, we have}\;
\big(i,\delta\big)(\alpha)=
\big(i,\delta\big)(\gamma\alpha).
\finisheqn
In particular, \eqref{eqn:deltapreservation} implies that
\startdisp
\big(i,\delta\big)\inv(2,0)\;\text{and}\;
\big(i,\delta\big)\inv(2,1)\;\text{are $\Xi_{12}$-spaces,}
\finishdisp
completing the verification of the part of \eqref{eqn:bothsetsxi12space}
referring to line \textbf{(A)}.
For line \textbf{(B)}, it follows immediately from
\eqref{sublem:finerdecomp} that each of
the two sets appearing the union of line \textbf{(B)}
are $\Xi_{12}$-orbits.  \textit{A fortiori}, then
\startdisp
\bigcup_{\epsilon=0,1}\hspace*{-0.45cm}\cdot\hspace*{.45cm}
\frac{1}{1+\vecti}
\Xi_{2}\alpha^{2\vecti}\hspace*{-0.7mm}
(\vecti,\vecti^{\epsilon})\;\text{and}\;\bigcup_{\epsilon=0,1}
\hspace*{-0.45cm}\cdot\hspace*{.45cm}
\frac{1}{\omega_8(1+\vecti)}
\Xi_{2}\alpha^{-2}\hspace*{-0.7mm}
(-1,\vecti^{\epsilon})\;\text{are $\Xi_{12}$-spaces,}
\finishdisp
which completes the verification of the part of
\eqref{eqn:bothsetsxi12space} referring to line \textbf{(B)}.

According to \eqref{eqn:bothsetsxi12space},
there are a total of $\Xi_{12}$-spaces in lines \textbf{(A)}
and \textbf{(B)}.  To complete the proof, it remains to show that the each $\Xi_{12}$-space
in line \textbf{(A)} equals the $\Xi_{12}$-space directly below
it in line \textbf{(B)}.

Since both $\conj\inv(\Gamma)$ and
$\conj\inv(\Gamma)\cap\SL{2}{\Int[\omega_8]}$ are groups,
we have
\startdisp\text{
$\conj\inv(\Gamma)-\conj\inv(\Gamma)\cap
\SL{2}{\Int[\omega_8]}$ is a $\conj\inv(\Gamma)\cap\SL{2}{\Int[\omega_8]}$-
space,
}
\finishdisp
under the action of left-multiplication.  Therefore, for
each $\alpha\in\conj\inv(\Gamma)\cap
\SL{2}{\Int[\omega_8]}$, we may define the operator
\startdisp
\mbox{\bmth $\ell$\ubmth}(\alpha):=\text{left-multiplication by
$\alpha$ on $\conj\inv(\Gamma)-\conj\inv(\Gamma)\cap
\SL{2}{\Int[\omega_8]}$}.
\finishdisp
The first line of \eqref{eqn:twopartsdescribed}
implies that $\frac{1}{\omega_8}\left(\begin{smallmatrix}
\vecti&0\\
0&1\end{smallmatrix}\right)\in\conj\inv(\Gamma)\cap
\SL{2}{\Int[\omega_8]}$.
Our next claim is that,
\starteqn\label{eqn:operatorisomorphism}
\begin{aligned}&&&\text{for each of lines \textbf{(A)} and
\textbf{(B)}, the operator
$\mbox{\bmth $\ell$\ubmth}\left(\frac{1}{\omega_8}\begin{pmatrix}
\vecti&0\\
0&1\end{pmatrix}\right)$}\\
&&&\hspace*{1cm}\text{
provides an bijection of the first $\Xi_{12}$-space in the union
to the second.}
\end{aligned}
\finisheqn
Since $\mbox{\bmth $\ell$\ubmth}\left(\frac{1}{\omega_8}\left(\begin{smallmatrix}
\vecti&0\\
0&1\end{smallmatrix}\right)\right)$ obviously has an inverse,
namely $\mbox{\bmth $\ell$\ubmth}\left(\left(\frac{1}{\omega_8}\left(\begin{smallmatrix}
\vecti&0\\
0&1\end{smallmatrix}\right)\right)\inv\right)$,
the only issue in proving \eqref{eqn:operatorisomorphism}
is verifying that the operator $\mbox{\bmth $\ell$\ubmth}\left(\frac{1}{\omega_8}\left(\begin{smallmatrix}
\vecti&0\\
0&1\end{smallmatrix}\right)\right)$ does indeed
map the set on the left in each line into the set on the right.
With regards to line \textbf{(A)}, it is readily verified from
the definition of $\textbf{HYP}(i,\delta)$ that for any $\alpha\in
\SL{2}{\Complex}$
satisfying, say $\textbf{HYP}(i_{\alpha},\delta_{\alpha})$,
\startdisp
\big(i,\delta\big)\left(\mbox{\bmth $\ell$\ubmth}\left(\frac{1}{\omega_8}\begin{pmatrix}
\vecti&0\\
0&1\end{pmatrix}\right)\alpha\right)=(i_{\alpha},\red_2(\delta_{\alpha}+1)).
\finishdisp
Applying this observation to the situation in line \textbf{(A)}, we have
\startdisp
\mbox{\bmth $\ell$\ubmth}\left(\frac{1}{\omega_8}\begin{pmatrix}
\vecti&0\\
0&1\end{pmatrix}\right)\Big(\big(i,\delta\big)\inv(2,0)\Big)=
\big(i,\delta\big)\inv(2,1)
\finishdisp
With regards to line \textbf{(B)}, the key observation is that,
as follows from the description of the three $\Xi$-subsets
give in \eqref{eqn:xisetdescription}, we have the commutation relation
\startdisp
\begin{pmatrix}
\vecti&0\\
0&1\end{pmatrix}\Xi=\Xi\begin{pmatrix}
\vecti&0\\
0&1\end{pmatrix},\;\text{for}\; \Xi=\Xi_1,\,\Xi_2,\, \Xi_{12}.
\finishdisp
From the commutation relation and the relevant definitions, it is easily
calculated that
\startdisp
\mbox{\bmth $\ell$\ubmth}\left(\frac{1}{\omega_8}\begin{pmatrix}
\vecti&0\\
0&1\end{pmatrix}\right)\bigcup_{\epsilon=0,1}\hspace*{-0.45cm}\cdot\hspace*{.45cm}
\frac{1}{1+\vecti}
\Xi_{2}\alpha^{2\vecti}\hspace*{-0.7mm}
(\vecti,\vecti^{\epsilon})=\bigcup_{\epsilon=0,1}
\hspace*{-0.45cm}\cdot\hspace*{.45cm}
\frac{1}{\omega_8(1+\vecti)}
\Xi_{2}\alpha^{-2}\hspace*{-0.7mm}
(-1,\vecti^{\epsilon}).
\finishdisp
Thus \eqref{eqn:operatorisomorphism} is verified.  Because
of \eqref{eqn:operatorisomorphism}, the verification
of the second line of \eqref{eqn:twopartsdescribed} is reduced
to showing that the first set in line \textbf{(A)} equals the
first set in line \textbf{(B)}, \textit{i.e.} that
\starteqn\label{eqn:i2delta0description}
\big(i,\delta\big)\inv(2,0)=
\bigcup_{\epsilon=0,1}\hspace*{-0.45cm}\cdot\hspace*{.45cm}
\frac{1}{1+\vecti}
\Xi_{2}\alpha^{2\vecti}\hspace*{-0.7mm}
(\vecti,\vecti^{\epsilon}).
\finisheqn
Since the left-side is the union of $\Xi_{12}$-orbits of
$\frac{1}{1+\vecti}M_2^{2\vecti}$, we have, by a combination
of Proposition \ref{prop:heckedecomp} and
Sublemma \ref{sublem:finerdecomp} that
\starteqn\label{eqn:xi12orbitsunion}
\big(i,\delta\big)\inv(2,0)=
\bigcup_{\left\{\stackrel{m\in\Int[\vecti]| \; m|2\vecti,}
{\frac{2\vecti}{m}\;\text{standard}}
\right\}}
\hspace*{-1.07cm}\cdot\hspace*{1cm}
\bigcup_{\Xi}
\hspace*{-0.3cm}\cdot\hspace*{0.3cm}\frac{1}{1+\vecti}\Xi \alpha{^\mathrm{2\vecti}}(m,x).
\finisheqn
On the right side, $\Xi$ in the union ranges over a subset,
possibly empty, of $\{\Xi_{12},\,\Xi_1,\,\Xi_2\}$, depending
on $m,x$.  We only have to determine the subset for each
of the finitely many possibilities of $m,x$.  In order
to facilitate this, let us note first that if an element
of the right side of \eqref{eqn:xi12orbitsunion} is written
in the form $\frac{1}{1+\vecti}
\xi\alpha^{\mathrm{2\vecti}}(m,x)$, with $\xi\in\Xi$,
then we have
\starteqn\label{eqn:xialphaprimedquad}
\begin{pmatrix}a&b\\c&d\end{pmatrix}
=\xi\alpha^{\mathrm{2\vecti}}(m,x).
\finisheqn
Supposing that
\startdisp
\xi=\begin{pmatrix}\xi_{11}&\xi_{12}\\ \xi_{21}&\xi_{22}\end{pmatrix},
\finishdisp
we have
\startdisp
(a,b,c,d)=\frac{1}{1+\vecti}\left(m\xi_{11},x\xi_{12}+\frac{N}{m}\xi_{12},
m\xi_{12},x\xi_{21}+\frac{N}{m}\xi_{22}\right).
\finishdisp
In the present context, we have
\startdisp
\big(i,\delta\big)\big((a,b,c,d)\big)=(2,0).
\finishdisp
Therefore,
\starteqn\label{eqn:primedquadruple20case}
(a',\,b',\,c',\,d')=(1+\vecti)(a,b,c,d)=\left(m\xi_{11},x\xi_{11}+\frac{N}{m}
\xi_{12},m\xi_{21},x\xi_{21}+\frac{N}{m}\xi_{22}\right).
\finisheqn
The key requirement for
\startdisp
\alpha=\begin{pmatrix}a&b\\c&d\end{pmatrix}
\finishdisp
to be an element of $\conj\inv(\Gamma)$ is that $(a',\,b',\,c',\,d')$
satisfies property \textbf{R0}$(4)$.  That is we have
\starteqn\label{eqn:R04interpreted}
\alpha\in\conj\inv(\Gamma)\;\text{implies}\;x'\equiv 1\mod 1+\vecti,
\finisheqn
for any permutation of $(a',\,b',\,c',\,d')$.
The reason for \eqref{eqn:R04interpreted} is that
if $\alpha\in\conj\inv(\Gamma)$
then Lemma \ref{lem:multiplesatsprops}, Part (a), applies
with $n=2$.  In particular,
we obtain $a',c'\equiv 1\mod 1+\vecti$.
By \eqref{eqn:primedquadruple20case}, therefore, if
$\alpha\in\conj\inv(\Gamma)$,
$1+\vecti$ does not divide $m$.  Since $m|2\vecti$, with $\frac{N}{m}$ standard,
there is only one choice for $m$ such that $\alpha$ can be in
$\conj\inv(\Gamma)$, namely $m=\vecti$.  Under the assumption
that $m=\vecti$, \eqref{eqn:R04interpreted} becomes
\starteqn\label{eqn:primedquadxiprimes}
(a',\,b',\,c',\,d')=(\vecti\xi_{11},x\xi_{11}+2\xi_{12},
\vecti\xi_{21},x\xi_{21}+2\xi_{22})
\finisheqn
From \eqref{eqn:R04interpreted}, we deduce that, in particular
$(b',\,d')\equiv (1,1)\mod 1+\vecti$.  From \eqref{eqn:primedquadxiprimes}
we therefore deduce that
\startdisp
x\xi_{i1}\equiv 1\mod 1+\vecti,\;\text{for}\; i=1,2.
\finishdisp
Thus, $x\equiv 1\mod 1+\vecti$ and
$\xi_{i1}\equiv 1\mod 1+\vecti$.  Since $x\in\Omega_{\frac{N}{m}}$
and
\startdisp
\Omega_2=\{0,1,\vecti,1+\vecti\},
\finishdisp
we obtain from $x\equiv 1\mod 1+\vecti$
that $x=\vecti^{\epsilon}$ for $\epsilon=0$ or $\epsilon=1$.
We obtain from $\xi_{i1}\equiv 1\mod 1+\vecti$
and the definition of the three $\Xi$ sets given in
\eqref{eqn:residuematrices} and \eqref{eqn:xi12defn}
that $\xi\in\Xi_2$.  Therefore,
we obtain
\starteqn\label{eqn:i2delta0descriptionalmost}
\big(i,\delta\big)(\alpha)\inv(2,0)=
\bigcup_{\epsilon}\hspace*{-0.3cm}\cdot\hspace*{.3cm}
\frac{1}{1+\vecti}
\Xi_{2}\alpha^{2\vecti}\hspace*{-0.7mm}
(\vecti,\vecti^{\epsilon}),
\finisheqn
where the $\epsilon$ in the union ranges over a subset of $\{0,1\}$.
In order to show that the $\epsilon$ actually ranges over
the entire subset $\{0,1\}$, it will suffice to choose a single
element $\xi\in\Xi_2$ and show that
\starteqn\label{eqn:directverification}
\conj\left(\frac{1}{1+\vecti}
\xi\alpha^{2\vecti}\hspace*{-0.7mm}(\vecti,\vecti^{\epsilon})\right)\in
\mathrm{Mat}_3(\Int[\vecti]),\;\text{for each $\epsilon\in\{0,1\}$}.
\finisheqn
We choose the particular element $\xi\in\Xi_2$ by setting
\startdisp
\xi=\begin{pmatrix}1&0\\
1&1\end{pmatrix}.
\finishdisp
Calculating $\frac{1}{1+\vecti}\xi\alpha^{2\vecti}\hspace*{-0.7mm}
(\vecti,\vecti^{\epsilon})$ directly, we obtain
\startdisp
\frac{1}{1+\vecti}\xi\alpha^{2\vecti}\hspace*{-0.7mm}
(\vecti,\vecti^{\epsilon})=\frac{1}{1+\vecti}\begin{pmatrix}
\vecti&\vecti^{\epsilon}\\ \vecti&2+\vecti^{\epsilon}
\end{pmatrix}.
\finishdisp
Using \eqref{eqn:imagematrixconj}, we calculate
\startdisp
\conj\left(\frac{1}{1+\vecti}
\xi\alpha^{2\vecti}\hspace*{-0.7mm}(\vecti,\vecti^{\epsilon})\right)=
\begin{pmatrix}2+2\vecti^{\epsilon}&2\vecti+2\vecti^{1+\epsilon}&
\vecti(2+\vecti^{\epsilon})-\vecti^{1+\epsilon}\\
3\vecti+(-1)^{\epsilon}\vecti+\vecti^{1+\epsilon}&
1+(-1)^{\epsilon}+2\vecti^{\epsilon}&
\vecti^{1+\epsilon}+\vecti(2+\vecti^{\epsilon})\\
1+\vecti^{\epsilon}(2+\vecti^{\epsilon})&
-\vecti+\vecti^{1+\epsilon}(2+\vecti^{\epsilon})&
\vecti(2+\vecti^{\epsilon})+\vecti^{1+\epsilon}
&\end{pmatrix},
\finishdisp
the right side of which belongs to $\mathrm{Mat}_3(\Int[\vecti])$,
for $\epsilon\in\{0,1\}$.
Therefore
in \eqref{eqn:i2delta0descriptionalmost}, $\epsilon$ ranges over
the entire set $\{0,1\}$.
So,
\startdisp
\big(i,\delta\big)\inv(2,0)=
\bigcup_{\epsilon=0,1}\hspace*{-0.45cm}\cdot\hspace*{.45cm}
\frac{1}{1+\vecti}
\Xi_{2}\alpha^{2\vecti}\hspace*{-0.7mm}
(\vecti,\vecti^{\epsilon}),
\finishdisp
which is \eqref{eqn:i2delta0description}.  By the comments
preceding \eqref{eqn:i2delta0description}, this completes
the proof of the second line of \eqref{eqn:twopartsdescribed},
and therefore of the Lemma.
\end{myproof}

To complete the proof of \ref{prop:inversematrixexplicitdescription},
write
\startdisp
\begin{aligned}
\conj\inv(\Gamma)&=&\conj\inv(\gamma)\cap\SL{2}{\Int[\vecti]}\bigcup
\hspace*{-.3cm}\cdot\hspace*{.3cm}\conj\inv(\Gamma)-\conj\inv(\Gamma)
\cap\SL{2}{\Int[\vecti]}
\end{aligned}
\finishdisp


\section{Explicit determination of the fundamental domain
for the action of $\SO{3}{\Int[\vecti]}$ on $\bbH^3$}
\label{sec:explicitfd}
We begin with the following definition, which is fundamental to everything
that follows.
\vspace*{0.3cm}

\noindent\textbf{Definition.}\hspace*{0.2cm}  Let $X$ be a topological space.
Suppose that $\Gamma$  is a group acting topologically on $X$,
\textit{i.e.}, $\Gamma\subseteq\Iso(X)$.
A subset $\scrF$ of $X$ is called an \textbf{exact fundamental
domain for the action of \boldmath $\Gamma$ on $X$\unboldmath} if the following
conditions are satisfied
\begin{itemize}
\item[\textbf{FD 1.}]  The $\Gamma$-translates of $\scrF$ cover $X$, \textit{i.e.},
\startdisp
X=\Gamma \scrF.
\finishdisp
\item[\textbf{FD 2.}]  Distinct $\Gamma$-translates of $\scrF$ intersect only
on their boundaries, \textit{i.e.},
\startdisp
\gamma_1,\gamma_2\in\Gamma,\, \gamma_1\neq\gamma_2\;\text{implies}\;
\gamma_1\scrF\cap
\gamma_2\scrF\subseteq\gamma_1\partial\scrF,\,\gamma_2\partial\scrF.
\finishdisp
\end{itemize}
Henceforth, we will drop the word \textbf{exact} and refer to such an $\scrF$
simply as a \textbf{fundamental domain}.

\subsection{The Grenier fundamental domain of a discrete
subgroup $\Gamma$ of $\mathrm{Aut}^+(\bbH^3)$}\label{subsec:grenierh3}

For the current section, \S\ref{subsec:grenierh3}, only,
$G$, instead of denoting $\SL{2}{\Complex}$, will denote $\SL{2}{\Complex}$.
Likewise, instead of denoting $\SO{3}{\Int[\vecti]}$ or
$\conj\inv\SO{3}{\Int[\vecti]}$,
$\Gamma$ will denote an arbitrary subgroup of $\SL{2}{\Complex}$,
satisfying certain conditions to be given below.  The main
examples to keep in mind are, first, $\Gamma=\SL{2}{\Int}$ the integer
subgroup of $\SL{2}{\Complex}$ and, second,
$\Gamma=\conj\inv(\SO{3}{\Int[\vecti]})$,
the inverse image of the integer subgroup of $\SO{3}{\Complex}$,
described explicitly as a group of fractional
linear transformations in Proposition
\ref{prop:inversematrixexplicitdescription}.

The main results of the present section, numbered
Theorem \ref{thm:grenierh3} and Theorem \ref{thm:grenierh3extended} below,
amount to an application of a general
result which is valid in a much wider context.  This wider context is
that of the integer subgroup $\mathbf{G}(\Int)$ of a
Chevalley group $\mathbf{G}$ acting on
the symmetric space $\mathbf{G}(\Complex)/K$.  In Chapter I of
\cite{brennerthesis} two general results, namely
Theorems I.2.6 and I.2.8 have already been given and these
theorems are general enough to imply Theorems \ref{thm:grenierh3}
and and \ref{thm:grenierh3extended}, below.
In lieu of a proof of Theorem \ref{thm:grenierh3}, we will
merely indicate how the relevant result of \cite{brennerthesis}
applies to the situation at hand to give Theorem \ref{thm:grenierh3}.

Theorems \ref{thm:grenierh3} and \ref{thm:grenierh3extended}
are introduced here in order to be applied to the case
$\Gamma=\conj\inv(\SO{3}{\Int[\vecti]})$, and so give Theorem
\ref{thm:gammagoodgrenier}.
A more complete exposition of the
theory of fundamental domains
in the context of Chevalley groups will appear in \cite{chevfd},
The results given in \cite{chevfd} will immediately imply Theorem
\ref{thm:gammagoodgrenier},
which will obviate the need of stating Theorems \ref{thm:grenierh3}
and \ref{thm:grenierh3extended} as an intermediate step in deducing
Theorem \ref{thm:gammagoodgrenier}.

\vspace*{0.3cm}
\noindent\bmth\textbf{Iwasawa decomposition of $\SL{2}{\Complex}$}.\ubmth
\hspace*{0.3cm}For the reader's convenience, we recall only those results
in the context of $\SL{2}{\Complex}$ which we need to proceed.
For proofs and the statements for $\SL{n}{\Complex}$, see
the ``Notation and Terminology" section of \cite{jol06}.  Let
\startdisp
\begin{aligned}
U&=&&\text{upper triangular unipotent matrices in $\SL{2}{\Complex}$, so
$U=\left\{\left.\begin{pmatrix}1&x\\0&1\end{pmatrix}\;\right|\;x\in\Complex
\right\}$},\\
A&=&&\text{diagonal elements of $\SL{2}{\Complex}$ with positive
diagonal entries, so $A=\left\{\left.\begin{pmatrix}y&0\\0&y\inv\end{pmatrix}
\;\right|\;y\in\Real_{+}\right\}$},\\
K&=&&\text{$\mathrm{SU}(2)$, so $K=\{k\in\SL{2}{\Complex}\;|\;kk^*=1\}$}.
\end{aligned}
\finishdisp
Here $x^*$ denotes the conjugate-transpose $\overline{x}^t$ of $x$.

\textit{We have the \textbf{Iwasawa decomposition}}
\startdisp
\SL{2}{\Complex}=UAK,
\finishdisp
\textit{and the product map $U\times A\times K\rightarrow UAK$ is a differential
isomorphism.}

The Iwasawa decomposition induces a system of coordinates $\phi$
on the symmetric space $\SL{2}{\Complex}/K$.  The mapping $\phi$
is a  diffeomorphism
between $\SL{2}{\Complex}/K$ and $\Real^3$.  The details are as follows.
The Iwasawa decomposition gives a uniquely determined
product decomposition of $gK\in\SL{2}{\Complex}/K$ as
\startdisp
gK=u(g)a(g)K,\;\text{where}\, u(g)\in U,\, a(g)\in A\;
\text{are uniquely determined by $gK$}
\finishdisp
Define the \textbf{Iwasawa coordinates} $x_{1}(g)$, $x_{2}(g)\in\Real$,
$y(g)\in\Real^+$ by the relations
\startdisp
u(g)=\begin{pmatrix}1&x_{1}(g)+\vecti x_{2}(g)\\
0&1\end{pmatrix}\,\quad a(g)=\begin{pmatrix}y(g)^{\half}&0\\
0&y(g)^{-\half}\end{pmatrix}.
\finishdisp
By the Iwasawa decomposition, the Iwasawa coordinates of $g$
are uniquely determined.  We emphasize that while $x_{1}(g)$ and $x_{2}(g)$
range over all the real numbers, $y(g)$ ranges over the positive
numbers.   As functions on $G$, $x_1$ $x_2$, and $y$ are invariant
under right-multiplication by $K$.  Thus $x_1$, $x_2$, and $y$
induce coordinates on $G/K$.
Now define the coordinate mappings
$\phi_i:\SL{2}{\Complex}/K\rightarrow \Real$, for $i=1,2,3$, by
\starteqn\label{eqn:phicoordsdefn}
\phi_1=-\log y,\;\phi_2=x_{1},\;\phi_3=x_{2},
\finisheqn
and set
\startdisp
\phi=(\phi_1,\phi_2,\phi_3): G/K\rightarrow\Real^3.
\finishdisp
The mapping $\phi$ is a diffeomorphism of $G/K$ onto $\Real^3$,
because the Iwasawa coordinate system is a diffeomorphism, as is $\log$.
Thus, there exists the inverse diffeomorphism
\startdisp
\phi\inv: \Real^3\rightarrow G/K.
\finishdisp
By \eqref{eqn:phicoordsdefn}, we can write, explicitly,
\starteqn\label{eqn:phiinverseexplicit}
\phi\inv(t_1,t_2,t_3)=t_2+t_3\vecti+e^{-t_1}\vectj,\forallindisp\;
t=(t_1,t_2,t_3)\in\Real^3.
\finisheqn

\vspace*{0.3cm}
\noindent\bmth\textbf{The quaternion model and the coordinate
system on $\SL{2}{\Complex}/K$}.\ubmth \hspace*{0.3cm}
We will use the model $G/K$ as the upper half-space $\bbH^3$,
defined as the following subset of the quaternions.
\starteqn\label{eqn:quatmodeldefn}
\bbH^3=\{x_{1}+x_{2}\vecti+y\vectj,
\;\text{where}\; x_{1}, \,x_{2}\in\Real,\; y\in\Real^+\}.
\finisheqn
Recall that $\SL{2}{\Complex}$ acts transitively on $\bbH^3$ by fractional
linear transformations.  See \S\Roman{ranrom}.0 of \cite{jol05} for
the details
of the action.  We note the relation
\starteqn\label{eqn:fraclintrans}
g\vectj=x_1(g)+x_2(g)\vecti+y(g)\vectj.
\finisheqn
As a result of \eqref{eqn:fraclintrans} and the Iwasawa decomposition,
we may identify $\SL{2}{\Complex}/K$ with $\bbH^3$.  So $\phi:G/K\rightarrow
\Real^3$ induces a diffeomorphism
\startdisp
\phi:\bbH^3\stackrel{\cong}{\longrightarrow}\Real^3.
\finishdisp
Because of \eqref{eqn:fraclintrans}, if $g$ is any element of $G$
such that $g\vectj=z$, then $\phi(g)=\phi(z)$.  Further, beause
of the way we set up the coordinates on $\bbH^3$,
$\phi: \bbH^3\rightarrow\Real^3$ is given explicitly by the same
formulas as \eqref{eqn:phicoordsdefn}.

As explained in, for example, \S\Roman{ranrom}.0 of \cite{jol05},
the kernel of the action of $\SL{2}{\Complex}$ on $\bbH^3$
is precisely the set $\{\pm I\}$, consisting of the identity
matrix and its negative.

For any oriented manifold $X$ equipped with a metric, use the notation
\startdisp
\mathrm{Aut}^+(X)=\;\text{group of
orientation-preserving isometric automorphisms
of $X$.}
\finishdisp
It is a fact that every element of $\mathrm{Aut}^+(X)$ is realized
by a fractional linear transformation in $\SL{2}{\Complex}$,
unique up to multiplication by $\pm 1$.  Therefore,
the action of $\SL{2}{\Complex}$ on $\bbH^3$
by fractional linear transformations induces an isomorphism
\starteqn\label{eqn:autometriesiso}
\SL{2}{\Complex}/\{\pm I\}\cong \mathrm{Aut}^+(\bbH^3).
\finisheqn

\vspace*{0.3cm}
\noindent\bmth\textbf{The stabilizer in $\Gamma$
of the first $j$ $\phi$-coordinates.}
\ubmth\hspace*{0.3cm}
In all that follows, if $i,j\in\Natural$, the notation $[i,j]$ is used to
denote the interval of \textit{integers} from $i$ to $j$,
inclusive.  The interval
$[i,j]$ is defined to be the empty set if $i>j$.
\begin{defn} \label{defn:projections} For $i,j\in\{1,2,3\}$, with $i\leq j$,
let \bmth$\phi_{[i,j]}$
be the \textbf{projection of $\bbH^3$ onto the $[i,j]$
factors of $\Real^3$}.  In other words,
we let
\startdisp
\phi_{[i,j]}=(\phi_i,\phi_{i+1},\ldots,\phi_j).
\finishdisp
\end{defn}
Since $\phi$ is a diffeomorphism of $\bbH^3$, $\phi_{[i,j]}$
is an smooth epimorphism of $\bbH^3$ onto $\Real^{i-j+1}$.

If $\scrK$ is any subset of $\{1,2,3\}$, of size $|\scrK|$, then we
can generalize in the obvious way to define the smooth epimorphism
\startdisp
\phi_{\scrK}: \bbH^3\rightarrow\Real^{|\scrK|}.
\finishdisp

Let $\Gamma$ be a group acting by diffeomorphisms
of $\bbH^3$.  For $\gamma\in\Gamma$ we also use $\gamma$ to denote the
diffeomorphism of $\bbH^3$ defined by the left action of $\gamma$
on $\bbH^3$.
Therefore, for $l\in\{1,\ldots 3\}$ the composition
$\phi_l\circ\gamma$ is the $\Real$-valued function on $\bbH^3$ defined by
\startdisp
\phi_l\circ\gamma(z)=\phi_l(\gamma z)\quad\text{for all}\; z\in \bbH^3.
\finishdisp

We use $\Gamma^{\phi_{[1,j]}}$ to denote the subgroup of $\Gamma$
whose action stabilizes the first $i$ coordinates.  In other words,
we set
\startdisp
\Gamma^{\phi_{[1,j]}}=\{\gamma\in\Gamma\;|\; \phi_{[1,j]}
=\phi_{[1,j]}\circ\gamma\}.
\finishdisp
We extend the definition of $\Gamma^{\phi_{[1,j]}}$ to $j=0,4$,
by adopting the conventions
\startdisp
\Gamma^{\phi_{[1,0]}}=\Gamma,\quad\text{and}\quad\Gamma^{\phi_{[1,4]}}=1.
\finishdisp
Note that, by definition, we have the descending sequence of
groups
\startdisp
\Gamma=\Gamma^{\phi_{[1,0]}}\geq
\Gamma^{\phi_1}\geq \Gamma^{\phi_{[1,2]}}\geq
 \Gamma^{\phi_{[1,3]}}\geq\Gamma^{\phi_{[1,4]}}=1.
\finishdisp
Note that the penultimate group
in this sequence, namely $\Gamma^{\phi_{[1,3]}}$,
equals, by definition, the kernel of the action of $\Gamma$ on $\bbH^3$.
Assuming that $\Gamma\subset\SL{2}{\Complex}$, \textit{i.e.} that
$\Gamma$ consists of fractional linear transformations, we always
have
\starteqn\label{eqn:kernelisverysmall}
\Gamma^{\phi_{[1,3]}}=\Gamma\cap\{\pm 1\}.
\finisheqn

Because the $\Gamma^{\phi_{[1,j]}}$ form a descending sequence, for
$k,j\in\{1,2,3\}$ with $k<j$,
we can consider
the left cosets of $\Gamma^{\phi_{[1,k]}}$ in $\Gamma^{\phi_{[1,j]}}$.
The left cosets are the sets of the form
$\Gamma^{\phi_{[1,j]}}\gamma_k$ for $\gamma_k\in\Gamma^{\phi_{[1,k]}}$.
Now let $i,j,k\in\{1,2,3\}$, $l\leq j$, $k<j$.
By the definition of $\Gamma^{\phi_{[1,j]}}$,
the function $\phi_l\circ\gamma_k$ depends only only on the left
$\Gamma^{\phi_{[1,j]}}$-coset to which $\gamma_k$ belongs.
Therefore, for fixed $z$
we may consider $\phi_l\circ\gamma_k(z)$ to be a well-defined function on
the set of left cosets $\Gamma^{\phi_{[1,j]}}\gamma_k$ of
$\Gamma^{\phi_{[1,k]}}$ in $\Gamma^{\phi_{[1,j]}}$.  We may
therefore, speak of the $\Real$-valued function
$\phi_l\circ\Gamma^{\phi_{[1,j]}}\gamma_k$.

In what follows we
will most often apply the immediately preceding paragraph when
$l=j$, and $k=j-1$.  For $\gamma\in\Gamma^{\phi_{[1,j-1]}}$
and $\Delta$ an arbitrary subset of $\Gamma^{\phi_{[1,j]}}$,
we have
\starteqn\label{eqn:Deltaonleftinphij}
\phi_j(\Delta\gamma z)=\{\phi_j(\gamma z)\}.
\finisheqn
therefore, by setting
\startdisp
\phi_j\circ\Gamma^{[1,j]}\gamma(z)=\phi_j(\gamma z),
\finishdisp
we obtain a well-defined function
\startdisp
\phi_j\circ\Gamma^{\phi_{[1,j]}}\gamma:\bbH^3\rightarrow\Real.
\finishdisp
The function $\phi_j\circ\Gamma^{\phi_{[1,j]}}\gamma$
depends only on the $\Gamma^{\phi_{[1,j]}}$-coset to which
$\gamma$ belongs.

For $\gamma\in\Gamma^{\phi_{[1,j-1]}}$, the $\Real$-valued function
$\phi_j\circ\Gamma^{\phi_{[1,j]}}\gamma$
gives the effect of the action of $\gamma\in\Gamma^{\phi_{[1,j-1]}}$
on the $j^{\rm th}$ coordinate of a point.  It
is clear from the definition that
\starteqn\label{eqn:identitycosetdefn}\text{
$\phi_j=\phi_j\circ\gamma$ if and only if $\Gamma^{\phi_{[1,j]}}\gamma$ is the identity left
coset of $\Gamma^{\phi_{[1,j]}}$ in $\Gamma^{\phi_{[1,j-1]}}$}.
\finisheqn

We now define the \bmth\textbf{difference function
$\Delta_{j,\gamma}$ associated to $\gamma\in\Gamma^{\phi_{[1,j-1]}}$},\ubmth
by setting
\starteqn\label{eqn:differencefunctiondefn}
\Delta_{j,\gamma}=\phi_j\phi\inv-
\phi_j\circ\gamma\phi\inv:\Real^3\rightarrow\Real.
\finisheqn
Then from \eqref{eqn:identitycosetdefn}, we deduce that
\startdisp\text{
$\Delta_{j,\gamma}$ is the (constant) 0-function
if and only if $\gamma\in\Gamma^{j}$.
}\finishdisp
Let $\mathbf{t}\in\Real^{3}$ with $z=\phi\inv(\mathbf{t})$ the corresponding
point in $\bbH^3$.  Then $\Delta_{j,\gamma}(\mathbf{t})$
measures the displacement
in the $j^{th}$ coordinate induced at $z$ by the action of $\gamma$
on $\bbH^3$.

We will in particular apply Definition \ref{defn:projections}
to the case when
$\Gamma$ is a discrete subgroup of $\SL{2}{\Complex}$.  Then,
by the isomorphism \eqref{eqn:autometriesiso},
$\Gamma/\{\pm I\}\subseteq \mathrm{Aut}^+(\bbH^3)$.  Thus,
$\Gamma/\{I\}$ is a (discrete)
subgroup of the group of metric automorphisms of $\bbH^3$.
So in particular, $\Gamma$ acts by diffeomorphisms
of $\bbH^3$.

It is an immediate consequence of the definitions that
for any group $\tilde{\Gamma}$ acting on $\bbH^3$
by diffeomorphisms, and any subgroup $\Gamma$
of $\tilde{\Gamma}$, we have, for $1\leq i\leq j\leq 3$,
\starteqn\label{eqn:gammastabilizerintersection}
\Gamma^{\phi_{[i,j]}}=(\tilde{\Gamma})^{\phi_{[i,j]}}\cap\Gamma.
\finisheqn
Applying \eqref{eqn:gammastabilizerintersection} to the case
of $\tilde{\Gamma}=\SL{2}{\Complex}$ and $i=1$, we deduce that
\starteqn\label{eqn:stabilizerintersection}
\Gamma^{\phi_{[1,j]}}=\Gamma\cap\SL{2}{\Complex}^{\phi_{[1,j]}},
\finisheqn
for any subgroup $\Gamma\subseteq\SL{2}{\Complex}$.
Because of \eqref{eqn:stabilizerintersection} it is very useful
to have an explicit expression for
$\SL{2}{\Complex}^{\phi_1}$.  We carry out the calculation
using the relations of \eqref{eqn:phicoordsdefn}.

Let $z\in\bbH^3$ with
\startdisp
z=x_1+x_2+y\vectj,
\finishdisp
as in \eqref{eqn:quatmodeldefn}.  Let
\startdisp
g\in\SL{2}{\Complex}\;\text{with}\;
g=\begin{pmatrix}a&b\\c&d\end{pmatrix}.
\finishdisp
Define
\starteqn\label{eqn:ycdzdefn}
y(c,d;z)=
\frac{y(z)}{||c z+d||^2},
\finisheqn
where in \eqref{eqn:ycdzdefn} and from now on,
for a quaternion $z$,
$||z||^2$ denotes the squared norm of a $z$, so that
$||z||^2=z\overline{z}$.
Then we have
\starteqn\label{eqn:ycdzrelation}
y(gz)=y(c,d;z).
\finisheqn
For the details of such calculations, see
\S\Roman{ranrom}.0 of \cite{jol05}.  Since
\startdisp
\phi_1:\bbH^3\rightarrow\Real\;\text{is defined
as}\;-\log y(\cdot),
\finishdisp
and $\log$ is injective, \eqref{eqn:ycdzrelation} implies
that
\starteqn\label{eqn:ginstabilizerprelim}\text{
$g\in\SL{2}{\Complex}^{\phi_1}$ if and only if
$y(c,d;z)=y(z)$ for all $z\in\bbH^3$.}
\finisheqn

By \eqref{eqn:ginstabilizerprelim} and \eqref{eqn:ycdzdefn}, we have
\starteqn\label{eqn:ginstabilizerprelim2}
g\in\SL{2}{\Complex}^{\phi_1}\;\text{if and only if}\;
||cz+d||^2=1\forallindisp z\in\bbH^3.
\finisheqn
Clearly, the condition $||cz+d||^2=1$ is satisfied
for all $z\in\bbH^3$ if and only if $c=0$ and $||d||=1$.  We therefore
deduce from \eqref{eqn:ginstabilizerprelim2} that
\starteqn\label{eqn:firstcoordfixing}
\SL{2}{\Complex}^{\phi_{1}}=\left\{\left.
\begin{pmatrix}\omega\inv&x\\0&\omega\end{pmatrix}\;\right|\;
x,\,\omega\in\Complex,
\,||\omega||=1\right\}.
\finisheqn

\vspace*{0.3cm}
\noindent\bmth\textbf{Axioms for the action of $\Gamma$.}
\ubmth\hspace*{0.3cm}
As before, suppose that $\Gamma$ is a group acting by diffeomorphisms
on $\bbH^3$, and let $\Gamma^{\phi_{[1,j]}}$ for $j\in\{1,2,3\}$
be defined
as above.  We will shortly state the four `\textbf{A}' axioms for the
action of $\Gamma$ on $\bbH^3$.  Before stating the axioms
we introduce a few pieces of terminology and make some easy
observations based on them.  For any subset $\scrK$ of
the interval of integers $[1,3]$, we let $\scrK^c=[1,3]-\scrK$
be the \textit{complement of $\scrK$ \textit in $[1,3]$}.

\begin{defn} \label{defn:independentof} Let $f$ be a real-valued function
\startdisp
f: \bbH^3\rightarrow \Real.
\finishdisp
Let $\scrK$ a subset of $[1,3]$.  We say that $f$ is \bmth\textbf{independent
of the $\scrK$ coordinates\hspace*{0.15cm}}\ubmth if for every $x,y\in\bbH^3$,
\startdisp
\phi_{\scrK^c}(x)=\phi_{\scrK^c}(y)\;\text{implies}
\;f(x)=f(y).
\finishdisp
\end{defn}

In other words, $f$ is independent of the coordinates in $\scrK$
if and only if $f$ is constant on the fibers of the projection
$\phi_{\scrK^c}$ onto the $\Real$-factors indexed by $\scrK^c$.
We will most often apply Definition \ref{defn:independentof} when
$\scrK$ meets one of the following two descriptions.
\startdisp
(1)\quad\scrK=\{1,\ldots,
i-1\}\qquad\text{or},\qquad(2)\quad\scrK=[i,j]^c.
\finishdisp
When Definition \ref{defn:independentof} applies in
case (1), we will say that $f$ is \bmth\textbf{independent of the
first $i-1$ coordinates}\ubmth.  When Definition \ref{defn:independentof}
applies in case (2), we will say that $f$ \bmth\textbf{depends only
on the $i^{\mathrm{th}}$ through $j^{\mathrm{th}}$ coordinates}\ubmth.
These definitions can be extended in the obvious
way from $\Real$-valued functions on $\bbH^3$ to functions on $\bbH^3$
taking values in any given set (for example $\Real^k$-valued functions).

For the next observation, we need to introduce the notion
of a section of a projection $\phi_{\scrK}$.  It will not really
matter which section we use, so for simplicity, we choose the zero section.
For a subinterval $[i,j]$ of $\{1,2,3\}$ of size $j-i+1$, define
\startdisp
\sigma_{[i,j]}^0: \Real^{j-1+1}\rightarrow\bbH^3
\finishdisp
by
\startdisp
\sigma_{[i,j]}^0(x_1,\ldots,x_{j-i+1})=(\underbrace{0,\ldots,0}_{i-1},
x_1,\ldots, x_{j-i+1},\underbrace{0,\ldots,0}_{3-j}).
\finishdisp
The map $\sigma_{[i,j]}^0$ is called the \bmth
\textbf{zero section of the projection $\phi_{[i,j]}$}\ubmth.
The terminology comes from the relation
\starteqn\label{eqn:sectionrelation}
\phi_{[i,j]}\sigma_{[i,j]}^0=\Id_{\Real^{i-j+1}},
\finisheqn
which is immediately verified.
The concept of the zero section of the projection can
 be generalized from the case of a projection
associated with an interval $[i,j]$ to that
of an arbitrary subset $\scrK$ of $\{1,2,3\}$, in the obvious way,
although we will not have any use for this generalization in the present
context.

By use of the zero section, we are able to
make a useful reformulation of the condition that $f:\bbH^3\rightarrow\Real$
is independent of the first $j-1$ coordinates.  Let $j\in\{2,3\}$
and $f$ a real values function on $\bbH^3$.  Then
\starteqn\label{eqn:independencereformulation}\text{
$f$ is independent of the first $j-1$ coordinates if and only if
$f\,\sigma^0_{[j,3]}\phi_{[j,3]}=
\sigma^0_{[j,3]}\phi_{[j,3]}\,f$.
}
\finisheqn

The reformulation \eqref{eqn:independencereformulation} allows
us to prove the following result.
\begin{lem}\label{lem:action1cons} Let $\Delta$ be a group
acting on $\bbH^3$, and for $j\in\{1,2,3\}$, let $\phi_{[j,3]}$
be the projection of $\bbH^3$ onto the last $3-j+1$-coordinates
and let $\sigma^0_{[j,3]}$ be the zero section of $\phi_{[j,3]}$.
Suppose that, for all $l\in[j,3]$ and $\delta\in\Delta$, the
functions $\phi_l\circ\delta$ are independent of the first
$j-1$ coordinates.  Then $\Delta$ has an induced action on $\Real^{3-j+1}$
defined by \starteqn\label{eqn:sectioninducedact}
\delta_{[j,3]}(\mathbf{t})=\phi_{[j,3]}(\delta\sigma^0_{[j,3]}(\mathbf{t}))
,\;\forallindisp
\;\mathbf{t}=(t_1,\ldots,t_{3-j+1})\in\Real^{3-j+1}. \finisheqn
\end{lem}
\begin{myproof}{Proof}
With $\delta_{[j,3]}$ defined as in \eqref{eqn:sectioninducedact},
we verify the relation
\starteqn
(\delta_1)_{[j,3]}(\delta_2)_{[j,3]}=(\delta_1\delta_2)_{[j,3]},
\forallindisp\; \delta_1,\,\delta_2\in\Delta
\finisheqn
In order to minimize the clutter, we drop all subscripts
$[j,3]$ from the $\phi$'s and $\sigma$'s and all superscripts
$0$ from the $\sigma$'s in the intermediate
steps of the calculation.  Applying,
successively, \eqref{eqn:sectioninducedact},
\eqref{eqn:independencereformulation}, \eqref{eqn:sectionrelation},
and \eqref{eqn:sectioninducedact}, we have
\startdisp
(\delta_1)_{[j,3]}(\delta_2)_{[j,3]}=
(\phi\delta_1\sigma)(\phi\delta_2\sigma)=
\phi(\delta_1\sigma\phi)\delta_2\sigma=
\phi(\sigma\phi\delta_1)\delta_2\sigma=
(\phi\sigma)(\phi(\delta_1\delta_2)\sigma)=
(\delta_1\delta_2)_{[j,3]}.
\finishdisp
This completes the proof of the Lemma.
\end{myproof}

We are now ready to state the \textbf{A} axioms.

In axiom \textbf{A 4}, we denote a direction
(unit vector) in $\Real^3$ by $\vectu$.  We denote the directional
derivative of a function $f:\Real^3\rightarrow\Real$ in the direction
$\vectu$ by $\mathrm{D}_{\vectu}f$.

\begin{itemize}
\item[\textbf{A 1.}]  For $j\in\{2,3\}$, $\gamma\in\Gamma^{\phi_{[1,j-1]}}$,
$\phi_j\circ\gamma$ is independent of the first $j-1$ coordinates.
\item[\textbf{A 2.}]  Let $\gamma\in\Gamma-\Gamma^{\phi_{1}}$.
Then the difference function $\Delta_{1,\gamma}$
has no critical zeros, \textit{i.e.}, no zeros which are also
critical points.
\item[\textbf{A 3.}]  Let $z\in\bbH^3$, $B\in\Real$ be given.
Then there is a neighborhood
$U$ of $z$ in $\bbH^3$ with the following property.
\starteqn\label{eqn:A3Ucond}\begin{gathered}\text{
There
are only finitely many left cosets
$\Gamma^{\phi_{1}}\gamma$ of $\Gamma^{\phi_{1}}$
in $\Gamma$}\\
\text{such that $\phi_1(\Gamma^{\phi_{1}}\gamma U)\cap(-\infty,B]$
is nonempty.}
\end{gathered}
\finisheqn
\item[\textbf{A 4.}]  Let $z\in\bbH^3$ be fixed.
Then there is a direction $\vectu$ depending only on $z$ such that we have
\starteqn\label{eqn:nocritpoints}
D_{\vectu}\Delta_{1,\gamma}(z)<0,\;\text{for every}\;
\gamma\in\Gamma-\Gamma^{\phi_{1}},
\finisheqn
\end{itemize}

Note that
\starteqn\label{eqn:a4impliesa2}\text{
Axiom \textbf{A 4} implies Axiom \textbf{A 2}.
}\finisheqn
The reason for \eqref{eqn:a4impliesa2} is that
\textbf{A 4} implies that $\Delta_{1,\gamma}$
has no critical points, so \textit{a fortiori}, no critical
zeros.  On the face of it, \textbf{A 4} is considerably
stronger than the statement that $\Delta_{1,\gamma}$
has no critical points, because the absence of critical
points for each $\Delta_{1,\gamma}$
would simply be a condition concerning the individual elements
$\gamma$ of $\Gamma-\Gamma^{\phi_{1}}$.  Because in
\textbf{A 4}, neither $\vectu$ nor the sign in \eqref{eqn:nocritpoints}
is allowed to depend on $\gamma$, \textbf{A 4} is essentially a statement
about the action of the group $\Gamma$ as a whole.  In the end
it does turn out that \textbf{A 4} is satisfied by all the discrete subgroups
$\Gamma$ of $\SL{2}{\Complex}$
that one would reasonably want to consider in this context.  The reason
for introducing \textbf{A 2} is that \textbf{A 2} isolates the
part of \textbf{A 4} necessary for the
conclusion of Theorem \ref{thm:grenierh3}, below,
to hold (under the assumption of \textbf{A 1} and \textbf{A 3}).
Therefore, using \textbf{A 2} to state the hypotheses
of Theorem \ref{thm:grenierh3}
helps clarify the reason Theorem \ref{thm:grenierh3}
holds true.

Applying Lemma \ref{lem:action1cons} to our situation, we
obtain the following result.
\begin{lem}\label{lem:action1consapplied}  Assume that the group
$\Gamma$ acts by diffeomorphisms on $\bbH^3$.  Further,
assume that the action of $\Gamma$ on $\bbH^3$ satisfies
axiom \textbf{A 1}.  Then $\Gamma^{\phi_1}$ has an action
on $\Real^2$ induced by \eqref{eqn:sectioninducedact},
\starteqn\label{eqn:sectioninducedactapplied}
\gamma_{[2,3]}(\mathbf{t})=\phi_{[2,3]}(\gamma\sigma^0_{[2,3]}(\mathbf{t}))
,\;\forallindisp\;\gamma\in\Gamma^{\phi_1},\;
\;\mathbf{t}=(t_1,t_2)\in\Real^{2}. \finisheqn
\end{lem}

We note that each \textbf{A} axiom is the version of
the \textbf{A} axiom from \S1.2 of \cite{brennerthesis},
adapted to the situation at hand.  Referring
to the notation of \cite{brennerthesis},
we are considering in this work the case of,
$X=\bbH^3$, $N=3$, $\scrM=\{1\}$, and
\startdisp
F^{\gamma,N-i+1}=F^{\gamma,3}=\phi_1\gamma\circ\phi\inv.
\finishdisp
Note that, in our situation, the second sentence
of \textbf{A 1} given in \cite{brennerthesis} becomes
redundant.  Because the axioms are the same, we can
apply Theorem 1.2.6, respectively Theorem 1.2.8, from \cite{brennerthesis}
to our situation to obtain Theorem \ref{thm:grenierh3}, respectively
Theorem \ref{thm:grenierh3extended}, below.

\begin{thm}\label{thm:grenierh3} Let $\bbH^3$,
$\phi$, $\Gamma$ be as above.
Assume that the action of $\Gamma$ on $\bbH^3$ satisfies axioms \textbf{A 1}
through \textbf{A 3} above.  Let $\scrG\subseteq\Real^{2}$ be a
fundamental domain for the induced action of
$\Gamma^{\phi_{[1,3]}}\backslash \Gamma^{\phi_{1}}$ on $\Real^{2}$. Assume
that $\scrG$ is closed.  Define
\startdisp
\scrF_1=\{z\in \bbH^3\;|\;
\phi_{1}(z)\leq \phi_{1}(\gamma z),\forallindisp
\gamma\in\Gamma^{\phi_1}\}.
\finishdisp
Set
\starteqn\label{eqn:scrFdefn}
\scrF=\phi_{[2,3]}\inv(\scrG)\cap \scrF_1. \finisheqn
Then we have
\begin{itemize}
\item[(a)] $\scrF$ is a
fundamental domain for the action of $\Gamma^{\phi_{[1,3]}}\backslash
\Gamma$ on
$\bbH^3$.
\item[(b)] We have
$\scrF_1$ closed, so that, by \eqref{eqn:scrFdefn} and the assumption
that $\scrG$ is closed,
$\scrF$ is closed.  Also,
\starteqn\label{eqn:intrrF1} \intrr\scrF_1=\{z\in \bbH^3\;|\;
\phi_{1}(z)< \phi_{1}(\gamma x),\forallindisp\;
\gamma\in\Gamma-\Gamma^{\phi_{1}}\}, \finisheqn
and
\starteqn\label{eqn:boundaryFp}
\partial\scrF_1=\{z\in\scrF_1\;|\; \phi_{1}(z)=\phi_{1}(\gamma z),\;\;
\text{for some} \;\gamma\in\Gamma-\Gamma^{\phi_1}\}.
\finisheqn
\end{itemize}
\end{thm}

\begin{thm}\label{thm:grenierh3extended}
Let $\bbH^3$,
$\phi$, $\Gamma$, $\scrG$ be as in Theorem \ref{thm:grenierh3}.
Suppose that the action of $\Gamma$ on
$\bbH^3$ satisfies \textbf{A 4} in addition to \textbf{A 1} through
\mbox{\textbf{A 3}}.  Suppose that
$\scrG=\overline{\intrr(\scrG)}$. Then we have the conclusions
of Theorem \ref{thm:grenierh3} and also
\starteqn\label{eqn:scrFp1ice}
\scrF_1=\overline{\intrr(\scrF_1)}.
\finisheqn
Further,
\starteqn\label{eqn:intrrF}
\intrr\scrF=\phi_{[2,3]}\inv(\intrr(\scrG))\cap
\intrr(\scrF_1),
\finisheqn
and
\starteqn\label{eqn:fisfintclosure}
\scrF=\overline{\intrr(\scrF)}.
\finisheqn
\end{thm}

Considering the coordinate system $\phi$ on $\bbH^3$ as fixed,
we may think of the fundamental domain $\scrF$ for
$\Gamma^{\phi_{[1,3]}}\backslash\Gamma$ to be a function of the fundamental
domain $\scrG$ for the induced action of $\Gamma^{\phi_1}$
on $\Real^2$.  When we wish to stress this dependence
of $\scrF$ on $\scrG$, we will write $\scrF(\scrG)$ instead of $\scrF$.

\begin{defn}  Suppose that the action of a group $\Gamma$ on $\bbH^3$
on $\scrF$ satisfies axioms \linebreak\textbf{A 1} through \textbf{A 4},
above.  Let $\scrG$ be a fundamental domain for the induced
action of $\Gamma^{\phi_{[1,3]}}\backslash
\Gamma^{\phi_1}$ on $\Real^2$ satisfying $\scrG=\overline{
\intrr(\scrG)}$.  Then the fundamental domain $\scrF(\scrG)$
for the action of $\Gamma^{\phi_{[1,3]}}\backslash\Gamma$
defined in \eqref{eqn:scrFdefn} is called the \bmth\textbf{good
Grenier fundamental domain for the action of $\Gamma$ on $\bbH^3$
associated to the fundamental domain $\scrG$}\ubmth.
\end{defn}
The reference to the fundamental domain $\scrG$ is often omitted in
practice.

Henceforth, we drop the explicit reference to $\Gamma^{\phi_{[1,3]}}$
and speak of a
\textit{fundamental domain of $\Gamma^{\phi_{[1,3]}}\backslash\Gamma$}
as a \textit{fundamental domain of $\Gamma$}.
By \eqref{eqn:kernelisverysmall},
$\Gamma$ is at worst a two-fold cover of
$\Gamma^{\phi_{[1,3]}}\backslash\Gamma$,
so this involves only a minor abuse of terminology.

\vspace*{0.3cm}
\noindent\bmth\textbf{The A axioms and discrete subgroups $\Gamma$
 of $\SL{2}{\Complex}$.}
\ubmth\hspace*{0.3cm}
 We now consider the problem of verifying
the \textbf{A} axioms for the examples of groups of diffeomorphisms
of $\bbH^3$ that arise in practice, namely discrete subgroups of
$\SL{2}{\Complex}$.  The following result says that all such
subgroups satisfy the first two axioms.
\begin{lem} \label{lem:a1a2} Let $\Gamma$
be a subgroup of $\SL{2}{\Complex}$,
acting on $\bbH^3$ on the left by fractional linear transformations.
Then $\Gamma$ satisfies Axiom \textbf{A 1}.
\end{lem}
\begin{myproof}{Proof}
  For \textbf{A 1}, note that, by
\eqref{eqn:stabilizerintersection} and \eqref{eqn:firstcoordfixing},
we have
\starteqn\label{eqn:gammastabilizerform}
\Gamma^{\phi_1}=\left\{\gamma\in\Gamma\;\left|\; \gamma=
\begin{pmatrix}\omega&b\\
0&\omega\inv
\end{pmatrix},\;||\omega||=1,\, b\in\Complex
\right.\right\}.
\finisheqn
It is easy to see that \textbf{A 1} is equivalent to the following condition.
\starteqn\label{eqn:a1confirmation}\text{For all $z_1, z_2\in\bbH^3$,
$\gamma\in\Gamma^{\phi_1}$,
$\phi_{[2,3]}(z_1)=\phi_{[2,3]}(z_2)$
implies $\phi_{[2,3]}(\gamma z_1)=\phi_{[2,3]}(\gamma z_1)$.
}\finisheqn
The hypothesis $\phi_{[2,3]}(z_1)=\phi_{[2,3]}(z_2)$ means,
by the definition of $\phi$ that, if
\startdisp
z_i=x_i+y\vectj,\; x\in\Complex,\, y\in\Real^+,\;\text{then}\; x_1=x_2.
\finishdisp
It is straightforward to calculate that if $\gamma$ is of the form
given in \eqref{eqn:gammastabilizerform}, then
\startdisp
\gamma z_i=\omega^2x(z_i)+\omega b+y(z_i)\vectj
\finishdisp
Thus,
\startdisp
x(\gamma z_1)=\omega^2x(z_1)+\omega b=\omega^2x(z_2)+\omega b=x(\gamma z_2),
\finishdisp
\textit{i.e.}, $\phi_{[2,3]}(\gamma z_1)=\phi_{[2,3]}(\gamma z_1)$.
This proves \eqref{eqn:a1confirmation} and therefore
the lemma.
\end{myproof}

Before proceeding to Axioms \textbf{A 2} and \textbf{A 4},
we calculate explicitly the difference function in the special
case at hand and fix some notation related to paths.
\begin{lem} \label{lem:differencefunctionexplicit} Let $\gamma\in\SL{2}{\Complex}$, given
explicitly as in \eqref{eqn:gammaa4lemform}.  Let $\phi$ be
the diffeomorphism of $\bbH^3$ onto $\Real^3$ defined in \eqref{eqn:phicoordsdefn},
and $\Delta_{1,\gamma}$ the difference function defined in \eqref{eqn:differencefunctiondefn}.
Use the notation $\vectt=\phi(z)$.  Then we have
\starteqn\label{eqn:differencefunctionexplicit}
\Delta_{1,\gamma}(t)=\Delta_{1,\gamma}(\phi(z))=-\log\left(||cz+d||^2\right).
\finisheqn
\end{lem}
\begin{myproof}{Proof}  We calculate $\Delta_{1,\gamma}$ by
applying the definition of $\Delta_{1,\gamma}$
\eqref{eqn:differencefunctiondefn}, \eqref{eqn:phicoordsdefn},
\eqref{eqn:ycdzrelation},
and \eqref{eqn:ycdzdefn} in turn.
\startdisp
\begin{aligned}
\Delta_{1,\gamma}(\phi(z))&=&&\phi_1(z)-\phi_1(\gamma z)\\
&=&&-\log y(z)+\log y(\gamma z)\\
&=&&-\log y(z)+\log y(c,d;z)\\
&=&&-\log y(z)+\log y(z)-\log\left(||cz+d||^2\right)\\
&=&&-\log\left(||cz+d||^2\right).
\end{aligned}
\finishdisp
This completes the proof of the lemma.
\end{myproof}

Let $\Gamma$
be a subgroup of $\SL{2}{\Complex}$,
acting on $\bbH^3$ on the left by fractional linear transformations.
Let $\gamma\in\Gamma$ be given by
\starteqn\label{eqn:gammaa4lemform}
\gamma=\begin{pmatrix}a&b\\c&d
\end{pmatrix},\;\text{with}\;a,b,c,d\in\Complex,\; ad-bc=1.
\finisheqn
Let
\starteqn\label{eqn:z0a4lem}
z_0=x_0+y_0\vectj\in\bbH^3,\;\text{with}\; x_0\in\Complex,\, y_0\in\Real^+.
\finisheqn
Let $z(s)$ be a smooth curve in $\bbH^3$ defined
on an interval $(-\epsilon,\epsilon)$ ($\epsilon>0$) and satisfying
\starteqn\label{eqn:pathconditions}
z(0)=z_0\;\text{and}\; z(s)=x_0+y(s)\vectj,\;\text{for}\; s\in
(-\epsilon,\epsilon).
\finisheqn
It follows immediately from \eqref{eqn:pathconditions}, the explicit
formulas for $\phi$ in \eqref{eqn:phicoordsdefn}
and the fact that $\phi$ is a diffeomorphism
of $\bbH^3$ with $\Real^3$ that $\phi z(s)=\phi(z(s))$ is a smooth curve
in $\Real^3$ satisfying
\starteqn\label{eqn:imagepathconditions}
\phi z(0)=\phi(z_0)\;\text{and}\;
\phi z (s)=(\Rept(x_0),\Impt(x_0),-\log(y(s))),\;\text{for}\; s\in
(-\epsilon,\epsilon).
\finisheqn
\begin{lem} \label{lem:a2}  Let $\Gamma$
be a subgroup of $\SL{2}{\Complex}$,
acting on $\bbH^3$ on the left by fractional linear transformations.
Let $\gamma\in\Gamma$ be as in \eqref{eqn:gammaa4lemform},
$z_0\in\bbH^3$ be as in \eqref{eqn:z0a4lem}.  Let $z(s)$
satisfy the conditions of \eqref{eqn:pathconditions}, so that $\phi z(s):=
\phi(z(s))$ satisfies the conditions of \eqref{eqn:imagepathconditions}.
\begin{itemize}
\item[(a)]  We have
\starteqn\label{eqn:derivativedifffunctcalc}
\left.\frac{\intd}{\intd s}\Delta_{1,\gamma}
(\phi z(s))\right|_{s=0}=\left(\frac{-2||c||^2y(0)}
{||cx_0+d||^2+||c||^2y(0)^2}\right)\,y'(0).
\finisheqn
\item[(b)]  The action of $\Gamma$ satisfies Axiom \textbf{A 4}.
\item[(c)]  The action of $\Gamma$ satisfies Axiom \textbf{A 2}.
\end{itemize}
\end{lem}
\begin{myproof}{Proof}  By Lemma \ref{lem:differencefunctionexplicit},
we have
\starteqn\label{eqn:explicitdifffunctform}
\Delta_{1,\gamma}(\vectt)=-\log||cz+d||^2.
\finisheqn
As usual, write
\startdisp
z=x(z)+y(z)\vectj,\;\text{with}\;x(z)\in\Complex, y(z)\in\Real^+.
\finishdisp
Since $d\in\Complex$, the definition of the squared
norm of a quaternion implies that
\startdisp
||cz+d||^2=||cx(z)+d||^2+||cy(z)||^2=||cx(z)+d||^2+||c||^2y(z)^2.
\finishdisp
Therefore, we can expand \eqref{eqn:explicitdifffunctform}, to
obtain the formula
\starteqn\label{eqn:explicitdifffunctformexpand}
\Delta_{1,\gamma}(\vectt)=-\log\left(||cx(z)+d||^2+||c||^2y(z)^2\right).
\finisheqn
We now apply \eqref{eqn:explicitdifffunctformexpand} to
the special case of the path $z(t)$ in $\bbH^3$.
For $s\in (-\epsilon,\epsilon)$, we calculate that
\startdisp
\Delta_{1,\gamma}(\phi z(s))=-\log\left(||cx_0+d||^2+||c||^2y(s)^2\right).
\finishdisp
Using freshman calculus we deduce that
\startdisp
\frac{\intd}{\intd s}\Delta_{1,\gamma}
(\phi z(s))=\frac{-2||c||^2y(s)y'(s)}{||cx_0+d||^2+||c||^2y(s)^2}
\finishdisp
Setting $s=0$, we obtain \eqref{eqn:derivativedifffunctcalc}.
This completes the proof of (a).

For (b) we are to show that given any $\vectt\in\Real^3$,
there is a direction $\vectu$,
independent of \linebreak $\gamma\in\Gamma-\Gamma^{\phi_1}$, satisfying
\eqref{eqn:nocritpoints}.  By
\eqref{eqn:gammastabilizerform} the assumption that
$\gamma\in\Gamma-\Gamma^{\phi_1}$
is equivalent to the assumption that $c\neq 0$, where $c$
is the lower-left entry of $\gamma$ as in \eqref{eqn:gammaa4lemform}.
Under the condition $c\neq 0$, it is easy to see that the
factor in front of $y'(0)$ on the left-hand side of
\eqref{eqn:derivativedifffunctcalc} is negative.  Therefore,
we deduce from \eqref{eqn:derivativedifffunctcalc} that
if $z(s)$ is as in part (a), then,
\starteqn\label{eqn:yprime0signimplication}
y'(0)>0\;\text{if and only if}\;\left.\frac{\intd}{\intd s}\Delta_{1,\gamma}
(\phi z(s))\right|_{s=0}<0.
\finisheqn
With $\vectt\in\Real^3$ given, let $z_0=\phi\inv(\vectt)$,
and
\startdisp
\text{$z(s)$ a smooth path in $\bbH^3$ satisfying the conditions
\eqref{eqn:pathconditions} and $y'(0)>0$}.
\finishdisp  Define $\vectu$ to be the unit
tangent vector to the path $\phi(z(s))$ at $s=0$.
Because $\phi$ is a diffeomorphism $\vectu$ is well-defined.
That $\phi$ is a diffeomorphism, together with the definition
of $\vectu$, implies that
\starteqn\label{eqn:directionalderivexpr}
D_{\vectu}\Delta_{1,\gamma}(\vectt)=\left.\frac{\intd}{\intd s}
\Delta_{1,\gamma}
(\phi z(s))\right|_{s=0}
\finisheqn
By \eqref{eqn:yprime0signimplication}, \eqref{eqn:directionalderivexpr},
and the choice of $z(s)$ satisfying $y'(0)>0$,
we have
\startdisp
D_{\vectu}\Delta_{1,\gamma}(\vectt)<0.
\finishdisp
Therefore, $\vectu$ satisfies the condition \eqref{eqn:nocritpoints}.
We conclude that the action of $\Gamma$ on $\bbH^3$
satisfies Axiom \textbf{A 4}.

Part (c) follows immediately from Part (b) and \eqref{eqn:a4impliesa2}.
\end{myproof}

The verification of Axiom \textbf{A 3}
requires more specialized techniques than the verification of \textbf{A 1},
\textbf{A 2} and \textbf{A 4}.  The reason is that unlike \textbf{A 1}
and \textbf{A 2}, \textbf{A 3} and \textbf{A 4} are essentially
statements about the action of the entire group $\Gamma$
on $\bbH^3$ in relation to the coordinate system, whereas \textbf{A 1}
and \textbf{A 2} merely concern the actions of individual elements
of $\Gamma$.  In particular, while Lemmas \ref{lem:a1a2}
and \ref{lem:a2} guarantee
that every group $\Gamma$ of fractional linear transformations
satisfies \textbf{A 1}, \textbf{A 2}, \linebreak \textbf{A 4}, it is easy to construct
examples of discrete groups $\Gamma$ of fractional linear
transformations which do not satisfy \textbf{A 3}.

We use a more indirect method to verify Axiom \textbf{A 3}.
The method for \textbf{A 3} amounts to verifying that all $\Gamma$
in a certain ``class" (the commensurability class, see
below Corollary \eqref{cor:a3} for the definition) either
satisfy \textbf{A 3} or fail to satisfy \textbf{A 3}.  We
then pick a specific $\Gamma_0$ in the class of $\Gamma=\conj\inv(\SO{3}{\Int[\vecti]})$
and show by direct calculations that $\Gamma_0$ satisfies
\textbf{A 3}, from which it follows via Corollary \eqref{cor:a3}
that $\Gamma$ satisfies \textbf{A 3}.  The class representative
$\Gamma_0$ is chosen so as to simplify the calculations to whatever
extent possible.

First, we recall the following definition from general topology
\begin{defn}  Let $X$ be a topological space, $x\in X$.
Then a \bmth\textbf{local base at $x$}\hspace*{0.15cm}\ubmth
is a collection of neighborhoods
$\scrB(x)$ of $x$ such that
\startdisp\text{
For every neighborhood $U$ of $x$ there is a neighborhood $B\in\scrB(x)$
such that $B\subset U$.
}
\finishdisp
\end{defn}

Now suppose that $(X,d)$ is a metric space with metric $d$.
For example, we may take $(X,d)=(\bbH^3,s^2)$, where
$s^2$ is the hyperbolic metric defined in terms of the above
coordinates by
\startdisp
\intd s^2=\frac{\intd x_1^2+\intd x_2^2+\intd y^2}{y^2}.
\finishdisp
By definition of the metric topology each $x\in X$ has a countable local base, consisting of
the open balls $B_{\frac{1}{n}}(x)$ centered at $x$, of radius
$1/n$, $n\in\Natural$.
We use the notation
\starteqn\label{eqn:localbasisspecific}
\scrB_{d}(x)=\left\{B_{\frac{1}{n}}(x)\;|\; n\in\Natural\right\},
\finisheqn
for the local base so described.  Because for each $z\in\bbH^3$,
$\scrB_{s^2}(z)$ is a local basis at $z$, we have a useful
reformulation of Axiom \textbf{A 3}.

\begin{lem} \label{lem:countablebasis}
Let $\Gamma$ be a group acting on $\bbH^3$ by diffeomorphisms,
as above.
\begin{itemize}
\item[(a)]  Let $\scrB(z)$ be a local base for the metric
topology of $\bbH^3$.  The action of $\Gamma$
on $\bbH^3$ satisfies \textbf{A 3} if and only if, for each $z\in\bbH^3$,
$B\in\Real$, there is an element $V\in\scrB(z)$, such that $V$
satisfies \eqref{eqn:A3Ucond}.
\item[(b)]  The action of $\Gamma$ on $\bbH^3$
satisfies \textbf{A 3} if and only if for each $z\in\bbH^3$
and $B\in\Real$, there exists $N(B,z)\in\Natural$ such that for
all $n>N$ the ball $B_{\frac{1}{n}}(z)$ satisfies \eqref{eqn:A3Ucond}.
\end{itemize}
\end{lem}
\begin{myproof}{Proof}  Part (a) follows from the fact
that if a neighborhood $U$ of $z$ satisfies \eqref{eqn:A3Ucond},
then so does any neighborhood $V$ of $z$ contained in $U$.
Given any neighborhood $U$ of $z$ satisfying \eqref{eqn:A3Ucond},
simply take $V\in\scrB(z)$ with $V\subseteq U$, which
is possible because $\scrB(z)$ is a basis of neighborhoods
at $z$.

In order to obtain part (b), apply part (a) to the special
case of  $\scrB(z)=\scrB_{s^2}(z)$, the basis of neighborhoods
for the metric topology on $\bbH^3$
defined in \eqref{eqn:localbasisspecific}.  Part (a)
gives an $N(B,z)$ such that $B_{N(B,z)\inv}$ satisfies
\eqref{eqn:A3Ucond}, and since for every $n>N(B,z)$,
$B_\frac{1}{n}\subseteq B_{N(B,z)\inv}$, $B_{\frac{1}{n}}$
also satisfies \eqref{eqn:A3Ucond}.
\end{myproof}

For the proof of Lemma \ref{lem:finiteindexA3}, below, we will
need some facts from elementary group theory, which we list
in Sublemma \ref{sublem:grouptheoryfacts}.
\begin{sublem}\label{sublem:grouptheoryfacts}  Let $\Gamma$ be a group,
with $\Gamma'$, $\Delta$ subgroups of $\Gamma$.  Define
\startdisp
\Delta'=\Gamma'\cap\Delta.
\finishdisp
\begin{itemize}
\item[(a)]  Assume for this part only that
\startdisp
[\Gamma:\Gamma']=n<\infty.
\finishdisp
Then there exist elements $\eta_1,\ldots \eta_{m-1}$,
each $\eta_i\neq 1$, such that
$\Gamma$ can be expressed as the disjoint union of $m$
$\Gamma'$ left cosets.  That is
\starteqn\label{eqn:indexcosets}
\Gamma=\Gamma'\cup\Gamma'\eta_1
\cup\Gamma'\eta_2\cup\cdots\cup\Gamma'\eta_{m-1}.
\finisheqn
\item[(b)]  Let $\gamma_1',\,\gamma_2'$ be arbitrary elements
of $\Gamma'$.  For $i=1,2$, consider the $\Delta'$ left cosets
$\Delta'\gamma_i'$ of $\gamma_i'$ in $\Gamma'$ and the $\Delta$
left cosets $\Delta\gamma_i'$ of $\gamma_i'$ in $\Gamma$.
Then we have
\starteqn\label{eqn:cosetequalityequivalence}
\Delta'\gamma_1'=\Delta'\gamma_2'\;\text{if and only if}\;
\Delta\gamma_1'=\Delta\gamma_2'.
\finisheqn
\end{itemize}
\end{sublem}
\begin{myproof}{Proof}  Part (a) is a standard part of elementary group
theory.

For part (b), we have $\Delta'\gamma_1'=\Delta'\gamma_2'$ if and
only if $\gamma'_1(\gamma'_2)\inv\in\Delta'$.  But
$\gamma'_1(\gamma'_2)\inv\in\Gamma'$, and also $\Delta'=\Gamma'\cap\Delta$.
Therefore, we have $\gamma'_1(\gamma'_2)\inv\in\Delta'$ if and only
if $\gamma'_1(\gamma'_2)\inv\in\Delta$.
Since, $\gamma'_1(\gamma'_2)\inv\in\Delta$
if and only if $\gamma'_1\Delta=\gamma'_2\Delta$, this completes
the proof of (b).
\end{myproof}

\begin{lem} \label{lem:finiteindexA3}
Let $\Gamma$ be a group of fractional linear
transformations acting on $\bbH^3$, and let $\Gamma'$ be a subgroup
of $\Gamma$.
\begin{itemize}
\item[(a)]    If $\Gamma$ satisfies Axiom \textbf{A 3}, then $\Gamma'$
also satisfies Axiom \textbf{A 3}.
\item[(b)] Suppose that $[\Gamma:\Gamma']<\infty$.   If $\Gamma'$
satisfies Axiom \textbf{A 3}, then $\Gamma$ also satisfies Axiom \textbf{A 3}.
\end{itemize}
\end{lem}
\begin{myproof}{Proof}
For (a), we suppose that $\Gamma'$ does not satisfy \textbf{A 3},
and we prove that in this case $\Gamma$ does not satisfy
\textbf{A 3}.  Then for
some $z\in\bbH^3$, $B\in\Real$ and for every neighborhood $U$ of
$z$ in $\bbH^3$, there is an infinite sequence
$\{\gamma_i'\}_{i=1}^{\infty}$ such that
\startdisp
\phi_1((\Gamma')^{\phi_1}\gamma_i'U)\cap(-\infty,B]\neq\emptyset\;
\text{and the}\;(\Gamma')^{\phi_1}\gamma_i'\;\text{are distinct.}
\finishdisp
By \eqref{eqn:Deltaonleftinphij}, we have
\startdisp
\phi_1((\Gamma')^{\phi_1}\gamma_i'U)=\phi_1(\Gamma^{\phi_1}\gamma_i'U).
\finishdisp
Therefore, $\{\gamma_i'\}_{i=1}^{\infty}$ is an infinite sequence
in $\Gamma$ such that
\startdisp
\phi_1(\Gamma^{\phi_1}\gamma_i'U)\cap(-\infty,B]\neq\emptyset.
\finishdisp
We may conclude that $\Gamma$ does not satisfy \textbf{A 3}
provided that we can prove that the cosets $\Gamma^{\phi_1}\gamma_i'$
are distinct.  However, we can prove that the cosets
$\Gamma^{\phi_1}\gamma_i'$ are distinct by applying
Sublemma \ref{sublem:grouptheoryfacts}, part (b).  Therefore,
if $\Gamma'$ fails to satisfy \textbf{A 3}, then
$\Gamma$ fails to satisfy \textbf{A 3}.  We have
proved the contrapositive of part (a), hence we have proved
part (a).

For part (b), suppose that $\Gamma'$ satisfies \textbf{A 3}.
Let $B\in\Real$ be fixed.
Set $m=[\Gamma:\Gamma']$.  By assumption, $m<\infty$.
Let $\eta_1,\ldots,\eta_{m-1}$ be the nontrivial coset representatives
of $\Gamma'$ in $\Gamma$, as in Sublemma \ref{sublem:grouptheoryfacts},
part (b).  By convention, set $\eta_0=\Id_{\Gamma}$.
For $i\in[0,m-1]$, let $N'(eta_i z, B)$ be the smallest
integer such that for all $n>N'(eta_i z, B)$, the ball
$B_{\frac{1}{n}}(\eta_i z)$ satisfies \eqref{eqn:A3Ucond},
with $\Gamma'$ in place of $\Gamma$.
Such $N'(\eta_i z, B)\in\Natural$ exist, by Lemma \ref{lem:countablebasis},
part (b), because $\Gamma'$
satisfies \textbf{A 3}.  Set
\starteqn\label{eqn:NofBzdefn}
N=\max\big(N'(\eta_0z),\ldots N'(\eta_{m-1}z)\big).
\finisheqn
We claim that $N=N(B,z)$, \textit{i.e.} that for all
$n>N$, the ball $B_{\frac{1}{n}}(z)$ satisfies
\eqref{eqn:A3Ucond} with respect to the action of $\Gamma$.
By Lemma \ref{lem:countablebasis}, part (b), the claim
implies that $\Gamma$ satisfies \textbf{A 3}, and
therefore, the proof of the claim will complete the proof
of the theorem.

In order to prove the claim, it clearly suffice
to prove that $B_{\frac{1}{N}}$ satisfies \eqref{eqn:A3Ucond}.
Assume the contrary, namely that there exists
an infinite sequence $\{\gamma_j\}_{j=1}^{\infty}$ of elements
of $\Gamma$ such that
\starteqn\label{eqn:a3proofcontradiction}\text{
the cosets $\Gamma^{\phi_1}\gamma_j$ are distinct and
$\phi_1(\gamma_jB_{\frac{1}{N}}(z))\cap(-\infty,B)\neq\emptyset$.}
\finisheqn
Infinitely many $\gamma_j$ in the infinite sequence
must belong to a single coset $\Gamma'\eta_i$ for $i\in[0,m-1]$.
By extraction of a subsequence of $\{\gamma_j\}_{j=1}^{\infty}$,
we may assume without loss of generality that each $\gamma_j$
is of the form $\gamma'_j\eta_i$ for a fixed $i\in[0,m-1]$.
By an application of Sublemma \ref{sublem:grouptheoryfacts},
we deduce from \eqref{eqn:a3proofcontradiction}
that the $(\Gamma')^{\phi_1}$-cosets $(\Gamma^{\phi_1}\gamma'_j$
are distinct.  Therefore, \eqref{eqn:a3proofcontradiction} that
\starteqn\label{eqn:a3proofcontradiction2}\text{
the cosets $(\Gamma')^{\phi_1}\gamma_j'$ are distinct and
$\phi_1(\gamma_j'\eta_i B_{\frac{1}{N}}(z))\cap(-\infty,B)\neq\emptyset$.}
\finisheqn
But because $\eta_i$ is an isometric automorphism of $\bbH^3$,
and by \eqref{eqn:NofBzdefn}, we have
\startdisp
\gamma_j'\eta_i B_{\frac{1}{N}}(z)=
\gamma_j' B_{\frac{1}{N}}(\eta_i z)\subseteq
\gamma_j' B_{1/N(\eta_iz,B)}(\eta_iz).
\finishdisp
Therefore, \eqref{eqn:a3proofcontradiction2} implies that
\startdisp\text{
the cosets $(\Gamma')^{\phi_1}\gamma_j'$ are distinct and
$\phi_1(\gamma_j'B_{1/N(\eta_iz,B)}(\eta_iz))\cap(-\infty,B)\neq\emptyset$,}
\finishdisp
which contradicts the choice of $N(\eta_iz,B)$.  Since
our assumption that there exists
an infinite sequence $\{\gamma_j\}_{j=1}^{\infty}$ of elements
of $\Gamma$ satisfying \eqref{eqn:a3proofcontradiction} lead
to a contradiction, we may conclude that no such sequence
exists.  This shows that $B_{\frac{1}{N}}$ satisfies \eqref{eqn:A3Ucond},
and completes the proof of the claim, and hence the proof of the lemma.
\end{myproof}

As an immediate consequence of Lemma \ref{lem:finiteindexA3}, we have
the following.
\begin{cor} \label{cor:a3} Suppose that $\Gamma_1$, $\Gamma_2$
are groups of fractional
linear transformations acting on $\bbH^3$ such that
\starteqn\label{eqn:commensurability}
[\Gamma_1:\Gamma_1\cap\Gamma_2],\;[\Gamma_2:\Gamma_1\cap\Gamma_2]<\infty.
\finisheqn
Then the action of $\Gamma_1$ satisfies \textbf{A 3} if and only if
the action of $\Gamma_2$ satisfies \textbf{A 3}.
\end{cor}

We call two subgroups $\Gamma_1, \Gamma_2$ of a group $G$ \textbf{commensurable}
if they satisfy the finite-index condition of \eqref{eqn:commensurability}.
Since commensurability is clearly an equivalence relation on the
set of all subgroup of $G$, we may speak of the \textbf{commensurability
class} of a subgroup $\Gamma$ of $G$.  Using this terminology,
we may restate the result of Corollary \ref{cor:a3} as follows.
\textit{For any commensurability class of subgroups of $\SL{2}{\Complex}$,
all of the subgroups in the class either satisfy \textbf{A 3} or fail
to satisfy \linebreak \textbf{A 3}.}

Inside the commensurability class of $\conj\inv(\SO{3}{\Int[\vecti]})$,
we find $\SL{2}{\Int[\vecti]}$.  See Lemma \ref{lem:indexlemma}, specifically
\eqref{eqn:index6}, for the commensurability result.  In light
of Corollary \ref{cor:a3}, in order to prove that $\SO{3}{\Int[\vecti]}$
satisfies \textbf{A 3}, we show that $\SL{2}{\Int[\vecti]}$ satisfies
\textbf{A 3}.  We now prepare the way for this demonstration.

\vspace*{0.3cm}
\noindent\bmth\textbf{Explicit description of the decomposition
of $\SL{2}{\Int[\vecti]}$ into left $\SL{2}{\Int[\vecti]}^{\phi_1}$
cosets.}\ubmth\hspace*{0.3cm}
By \eqref{eqn:ginstabilizerprelim} and
\eqref{eqn:gammastabilizerintersection}, applied to
$\Gamma=\SL{2}{\Int[\vecti]}$, we have
\starteqn\label{eqn:sl2intstabilizer}
\SL{2}{\Int[\vecti]}^{\phi_1}=\left\{\left.
\begin{pmatrix}\vecti^{-\epsilon}&x\\0&\vecti^{\epsilon}\end{pmatrix}\;\right|\;
x\in\Int[\vecti],
\,\epsilon\in\{0,1,2,3\}\right\}.
\finisheqn
Let
\startdisp
\gamma=\begin{pmatrix}a&b\\c&d\end{pmatrix}\in\SL{2}{\Int}
\finishdisp
be an arbitrary element of $\SL{2}{\Int[\vecti]}$.  It follows directly from
\eqref{eqn:sl2intstabilizer} that we have the following
description of the coset containing $\gamma$.
\starteqn\label{eqn:sl2intcosetrough}
\SL{2}{\Int[\vecti]}^{\phi_1}\gamma=
\left\{\left.
\begin{pmatrix}
\vecti^{-\epsilon}a+xc&\vecti^{-\epsilon}b+xd\\
\vecti^{\epsilon}c&\vecti^{\epsilon}d
\end{pmatrix}\;\right|\;
x\in\Int[\vecti],
\,\epsilon\in\{0,1,2,3\}\right\}.
\finisheqn
One can verify using elementary number theory that for fixed
$\epsilon\in\{0,1,2,3\}$,
\startdisp\text{The set of elements of $\SL{2}{\Complex}$
with second row
$\begin{pmatrix}\vecti^{\epsilon}c&\vecti^{\epsilon}d\end{pmatrix}$
is
$\left\{\left.
\begin{pmatrix}
\vecti^{-\epsilon}a+xc&\vecti^{-\epsilon}b+xd\\
\vecti^{\epsilon}c&\vecti^{\epsilon}d
\end{pmatrix}\;\right|\;
x\in\Int[\vecti]\right\}$}.
\finishdisp
Therefore \eqref{eqn:sl2intcosetrough} implies that
\starteqn\label{eqn:sl2intcoset}
\SL{2}{\Int[\vecti]}^{\phi_1}\begin{pmatrix}a&b\\c&d
\end{pmatrix}=\left\{\left.\begin{pmatrix}
*&*\\
\vecti^{\epsilon}c&\vecti^{\epsilon}d
\end{pmatrix}\in\SL{2}{\Int[\vecti]}\;\right|\; \epsilon\in\{0,1,2,3\}\right\}.
\finisheqn
Also, by elementary number theory, the row vector
$\begin{pmatrix}c&d\end{pmatrix}$
can be completed to an element of $\SL{2}{\Int[\vecti]}$ if and only
if $\mathrm{GCD}(c,d)=1$.
As a consequence of \eqref{eqn:sl2intcoset} we have the bijection
\startdisp
\{\SL{2}{\Int[\vecti]}^{\phi_1}\gamma\;|\;\gamma\in\SL{2}{\Int[\vecti]}\}
\longleftrightarrow
\{\begin{pmatrix}c&d
\end{pmatrix}\;|\;\mathrm{GCD}(c,d)=1\}/\{\times\vecti^\epsilon\;|\;\epsilon=
0,1,2,3\},
\finishdisp
which is described explicitly by \eqref{eqn:sl2intcoset}.

\begin{prop} \label{prop:sl2zia3} The action of
$\SL{2}{\Int[\vecti]}$ on $\bbH^3$
by fractional linear transformations satisfies Axiom \textbf{A 3}.
\end{prop}
\begin{myproof}{Proof}   Let $z\in\bbH^3$,
$B\in\Real$.  For any $w\in\bbH^3$, we write, as usual,
\starteqn\label{eqn:wcoordinates}
w=x(w)+y(w)\vectj,\;\text{for}\;x(w)\in\Complex,\; y(w)\in\Real^+.
\finisheqn
Define
\startdisp
U_z=\{w\in\bbH^3\;|\; 2y(z)>y(w)>\half y(z),\, ||x(w)-x(z)||\leq 1 \}.
\finishdisp
It is clear that $U_z$ is a neighborhood of $z$.  In order to complete
the verification of \textbf{A 3}, it will
therefore suffice to show that there are only finitely
many cosets $\SL{2}{\Int[\vecti]}^{\phi_1}\gamma$ of
$\SL{2}{\Int[\vecti]}^{\phi_1}$ in $\SL{2}{\Int[\vecti]}$ such that
\startdisp
\phi_1(\SL{2}{\Complex}^{\phi_1}\gamma U_z)\cap(-\infty,B]\neq\emptyset.
\finishdisp
Let $w\in U_z$ be arbitrary.
We have by the definition of $\Delta_{1,\gamma}$ and
by Lemma \ref{lem:differencefunctionexplicit} that
\starteqn\label{eqn:phi1gammaexplicit}
\phi_1(\SL{2}{\Int[\vecti]}^{\phi_1}\gamma w)
=\phi_1(w)-\Delta_{1,\gamma}\circ\phi(w)=\phi_1(w)+
\log\left(||cz+d||^2\right).
\finisheqn
Using the description of the left
$\SL{2}{\Int[\vecti]}^{\phi_1}$ cosets in $\SL{2}{\Int[\vecti]}$
given in \eqref{eqn:sl2intcoset},
and \eqref{eqn:phi1gammaexplicit}, we see that it will
suffice to show that there are finitely many pairs $(c,d)\in\Int[\vecti]^2$
with $\mathrm{GCD}(c,d)=1$, such that
\starteqn\label{eqn:A3sl2inequality}
\phi_1(w)+2\log\left(||cz+d||^2\right)<B,\;\text{for any $w\in U_z$}.
\finisheqn
The condition $w\in U_{z}$ implies, by the definition
of $U_z$, that $-\phi_1(w)<-\phi_1(z)+\log 2$.
Therefore, the inequality in \eqref{eqn:A3sl2inequality} implies that
\startdisp
\log(\left(||cw+d||^2\right)<B-\phi_1(z)+\log 2=B',
\finishdisp
where $B'$ is defined as $B'=B-\phi_1(z)+\log 2$.
(The real number $B'$ is fixed throughout the argument
because both $B$ and $z$ are fixed.)
Since the $\exp$ function is monotone increasing,
our task is reduced to showing that there are only finitely many
pairs $(c,d)\in\Int[\vecti]^2$ with $\mathrm{GCD}(c,d)=1$ such that
\starteqn\label{eqn:minimizationfiniteness}
||cw+d||^2<\exp(B').
\finisheqn
Using \eqref{eqn:wcoordinates}, we calculate that
\starteqn\label{eqn:minimizationfiniteness1}
||cw+d||^2=||cx(w)+d||^2+||c||^2y(w)^2\geq ||c||^2y(w)^2.
\finisheqn
However, because we are assuming $w\in U_z$, we have $y(w)>\half y(z)^2$,
so that
\starteqn\label{eqn:Uconditionlowerbound}
||c||^2y(w)^2\geq \frac{1}{4}||c||^2y(z)^2.
\finisheqn
Applying \eqref{eqn:Uconditionlowerbound}
to \eqref{eqn:minimizationfiniteness1}, yields
\startdisp
||cw+d||^2>\frac{1}{4}||c||^2y(z)^2.
\finishdisp
Therefore, only those relatively prime pairs
$(c,d)\in\Int[\vecti]$ such that
$||c||^2<4\exp(B')/y(z)^2$ can satisfy \eqref{eqn:minimizationfiniteness}.
This demonstrates that there are only finitely many allowable
$c\in\Int[\vecti]$ such that at least one
$(c,d)\in\Int[\vecti]^2$ satisfies \eqref{eqn:minimizationfiniteness}.

Now, let such a $c\in\Int[\vecti]$ be fixed.
Then \eqref{eqn:minimizationfiniteness} implies that
\starteqn\label{eqn:cxwplusupperest}
\begin{aligned}
||cx(w)+d||^2&<&&\exp B'-||c||^2y(w)^2\\
&<&&\exp B'-\frac{1}{4}||c||^2y(z)^2\\
&=&&B'',
\end{aligned}
\finisheqn
where $B''\in\Real$ is defined as $\exp B'-\frac{1}{4}||c||^2y(z)^2$,
and is fixed throughout this argument because $B,c$ and $z$ are fixed.
On the other hand, we have
\starteqn\label{eqn:cxzplusdest}
\begin{aligned}
||cx(z)+d||&\geq&& ||cx(z)+d||-||cx(z)-cx(w)||\\
&\geq&&||cx(z)+d||-||c||,
\end{aligned}
\finisheqn
where the last inequality follows from the fact that $w\in U$.
Clearly, there are only $d\in\Int[\vecti]$ such that
\startdisp
||cx(z)+d||\leq \sqrt{|B''|}+||c||.
\finishdisp
Therefore, by \eqref{eqn:cxzplusdest}, there are only
finitely many $d$ for which
\startdisp
||cx(z)+d||^2\leq |B''|.
\finishdisp
Consequently, for a fixed $c$
there are only finitely many $d$ such that $(c,d)\in\Int[\vecti]^2$
satisfies \eqref{eqn:minimizationfiniteness}.  Since it has
already been established that there are only finitely
many $c\in\Int[\vecti]$ such that there exists any $d\in\Int[\vecti]$
with $(c,d)$ satisfying \eqref{eqn:minimizationfiniteness}, there
are only finitely many $(c,d)\in\Int[\vecti]^2$ satisfying
\eqref{eqn:minimizationfiniteness}.  By the comments
preceding \eqref{eqn:minimizationfiniteness}, this completes
the proof of the lemma.
\end{myproof}

\noindent\bmth\textbf{Explicit description of $\scrF_1$.}
\ubmth\hspace*{0.3cm}
We now give a collection of general results that will
be useful, for various concrete cases of $\Gamma$
in establishing a collection of
explicit inequalities defining $\scrF_1$.

Note that with $\gamma$ as in \eqref{eqn:gammaa4lemform}, we have
\starteqn\label{eqn:secondrowmultdescription}
\begin{pmatrix}0&1\end{pmatrix}\gamma=\begin{pmatrix}c&d\end{pmatrix},
\finisheqn
from which we deduce the following
\begin{cor} \label{cor:scrF1altdescr} Let $\scrF_1$ be as defined in Theorem \ref{thm:grenierh3}.
Then we have the following alternate descriptions of $\scrF_1$
and its interior.
\starteqn\label{eqn:scrF1altdescr}
\scrF_1=\left\{z\in\bbH^3\;\left|\; \left\|\begin{pmatrix}0&1\end{pmatrix}\gamma
\begin{pmatrix}z\\1\end{pmatrix}\right\|^2\geq 1,\;\text{for all}\;
\gamma\in\Gamma\right.\right\},
\finisheqn
and
\starteqn\label{eqn:scrF1intrraltdescr}
\intrr(\scrF_1)=\left\{z\in\bbH^3\;\left|\;
\left\|\begin{pmatrix}0&1\end{pmatrix}\gamma
\begin{pmatrix}z\\1\end{pmatrix}\right\|^2>1,\;\text{for all}\;
\gamma\in\Gamma-\Gamma^{\phi_1}\right.\right\}.
\finisheqn
\end{cor}
\begin{myproof}{Proof}
Using the definition of $\scrF_1$ and $\Delta_{1,\gamma}$, we have
\startdisp
\scrF_1=\{z\in\bbH^3\;|\; \Delta_{1,\gamma}(z)\leq 0\;\text{for all}\;
\gamma\in\Gamma\}.
\finishdisp
Using the expression for $\Delta_{1,\gamma}$
in \eqref{eqn:differencefunctionexplicit}, and
\eqref{eqn:secondrowmultdescription}, we obtain \eqref{eqn:scrF1altdescr}.

In order to obtain \eqref{eqn:scrF1intrraltdescr}, we
apply the same argument, but with strict equalities
replacing the non-strict inequalities and $\gamma$
ranging over $\Gamma-\Gamma^{\phi_1}$ instead of all of $\Gamma$.
\end{myproof}

The main technical result needed to calculate explicit
equalities for $\scrF_1$ in the case
of specific $\Gamma$ is Lemma \ref{lem:latticeandspheres}.  For
the lemma, we need the following elementary notion.
\begin{defn} \label{defn:epsilondense} Let $(X,d)$ be a metric space.  Let $\Xi$
be a subset of $X$, and let $\kappa>0$.  Then
$\Xi$ is said to be \bmth\textbf{$\kappa$-dense
in $X$}\ubmth\hspace*{0.05cm} if
\startdisp
\text{for all}\; x\in X,\;\text{there exists}\; \xi\in\Xi\;\text{such that}\,
d(x,\xi)\leq \kappa,
\finishdisp
in other words, when every open ball in $X$ of radius $\kappa$
intersects $\Xi$.
\end{defn}
We will apply Definition \ref{defn:epsilondense} most often
in the case when $(X,d)=(\Complex,||\cdot ||)$, the complex
numbers with the usual norm, and when $\Omega$ is a lattice
in $\Complex$, a translate of a lattice, or a similar discrete set.
We have the following elementary properties of $\kappa$-density
in $\Complex$.

\vspace*{0.3cm}
\noindent\begin{tabular}{p{4.3cm}p{10.2cm}}
\textbf{Translation Invariance}&  Let $x\in\Complex$.
The set $\Omega\subset \Complex$
is $\kappa$-dense in $\Complex$ if and only if $x+\Omega$ is
$\kappa$-dense in $\Complex$. \vspace*{0.1cm} \\
\textbf{Behavior under Dilation}&Let
$c\in\Complex-\{0\}$.  The set $\Omega\subset \Complex$
is $\kappa$-dense in $\Complex$ if and only if
then $c\Omega$ is $||c||\kappa$-dense in $\Complex$. \vspace*{0.1cm} \\
\textbf{Standard Lattice}&The standard lattice $\Int[\vecti]$
is $\frac{\sqrt{2}}{2}$-dense in $\Complex$.
\end{tabular}
\hspace*{0.3cm}

Before proceeding, we list some elementary properties of the
subsets of the type used to describe $\scrF_1$ in Corollary \ref{cor:scrF1altdescr},
which will be used repeatedly in manipulating such sets.
None of the elementary properties require proof.  Throughout
the properties, $\rho$ represents a nonzero real number unless stated
otherwise.

\vspace*{0.3cm}
\noindent\bmth
\textbf{Elementary Properties of the sets $\left\{z\in\bbH^3\;|\; \left\|
\begin{pmatrix}0&1\end{pmatrix}\gamma
\begin{pmatrix}z\\1\end{pmatrix}\right\|^2\geq\rho,\;\text{for all}\;
\gamma\in\Gamma\right\}$.}
\ubmth
\begin{itemize}
\item[\textbf{1.}]$\gamma=\begin{pmatrix}a&b\\c&d\end{pmatrix}\;\text{implies}\;
\begin{pmatrix}0&1\end{pmatrix}\gamma
\begin{pmatrix}z\\1\end{pmatrix}=\begin{pmatrix}1&d\end{pmatrix}
\begin{pmatrix}cz\\1\end{pmatrix}.$
\item[\textbf{2.}]$\left\{z\in\bbH^3\;\left|\; \left\|
\begin{pmatrix}0&1\end{pmatrix}\gamma
\begin{pmatrix}z\\1\end{pmatrix}\right\|\right.^2\geq\rho^2,\;\text{for all}\;\gamma\in
\Gamma\right\}=$\vspace*{2mm}\newline
\hspace*{3cm}
$\left\{z\in\bbH^3\;\left|\; \left\|
\begin{pmatrix}1&d\end{pmatrix}
\begin{pmatrix}cz\\1\end{pmatrix}\right\|\right.^2\geq\rho^2\;\text{for all}\;
\begin{pmatrix}c&d\end{pmatrix}\in\begin{pmatrix}0&1\end{pmatrix}\Gamma\right\}.$
\item[\textbf{3.}]  \textbf{Inclusion Reversal.}
If $\Omega_1\subseteq\Omega_2
\subseteq\Complex^2$, then
\begin{multline}
\left\{z\in\bbH^3\;\left|\; \left\|
\begin{pmatrix}1&d\end{pmatrix}
\begin{pmatrix}cz\\1\end{pmatrix}\right\|\right.^2\geq\rho^2,\;\text{for all}
\begin{pmatrix}c&d\end{pmatrix}\in\Omega_2
\right\}\vspace*{2mm}\notag\\
\hspace*{3cm}
\bigsubseteq\left\{z\in\bbH^3\;\left|\; \left\|
\begin{pmatrix}1&d\end{pmatrix}
\begin{pmatrix}cz\\1\end{pmatrix}\right\|\right.^2\geq\rho^2,\;\text{for all}
\begin{pmatrix}c&d\end{pmatrix}\in\Omega_1
\right\}.
\end{multline}
\item[\textbf{4.}]  \textbf{Constant Scalar Multiple.}  Let
$\rho\in\Complex^{\times}$.  Then
\startdisp
\left\|\frac{1}{\rho}\begin{pmatrix}1&d\end{pmatrix}
\begin{pmatrix}cz\\1\end{pmatrix}\right\|^2\geq 1\;\text{if and only if}
\left\|\begin{pmatrix}1&d\end{pmatrix}
\begin{pmatrix}cz\\1\end{pmatrix}\right\|^2\geq ||\rho||^2.
\finishdisp
\item[\textbf{5.}] \textbf{Scalar Function Multiple.}  Let $\rho(\cdot)$
be a non-vanishing complex-valued function on $\Gamma$.  Then
\startdisp
\left\|\frac{1}{\rho(\gamma)}\begin{pmatrix}1&d\end{pmatrix}
\begin{pmatrix}cz\\1\end{pmatrix}\right\|^2=
\frac{1}{||\rho||^2}\left\|\begin{pmatrix}1&d\end{pmatrix}
\begin{pmatrix}cz\\1\end{pmatrix}\right\|^2,
\finishdisp
so that
\begin{multline*}
\left\{z\in\bbH^3\;\left|\; \left\|\frac{1}{\rho(\gamma)}
\begin{pmatrix}1&0\end{pmatrix}\gamma
\begin{pmatrix}z\\1\end{pmatrix}\right\|\right.^2\geq 1,\;\text{for all}\;
\gamma\in\Gamma
\right\}\vspace*{2mm}\\
\hspace*{3cm}
=\left\{z\in\bbH^3\;\left|\; \left\|
\begin{pmatrix}1&0\end{pmatrix}\gamma
\begin{pmatrix}z\\1\end{pmatrix}\right\|\right.^2\geq||\rho(\gamma)||^2
,\;\text{for all}\;
\gamma\in\Gamma
\right\}.
\end{multline*}
\end{itemize}

In contrast to the elementary properties,
the property in Lemma \ref{lem:latticeandspheres}
does deserve a short proof.

\begin{lem}\label{lem:latticeandspheres}
Let $\Omega\subseteq\Complex$.
Assume that $\Omega$ is $\kappa$-dense in $\Complex$ under
the usual norm.  Let $\rho_1, \rho_2>0$ and let $c\in\Complex-\{0\}$.
Further assume that the following condition is satisfied,
\starteqn\label{eqn:relationrhois}
\rho_1^2-\kappa^2\geq\frac{\rho_2^2}{||c||^2}.
\finisheqn\label{eqn:latticeandspheres}
Then we have the containment
\begin{multline}\label{eqn:latticeandspheresstatement}
\left\{z\in\bbH^3\;\left|\;\left\|\begin{pmatrix}1&d
\end{pmatrix}\begin{pmatrix}z\\1
\end{pmatrix}\right\|^2\geq\rho_1^2\;\text{for all}\;
 d\in\Omega\right.\right\}\vspace*{2mm}\\
\hspace*{3cm}\bigsubseteq
\left\{z\in\bbH^3\;\left|\;\left\|\begin{pmatrix}1&\xi
\end{pmatrix}\begin{pmatrix}cz\\1
\end{pmatrix}\right\|^2\geq\rho_2^2\;\text{for all}\; \xi\in
\Complex\right.\right\}.
\end{multline}
\end{lem}
\begin{myproof}{Proof}
Let $z\in\bbH^3$, as in \eqref{eqn:z0a4lem}, be an element
of the left-hand side of \eqref{eqn:latticeandspheresstatement}.
Since $\Omega$ is $\kappa$-dense, there exists $d\in\Omega$,
such that
\starteqn\label{eqn:epsilondenseapplication}
||x(z)-d||^2\leq \kappa^2.
\finisheqn
On the other hand, because $z$ belongs to the left-hand side of
\eqref{eqn:latticeandspheresstatement}, we have
\starteqn\label{eqn:latticeandspheresapplication}
||z+d||^2=||x(z)+d||^2+y(z)^2\geq\rho_1^2.
\finisheqn
Combining the inequalities \eqref{eqn:epsilondenseapplication} and
\eqref{eqn:latticeandspheresapplication}, we deduce that
\startdisp
y(z)^2\geq\rho_1^2-\kappa^2.
\finishdisp
Therefore, by \eqref{eqn:relationrhois}, we have
\starteqn\label{eqn:latticeandspheresinter}
y(z)^2\geq\frac{\rho_2^2}{||c||^2}.
\finisheqn
Now, let $\xi\in\Complex$.  We have
\startdisp
\left\|\begin{pmatrix}1&\xi
\end{pmatrix}\begin{pmatrix}cz\\1
\end{pmatrix}\right\|^2=||cz+\xi||^2=||cx(z)+\xi||^2+||c||^2y(z)^2\geq
||c||^2y(z)^2\geq\rho_2^2,
\finishdisp
where we have applied \eqref{eqn:latticeandspheresinter}.
We have shown that $z$ belongs to the right side of
\eqref{eqn:latticeandspheresstatement},
and this completes the proof of the
containment \eqref{eqn:latticeandspheresstatement}.
\end{myproof}

We next explicitly state, in the form of
Corollaries \ref{cor:spheresandcspheres} and
\ref{cor:spheresandcspherescont}, the particular applications
 of Lemma \ref{lem:latticeandspheres}
that we will need when treating the particular
cases of $\Gamma=\SL{2}{\Int[\vecti]}$
and $\Gamma=\conj\inv(\SO{3}{\Int[\vecti]})$.
\begin{cor}\label{cor:spheresandcspheres}
Let $c\in\Int[\vecti]-\{0\}$.  Then we have
the following containment
\begin{multline}\label{eqn:spheresandcspheres}
\left\{z\in\bbH^3\;\left|\;
\left\|\begin{pmatrix}1&d\end{pmatrix}\begin{pmatrix}z\\
1\end{pmatrix}\right\|\geq 1,\text{for all}\;d\in\Int[\vecti]|\right.\right\}
\vspace*{0.2cm}\\
\hspace*{3cm}\bigsubseteq
\left\{z\in\bbH^3\;\left|\;
\left\|\begin{pmatrix}1&d\end{pmatrix}\begin{pmatrix}cz\\
1\end{pmatrix}\right\|\geq 1,\text{for all}\;d\in\Int[\vecti]\right.\right\}.
\end{multline}
\end{cor}
\begin{myproof}{Proof}
Since $c\in\Int[\vecti]-\{0\}$, we have two cases, namely $||c||=1$,
and $||c||^2\geq 2$.\vspace*{3mm}

\noindent\textit{Case 1: $||c||=1$.}  We have $c$ a unit,
so that $c\inv\in\Int[\vecti]$ and $c\inv\Int[\vecti]=\Int[\vecti]$.
Moreover, for $d\in\Int[\vecti]$, we have
\startdisp
\begin{pmatrix}1&d\end{pmatrix}\begin{pmatrix}cz\\ 1\end{pmatrix}=
\frac{1}{c\inv}
\begin{pmatrix}1&dc\inv\end{pmatrix}\begin{pmatrix}z\\ 1\end{pmatrix}
\finishdisp
We can apply Elementary Property \textbf{4, Constant Scalar Multiple},
with $\rho=c$, to the right side of \eqref{eqn:spheresandcspheres}
to conclude that we have equality in \eqref{eqn:spheresandcspheres}.

\vspace*{3mm}
\noindent\textit{Case 2: $||c||\geq 2$.}  The inequality
\eqref{eqn:relationrhois} is verified with $\rho_1=\rho_2=1$,
$\kappa^2=\half$ .  Therefore, by
the property \textbf{density of the standard
lattice}, we see that the hypotheses
of Lemma \eqref{lem:latticeandspheres}
are satisfied with $\Omega=\Int[\vecti]$.
Lemma \ref{lem:latticeandspheres} and the Elementary Property
\textbf{3, Inclusion Reversal} yield
\startdisp
\begin{aligned}
\left\{z\in\bbH^3\;\left|\;
\left\|\begin{pmatrix}1&d\end{pmatrix}\begin{pmatrix}z\\
1\end{pmatrix}\right\|\geq 1,\text{for all}\;d\in\Int[\vecti]\right.\right\}
&\bigsubseteq&&
\left\{z\in\bbH^3\;\left|\;
\left\|\begin{pmatrix}1&\xi\end{pmatrix}\begin{pmatrix}cz\\
1\end{pmatrix}\right\|\geq 1,\text{for all}\;\xi\in\Complex\right.\right\}\\
&\bigsubseteq&&\left\{z\in\bbH^3\;\left|\;
\left\|\begin{pmatrix}1&d\end{pmatrix}\begin{pmatrix}cz\\
1\end{pmatrix}\right\|\geq 1,\text{for all}\;d\in\Int[\vecti]\right.\right\}.
\end{aligned}
\finishdisp
\end{myproof}

For Corollary \ref{cor:spheresandcspherescont} define the following
subsets of $\bbH^3$.
\begin{gather}\label{eqn:scrFuppersdefn}
\scrF_1^{(1)}=\left\{z\in\bbH^3\;\left|\;\left\|
\begin{pmatrix}1&d\end{pmatrix}
\begin{pmatrix}cz\\ 1\end{pmatrix}\right\|^2\geq 1\;\text{for}
\; c\neq 0,\,d\in\Int[\vecti]
\right.\right\}.
\\
\scrF_1^{(2)}=\left\{z\in\bbH^3\;\left|\;\left\|
\begin{pmatrix}1&d\end{pmatrix}
\begin{pmatrix}cz\\ 1\end{pmatrix}\right\|^2\geq 2\;\text{for}
\;\begin{pmatrix}c&d\end{pmatrix}\hspace*{-1mm}\in\hspace*{-2mm}
\bigcup_{\stackrel{\omega\in\Int[\vecti]}{||\omega||=1}}
\hspace*{-5.5mm}\cdot\hspace*{5.5mm}
\hspace*{-2mm}
\begin{pmatrix}\omega&1+(1+\vecti)\Int[\vecti]\end{pmatrix}\bigcup
\hspace*{-3mm}\cdot
\bigcup_{\stackrel{c\in\Int[\vecti]}{||c||>1}}
\hspace*{-5.25mm}\cdot\hspace*{5.25mm}\hspace*{-2mm}
\begin{pmatrix}c&\Int[\vecti]\end{pmatrix}\right.\right\}.\notag
\\ \notag\\
\scrF_1^{(2,1)}=\left\{z\in\bbH^3\;\left|\;\left\|
\begin{pmatrix}1&d\end{pmatrix}
\begin{pmatrix}z\\ 1\end{pmatrix}\right\|^2\geq 2\;\text{for all}
\; d\in 1+(1+\vecti)\Int[\vecti]\right.\right\}.\notag
\end{gather}

\begin{cor}\label{cor:spheresandcspherescont}  Let
$\scrF_1^{(1)}, \scrF_1^{(2)}$ and $\scrF_1^{(2,1)}$,
be as defined in \eqref{eqn:scrFuppersdefn}.  Then
we have the inclusion
\startdisp
\scrF_1^{(2,1)}\subseteq \scrF_1^{(1)}\cap \scrF_1^{(2)}.
\finishdisp
\end{cor}
\begin{myproof}{Proof}
In order to prove that $\scrF_{1}^{(1,2)}\subseteq\scrF_1^{(1)}$,
apply Lemma \ref{lem:latticeandspheres} with $\rho_1^2=2$,
$\rho_2^2=1$, and $c\in\Int[\vecti]$, and $\kappa=1$.  Note
the set $1+(1+\vecti)\Int[\vecti]$ is a translation of
the dilation of the standard integer lattice by $1+\vecti$.  Therefore,
we can use the properties of $\kappa$-density to conclude
that $1+(1+\vecti)\Int[\vecti]$ is $1$-dense in $\Complex$.
Therefore, the hypotheses of Lemma \ref{lem:latticeandspheres}
are satisfied, and Lemma \ref{lem:latticeandspheres}
gives $\scrF_{1}^{(1,2)}\subseteq\scrF_1^{(1)}$.

In order
to prove that $\scrF_1^{1,2}\subseteq\scrF_1^{(2)}$,
apply Lemma \ref{lem:latticeandspheres} with $\rho_1^2=\rho_2^2=2$,
and $c\in\Int[\vecti]$ satisfying $||c||^2\geq 2$.  Lemma
\ref{lem:latticeandspheres} and the Elementary Property
\textbf{3, Inclusion reversal}, imply that
\starteqn\label{eqn:scrF1scontainment1}
\begin{aligned}
\scrF_1^{2,1}&\bigsubseteq&&\left\{z\in\bbH^3\;\left|\;\left\|
\begin{pmatrix}1&d\end{pmatrix}
\begin{pmatrix}cz\\ 1\end{pmatrix}\right\|^2\geq 2\;\text{for}
\;\begin{pmatrix}c&d\end{pmatrix}\hspace*{-1mm}\in
\bigcup_{\stackrel{c\in\Int[\vecti]}{||c||>1}}
\hspace*{-5.25mm}\cdot\hspace*{5.25mm}\hspace*{-2mm}
\begin{pmatrix}c&\Complex\end{pmatrix}\right.\right\}\\
&\bigsubseteq&&\left\{z\in\bbH^3\;\left|\;\left\|
\begin{pmatrix}1&d\end{pmatrix}
\begin{pmatrix}cz\\ 1\end{pmatrix}\right\|^2\geq 2\;\text{for}
\;\begin{pmatrix}c&d\end{pmatrix}\hspace*{-1mm}\in
\bigcup_{\stackrel{c\in\Int[\vecti]}{||c||>1}}
\hspace*{-5.25mm}\cdot\hspace*{5.25mm}\hspace*{-2mm}
\begin{pmatrix}c&\Int[\vecti]\end{pmatrix}\right.\right\}.
\end{aligned}
\finisheqn
Now, let $\omega\in\Int[\vecti]$, $||\omega||=1$, so that
$\omega\inv\in\Int[\vecti]$.
The Elementary Property \textbf{4, Constant Scalar Multiple},
applied with $\rho=\omega\inv$ yields
\startdisp
\left\|
\begin{pmatrix}1&d\end{pmatrix}
\begin{pmatrix}\omega z\\ 1\end{pmatrix}\right\|^2=
\left\|
\begin{pmatrix}1&\omega\inv d\end{pmatrix}
\begin{pmatrix}z\\ 1\end{pmatrix}\right\|^2.
\finishdisp
Thus, we have
\starteqn\label{eqn:scrF1scontainment2}
\begin{aligned}
\scrF_1^{(2,1)}&=&&
\left\{z\in\bbH^3\;\left|\;\left\|
\begin{pmatrix}1&d\end{pmatrix}
\begin{pmatrix}cz\\ 1\end{pmatrix}\right\|^2\geq 2\;\text{for}
\;\begin{pmatrix}c&d\end{pmatrix}\hspace*{-1mm}\in\hspace*{-2mm}
\bigcup_{\stackrel{\omega\in\Int[\vecti]}{||\omega||=1}}
\hspace*{-5.5mm}\cdot\hspace*{5.5mm}
\hspace*{-2mm}
\begin{pmatrix}\omega&\omega\inv\big(1+(1+\vecti)\Int[\vecti]\big)
\end{pmatrix}\right.\right\}\\
&=&&
\left\{z\in\bbH^3\;\left|\;\left\|
\begin{pmatrix}1&d\end{pmatrix}
\begin{pmatrix}cz\\ 1\end{pmatrix}\right\|^2\geq 2\;\text{for}
\;\begin{pmatrix}c&d\end{pmatrix}\hspace*{-1mm}\in\hspace*{-2mm}
\bigcup_{\stackrel{\omega\in\Int[\vecti]}{||\omega||=1}}
\hspace*{-5.5mm}\cdot\hspace*{5.5mm}
\hspace*{-2mm}
\begin{pmatrix}\omega&1+(1+\vecti)\Int[\vecti]
\end{pmatrix}\right.\right\},
\end{aligned}
\finisheqn
where the substitution in the latter line follows from the
fact that
\startdisp
\red_{1+\vecti}(\omega)=\red_{1+\vecti}(\omega\inv)=1.
\finishdisp
By combining \eqref{eqn:scrF1scontainment1} and \eqref{eqn:scrF1scontainment2},
we have $\scrF_{1}^{(2,1)}\subseteq\scrF_{1}^{(2)}$.
Since we have already shown that $\scrF_1^{(2,1)}\subseteq\scrF_1^{(1)}$,
this completes the proof of the corollary.
\end{myproof}

\noindent\bmth\textbf{Example: The Picard domain $\scrF$
for $\SL{2}{\Int[\vecti]}$.}
\ubmth\hspace*{0.3cm}
Define the following rectangle in $\Real^2$:
\startdisp
\scrG_{ \SL{2}{\Int[\vecti]}^{\phi_1}}
=\left\{(t_1,t_2)\in\Real^2\;\left|\; t_1\in\left[-\half,\half\right],\,
t_2\in\left[0,\half\right]\right.\right\}.
\finishdisp
It is easy to verify, from the explicit description of
$\SL{2}{\Int[\vecti]}^{\phi_1}$, given in \eqref{eqn:sl2intstabilizer}
that $\scrG_{ \SL{2}{\Int[\vecti]}^{\phi_1}}$ is a fundamental
domain for the action of $\SL{2}{\Int[\vecti]^{\phi_1}}/\{\pm 1\}$.
Further, it is obvious that
\startdisp
\scrG_{ \SL{2}{\Int[\vecti]}^{\phi_1}}=
\overline{\intrr(\scrG_{ \SL{2}{\Int[\vecti]}^{\phi_1}})}.
\finishdisp
By putting together Lemmas \eqref{lem:a1a2}, and \eqref{lem:a2}, and
Proposition \ref{prop:sl2zia3}, we see that the
action of $\SL{2}{\Int[\vecti]}$ satisfies the \textbf{A} axioms.
Therefore, Theorems \ref{thm:grenierh3} and \ref{thm:grenierh3extended}
apply.  We deduce that, with $\scrF_1$,
$\scrF(\scrG_{ \SL{2}{\Int[\vecti]}^{\phi_1}})$ defined as in Theorem
\ref{thm:grenierh3}, we have
\startdisp\text{
$\scrF:=\scrF(\scrG_{ \SL{2}{\Int[\vecti]}^{\phi_1}})$ is a good
Grenier fundamental domain for $\SL{2}{\Int[\vecti]}$.
}
\finishdisp
The fundamental domain $\scrF$ is defined
in \S\Roman{ranrom}.1 of \cite{jol05}, where, in keeping with
classical terminology, $\scrF$ is called the \textbf{Picard domain}.

In order to complete the example, we now give an
explicit description of the set $\scrF_1$, which will
allow the reader to see that ``our" $\scrF$ is exactly
the same as the Picard domain.  We claim that
$\scrF_1$ is the subset of $\Real^3$ whose image under
the diffeomorphism $\phi\inv$ is given as follows.
\starteqn\label{eqn:scrF1sl2zidescript}
\phi\inv(\scrF_1)=\{z\in\bbH^3\;|\; ||z-m||\geq 1,\;\text{for all}\,
m\in\Int[\vecti]\}.
\finisheqn
In order to verify the claim, we use
the expression for $\phi\inv(\scrF_1)$ given in Corollary \ref{cor:scrF1altdescr}.
Because
\startdisp
\begin{pmatrix}0&1\end{pmatrix}\SL{2}{\Int[\vecti]}=\{\begin{pmatrix}c&d
\end{pmatrix}\;|\;c,d,\in\Int[\vecti],\; \mathrm{GCD}(c,d)=1\},
\finishdisp
we can apply the Elementary Property \textbf{2} to obtain
\starteqn\label{eqn:sl2ziexample1}
\phi\inv(\scrF_1)=\left\{z\in\bbH^3\;\left|\; \left\|
\begin{pmatrix}1&d\end{pmatrix}
\begin{pmatrix}cz\\1\end{pmatrix}\right\|^2\geq 1,\;\text{for all}\;
(c,d)\in\Int[\vecti]^2,\; \text{with}\;\mathrm{GCD}(c,d)=1\right.\right\}.
\finisheqn
We have
\startdisp
(1,\Int[\vecti])\subseteq\{(c,d)\in\Int[\vecti]\;|\;(c,d)=1\}.
\finishdisp
So by Elementary Property \textbf{3, Inclusion Reversal},
the right-hand side of \eqref{eqn:sl2ziexample1} is
contained in the set
\starteqn\label{eqn:sl2ziexample2}
\left\{z\in\bbH^3\;\left|\; \left\|
\begin{pmatrix}1&d\end{pmatrix}
\begin{pmatrix}z\\1\end{pmatrix}\right\|^2\geq 1,\;\text{for all}\;
d\in\Int[\vecti]\right.\right\}
\finisheqn
On the other hand, we apply Corollary \ref{cor:spheresandcspheres}
to see that the set in \eqref{eqn:sl2ziexample2}
is contained in the right-hand side of \eqref{eqn:sl2ziexample1}.
Therefore, we can substitute
the set in \eqref{eqn:sl2ziexample2} for
the right-hand side of
\eqref{eqn:sl2ziexample1}.
The substitution yields \eqref{eqn:scrF1sl2zidescript}.

Of the infinite set of
inequalities defining $\scrF_1$, all except the one
with $d=0$, \textit{i.e.} $||z||^2\geq 1$, are
trivially satisfied on
$\phi_{[2,3]}\inv\left(\scrG_{ \SL{2}{\Int[\vecti]}^{\phi_1}}\right)$.
Thus, from \eqref{eqn:scrF1sl2zidescript} and \eqref{eqn:scrFdefn},
we recover the description of the Picard domain by
finitely many inequalities given in \S\Roman{ranrom}.1 of \cite{jol05}.

\subsection{The good Grenier fundamental domain for $\Gamma=
\conj\inv(\SO{3}{\Int[\vecti]})$}\label{subsec:gammafd}

We now proceed to consider the special case of
$\conj\inv(\SO{3}{\Int[\vecti]})$ in Theorems
\ref{thm:grenierh3} and \ref{thm:grenierh3extended} above.
In keeping with the general practice of this chapter,
we will go back to using $G$ to denote
$\SO{3}{\Complex}$ and $\Gamma$ to denote
$\SO{3}{\Int[\vecti]}$, exclusively.  Since we are always
in this section in the setting of subgroups of $\SL{2}{\Complex}$,
we will abuse notation slightly and use $\Gamma$ to denote
the isomorphic inverse image $\conj\inv(\Gamma)$ of
$\Gamma=\SO{3}{\Int[\vecti]}$ in $\SL{2}{\Complex}$.

Also, we treat $\Real^2$, the image of the projection
$\phi_{[2,3]}$, as $\Complex$, by identifying the point $(t_1,t_2)\in\Real^2$
with $t_1+\vecti t_2$.  Thus, our ``new" $\phi_{[2,3]}$ is defined
in terms of the ``old" $\phi$-coordinates by
\starteqn\label{eqn:phi23new}
\phi_{[2,3]}(z)=\phi_2(z)+\vecti\phi_3(z).
\finisheqn

We will give an expression for a good Grenier fundamental
domain $\scrF(\scrG)$ in terms of explicit inequalities,
in \eqref{eqn:scrFscrGfirstform}, and again as a
convex polytope in $\bbH^3$, in Proposition \ref{prop:gammafunddom},
below.

\vspace*{3mm}
\noindent\textbf{Statement of main results.}
\begin{prop}\label{prop:gammaaaxioms}
Let $\Gamma=\conj\inv(\SO{3}{\Int[\vecti]})$,
given as a set of fractional linear transformations
explicitly in Proposition \ref{prop:inversematrixexplicitdescription}.
Then $\Gamma$ satisfies the \textbf{A} axioms.
\end{prop}
\begin{myproof} {Proof} Since $\Gamma$ is a subgroup of $\SL{2}{\Int[\vecti]}$,
Lemma \ref{lem:a1a2} implies that $\Gamma$ satisfies Axioms\linebreak \textbf{A 1},
\textbf{A 2} and \textbf{A 4}.  By Lemma \ref{lem:indexlemma}, we have that
\startdisp
[\Gamma:\Gamma\cap\SL{2}{\Int[\vecti]}],\;[\SL{2}{\Int[\vecti]}:
\Gamma\cap\SL{2}{\Int[\vecti]}]<\infty,
\finishdisp
\textit{i.e.}, $\Gamma$ and $\SL{2}{\Int[\vecti]}$ are in the same
commensurability class.  We may apply Corollary \ref{cor:a3}
and \ref{prop:sl2zia3} to conclude that $\Gamma$ satisfies
\textbf{A 3}.
\end{myproof}

As a result of Proposition \ref{prop:gammaaaxioms},
we can apply Theorems \ref{thm:grenierh3} and \ref{thm:grenierh3extended}
to $\Gamma$.  The combined results are as follows.

\begin{thm}\label{thm:gammagoodgrenier}
Let $\Gamma=\conj\inv(\SO{3}{\Int[\vecti]})$, acting on $\bbH^3$
on the left by fractional linear transformations.  Otherwise,
the notation is as in Theorems \ref{thm:grenierh3}
and \ref{thm:grenierh3extended}.  Let $\scrG$ be a fundamental domain for the induced action of
$\Gamma^{\phi_{[2,3]}}/\{\pm 1\}$ on $\Real^2$.  Assume further
that $\scrG=\overline{\intrr(\scrG)}$.  Define
\starteqn\label{eqn:scrF1gammadefn}
\scrF_1=\{z\in\bbH^3\;|\; \phi_1(z)\leq \phi_1(\gamma z),\;
\text{for all}\; \gamma\in\Gamma\}.
\finisheqn
Set
\startdisp
\scrF(\scrG)=\phi\inv_{{[2,3]}}(\scrG)\cap\scrF_1.
\finishdisp
\begin{itemize}
\item[(a)]  We have $\scrF(\scrG)$ a fundamental domain for the action
of $\Gamma/\{\pm 1\}$ on $\bbH^3$.
\item[(b)]  We have
\starteqn\label{eqn:scrFtopology}
\scrF(\scrG)=\overline{\intrr\big(\scrF(\scrG)\big)}.
\finisheqn
\item[(c)]
Further, $\intrr\big(\scrF_1\big)$ and $\intrr\big(\scrF(\scrG)\big)$
have explicit descriptions as follows.
\starteqn\label{eqn:scrF1intrrdescription}
\intrr\big(\scrF_1\big)=
\{z\in\bbH^3\;|\; \phi_1(z)< \phi_1(\gamma z),\;
\text{for all}\; \gamma\in\Gamma-\Gamma^{\phi_1}\},
\finisheqn
and
\starteqn\label{eqn:scrFintrr}
\intrr\big(\scrF(\scrG)\big)=\phi_{[2,3]}\inv\big(\intrr(\scrG)\big)\cap
\intrr(\scrF_1),
\finisheqn
\end{itemize}
\end{thm}

In order to complete the concrete description of $\scrF(\scrG)$,
the principal remaining task is to describe $\scrF_1$ by
a set of explicit inequalities.  The other task, namely
giving an example of a suitable fundamental
domain $\scrG$ for the induced action of
$\Gamma^{\phi_{[2,3]}}/\{\pm \Id\}$, is much easier and will
be left for later.  Here, then, is the result concerning
$\scrF_1$ which we are now aiming for.

\begin{prop}\label{prop:scrF1forgamma}\bmth
\textbf{First form of $\scrF_1$.}\ubmth \hspace*{1mm} Let $\scrF_1$ be as
defined in \eqref{eqn:scrF1gammadefn}.  All other
notation has the same meaning as in Theorem \ref{thm:gammagoodgrenier}.
Then we have
\starteqn\label{eqn:scrF1gammaexplicit}
\scrF_1=\{z=x(z)+y(z)\vectj\in\bbH^3\;
|\; ||x(z)-d||^2+y(z)^2\geq 2,\;\text{for}\;
d\in1+(1+\vecti)\Int[\vecti]\},
\finisheqn
and $\intrr(\scrF_1)$ is the same as in \eqref{eqn:scrF1gammaexplicit},
but with strict inequality instead of nonstrict inequality.
\end{prop}
The proof of Proposition \ref{prop:scrF1forgamma} depends
on the lemmas in the previous section, and a sequence of elementary
lemmas, which we now give.
\begin{lem}\label{lem:xi2lem}  Let $r,s\in\Int[\vecti]$ with $\mathrm{GCD}(r,s)=1$, satisfying
\startdisp
\begin{pmatrix}r&s\end{pmatrix}\equiv \begin{pmatrix}1&0\end{pmatrix}\mod
1+\vecti,\quad \text{resp.,}\begin{pmatrix}1&1\end{pmatrix}\mod 1+\vecti.
\finishdisp
Then there exists a pair of integers $p,q\in\Int[\vecti]$ such that
\starteqn\label{eqn:determinantcompletion}
\begin{pmatrix}p&q\\r&s\end{pmatrix}\in\SL{2}{\Int[\vecti]},
\finisheqn
and satisfying
\starteqn\label{eqn:congruencecompletion}
\begin{pmatrix}p&q\\r&s\end{pmatrix}\equiv \begin{pmatrix}1&1\\
1&0\end{pmatrix}\mod
1+\vecti,\quad \text{resp.,}
\begin{pmatrix}1&0\\1&1\end{pmatrix}\mod 1+\vecti.
\finisheqn
\end{lem}
\begin{myproof}{Proof}  Since $\mathrm{GCD}(r,s)=1$,
there exist $p,q\in\Int[\vecti]$
satisfying \eqref{eqn:determinantcompletion}.  Then either $p,q,r,s$
satisfy \eqref{eqn:congruencecompletion}, in which case
we are done; or, $p,q,r,s$ satisfy
\starteqn\label{eqn:congruencecompletionnegation}
\begin{pmatrix}p&q\\r&s\end{pmatrix}\equiv \begin{pmatrix}0&1\\
1&0\end{pmatrix}\mod
1+\vecti,\quad \text{resp.,}\begin{pmatrix}0&1\\1&1\end{pmatrix}\mod 1+\vecti.
\finisheqn
If $p,q,r,s$ satisfy \eqref{eqn:congruencecompletionnegation}, then we
perform the elementary row operation of adding the second
row of $\left(\begin{smallmatrix}p&q\\r&s\end{smallmatrix}\right)$
to the first row in order to produce new $p,q$.  The
$p,q,r,s$ produced by the elementary row operation
satisfy both \eqref{eqn:determinantcompletion}
and \eqref{eqn:congruencecompletion}.
\end{myproof}

Using the definition of $\Xi_2$ in \eqref{eqn:residuematrices},
we deduce from Lemma \ref{lem:xi2lem} that
\starteqn\label{eqn:xi2secondrow}
\begin{pmatrix}0&1\end{pmatrix}\Xi_{2}=\red_{1+\vecti}\inv
\left(\left\{\begin{pmatrix}1&0\end{pmatrix},\;
\begin{pmatrix}1&1\end{pmatrix}\right\}\right)=
\begin{pmatrix}1+(1+\vecti)\Int[\vecti]& \Int[\vecti]\end{pmatrix}.
\finisheqn

\begin{lem}\label{lem:rowofgammadescrip}  Let the matrices
$\alpha^{\mathrm{N}}(m,x)$ be as in \eqref{eqn:MofNmxdefn},
the $\Xi$-subsets of $\SL{2}{\Int[\vecti]}$ be as in \eqref{eqn:xi12defn}
and \eqref{eqn:residuematrices}.
\begin{itemize}
\item[(a)]  We have the containment
\starteqn\label{eqn:nonintegersetrow}
\begin{pmatrix}1&1+(1+\vecti)
\Int[\vecti]\end{pmatrix}\hspace*{2mm}\bigsubseteq
\hspace*{2mm}
\begin{pmatrix}0&1\end{pmatrix}
\bigcup_{\epsilon=0,1}
\hspace*{-0.45cm}\cdot\hspace*{.3cm}
\Xi_2\alpha^{2\vecti}\hspace*{-0.7mm}
(\vecti,\vecti^{\epsilon})
\finisheqn
\item[(b)]  We have the containment
\begin{multline}\label{eqn:lastrowcontainmentb}
\begin{pmatrix}0&1\end{pmatrix}
\left(
\bigcup_{\delta=0,1}
\hspace*{-0.5cm}\cdot\hspace*{.3cm}
\left(\Xi_{12}\alpha^{\vecti^{\delta}}\hspace{-0.5mm}
(\vecti^{\delta},0)\bigcup\hspace*{-0.32cm}
\cdot\hspace*{0.12cm}
\left(\bigcup_{\epsilon=0,1}\hspace*{-0.45cm}\cdot\hspace*{.25cm}
\frac{1}{1+\vecti}
\Xi_{2}\alpha^{2\vecti^{1+\delta}}\hspace*{-0.7mm}
(\vecti^{1+\delta},\vecti^{\epsilon})
\right)\right)\right)\\
\hspace*{3cm}\bigsubseteq
\begin{pmatrix}\Int[\vecti]&\Int[\vecti]\end{pmatrix}
\bigcup\hspace*{-3.5mm}\cdot\hspace*{3mm}\frac{1}{1+\vecti}
\left(\bigcup_{\stackrel{\omega\in\Int[\vecti]}{||\omega||=1}}
\hspace*{-.55cm}\cdot\hspace*{.35cm}
\begin{pmatrix}\omega&1+(1+\vecti)\Int[\vecti]\end{pmatrix}
\bigcup\hspace{-3mm}\cdot\hspace*{1.5mm}
\bigcup_{\stackrel{c\in\Int[\vecti]}{||c||>1}}
\hspace*{-.55cm}\cdot\hspace*{.35cm}
\begin{pmatrix}c&\Int[\vecti]\end{pmatrix}
\right).
\end{multline}
\end{itemize}
\end{lem}
\begin{myproof}{Proof}
For (a), since $\red_{1+\vecti}(-\vecti)=1$, \eqref{eqn:xi2secondrow} implies
that
\startdisp
\begin{pmatrix}-\vecti&\Int[\vecti]\end{pmatrix}\subseteq
\begin{pmatrix}0&1\end{pmatrix}\Xi_{2}.
\finishdisp
Thus,
\starteqn\label{eqn:lastrowcontainmentinter}
\bigcup_{\epsilon=0,1}
\begin{pmatrix}-\vecti&\Int[\vecti]\end{pmatrix}
\alpha^{2\vecti}\hspace*{-0.7mm}
(\vecti,\vecti^{\epsilon})
\bigsubseteq
\begin{pmatrix}0&1\end{pmatrix}
\bigcup_{\epsilon=0,1}
\Xi_2\alpha^{2\vecti}\hspace*{-0.7mm}
(\vecti,\vecti^{\epsilon}).
\finisheqn
Using \eqref{eqn:MofNmxdefn},
we calculate
the left-hand side of \eqref{eqn:lastrowcontainmentinter},
and we obtain
\startdisp
\bigcup_{\epsilon=0,1}
\begin{pmatrix}-\vecti&\Int[\vecti]\end{pmatrix}
\alpha^{2\vecti}\hspace*{-0.7mm}
(\vecti,\vecti^{\epsilon})=\bigcup_{\epsilon=0,1}
\begin{pmatrix}1&\vecti^{\epsilon-1}+2\Int[\vecti]
\end{pmatrix}=\begin{pmatrix}1&\bigcup_{\epsilon=0,1}
\vecti^{\epsilon-1}+2\Int[\vecti]\end{pmatrix}.
\finishdisp
So we have
\starteqn\label{eqn:lastrowcontainmentinter2}
\begin{pmatrix}1&\bigcup_{\epsilon=0,1}
\vecti^{\epsilon-1}+2\Int[\vecti]\end{pmatrix}\;\text{is contained
in the right side of \eqref{eqn:lastrowcontainmentinter}.}
\finisheqn
It is easily verified that
\startdisp
\bigcup_{\epsilon=0,1}\vecti^{\epsilon}+2\Int[\vecti]=
1+(1+\vecti)\Int[\vecti].
\finishdisp
Therefore, \eqref{eqn:lastrowcontainmentinter2} implies that
\startdisp
\begin{pmatrix}
1&1+(1+\vecti)\Int[\vecti]\end{pmatrix}\;\text{is contained
in the right side of \eqref{eqn:lastrowcontainmentinter}.}
\finishdisp
This completes the proof of part (a) of the lemma.

We now verify part (b).  We have by the definitions that
$\alpha^{\mathrm{N}}\hspace*{-.7mm}(m,x)$ and the subsets
$\Xi_1,\,\Xi_2,\,\Xi_{12}$
of $\SL{2}{\Int[\vecti]}$
are subsets of $\mathrm{Mat}_{2}(\Int[\vecti])$.
Therefore,
\starteqn\label{eqn:lastrowcontainmentb1inter}
\Xi\alpha^{\mathrm{N}}\hspace*{-.7mm}(m,x)
\subseteq\mathrm{Mat}_{2}(\Int[\vecti])\;\text{for each $\Xi$-subset,
$m,\,N\in\Int[\vecti]$, $m|N$, $x\in\Omega_{\frac{N}{m}}$}.
\finisheqn
We clearly have
\startdisp
\begin{pmatrix}0&1\end{pmatrix}\mathrm{Mat}_2(\Int[\vecti])=
\begin{pmatrix}\Int[\vecti]&\Int[\vecti]\end{pmatrix},
\finishdisp
Therefore the containment \eqref{eqn:lastrowcontainmentb1inter} implies that
\starteqn\label{eqn:lastrowcontainmentb1}\begin{pmatrix}0&1\end{pmatrix}
\bigcup_{\delta,=0,1}
\hspace*{-0.5cm}\cdot\hspace*{.5cm}
\Xi_{12}\alpha^{\vecti^{\delta}}\hspace{-0.5mm}(\vecti^{\delta},0)
\subseteq
\begin{pmatrix}\Int[\vecti]&\Int[\vecti]\end{pmatrix},
\finisheqn
and
\startdisp
\begin{pmatrix}0&1\end{pmatrix}
\bigcup_{\delta,\epsilon=0,1}\hspace*{-0.55cm}\cdot\hspace*{.55cm}
\frac{1}{1+\vecti}
\Xi_{2}\alpha^{2\vecti^{1+\delta}}\hspace*{-0.7mm}
(\vecti^{1+\delta},\vecti^{\epsilon})\bigsubseteq\frac{1}{1+\vecti}
\begin{pmatrix}\Int[\vecti]&\Int[\vecti]\end{pmatrix}
\finishdisp
Decomposing $\Int[\vecti]$ into disjoint sets of square norm $0,1$
and $2$ or greater, we may rewrite this as
\starteqn
\label{eqn:lastrowcontainmentb2}
\bigcup_{\delta,\epsilon=0,1}\hspace*{-0.55cm}\cdot\hspace*{.55cm}
\frac{1}{1+\vecti}
\Xi_{2}\alpha^{2\vecti^{1+\delta}}\hspace*{-0.7mm}
(\vecti^{1+\delta},\vecti^{\epsilon})\bigsubseteq
\frac{1}{1+\vecti}\left(\begin{pmatrix}0&\Int[\vecti]\end{pmatrix}
\bigcup\hspace*{-3mm}\cdot
\bigcup_{\stackrel{\omega\in\Int[\vecti]}{||\omega||=1}}
\hspace*{-5.5mm}\cdot\hspace*{2mm}
\begin{pmatrix}\omega&\Int[\vecti]\end{pmatrix}\bigcup
\hspace*{-3mm}\cdot
\bigcup_{\stackrel{c\in\Int[\vecti]}{||c||>1}}
\hspace*{-5mm}\cdot\hspace*{2mm}
\begin{pmatrix}c&\Int[\vecti]\end{pmatrix}\right).
\finisheqn
We next claim that, for $\epsilon,\,\delta\in\{0,1\}$,
\starteqn\label{eqn:lastrowclaim}
\left\{\begin{pmatrix}c&d\end{pmatrix}
\in\begin{pmatrix}0&1\end{pmatrix}\Xi_{2}
\alpha^{2\vecti^{1+\delta}}\hspace*{-0.7mm}
(\vecti^{1+\delta},\vecti^{\epsilon})\;|\;||c||\leq 1\right\}\subseteq
\bigcup_{\stackrel{\omega\in\Int[\vecti]}{||\omega||=1}}\hspace*{-1.5mm}
\begin{pmatrix}\omega&1+(1+\vecti)\Int[\vecti]\end{pmatrix}.
\finisheqn
Together with \eqref{eqn:lastrowcontainmentb2}, \eqref{eqn:lastrowclaim}
implies that
\starteqn\label{eqn:lastrowcontainmentb3}
\bigcup_{\delta,\epsilon=0,1}\hspace*{-0.55cm}\cdot\hspace*{.55cm}
\frac{1}{1+\vecti}
\Xi_{2}\alpha^{2\vecti^{1+\delta}}\hspace*{-0.7mm}
(\vecti^{1+\delta},\vecti^{\epsilon})\bigsubseteq
\frac{1}{1+\vecti}\left(
\bigcup_{\stackrel{\omega\in\Int[\vecti]}{||\omega||=1}}
\hspace*{-5.5mm}\cdot\hspace*{2mm}
\begin{pmatrix}\omega&1+(1+\vecti)\Int[\vecti]\end{pmatrix}\bigcup
\hspace*{-3mm}\cdot
\bigcup_{\stackrel{c\in\Int[\vecti]}{||c||>1}}
\hspace*{-5mm}\cdot\hspace*{2mm}
\begin{pmatrix}c&\Int[\vecti]\end{pmatrix}\right).
\finisheqn
Together, \eqref{eqn:lastrowcontainmentb1} and \eqref{eqn:lastrowcontainmentb3}
imply \eqref{eqn:lastrowcontainmentb}.  Therefore, we
have reduced the proof of part (b) of the lemma to the proof of
\eqref{eqn:lastrowclaim}.

In order to verify the claim
\eqref{eqn:lastrowclaim}, we first use \eqref{eqn:xi2secondrow}
to calculate that
\starteqn\label{eqn:lastrowcontainmentb4}
\begin{aligned}
\begin{pmatrix}0&1\end{pmatrix}
\Xi_{2}\alpha^{2\vecti^{1+\delta}}\hspace*{-0.7mm}(\vecti^{1+\delta}
,\vecti^{\epsilon})
&=&&\begin{pmatrix}1+(1+\vecti)\Int[\vecti]&\Int[\vecti]\end{pmatrix}
\alpha^{2\vecti^{1+\delta}}\hspace*{-0.7mm}
(\vecti^{1+\delta},\vecti^{\epsilon})\\
&=&&\big(1+(1+\vecti)\Int[\vecti]\big)\begin{pmatrix}\vecti^{1+\delta}&\vecti^{\epsilon}
\end{pmatrix}+\begin{pmatrix}0&2\Int[\vecti]\end{pmatrix}.
\end{aligned}
\finisheqn
So let $\begin{pmatrix}c&d\end{pmatrix}$ be as in the left-hand side of
\eqref{eqn:lastrowclaim}.  By \eqref{eqn:lastrowcontainmentb4},
we have
\starteqn\label{eqn:lastrowcontainmentb5}
\begin{pmatrix}c&d\end{pmatrix}\in
\begin{pmatrix}\omega'\vecti^{1+\delta}&\omega'\vecti^{\epsilon}+
2\Int[\vecti]\end{pmatrix}\;\text{for}\;\omega'\in1+(1+\vecti)\Int[\vecti].
\finisheqn
So $c=\omega'\vecti^{1+\delta}$.
Since $||c||\leq1$, we deduce that $||\omega'||\leq1$.  But we also
have $\omega\in 1+(1+\vecti)\Int[\vecti]$, so that
$\omega\neq 0$.  Thus,
$||\omega'||=1$.  Set $\omega=\omega'\vecti^{1+\delta}$,
so that $c=\omega$.  Note that
\startdisp
d\in\omega\vecti^{\epsilon-1-\delta}+
2\Int[\vecti]=1+2\Int[\vecti].
\finishdisp
Summing up these observations, we see that
\eqref{eqn:lastrowcontainmentb5} can be rewritten as
\startdisp
\begin{pmatrix}c&d\end{pmatrix}\in
\begin{pmatrix}\omega&1+2\Int[\vecti]\end{pmatrix}
\text{with}\;\omega\in\Int[\vecti],\,||\omega||=1.
\finishdisp
Since $\begin{pmatrix}c&d\end{pmatrix}$ was an arbitrary
element of the left-hand side of \eqref{eqn:lastrowclaim},
this completes the verification of \eqref{eqn:lastrowclaim}.
By the comments immediately following \eqref{eqn:lastrowcontainmentb3},
this completes the proof of part (b) of the lemma.
\end{myproof}

For the proof of Lemma \ref{lem:scrF1determination}, below,
we will need to introduce the function
\startdisp
\omega_8^{\delta}:\Gamma\rightarrow\bbS^1.
\finishdisp
By the description of $\Gamma$ given in
\eqref{eqn:inversematrixexplicitdescription},
there is a well-defined function $\omega_8^{\delta}$
on $\Gamma$ defined by
\starteqn\label{eqn:omegaeightdeltadefn}
\omega_8^{\delta}(\gamma)=\begin{cases}1&\text{for}\;\gamma\in\mathrm{Mat}_2
(\Int[\vecti]),
\\
\omega_8&\text{for}\;\gamma\in\omega_8\mathrm{Mat}_2
(\Int[\vecti]).
\end{cases}
\finisheqn
Recalling the function $\big(i,\delta\big)$ on $\Gamma$
defined at \eqref{eqn:ideltadefn}, we see that $\omega_8^{\delta}$
may also be defined by
\startdisp
\omega_8^{\delta}:=\omega_8^{\mathrm{pr}_2\big(i,\delta\big)},
\finishdisp
where $\mathrm{pr}_2$ denotes projection onto the second factor.

Applying Elementary Property \textbf{5, Scalar Function Multiple},
with $\rho=\omega_8^{\delta}$, so that $||\omega_8^{\delta}=1$,
we deduce that
\begin{multline}\label{eqn:scalarfunctionpropapplied}
\left\{z\in\bbH^3\;\left|\; \left\|\frac{1}{\omega_8^{\delta}(\gamma)}
\begin{pmatrix}1&0\end{pmatrix}\gamma
\begin{pmatrix}z\\1\end{pmatrix}\right\|\right.^2\geq 1,\;\text{for all}\;
\gamma\in\Gamma
\right\}\vspace*{2mm}\\
\hspace*{3cm}
=\left\{z\in\bbH^3\;\left|\; \left\|
\begin{pmatrix}1&0\end{pmatrix}\gamma
\begin{pmatrix}z\\1\end{pmatrix}\right\|\right.^2\geq1
,\;\text{for all}\;
\gamma\in\Gamma
\right\}.
\end{multline}

Recall from \eqref{eqn:scrFuppersdefn} the subsets $\scrF_1^{(1)},
\,\scrF_1^{(2)},\,\scrF_1^{(1,2)}$ of $\bbH^3$,
which we we now relate to $\scrF_1$.
\begin{lem}\label{lem:scrF1determination}
\begin{itemize}
\item[(a)]  We have the inclusions
\starteqn\label{eqn:scrF1inclusionschain}
\scrF_1^{(2)}\cap\scrF_1^{(1)}\subseteq\scrF_1\subseteq\scrF_1^{(2,1)}.
\finisheqn
\item[(b)]  Each inclusion in \eqref{eqn:scrF1inclusionschain}
can be replaced by an equality, so that in particular, we have
\startdisp
\scrF_1=\scrF_1^{(2,1)}.
\finishdisp
\end{itemize}
\end{lem}
\begin{myproof}{Proof}  For (a),
by \eqref{eqn:inversematrixexplicitdescription} and \eqref{eqn:omegaeightdeltadefn}
we have
\starteqn\label{eqn:omegadeltaGamma}
\left\{\left.
\frac{1}{\omega_8^{\delta}(\gamma)}\gamma\;\right|\;\gamma\in\Gamma\right\}
=\bigcup_{\delta,=0,1}
\hspace*{-0.5cm}\cdot\hspace*{.3cm}
\left(\Xi_{12}\alpha^{\vecti^{\delta}}\hspace{-0.5mm}
(\vecti^{\delta},0)\bigcup\hspace*{-0.32cm}
\cdot\hspace*{0.18cm}
\left(\bigcup_{\epsilon=0,1}\hspace*{-0.45cm}\cdot\hspace*{.25cm}
\frac{1}{1+\vecti}
\Xi_{2}\alpha^{2\vecti^{1+\delta}}\hspace*{-0.7mm}
(\vecti^{1+\delta},\vecti^{\epsilon})
\right)\right).
\finisheqn

By applying \eqref{eqn:scalarfunctionpropapplied}
to the right-hand side of \eqref{eqn:scrF1altdescr},
we have
\starteqn\label{eqn:scrF1incaseofgamma}
\scrF_1=\left\{z\in\bbH^3\;\left|\; \left\|
\begin{pmatrix}0&1\end{pmatrix}\frac{1}{\omega_8^{\delta}(\gamma)}\gamma
\begin{pmatrix}z\\1\end{pmatrix}\right\|^2\geq 1,\;\text{for}\;
\frac{1}{\omega_8^{\delta}(\gamma)}\gamma\;\text{of the form in \eqref{eqn:omegadeltaGamma}}
\right.\right\}.
\finisheqn
The containments
\eqref{eqn:nonintegersetrow} and \eqref{eqn:lastrowcontainmentb}
applied to \eqref{eqn:omegadeltaGamma}
imply that
\begin{multline}\label{eqn:multicontainments}
\begin{pmatrix}1&1+(1+\vecti)
\Int[\vecti]\end{pmatrix}\hspace*{2mm}\\ \bigsubseteq
\left\{\left.\begin{pmatrix}0&1\end{pmatrix}
\frac{1}{\omega_8^{\delta}(\gamma)}\gamma
\;\right|\;
\frac{1}{\omega_8^{\delta}(\gamma)}\gamma\;
\text{of the form in \eqref{eqn:omegadeltaGamma}}\right\}\\
\bigsubseteq
\begin{pmatrix}\Int[\vecti]&\Int[\vecti]\end{pmatrix}
\bigcup\hspace*{-3.25mm}\cdot\hspace*{2.75mm}\frac{1}{1+\vecti}
\left(\bigcup_{\stackrel{\omega\in\Int[\vecti]}{||\omega||=1}}
\hspace*{-.55cm}\cdot\hspace*{.35cm}
\begin{pmatrix}\omega&1+(1+\vecti)\Int[\vecti]\end{pmatrix}
\bigcup\hspace{-3mm}\cdot\hspace*{1.5mm}
\bigcup_{\stackrel{c\in\Int[\vecti]}{||c||>1}}
\hspace*{-.55cm}\cdot\hspace*{.35cm}
\begin{pmatrix}c&\Int[\vecti]\end{pmatrix}
\right),
\end{multline}
which, by \eqref{eqn:scrFuppersdefn}, is part (a).

Part (b) is a restatement of Corollary \ref{cor:spheresandcspherescont}.
\end{myproof}

\begin{myproof}{Completion of Proof of Proposition \ref{prop:scrF1forgamma}}
By Lemma \ref{lem:scrF1determination},
$\scrF_1=\scrF_1^{2,1}$.  By applying the relation
\startdisp
\left\|\begin{pmatrix}1&d\end{pmatrix}
\begin{pmatrix}cz\\1\end{pmatrix}\right\|^2=
||cz(x)+d||^2+y(z)^2,
\finishdisp
we can readily rewrite $\scrF_1^{2,1}$ in the form given on the right-hand
side of \eqref{eqn:scrF1gammaexplicit}.
\end{myproof}
\vspace*{3mm}

\noindent\bmth\textbf{Fundamental domain $\scrG$
for $\Gamma^{\phi_1}$.}\ubmth\hspace*{3mm}
In order to complete the explicit determination of a good
Grenier fundamental domain $\scrF$ for $\Gamma$, it remains
to give describe a suitable fundamental domain $\scrG$ for $\Gamma^{\phi_1}$.
Using \eqref{eqn:stabilizerintersection}, \eqref{eqn:firstcoordfixing},
and the description of $\Gamma$ in
\eqref{eqn:inversematrixexplicitdescription}
we deduce that
\starteqn\label{eqn:gammaphi1explicit}
\Gamma^{\phi_1}=\left\{\left.
\begin{pmatrix}\omega_8^{\delta}&\omega_8^{\delta}b\\
0&\omega_8^{-\delta}\end{pmatrix}\;\right|\; b\in(1+\vecti)\Int[\vecti],\;
\delta\in\{0,1\}\right\}.
\finisheqn
It follows from \eqref{eqn:gammaphi1explicit} that the subgroup
of unipotent elements of $\Gamma^{\phi_1}$ is
\starteqn\label{eqn:gammaphi1Uexplicit}
(\Gamma^{\phi_1})_U=\begin{pmatrix}1&(1+\vecti)\Int[\vecti]\\0&1\end{pmatrix}.
\finisheqn
We make note of certain group-theoretic properties of $\Gamma^{\phi_1}$
and $(\Gamma^{\phi_1})_U$ that will be used in determining
the fundamental domains.
First, we define the following generating elements:
\starteqn\label{eqn:gammagensdefn}
R_{\frac{\pi}{2}}=\begin{pmatrix}\omega_8&0\\0&
\omega_8\inv\end{pmatrix},\;
T_{1+\vecti}=\begin{pmatrix}1&1+\vecti\\0&1\end{pmatrix},\;\text{and}\;
T_{1-\vecti}=\begin{pmatrix}1&1-\vecti\\0&1\end{pmatrix}.
\finisheqn
It is easily verified, using \eqref{eqn:gammaphi1explicit}
and \eqref{eqn:gammaphi1Uexplicit}, that
\starteqn\label{eqn:generatingsets}
(\Gamma^{\phi_1})_U=\langle T_{1+\vecti},\, T_{1-\vecti}\rangle,\quad
\Gamma^{\phi_1}=\langle R_{\frac{\pi}{2}}, \,T_{1+\vecti},
\,T_{1-\vecti}\rangle.
\finisheqn
We calculate, from the definition of $R_{\frac{\pi}{2}}$
and \eqref{eqn:generatingsets}, that
\startdisp
\conj(R_{\frac{\pi}{2}})(\Gamma^{\phi_1})_U=(\Gamma^{\phi_1})_U.
\finishdisp
Since $\Gamma^{\phi_1}$ is generated by $(\Gamma^{\phi_1})_U$
and $R_{\frac{\pi}{2}}$, and $R_{\frac{\pi}{2}}$ has order
$4$, we deduce that
\starteqn\label{eqn:Unormality}\text{
$(\Gamma^{\phi_1})_U$ is normal in $\Gamma^{\phi_1}$ with
$[\Gamma^{\phi_1}:(\Gamma^{\phi_1})_U]=4$.}
\finisheqn
Let $T$ be any element of $(\Gamma^{\phi_1})_U$.
Then we have a more precise version of \eqref{eqn:Unormality},
\starteqn\label{eqn:cyclicgpcosetreps}\text{The
group $\langle TR_{\frac{\pi}{2}}\rangle$ of order $4$
is a set of representatives for the coset group $\Gamma^{\phi_1}/(\Gamma^{\phi_1})_U$.}
\finisheqn
Applying \eqref{eqn:cyclicgpcosetreps} to the case $T=T_{1-\vecti}$, we have
\starteqn\label{eqn:cyclicgpcosetrepsT1plusi}\text{The
group $\langle T_{1-\vecti}R_{\frac{\pi}{2}}\rangle$ of order $4$
is a set of representatives for the coset group $\Gamma^{\phi_1}/(\Gamma^{\phi_1})_U$.}
\finisheqn

It is easily verified that the action of $R_{\frac{\pi}{2}}$
on $\Complex$ is rotation by an angle $\pi/2$ about the
fixed point $0$.  Furthermore, calculate from \eqref{eqn:gammagensdefn}
that
\startdisp
T_{1-\vecti}R_{\frac{\pi}{2}}=\conj(T_{1})R_{\frac{\pi}{2}}.
\finishdisp
Therefore,
\starteqn\label{eqn:rotationdescription}\text{The action
of $T_{1-\vecti}R_{\frac{\pi}{2}}$ on $\Complex$
is rotation by $\pi/2$ about $1$.}
\finisheqn

By \eqref{eqn:Unormality} and \eqref{eqn:cyclicgpcosetrepsT1plusi},
we are in the situation of Lemma \ref{lem:finiteindexfunddom},
below, if we take $X=\bbH^3$,
$\Gamma=(\Gamma^{\phi_1})_U$
and $\tilde{\Gamma}=\Gamma^{\phi_1}$,
$\Delta=\langle T_{1-\vecti}R_{\frac{\pi}{2}}\rangle.$
\begin{lem}\label{lem:finiteindexfunddom}  Let $\Gamma\subseteq\tilde{\Gamma}$
be discrete subgroups acting topologically on a manifold $X$,
such that $\Gamma$ is normal in $\tilde{\Gamma}$, and the action of $\Gamma$
has a fundamental domain $\scrF$.  Suppose further that $\Delta$
is a set of coset representatives of $\tilde{\Gamma}/\Gamma$ in $\Gamma$
such that $\Delta$ forms a group and such that the action
of $\Delta$ on $X$ preserves $\scrF$.  Let $\tilde{\scrF}$ be a
fundamental domain for the action of $\Delta$ on $\scrF$.
Then $\tilde{\scrF}$ is a fundamental domain for the action of $\tilde{\Gamma}$
on $X$.
\end{lem}
For the proof of Lemma \ref{lem:finiteindexfunddom}, which
is not difficult, see the proof
of Lemma 2.1.13 in \cite{brennerthesis}.  Applied to the case at hand,
Lemma \ref{lem:finiteindexfunddom} yields the following statement.
\begin{cor}\label{cor:finiteindexfunddom}
Let $\scrG_U$ be a
fundamental domain for the action of $(\Gamma^{\phi_1})_{U}$
on $\bbH^3$, satisfying
\startdisp
T_{1+\vecti}R_{\frac{\pi}{2}}(\scrG_U)=\scrG_U.
\finishdisp
Let $\scrG$ be a fundamental domain
for the action of $\langle T_{1+\vecti}R_{\frac{\pi}{2}}\rangle$
on $\scrG$.
Then $\scrG$ a fundamental domain for the action of
$\Gamma^{\phi_1}$ on $\bbH^3$.
\end{cor}

In order to define and work with the sets $\scrG_U$ and $\scrG$
which will be fundamental domains for the action of
$\Gamma^{\phi_1}_U$ and $\Gamma^{\phi_1}$, it is useful to introduce
the notion of a convex hull in a totally geodesic metric space.

A metric space $(X,d)$ will be called \textbf{totally geodesic} if
for every pair of points $p_1,p_2\in X$, $p_1\neq p_2$ there is a unique
geodesic segment connecting $p_1,p_2$.  In this situation,
the (closed) geodesic segment connecting
$p_1,p_2$ will be denoted $[p_1,p_2]_{d}$.  A point $x\in X$
is said to lie \bmth\textbf{between $p_1$ and $p_2$}\ubmth\hspace*{0.5mm}
when $x$ lies on $[p_1,p_2]_{d}$.  We then say that
$\scrS\subset X$ is \textbf{convex} when $p_1,\,p_2\in\scrC$
and $p_3$ between $p_1$ and $p_2$ implies
that $p_3\in\scrS$.  Let $p_1,\ldots,p_r$ be $r$ points
in $X$.  The points determine a set
\startdisp
\scrC_{d}(p_1,\ldots, p_r)
\finishdisp
called the \textbf{convex closure} of $p_1,\ldots,p_r$,
described as the smallest convex subset of $X$ containing the set
$\{p_1,\ldots,p_r\}$.

Obviously, we can apply the notion of convex hull to any set $\scrS$,
rather than a finite set of points.  The definition remains the same,
namely that $\scrC_d(\scrS)$ is the smallest convex subset of $X$
containing $\scrS$.  In general we will use the notation
\startdisp
\scrC_d(\scrS_1,\ldots,\scrS_r)=\scrC_d\left(\bigcup_{i=1,\ldots r}
\scrS_i\right).
\finishdisp

We record the following elementary properties of convex hulls
\begin{itemize}
\item[\textbf{CH 1}]  $\scrS_1\subseteq\scrS_2$ implies that $\scrC_d(\scrS_1)
\subseteq\scrC_d(\scrS_2)$.
\item[\textbf{CH 2}]  $\scrS\subseteq\scrC_d(\scrS)$; $\scrC_d(\scrC_d(\scrS))=
\scrC_d(\scrS)$; $\scrS_1\subseteq\scrC_d(\scrS_2)$ implies
$\scrC_d(\scrS_1)\subseteq\scrC_d(\scrS_2)$.
\item[\textbf{CH 3}]  $\scrC_d(\scrC_d(\scrS_1),\scrS_2)=\scrC_d(\scrS_1,\
\scrS_2)$
\end{itemize}
\begin{myproof}{Proof}  We prove only \textbf{CH 3}, which is the only
one not following immediately from the definitions.
The inclusion of the right-hand side in the left-hand
side follows from \textbf{CH 1}.  In order to prove the reverse
inclusion, it will suffice, also by \textbf{CH 1} to prove that
$\scrC_d(\scrS_1), \scrS_2\subseteq\scrC_d(\scrS_1,\scrS_2)$.
By \textbf{CH 2}, we have
\startdisp
\scrS_1,\scrS_2\subseteq\scrC_d(\scrS_1,\scrS_2).
\finishdisp
Then by \textbf{CH 1} and \textbf{CH 2}, applied successively, we have
\startdisp
\scrC_d(\scrS_1)\subseteq\scrC_d(\scrC_d(\scrS_1,\scrS_2))=\scrC_d(\scrS_1,
\scrS_2).
\finishdisp
This completes the proof of \textbf{CH 3}.
\end{myproof}

Note that \textbf{CH 3} has obvious generalizations to more than $2$ sets,
$\scrS_i$, which we do not need here.  Further properties of
the convex hull in the specific context of hyperbolic $3$-space,
including an explicit description in coordeinatess
of certain non-Euclidean convex hulls,
will be given in Properties \textbf{HCH 1} through
\textbf{HCH 4} below.

In particular, if we apply these notions to $X=\Real^2$
with the ordinary Euclidean metric $\mathrm{Euc}$, then
the geodesic segment $[p_1,p_2]_{\mathrm{Euc}}$ is
just the line-segment joining $p_1,p_2$.  Further,
provided that not all the $p_i$ are collinear, $\scrC(p_1,\ldots p_r)$
is a closed convex polygon whose vertices are located at a subset of
$\{p_1,\ldots,p_r\}$.

We first use the notion of convex closure to record an
elementary facts concerning the fundamental domains
of groups of translations acting on $\Real^2$, identified with $\Complex$
in the usual way.  Let $\omega_1,\omega_2\in\Complex$
be linearly independent over $\Real$.  Then
$\Int\omega_1+\Int\omega_2$ is a lattice in $\Complex$,
and it is well known that all lattices in $\Complex$
are of this form for suitable $\omega_1, \omega_2$
Let $T$ denote the group of translations by elements
of $\Int\omega_1+\Int\omega_2$ acting on $\Complex$.  Then we have
\starteqn\label{eqn:latticetransF}
\scrC(0,\omega_1,\omega_2,\omega_1+\omega_2)\;\text{is a fundamental
domain for the action of $\Int\omega_1+\Int\omega_2$ on $\Complex$}.
\finisheqn
Now we define the following polygons in $\Complex\cong\Real^2$.  Let
\startdisp
\scrG_U=\scrC_{\mathrm{Euc}}(0,2,1+\vecti,1-\vecti),
\finishdisp
and let
\starteqn\label{eqn:scrGdefn}
\scrG=\scrC_{\mathrm{Euc}}(1,2,1+\vecti).
\finisheqn
The relation between the polygons is that $\scrG_U$ is a square
centered at $1$, while $\scrG$ is an isosceles right triangle
inside $\scrG_U$,
with vertices at the center of $\scrG_U$ and two of the corners
of $\scrG_U$.  Therefore, it follows from \eqref{eqn:rotationdescription}
that we have
\starteqn\label{eqn:scrGfunddominscrGU}
\scrG_U=\bigcup_{i=0,1,2,3}(T_{1+\vecti}R_{\frac{\pi}{2}})^i\scrG,\;\text{with}
\;(T_{1+\vecti}R_{\frac{\pi}{2}})^i\scrG\cap\scrG
\subseteq\partial\scrG,\;\text{for}\; i\not\equiv0\mod 4.
\finisheqn

\begin{lem}\label{lem:scrGfunddom}  Let $\Gamma^{\phi_1}$
be as given in \eqref{eqn:gammaphi1explicit} and
$(\Gamma^{\phi_1})_U$ as given in
\eqref{eqn:gammaphi1Uexplicit}.
\begin{itemize}
\item[(a)]  The set $\scrG_U$ is a fundamental domain for the
induced action of $(\Gamma^{\phi_1})_U$ on $\Complex\cong\Real^2$.
\item[(b)]  $\scrG$ is a fundamental domain for the induced
action of $\langle T_{1+\vecti}R_{\frac{\pi}{2}}\rangle$
on $\scrG_U$.
\item[(c)]  The set $\scrG$ is a fundamental domain for the
induced action of $\Gamma^{\phi_1}$ on $\Complex\cong\Real^2$.
\end{itemize}
\end{lem}
\begin{myproof}{Proof}  The action of $\Gamma^{\phi_1}_U$
on $\Complex$ is by translations by the elements of
$(1+\vecti)\Int[\vecti]$.  Since
\startdisp
(1+\vecti)\Int[\vecti]=\Int(1+\vecti)+\Int(1-\vecti),
\finishdisp
we can apply \eqref{eqn:latticetransF} with $\omega_1=1+\vecti$
and $\omega_2=1-\vecti$ to conclude that
\startdisp
\scrC_{\mathrm{Euc}}(0,1+\vecti,1-\vecti,(1+\vecti)+(1-\vecti))
\finishdisp
is a fundamental domain for the action of $(\Gamma^{\phi_1})_U$
on $\Complex=\Real^2$.  Thus we obtain (a).

Part (b) is a restatement of \eqref{eqn:scrGfunddominscrGU}.

In view of Parts (a) and (b), Part (c)
is an application of Corollary \ref{cor:finiteindexfunddom}.
\end{myproof}
\vspace*{3mm}

\bmth
\noindent\textbf{Form of $\scrF$ in terms of explicit inequalities.}
\ubmth\hspace*{3mm}
Combining Part (c) of Lemma \ref{lem:scrGfunddom},
Proposition \ref{prop:scrF1forgamma}, and \eqref{eqn:scrFdefn},
we deduce that
\startdisp
\scrF(\scrG)=\{z\in\bbH^3\;|\; \phi_{[2,3]}(z)\in\scrC_{\mathrm{Euc}}(1,2,1+\vecti),\;
||x(z)-m||^2+y(z)^2\geq 2,\text{for}\; m\in1+(1+\vecti)\Int[\vecti]\}.
\finishdisp
By \eqref{eqn:phi23new}, the first
condition in the description of $\scrF(\scrG)$ above may be replaced
by
\starteqn\label{eqn:trianglecondition}
x(z)\in\scrC_{\mathrm{Euc}}(1,2,1+\vecti)
\finisheqn
Let $z\in\Complex$ satisfying \eqref{eqn:trianglecondition}.
The element $m=1$ is the element
 of $1+(1+\vecti)\Int[\vecti]$ closest to $x(z)$.
Therefore, for $z$ satisfying
\eqref{eqn:trianglecondition},
the condition
\startdisp
||x(z)-m||^2+y(z)^2\geq 2,\text{for all}\; m\in1+(1+\vecti)\Int[\vecti]
\finishdisp
reduces to $||x(z)-1||^2+y(z)^2\geq 2$.  So we may rewrite
the description of $\scrF(\scrG)$ in the form
\starteqn\label{eqn:scrFscrGfirstform}
\scrF(\scrG)=\{z\in\bbH^3\; x(z)\in\scrC_{\mathrm{Euc}}(1,2,1+\vecti),\;
||x(z)-1||^2+y(z)^2\geq 2\}.
\finisheqn
\vspace*{3mm}

\noindent\bmth
\textbf{Additional facts regarding convex hulls and totally geodesic
hypersurfaces in $\overline{\bbH^3}$.}
\ubmth\hspace*{3mm}
We now extend our ``geodesic hull" treatment of $\scrF$ from the boundary
into the interior of $\bbH^3$.
We first recall certain additional
facts regarding convex hulls and totally geodesic hypersurfaces
in $\bbH^3$.

The description of the geodesics in $\bbH^2$ is well known, but
the corresponding description of the totally geodesic surfaces in
$\bbH^3$ perhaps not as well known, so we recall it here.
Henceforth we abbreviate
``totally geodesic'' by t.g.  Although all t.g. surfaces are
related by isometries, in our model they have two basic types. The
first type is a vertical upper half-plane passing through the
origin with angle $\theta$ measured counterclockwise from the real
axis, which we denote by $\bbH^2(\theta)$.  The second type is an
upper hemisphere centered at the origin with radius $r$, which we
will denote by $\bbS^+_r(0)$.  The t.g. surfaces of $\bbH^3$ are
the $\bbH^2(\theta)$, the $\bbS^+_r(0)$, and their translates by
elements of $\Complex$.  For each of the basic t.g. surfaces, we
produce an isometry $g\in\text{Aut}(\bbH^3)$, necessarily
orientation-reversing, such that $\Fix(g)$ is precisely
the surface in question.  The existence of such a
$g$ shows that the surface is a t.g.
surface.

We define
\startdisp
\overline{\bbH^3}=\bbH^3\cup\Complex\cup\infty
\finishdisp
to be the usual closure of $\bbH^3$ and extend the action of fractional
linear transformations and the notion of the convex hull in
the usual way.  For any subset $\scrS$ of $\bbH^3$, $\overline{\scrS}$
will denote the closure in $\overline{\bbH^3}$.  For $g\in\mathrm{Aut}
(\bbH^3)$, we will likewise use $g$ to denote the extension of $g$
to the closure $\overline{\bbH^3}$.
Henceforth, we will work exclusively in the setting of the closure
$\overline{\bbH^3}$ of $\bbH^3$.  Thus, we will
actually identify the closures of the t.g. surfaces.

The basic orientation-reversing isometry of $\overline{\bbH^3}$
may be denoted $R^*$.
With $x_1+x_2\vecti+y\vectj \in\overline{\bbH^3}$, we have
\startdisp R^*(x_1+x_2\vecti+y\vectj)=x_1-x_2\vecti+y\vectj.
\finishdisp Clearly, we have $\Fix(R^*)=\overline{\bbH^2(0)}$.  To obtain
isometries corresponding to the other vertical planes, let
\startdisp
R_{\theta}=\begin{pmatrix}e^{\vecti\theta/2}&0\\
0&e^{-\vecti\theta/2}\end{pmatrix}. \finishdisp Because
$R_{\theta}\overline{\bbH^2(0)}=\overline{\bbH^2(\theta)}$,
we have \startdisp
\Fix(\conj(R_{\theta})R^*)=\overline{\bbH^2(\theta)}.
\finishdisp

To define the isometry $I$ such that $\Fix(I)$ is the basic
hemisphere $\overline{\bbS_0^+(1)}$, let $\overline{z}$ denote the conjugate
of the quaternion $z$, \textit{i.e.} if $z=x_1+x_2\vecti+y\vectj$
then $\overline{z}=x_1-x_2\vecti-y\vectj$.  For $z\in\overline{\bbH^3}$, set
\startdisp I(z)=1/\overline{z}. \finishdisp We have the equality
$z/I(z)=||z||^2$.  Observe that $\overline{\bbS^+_1(0)}$ is precisely the set
of quaternions in $\overline{\bbH^3}$ of norm one.   Thus,
 $\Fix(I)=\overline{\bbS^+_1(0)}$. For the more general
hemispheres $\overline{\bbS^+_r(0)}$, set \startdisp A(r)=\begin{pmatrix}
\sqrt{r}&0\\
0&\frac{1}{\sqrt{r}}
\end{pmatrix}.
\finishdisp Then, since $A(r)\overline{\bbS^+_1(0)}=
\overline{\bbS^+_r(0)}$, we have
$\Fix(\conj(A(r))I)=\overline{\bbS_r^+(0)}$.

In order to denote the convex hull in $\overline{\bbH^3}$, we use
the notation $\coh$.  Therefore, if $\intd s^2$ is the hyperbolic
metric on $\overline{\bbH^3}$, we have
\startdisp
\coh(p_1,\ldots, p_r)=\mathscr{C}_{\intd s^2}(p_1,\ldots, p_r),
\finishdisp
in terms of our original notational conventions.

If $\scrS$ is a subset of $\overline{\bbH^i}$, $i=2$ or $3$
contained in a unique geodesic, we use $\bbH^1(\scrS)$ to denote
the unique geodesic containing it.  Whenever $\bbH^1(\scrS)$
exists for $\scrS$ a finite set of points $\{p_1,\ldots p_r\}$, we
have $\coh(p_1,\ldots p_r)$ a geodesic segment with endpoints at
two of the $p_i$.

If $\scrS$ is a subset of $\overline{\bbH^3}$ contained in a
unique t.g. surface, then we use $\bbH^2(\scrS)$ to denote that
surface. According to the description of t.g. surfaces given
above, $\bbH^2(\scrS)$ is a vertical half-plane or an
upper-hemisphere centered at a point of $\Complex$. Whenever
$\bbH^2(\scrS)$ exists for $\scrS$ a set of points $\{p_1,\ldots
p_r\}$, we have $\coh(p_1,\ldots, p_r)$ a region in
$\bbH^2(\scrS)$.

Let $p_1,\ldots p_r\in\overline{\bbH^3}$,
for $r>3$ not lying on the same totally
geodesic surface, such that, for each $i$, $1\leq i\leq r$,
\[
p_i\notin\coh(p_1,\ldots,p_{i-1},p_{i+1},\ldots,p_r).
\]
Then the set $\coh(p_1,\ldots, p_r)$ will be called the
\textbf{solid convex polytope with vertices at \boldmath$p_1,\ldots,
p_r$\unboldmath}.  It is clear that for any $p_1,\ldots
p_r\in\overline{\bbH^3}$ not lying in the same totally geodesic
surface, $\coh(p_1,\ldots, p_r)$ is a solid convex polytope with vertices
consisting of some subset of the $r$ points.
\vspace*{3mm}

We now collect some useful properties of the convex hull specific
to the setting of $\overline{\bbH^3}$.
\begin{itemize}
\item[\textbf{HCH 1}]  Let $p_1,\ldots,p_r\in\overline{\bbH^3}-\{
\infty\}$.
Then
\startdisp
\pi_{[2,3]}\scrC_{\bf H}(p_1,\ldots,p_r)=\scrC_{\rm Euc}(\pi_{[2,3]}p_1,\ldots,
\pi_{[2,3]}p_r)
\finishdisp
\item[\textbf{HCH 2}]  Assume that $\bbH^2(p_1,\ldots,p_r)=\bbS_r(p)$,
for some $r\in\Real_+$ and $p\in\Complex$.  Then we have
\startdisp
\scrC_{\bf H}(p_1,\ldots p_r)=\pi_{[2,3]}\inv(\scrC_{\rm Euc}(\pi_{[2,3]}p_1,\ldots,
\pi_{[2,3]}p_r))\bigcap\bbS_r(p).
\finishdisp
\item[\textbf{HCH 3}]  Assume that $p_1,p_2\in\overline{\bbH^2}-\{\infty\}$
satisfy $x(p_1)\neq x(p_2)$ in $\Real$.  Assume without loss of generality
that $x(p_1)<x(p_2)$.
Then we have the following description of the convex hull
of $\{p_1,p_2,\infty\}$, namely
\starteqn\label{eqn:ch2eqn}
\coh(p_1,p_2,\infty)=\big\{p\in\overline{\bbH^2}\;|\;x(p_1)\leq
x(p)\leq x(p_2),\,p\;\text{lies above}\; [p_1,p_2]_{\bf H} \big\}.
\finisheqn
\item[\textbf{HCH 4}]  Assume that $p_1,p_2\in\overline{\bbH^3}-\{\infty\}$
satisfy $x(p_1)\neq x(p_2)$ in $\Complex$.  Then
we have the following description of the convex hull of $\{p_1,p_2,\infty\}$,
namely
\starteqn\label{eqn:ch3eqn}
\coh(p_1,p_2,\infty)=\big\{p\in
\overline{\bbH^2\big(\{p_1,p_2,\infty\}\big)}\;|\; x(p)\in [x(p_1),
x(p_2)]_{\mathrm{Euc}},\; p\;\text{lies above}\; [p_1,p_2]_{\bf H}\big\}.
\finisheqn
\end{itemize}
\begin{myproof}{Proof of the \textbf{HCH} Properties}
Property \textbf{HCH 1} follows easily from the fact that
 $\pi_{[2,3]}=x(\cdot)$
projects geodesic  segments in $\overline{\bbH^3}$ to straight-line
segments in $\Complex$, \textit{i.e.}, to Euclidean geodesics.

For \textbf{HCH 2}, the inclusion
\starteqn\label{eqn:hch2lhs1}
\scrC_{\bf H}(p_1,\ldots,p_r)\subseteq\pi\inv_{[2,3]}(\scrC_{\rm Euc}
(\pi_{[2,3]}p_1,\ldots,\pi_{[2,3]}p_r))
\finisheqn
follows from \textbf{HCH 1}.  Since $\bbS_r(p)$ is by the assumption, a
a convex set containing $p_1,\ldots p_r$,
\starteqn\label{eqn:hch2lhs2}
\scrC_{\bf H}(p_1,\ldots, p_r)\subseteq \bbS_r(p).
\finisheqn
by \eqref{eqn:hch2lhs1} and \eqref{eqn:hch2lhs2} the left-hand
side of \textbf{HCH 2} is contained in the right-hand side.

Now, let
\starteqn\label{eqn:hch2rhs0}
q\in\bbS_r(p)\cap\pi\inv_{[2,3]}(\scrC_{\rm Euc}
(\pi_{[2,3]}p_1,\ldots,\pi_{[2,3]}p_r)),
\finisheqn
Then
\startdisp
\pi_{[2,3]}(q)\in\scrC_{\rm Euc}(\pi_{[2,3]}p_1,\ldots,\pi_{[2,3]}p_r).
\finishdisp
Applying the equality of \textbf{HCH 1}, we have
\startdisp
\pi_{[2,3]}q\in\pi_{[2,3]}\coh(p_1,\ldots,p_r).
\finishdisp
Thus, there exists $q'\in\overline{\bbH^3}$ satisfying the conditions
\starteqn\label{eqn:hch2rhs1}
\pi_{[2,3]}q=\pi_{[2,3]}q'\;\text{and}\; q'\in\coh(p_1,\ldots, p_r)
\finisheqn
Now consider the restriction of $\pi_{[2,3]}|_{\bbS_r(p)}$ of $\pi_{[2,3]}$
to the upper-half sphere.  It is clear that the restriction is
injective.  By the assumption \eqref{eqn:hch2rhs0}, we have that
$q$ is in the domain of the restriction.  Since by the hypothesis
of \textbf{HCH 2},
\startdisp
\coh(p_1,\ldots, p_r)\subseteq\bbS_r(p),
\finishdisp
the second statement of \eqref{eqn:hch2rhs1} implies that
$q'$ is in the domain of the restriction.  Therefore, the
first statement of \eqref{eqn:hch2rhs1} can be rewritten
\startdisp
\pi_{[2,3]}|_{\bbS_r(p)}q=\pi_{[2,3]}|_{\bbS_r(p)}q'.
\finishdisp
By the injectivity of the restriction $q=q'$.  So
by the second statement of \eqref{eqn:hch2rhs1}, we have
$q\in\coh(p_1,\ldots p_r)$.  Since $q$ was an arbitrary
point satisfying \eqref{eqn:hch2rhs0}, this completes
the proof of the inclusion of the right-hand side in the left-hand side.

For Property \textbf{HCH 3}, since we are working
in the context of $\overline{\bbH^2}$ the assumption that $x(p_1)<x(p_2)$
implies that
\starteqn\label{eqn:arcdescription}\text{
$[p_1,p_2]_{\bf H}$ is the arc connecting
$p_1,\, p_2$ on a circle centered at a point of $\Real$.}
\finisheqn
The set on the right-hand
side of \eqref{eqn:ch2eqn} is merely the hyperbolic triangle
whose sides are the arc $[p_1,p_2]_{\bf H}$ and the semi-infinite
vertical lines beginning at $p_1$, $p_2$.  It is well-known
that this hyperbolic triangle is convex in the hyperbolic metric.
Therefore, $\coh(p_1,p_2,\infty)$ is contained in the right-hand
side of \eqref{eqn:ch2eqn}.  For the reverse inclusion,
we have, by the definition of the convex hull,
\starteqn\label{eqn:lowerarcinconvexhull}
[p_1,p_2]_{\bf H}\subseteq\coh(p_1,p_2,\infty).
\finisheqn
Let $p$ be an element of the right-hand side of \eqref{eqn:ch2eqn}.
By \eqref{eqn:arcdescription} and the assumption that $x(p_1)\leq x(p)\leq
x(p_2)$,
there is a unique point in $\overline{\bbH^3}$,
\startdisp
p'\in [p,x(p)]_{\bf H}\cap[p_1,p_2]_{\bf H}.
\finishdisp
By definition, $p'$ is the point
directly below $p$ and on the arc $[p_1,p_2]_{\bf H}$.
Since $p'\in[p_1,p_2]_{\bf H}$, \eqref{eqn:lowerarcinconvexhull}
implies that $p'\in\coh(p_1,p_2,\infty)$.  Taking
into account that $\coh(p_1,p_2,\infty)$ is by definition
a convex set containing $\infty$, we deduce that
\starteqn\label{eqn:convexityverticalhalflines}\text{
$p',\infty\in\coh(p_1,p_2,\infty)$, and thus $[p',\infty]_{\bf H}\subseteq
\coh(p_1,p_2,\infty)$.}
\finisheqn
Since $[p',\infty]_{\bf H}$ is the semi-infinite vertical line
beginning at $p'$ and $p'\in [p,x(p)]_{\bf H}$, the vertical
segment connecting $p$ to the point $x(p)$ directly below $p$
on $\Real$, we have $p\in[p',\infty]_{\bf H}$.  So, by
\eqref{eqn:convexityverticalhalflines}, we have
\startdisp
p\in\coh(p_1,p_2,\infty).
\finishdisp
Since $p$ was an arbitrary point of the right-hand side
of \eqref{eqn:ch2eqn}, this completes the proof of the inclusion
of the right-hand side of \eqref{eqn:ch2eqn} into the left.

For \textbf{HCH 4}, note that the assumption that $x(p_1)\neq x(p_2)$
is equivalent to saying that $\infty$ does not lie on the geodesic
$\bbH^1(p_1,p_2)$.  Therefore, the assumption implies that
$\bbH^2(p_1,p_2,\infty)$ is well-defined.
The t.g. surfaces in $\overline{\bbH^3}$ containing $\infty$
are the vertical planes, so that $\bbH^2(p_1,p_2,\infty)$
is a vertical plane.
It is possible to choose a
\starteqn\label{eqn:translationrotation}
TR\in T_{\Complex}R_{\Real/2\pi\Int}\subseteq \mathrm{Aut}^+(\overline{\bbH^3})\;
\text{such that}\; TR\bbH^2(p_1,p_2,\infty)=\overline{\bbH^2(0)}
\finisheqn
Clearly $TR$ fixes $\infty$ and maps
convex hulls to convex hulls.  Further
$\overline{\bbH^2(0)}\cong\overline{\bbH^2}$.
Therefore, we have reduced \textbf{HCH 4}
to \textbf{HCH 3}, provided that we verify that $TR$
as in \eqref{eqn:translationrotation} maps the
set on the right-hand side of \eqref{eqn:ch3eqn} to
the set on the right-hand side of \eqref{eqn:ch2eqn},
with $TRp_i$ in place of $p_i$ for $i=1,2$.
But note that the right hand side of \eqref{eqn:ch2eqn}
can be rewritten as
\starteqn\label{eqn:ch2eqnsillyrewrite}
\big\{p\in\overline{\bbH^2}\;|\;x(p)\in
[x(p_1),x(p_2)]_{\rm Euc},\,p\;\text{lies above}\; [p_1,p_2]_{\bf H} \big\},
\finisheqn
since in the context of
\eqref{eqn:ch2eqn}, $[x(p_1),x(p_2)]_{\rm Euc}$ is just a fancy way of writing
the real interval $[x(p_1),x(p_2)]$.  Since
\startdisp
TR\circ y(\cdot)=y(\cdot),
\finishdisp
$TR$ maps the set lying above an arc to the set lying above
the image of the arc.  Therefore, \eqref{eqn:ch2eqnsillyrewrite} exactly
matches the right-hand side of \eqref{eqn:ch3eqn}.
This completes the proof of \textbf{HCH 4}.
\end{myproof}

\noindent\bmth\textbf{Description of $\scrF(\scrG)$ as a convex
polytope.}\hspace*{3mm}
\ubmth
We have already given one form of $\scrF(\scrG)$ in
\eqref{eqn:scrFscrGfirstform}.  In order to express the set
on the right-hand side of \eqref{eqn:scrFscrGfirstform}
as a convex polytope, we record the following general
observations.  Let $\Omega\subseteq\Complex$, $\Omega\neq\emptyset$
and let $f$ be any function
\starteqn\label{eqn:fforgraph}
f: \Omega\rightarrow\Real_{\geq0}\cup\infty.
\finisheqn
Then the graph of $f$ is defined as subset of $\overline{\bbH^3}$, namely,
\startdisp
\Graph(f)=
\{w+f(w)\vectj\;|\;w\in\Omega\}.
\finishdisp
We also define the \bmth\textbf{subset $\scrS_f$ of
of $\overline{\bbH^3}$ lying on or above
the graph of $f$}\ubmth, by setting
\starteqn\label{eqn:abovegraph}
\scrS_f:=
\{z\in\overline{\bbH^3}\;|\; x(z)\in\Omega,\, y(z)\geq f(x(z))\}.
\finisheqn
By convention, we always include $\infty$ in $\scrS_f$.

The following proposition, Lemma \ref{lem:graphconvexhull},
may be thought of as relating a convexity property of $f$
to the convexity of the set $\scrS_f$.  Assume that $\Omega\subset\Complex$,
\textit{i.e.} the domain of $f$, is a convex subset of $\Complex$.  We
\textit{define} $f$ as in \eqref{eqn:fforgraph} to
be a \textbf{convex function} if and only if it satisfies the property
\starteqn\label{eqn:fconvex}
x_1, x_2 \in\Omega,\; p\in[x_1+f(x_1)\vectj,x_2+f(x_2)\vectj]_{\mathbf{H}}
\;\text{implies}\; y(p)\geq f(p),
\finisheqn
in other words if and only if the geodesic segment between
any two points on $\Graph(f)$ lies above the graph of $f$ (in $\scrS_f$).

\begin{lem}\label{lem:graphconvexhull}  Let $\Omega\subset\Complex$,
$\Omega\neq\emptyset$ be convex.  Let $f$ a function as
in \eqref{eqn:fforgraph},
$\scrS_f$ the subset of $\bbH^3$ lying on or above $\Graph(f)$
as in \eqref{eqn:abovegraph}.  Then $f$ is convex in the sense
of \eqref{eqn:fconvex} if and only if $\scrS_f$ is a convex
subset of $\bbH^3$, and in that case we have
\starteqn\label{eqn:scrSfdescription}
\scrS_f=\coh(\Graph(f),\infty).
\finisheqn
\end{lem}
\begin{myproof}{Proof}  First, note that if $\scrS_f$ is convex,
then it immediately follows that the function $f$ is convex.
The reason is that $x_i+f(x_i)\vectj$ for $i=1,2$ are points
on $\Graph(f)$, therefore in $\scrS_f$, and the convexity
of $\scrS_f$ implies that the entire geodesic segment
\startdisp
[x_1+f(x_1)\vectj,x_2+f(x_2)\vectj]_{\mathbf{H}}\subseteq\scrS_f.
\finishdisp
For $p$ a point in the geodesic segment, the condition
$y(p)\geq f(p)$ then follows from the definition of $\scrS_f$.

For the converse, suppose that $f$ is a convex function.
We first claim that for arbitrary $p_1, p_2$ in $\scrS_f$,
we have $\{p_1,p_2,\infty\}\subseteq\coh(x_1+f(x_1)\vectj,x_2+f(x_2)\vectj,\infty)$.
From this claim and \textbf{CH 1} it follows that
\startdisp
\coh(p_1,p_2,\infty)\subseteq\coh(x_1+f(x_1)\vectj,x_2+f(x_2)\vectj,\infty).
\finishdisp
The assumption of convexity  $f$ means that we have the containment
\starteqn\label{eqn:graphinscrSf}
[x_1+f(x_1)\vectj,x_2+f(x_2)\vectj]\subseteq\scrS_f.
\finisheqn
For $\infty$, the containment is obvious.  Since $p_i\in\scrS)f$,
$y(p_i)\geq f(x_i)$.  Thus,
\startdisp
p_i\in[x_i+y_i,\infty)_{\bbH}
\finishdisp
By the description of $\coh(x_1+f(x_1),x_2+f(x_2),\infty)$
given in \textbf{CH 3} and the definition of $\scrS_f$,
\eqref{eqn:graphinscrSf} yields
\starteqn\label{eqn:graphinscrSf2}
\coh(x_1+f(x_1)\vectj,x_2+f(x_2)\vectj,\infty)\subseteq\scrS_f.
\finisheqn
Combining \eqref{eqn:graphinscrSf2} with the above claim, we have
\startdisp
\coh(p_1,p_2,\infty)\subseteq\scrS_f,
\finishdisp
so in particular $\scrS_f$ is convex.  From \eqref{eqn:graphinscrSf2}
and the description of $\coh(x_1+f(x_1),x_2+f(x_2),\infty)$
given in \textbf{CH 3}, we obtain the description of $\scrS_f$
in the ``convex case".
\end{myproof}

In light of the fact that $\bbS_r(p)$ is the graph of the function
\startdisp
f(x)=\sqrt{r^2-||x-p||^2},
\finishdisp
we can rewrite property \textbf{HCH 2} in the following form
\begin{itemize}
\item[\textbf{HCH 2$'$}]  Assume that $\bbH^2(p_1,\ldots,p_r)=\bbS_r(p)$,
for some $r\in\Real_+$ and $p\in\Complex$.  Then
$\scrC_{\bf H}(p_1,\ldots p_r)$ is the portion of the graph
of $f(x)=\sqrt{r^2-||x-p||}$ in $\bbH^3$ lying above
$\scrC_{\rm Euc}(\pi_{[2,3]}p_1,\ldots,\pi_{[2,3]}p_r)$.
\end{itemize}

\begin{prop}\label{prop:gammafunddom}
The solid convex polytope with four vertices given by
\starteqn\label{eqn:scrFscrGsecondform}
\scrF(\scrG)=\scrC_{\bbH}(1+\sqrt{2}\vectj,2+\vectj,1+\vecti+\vectj,\infty)
\finisheqn
is a good Grenier fundamental domain for the action of
$\Gamma=\conj\inv(\SO{3}{\Int[\vecti]})$ on $\bbH^3$.
\end{prop}
\begin{myproof}{Proof}  In \eqref{eqn:scrFscrGfirstform}, we have described the fundamental
domain as the set $\scrS_f$ lying above the graph of the function
\startdisp
f(x)=\sqrt{2-||x-1||^2},
\finishdisp
defined on the domain
\startdisp
\Omega=\scrC_{\rm Euc}(1,2,1+\vecti).
\finishdisp
Since the graph of $f$ is a t.g. surface (sphere of radius $\sqrt{2}$
centered at the point $x=1$) in $\overline{\bbH^3}$, it
is clear that $f$ is convex.  Therefore, $\scrS_f$
is convex and has the description given in Lemma \ref{lem:graphconvexhull}.
Therefore,
\starteqn\label{eqn:scrFscrGsecondformprelim}
\begin{aligned}
\scrF(\scrG)&=&&\coh(\text{Graph}f,\infty)\\
&=&&\coh(\scrC_{\bf H}(1+\sqrt{2}\vectj,2+\vectj,1+\vecti+\vectj),\infty)\\
&=&&\scrC_{\bf H}(1+\sqrt{2}\vectj,2+\vectj,1+\vecti+\vectj,\infty),
\end{aligned}
\finisheqn
where we have applied $\textbf{HCH 2$'$}$ and $\textbf{CH 3}$
in the second and third lines of
\eqref{eqn:scrFscrGsecondformprelim}.  This completes the proof of
\eqref{eqn:scrFscrGsecondform}.
\end{myproof}

\section{The split real form $\conj(\SL{2}{\Real})$ and its
intersection $\SOM{2}{1}_{\Int}$ with $\SO{3}{\Int[\vecti]}$}
\label{sec:realform}
We will now use the results of \S\ref{sec:lattice} and \S\ref{sec:explicitfd}
to deduce a realization of $\Gamma_{\Int}=\SOM{2}{1}_{\Int}$ as a group
of fractional linear transformations, as well as a description
of a fundamental domain for $\Gamma_{\Int}$ acting on $\bbH^2$
that is in some sense (to be explained precisely below) compatible with the
fundamental domain of $\Gamma$ acting on $\bbH^3$.

\subsection{$\SOM{2}{1}_{\Int}$ as a group of fractional
linear transformations}
\label{subsec:realformflts}
We maintain to the notational conventions established in
\S\ref{subsec:flts}.  In particular, $G=\SO{3}{\Complex}$
and $\Gamma=\SO{3}{\Int[\vecti]}$.  It is crucial, for the moment,
that we observe the distinction between $G,\Gamma$ and their
isomorphic images under $\conj\inv$.
\begin{defn} \label{defn:gammaz} Set
\starteqn\label{eqn:gammazdefn}
\Gamma_{\Int}=\conj(\SL{2}{\Real}\cap\conj\inv(\Gamma)).
\finisheqn
\end{defn}
\begin{rem}
Note that the elements of $\Gamma_{\Int}$ do not have real entries!
The na\"{i}ve approach to the definition of $\Gamma_{\Int}$
would be to take
the elements of $\Gamma$ with real entries, as in the case of
$\SL{2}{\Int[\vecti]}$ and $\SL{2}{\Int}$.  However, this clearly
cannot be the right definition
because the resulting discrete group would be contained
in $\SOR{3}$, hence compact, and hence finite.
The justification for Definition \ref{defn:gammaz}
is contained in Proposition \ref{prop:gammaziso}, below.
\end{rem}

Recall the orthonormal basis $\beta$ for $\mathrm{Lie}(\SL{2}{\Complex})$
defined at $\eqref{eqn:beta}$.  Define a new basis
$\eta$ by specifying the change-of-basis matrix
\starteqn\label{eqn:etadefn}
\alpha^{\beta\mapsto\eta}=\mathrm{diag}(1,-\vecti,1).
\finisheqn

Let $V_{\Real}$ be a \textit{real} vector space of dimension 3.
Let $\SOM{2}{1}$ denote the group of unimodular linear automorphisms
of $V_{\Real}$
preserving a form $B_{\Real}$ on $V_{\Real}$ of bilinear signature $(2,1)$.
For definiteness,
we will take
\startdisp
V_{\Real}=\Real\text{-span}(\eta)\subseteq\mathrm{Lie}(\SL{2}{\Complex}),\quad
B_{\Real}=B|_{V_{\Real}},
\finishdisp
where $\beta'$ is the basis of $\mathrm{Lie}(\SL{2}{\Complex})$
defined at
\eqref{eqn:betaprime}, and $B$ is as usual the Killing form
on $\mathrm{Lie}(\SL{2}{\Complex})$.  From the fact that $\beta$
is an orthonormal set under $B$ and from \eqref{eqn:etadefn},
it is immediately verified that $B|_{\Real}$ has signature $(2,1)$.
Note also that
\startdisp
V:=V_{\Real}\otimes\Complex=\mathrm{Lie}(\SL{2}{\Complex}).
\finishdisp

By considering $\SOM{2}{1}$ as a subset of $\GL{3}{\Real}$
we obtain the \bmth\textbf{standard representation of $\SOM{2}{1}$}\ubmth.
We define $\SOM{2}{1}_{\Int}$ to be the matrices with integer
coefficients in the standard representation of $\SOM{2}{1}$.

Recall from \eqref{eqn:conjwithresptobasis} the definition of
the morphism
\startdisp
\conj_{\eta}:=\conj_{V,\eta}:
\SL{2}{\Complex}\rightarrow\SL{3}{\Real}.
\finishdisp

\begin{prop} \label{prop:gammaziso} Let $\Gamma_{\Int}$ as defined
in \eqref{defn:gammaz}.  Then the restriction of $\conj_{\eta}$
to $V_{\Real}$ provides an isomorphism
\starteqn\label{eqn:liegroupsiso}
\conj_{\eta}:\SL{2}{\Real}/\{\pm I\}\rightarrow\SOM{2}{1}^0
\finisheqn
of Lie groups.  The isomorphism of \eqref{eqn:liegroupsiso}
further restricts to an isomorphism of discrete subgroups
\starteqn\label{eqn:discretegroupsiso}
\conj_{\eta}: \conj\inv(\Gamma_{\Int})\rightarrow\SOM{2}{1}_{\Int}.
\finisheqn
As a result, $\conj_{\eta}\conj\inv$ exhibits an isomorphism
\starteqn\label{eqn:gammaintjustification}
\Gamma_{\Int}\cong\SOM{2}{1}_{\Int}.
\finisheqn
\end{prop}
\begin{myproof}{Proof}
The image of $\SL{2}{\Real}$ under $\conj_{\eta}$ preserves
$B_{\Real}$ because of property \textbf{B3} from \S\ref{subsec:flts}.
It is clear that the $\ker(\conj_{\eta})$ is just the center
of $\SL{2}{\Real}$, \textit{i.e.}, $\{\pm I\}$.  A comparison
of dimensions completes the proof that
$\conj_{\eta}$ in \eqref{eqn:liegroupsiso}
is an isomorphism.

For the remaining statements, first observe from Definition
\ref{defn:gammaz} that
\startdisp
\conj\inv(\Gamma)\cap
\SL{2}{\Real}=\{\alpha\in\SL{2}{\Real}\;|\;
\conj(\alpha)\in\Mats_3(\Int[\vecti])\}.
\finishdisp
By the definition of $\Gamma_{\Int}$ in \eqref{eqn:gammazdefn}, we
therefore have
\starteqn\label{eqn:conjetadiscreteimage}
\begin{aligned}
\conj_{\eta}|_{\SL{2}{\Real}}(\conj\inv(\Gamma_{\Int}))
&=&&\conj_{\eta}(\conj\inv(\Gamma)\cap
\SL{2}{\Real})\\
&=&&\conj_{\eta}\left(\{\alpha\in\SL{2}{\Real}\;|\;
\conj(\alpha)\in\Mats_3(\Int[\vecti])\}\right)\\
&=&&\{\conj_{\eta}(\alpha)\;|\; (\conj\circ\conj_{\eta}\inv)
(\conj_{\eta}(\alpha))\in\Mats_3(\Int[\vecti])\}\\
&=&&\left\{\conj_{\eta}(\alpha)\;|\; \conj_{\SL{3}{\Complex}}\left(\alpha^{\beta
\mapsto\eta}\right)
(\conj_{\eta}(\alpha))\in\Mats_3(\Int[\vecti])\right\}\\
&=&&\left\{\conj_{\eta}(\alpha)\;|\;
\conj_{\SL{3}{\Complex}}\left(\mathrm{diag}(1,-\vecti,1)\right)
(\conj_{\eta}(\alpha))\in\Mats_3(\Int[\vecti])\right\}\\
&=&&\{\conj_{\eta}(\alpha)\;|\;
(\conj_{\eta}(\alpha))\in\Mats_3(\Int[\vecti])\},
\end{aligned}
\finisheqn
where in the antepenultimate and penultimate lines we have used relations
\eqref{eqn:changeofbasisconj} and \eqref{eqn:etadefn}, respectively.
Comparing the first and last expressions in the string
of equalities in \eqref{eqn:conjetadiscreteimage}, we
we see that the isomorphism $\conj_{\eta}|_{\SL{2}{\Real}}$
restricts to an isomorphism of $\conj\inv(\Gamma_{\Int})$
onto $\SOM{2}{1}_{\Int}$.  This completes the proof of the
second statement of the Proposition.

The isomorphism in \eqref{eqn:gammaintjustification} is an
immediate consequence of the second statement of the Proposition.
\end{myproof}

The next Proposition, \ref{prop:gammaintexplicitdescrip}, is the analogue
of Proposition \ref{prop:inversematrixexplicitdescription}
for the real form of the complex group.
Proposition \ref{prop:gammaintexplicitdescrip} below
is, in contrast, almost a triviality to prove at this point,
since it can be deduced rather readily from Proposition
\ref{prop:inversematrixexplicitdescription} and a few elementary
preliminary results, which we now state.

We begin by establishing a few conventions on the notation
\text{arg}.  First, \text{arg} is for us an \textit{even} function on $\Complex^\times$
taking values in the interval $[0,\pi)$.  It is defined by
\startdisp
\arg(z)=\theta,\;\text{where}\;\theta\in[0,\pi)\;\text{is such that}\;
z=re^{\vecti\theta},\;\text{with}\;r\in\Real.
\finishdisp
Note that in this definition, $r$ is not assumed to be positive.
Second, the function $\arg$ is extended to $n$-tuples of nonzero
complex numbers, but only those of the form
\startdisp
(r_1e^{\vecti\theta},\ldots,r_ne^{\vecti\theta}),\; r_i\in\Real-\{0\},
\;\text{for}\; i=1,\ldots, n,
\finishdisp
\textit{i.e.}, only to those for which the argument of all the entries
is the same.  Accordingly, any time $\arg$ of a $n$-tuple
appears in an equality, the equality must be read as containing
the implicit assertion that $\arg$ of the $n$-tuple exists.
\begin{lem}\label{lem:argumentcondition}  Let $x\in\Complex^{\times}$ We have
\startdisp
\arg(x)=\pi/4\;\text{or}\;3\pi/4,\;x\in\Int[\vecti]\;\text{implies that}\;
\ord_{1+\vecti}(x)> 0.
\finishdisp
\end{lem}
\begin{myproof}{Proof}  For definiteness, assume that
$\arg(x)=\pi/4$, as the proof in the case $3\pi/4$ is identical.
Let $L_{\pi/4}$ be the line passing through the origin of $\Complex$
at angle $pi/4$ from the real axis, and consider the additive
subgroup
\startdisp
L_{\pi/4}\cap \Int[\vecti],
\finishdisp
containing $x$ by assumption.
This additive subgroup is isomorphic to an additive subgroup
of $\Real$, and so is singly generated,
namely by an element of minimum modulus.  It is easy to
see however that an element of minimum modulus in the group
is $1+\vecti$.  Thus
\startdisp
L_{\pi/4}\cap\Int[\vecti]=\Int\,\cdot\,1+\vecti,
\finishdisp
so that $1+\vecti|x$.  Therefore, $\ord_{1+\vecti}> 0$.
\end{myproof}

For Proposition \ref{prop:gammaintexplicitdescrip},
it is necessary to recall the group $\Xi$-subgroups of
defined in \eqref{eqn:xi12defn} and \eqref{eqn:residuematrices}.  For
each of the three $\Xi$-subgroups, we define
\starteqn\label{eqn:xi12intdefn}
(\Xi)_{\Int}=\Xi\cap\SL{2}{\Real}.
\finisheqn
The following result both justifies this notation and clarifies
the meaning of Proposition \ref{prop:gammaintexplicitdescrip}, below.
\begin{lem}
Each $(\Xi)_{\Int}$-group can be given the following description.
\starteqn\label{eqn:xiintsetdescription}
\begin{gathered}
\text{For fixed}\;
\begin{pmatrix}\overline{p}&\overline{q}\end{pmatrix},
\begin{pmatrix}\overline{r}&\overline{s}\end{pmatrix}\in
\left\{
\begin{array}{l}
\begin{pmatrix}1&1\end{pmatrix},\vspace*{0.15cm}\\
\begin{pmatrix}1&0\end{pmatrix},\vspace*{0.15cm}\\
\begin{pmatrix}0&1\end{pmatrix}
\end{array}
\right\}\subset(\SL{2}{\Int[\vecti]/(2)})^2,\\
\Xi=\red_{2}\inv\left(
\left\{\begin{pmatrix}\overline{p}&\overline{q}\\
\overline{r}&\overline{s}\end{pmatrix},\,
\begin{pmatrix}\overline{r}&\overline{s}\\
\overline{p}&\overline{q}
\end{pmatrix}\right\}\right).
\end{gathered}
\finisheqn
In order to obtain $\Xi_{12}$, we may take, in \eqref{eqn:xiintsetdescription},
\startdisp
\begin{pmatrix}\overline{p}&\overline{q}\end{pmatrix}=
\begin{pmatrix}1&0\end{pmatrix}\;\text{and}
\begin{pmatrix}\overline{r}&\overline{s}\end{pmatrix}
=\begin{pmatrix}0&1\end{pmatrix}
\finishdisp
Further, we may take
\startdisp
\begin{pmatrix}\overline{p}&\overline{q}\end{pmatrix}=
\begin{pmatrix}1&1\end{pmatrix},\;\text{in order
to obtain both $\Xi_{1}$ and $\Xi_{2}$,}
\finishdisp
and
\startdisp\begin{gathered}
\begin{pmatrix}\overline{r}&\overline{s}\end{pmatrix}
=\begin{pmatrix}0&1\end{pmatrix},\,\text{in order to obtain $\Xi_1$},\\
\begin{pmatrix}\overline{r}&\overline{s}\end{pmatrix}
=\begin{pmatrix}1&0\end{pmatrix},\,\text{in order to obtain $\Xi_2$}.
\end{gathered}
\finishdisp
\end{lem}
\begin{myproof}{Proof}  Since the factorization of $2$ in $\Int[\vecti]$
is
\startdisp
2=\vecti^3(1+\vecti)^2
\finishdisp
we have
\startdisp
\res_{1+\vecti}=\res_2\;\text{on}\; \Int\subseteq\Int[\vecti].
\finishdisp
Therefore, \eqref{eqn:xiintsetdescription} is just a restatement
of \eqref{eqn:xisetdescription}.  The explicit specifications
of the matrices $\begin{pmatrix}\overline{p}&\overline{q}\end{pmatrix}$
and $\begin{pmatrix}\overline{p}&\overline{q}\end{pmatrix}$
are immediately deduced from \eqref{eqn:xi12defn} and
\eqref{eqn:residuematrices}.
\end{myproof}

\begin{prop}\label{prop:gammaintexplicitdescrip}  With
$\Gamma_{\Int}$ defined as in \eqref{eqn:gammazdefn},
we have
\starteqn\label{eqn:gammaintdescription}
\conj\inv(\Gamma_{\Int})=(\Xi_{12})_{\Int}\bigcup
\frac{1}{\sqrt{2}}({\Xi_{2}})_{\Int}\begin{pmatrix}1&-1\\0&2
\end{pmatrix}.
\finisheqn
\end{prop}
\begin{myproof}{Proof}  We briefly recall certain notations from Section
\ref{subsec:inverseimagedescription}.  In particular,
we represent $\alpha\in\conj\inv(\Gamma)$ in the form
\startdisp
\alpha=\begin{pmatrix}a&b\\c&d\end{pmatrix}.
\finishdisp
We recall the notations $(x,y,z,w)$ for an arbitrary
permutation of $(a,b,c,d)$ with corresponding permutation
$(x',y',z',w')$ of $(a',b',c',d')$.
Further, we set
\startdisp
(i,\delta):=\big(i,\delta\big)(\alpha),
\finishdisp
so that $(i,\delta)\in\{0,2\}\times\{0,1\}$, and
\starteqn\label{eqn:alphaprimerecalled}
\alpha'=\begin{pmatrix}a'&b'\\c'&d'\end{pmatrix}\in\mathrm{Mat}_2
(\Int[\vecti]),\;\ord_{1+\vecti}(x')=0\;\text{where}\;\alpha'
:=(1+\vecti)^{i/2}\omega_8^{\delta}\alpha.
\finisheqn
An elementary computation using \eqref{eqn:omegaeightdefn} gives
\starteqn\label{eqn:argumentcomp}
\arg\left((1+\vecti)^{i/2}\omega_8^{\delta}\right)=
(i/2+\delta)\frac{\pi}{4}.
\finisheqn
By \eqref{eqn:alphaprimerecalled} and \eqref{eqn:argumentcomp},
we have the equivalence of conditions
\starteqn\label{eqn:alphaprimeargcond}
\alpha\in\SL{2}{\Real}\;\text{if and only if}\;\arg(x')=
(i/2+\delta)\frac{\pi}{4}.
\finisheqn
Suppose that $i/2+\delta$ is odd.  Then,
by \eqref{eqn:alphaprimeargcond}, Lemma \ref{lem:argumentcondition}
applies to $x'$ and we obtain
\startdisp
\ord_{1+\vecti}(x')> 0.
\finishdisp
This contradicts \eqref{eqn:alphaprimerecalled}.  Therefore,
we may assume without loss of generality that $i/2+\delta$
is even.  Since $(i,\delta)\in\{0,2\}\times\{0,1\}$,
we may assume, equivalently, that $(i,\delta)$ is either
$(0,0)$ or $(2,1)$.  We now examine these two cases.

\vspace*{3mm}
\noindent\textit{Case 1: $(i,\delta)=(0,0)$.}  By the
first line of \eqref{eqn:twopartsdescribed}, we have
\startdisp
\big(i,\delta\big)\inv((0,0))=\Xi_{12}\alpha^{\vecti^{0}}(\vecti^{0},0)=
\Xi_{12}.
\finishdisp
Thus,
\starteqn\label{eqn:gammaint00description}
\big(i,\delta\big)\inv((0,0))\cap\SL{2}{\Real}=\Xi_{1,2}\cap\Real=
(\Xi_{1,2})_{\Int},
\finisheqn
using the definition \eqref{eqn:xi12intdefn} of $(\Xi_{12})_{\Int}$.

\vspace*{3mm}
\noindent\textit{Case 1: $(i,\delta)=(2,1)$.}  By the second
line of \eqref{eqn:twopartsdescribed}, we have
\startdisp
\big(i,\delta\big)\inv((2,1))=
\bigcup_{\epsilon=0,1}\hspace*{-0.45cm}\cdot\hspace*{.55cm}
\frac{1}{\omega_8(1+\vecti)}
\Xi_{2}\alpha^{-2}\hspace*{-0.7mm}
(-1,\vecti^{\epsilon}).
\finishdisp
By \eqref{eqn:alphaprimerecalled}, therefore,
\starteqn\label{eqn:gammaint21descriptioninter}
\{\alpha'\;|\;\alpha\in\conj\inv(\Gamma),\,
\big(i,\delta\big)(\alpha)=(2,1)\}=
\bigcup_{\epsilon=0,1}\hspace*{-0.45cm}\cdot\hspace*{.55cm}
\Xi_{2}\alpha^{-2}\hspace*{-0.7mm}
(-1,\vecti^{\epsilon}).
\finisheqn
Write $\alpha$ as in \eqref{eqn:gammaint21descriptioninter}, so that
\startdisp
\alpha'=\begin{pmatrix}p&q\\r&s\end{pmatrix}
\alpha^{-2}\hspace*{-0.7mm}
(-1,\vecti^{\epsilon}),\;\text{where}
\;\begin{pmatrix}p&q\\r&s\end{pmatrix}
\in\Xi_{12}.
\finishdisp
It is easily computed that such an $\alpha'$ can be expressed as
\starteqn\label{eqn:alphaprimemultform}
\alpha'=\begin{pmatrix}a'&b'\\c'&d'\end{pmatrix},\;\text{with}\;
(a',b',c',d')=(-p,p\vecti^\epsilon+2q,-r,r\vecti^\epsilon+2s).
\finisheqn
By \eqref{eqn:gammaint21descriptioninter},
\eqref{eqn:alphaprimeargcond}, and \eqref{eqn:alphaprimemultform} we have
\begin{multline}\label{eqn:gammaint21descriptioninter2}
\{\alpha'\;|\;\alpha\in\conj\inv(\Gamma)\cap\SL{2}{\Real},\,
\big(i,\delta\big)(\alpha)=(2,1)\}=\\
\bigcup_{\epsilon=0,1}\hspace*{-0.45cm}\cdot\hspace*{.55cm}
\{\alpha'\;\text{of the form \eqref{eqn:alphaprimemultform}}\;
\text{with}\;\arg(-p,p\vecti^\epsilon+2q,-r,r\vecti^\epsilon+2s)=\pi/2\}.
\end{multline}
We now claim that
\starteqn\label{eqn:argumentclaim}\text{
For $\begin{pmatrix}p&q\\r&s\end{pmatrix}\in\Xi_{12}$,
$\arg(-p,p\vecti^\epsilon+2q,-r,r\vecti^\epsilon+2s)=\pi/2$
if and only if $\epsilon=0$, $\arg(p,q,r,s)=\pi/2$.}
\finisheqn
In order to prove the claim, first suppose that $\epsilon=1$.
We will derive a contradiction.  Since $\arg(-p,-r)=\arg(p,r)=\pi/2$,
we have
\startdisp
\text{$p,r$ are of the form $p=\Impt(p)\vecti,\, r=\Impt(r)\vecti$, for
$\Impt(p), \Impt(r)\in\Int$}.
\finishdisp
Thus we have $\arg(-\Impt(p)+2q,-\Impt(r)+2s)=\pi/2$, from which
we deduce that
\startdisp
\Impt(p)=2\Rept(q),\;\Impt(r)=2\Rept(s).
\finishdisp
However, this implies that
\startdisp
\begin{pmatrix}\overline{p}\\\overline{r}\end{pmatrix}=
\begin{pmatrix}\overline{0}\\\overline{0}\end{pmatrix},
\finishdisp
which, by \eqref{eqn:xi12defn}, contradicts the assumption that
\startdisp
\begin{pmatrix}p&q\\r&s\end{pmatrix}\in\Xi_{12}.
\finishdisp
Thus we have $\epsilon=0$.
Assuming that $epsilon=0$, we clearly have
\startdisp
\arg(-p,p\vecti^\epsilon+2q,-r,r\vecti^\epsilon+2s)=
\arg(p,q,r,s).
\finishdisp
Therefore,
we have the equivalence,
\startdisp
\arg(-p,p\vecti^\epsilon+2q,-r,r\vecti^\epsilon+2s)=\pi/2\text{if and
only if}\;\epsilon=0,\;\arg(p,q,r,s)=\pi/2,
\finishdisp
This completes the proof of the claim \eqref{eqn:argumentclaim}

Now note that
$\arg(p,q,r,s)=\pi/2$ is equivalent to
\startdisp
\begin{pmatrix}0&-\vecti\\-\vecti&0\end{pmatrix}
\begin{pmatrix}p&q\\r&s\end{pmatrix}\in\Xi_{2}\cap\SL{2}{\Real}
:=(\Xi_2)_{\Int}.
\finishdisp
Therefore, using \eqref{eqn:argumentclaim} and \eqref{eqn:gammaint21descriptioninter2}
we have
\begin{multline*}
\{\alpha'\;|\;\alpha\in\conj\inv(\Gamma)\cap\SL{2}{\Real},\,
\big(i,\delta\big)(\alpha)=(2,1)\}=\\
\bigcup_{\epsilon=0,1}\hspace*{-0.45cm}\cdot\hspace*{.55cm}
\{\alpha'\;\text{of the form \eqref{eqn:alphaprimemultform}}\;
\text{with}\;\arg(p,q,r,s)=\pi/2\}
\end{multline*}

Applying these observations to \eqref{eqn:gammaint21descriptioninter2},
we obtain
\startdisp
\{\alpha'\;|\;\alpha\in\conj\inv(\Gamma_{\Int}),\,
\big(i,\delta\big)(\alpha)=(2,1)\}=
\begin{pmatrix}0&\vecti\\ \vecti&0\end{pmatrix}(\Xi_2)_{\Int}
\alpha^{-2}\hspace*{-0.7mm}
(-1,1).
\finishdisp
Therefore,
\startdisp
\begin{aligned}
\big(i,\delta\big)\inv(2,1)\cap\SL{2}{\Real}&=&&\frac{1}{\sqrt{2}\vecti}
\begin{pmatrix}0&\vecti\\ \vecti&0\end{pmatrix}(\Xi_2)_{\Int}
\alpha^{-2}\hspace*{-0.7mm}(-1,1)\\
&=&&\frac{1}{\sqrt{2}}
\begin{pmatrix}0&1\\ 1&0\end{pmatrix}(\Xi_2)_{\Int}
\alpha^{-2}\hspace*{-0.7mm}(-1,1)\\
&=&&\frac{1}{\sqrt{2}}
\begin{pmatrix}0&1\\ 1&0\end{pmatrix}\conj(-I)(\Xi_2)_{\Int}
\alpha^{-2}\hspace*{-0.7mm}(-1,1)\\
&=&&\frac{1}{\sqrt{2}}
\begin{pmatrix}0&-1\\ -1&0\end{pmatrix}(\Xi_2)_{\Int}
\begin{pmatrix}-1&0\\0&-1\end{pmatrix}
\alpha^{-2}\hspace*{-0.7mm}(-1,1)\\
&=&&\frac{1}{\sqrt{2}}
(\Xi_2)_{\Int}
\alpha^{2}\hspace*{-0.7mm}(1,-1),
\end{aligned}
\finishdisp
where to obtain the last line we have used the fact established in
\S\ref{subsec:flts} that $\Xi_2$ is preserved under the action
of $\Xi_{12}$ by left-multiplication.  Using the definition
of $\alpha^{2}\hspace*{-0.7mm}(1,-1)$, we have
\starteqn\label{eqn:gammaint21description}
\big(i,\delta\big)\inv(2,1)\cap\SL{2}{\Real}=
\frac{1}{\sqrt{2}}
(\Xi_2)_{\Int}\begin{pmatrix}1&-1\\
0&2\end{pmatrix}.
\finisheqn

Gathering together \eqref{eqn:gammaint00description} and
\eqref{eqn:gammaint21description} we deduce \eqref{eqn:gammaintdescription}
\end{myproof}

From \eqref{eqn:gammaintdescription}, we deduce the analogue of
Lemma \ref{lem:indexlemma}
\begin{lem}\label{lem:indexlemmarcase}  Let $\conj\inv(\Gamma_{\Int})$
be the discrete subgroup of $\SL{2}{\Real}$ defined in \ref{eqn:gammazdefn},
and given explicitly in matrix form
in \eqref{eqn:gammaintdescription}.  All the
other notation is also as in
Proposition \ref{prop:gammaintexplicitdescrip}.
\begin{itemize}\item[(a)]  We have
\startdisp
\conj\inv(\Gamma_{\Int})\cap\SL{2}{\Int}=(\Xi_{12})_{\Int}.
\finishdisp
\item[(b)]  We have
\starteqn
[\conj\inv(\Gamma_{\Int}):(\Xi_{12})_{\Int}]=2,\quad [\SL{2}{\Int}:\Xi_{12}]=3.
\finisheqn
Explicitly, a representative of the unique non-identity right coset of
$(\Xi_{12})_{\Int}$ in $\conj\inv(\Gamma)$ is
\startdisp
\frac{1}{\sqrt{2}}
\begin{pmatrix}1&0\\1&1\end{pmatrix}\begin{pmatrix}1&-1\\0&2\end{pmatrix}.
\finishdisp
\end{itemize}
\end{lem}

\subsection{Explicit determination of fundamental
domain for $\SOM{2}{1}_{\Int}$ acting on $\bbH^2$.}
The main point of this section is that, provided the fundamental
domain $\scrG_{\Real}$ of the the standard unipotent subgroup
of $\conj\inv(\Gamma_{\Int})$ is chosen in a way that is compatible with the
choice of $\scrG$ in \eqref{eqn:scrGdefn}, then the good Grenier
fundamental domain $\scrF_{\Real}(\scrG_{\Real})$ for
$\conj\inv(\Gamma_{\Int})$ corresponding to $\scrG_{\Real}$
will have a close geometric relationship to $\scrF(\scrG)$.  Based
on the classical example of Dirichlet's fundamental domain
for $\SL{2}{\Int}$ acting on $\bbH^2$ and the Picard domain,
one might guess that we would have the equality
\startdisp
\scrF_{\Real}(\scrG_{\Real})=\scrF(\scrG)\cap\bbH^2.
\finishdisp
In fact, this intersection property cannot hold, because of the presence
of additional torsion elements (the powers of $\omega_8I_2$) in
$\conj\inv(\Gamma)$.  However, in a sense which will be made precise
in Proposition \ref{prop:funddomsrelation}, below, the next best thing holds.
Namely, the intersection
of the set consisting of \textit{two} $\Gamma$-translates of $\scrF(\scrG)$
with $\bbH^2$ equals $\scrF_{\Real}(\scrG_{\Real})$, for the choice
of $\scrG_{\Real}$ in \eqref{eqn:scrGRdefn}, below.
\vspace*{0.3cm}

\noindent \bmth\textbf{Grenier domains for discrete
subgroups of $\SL{2}{\Real}$.}\ubmth\hspace*{2mm}
We begin by sketching without proof the foundations of the subject in terms
consistent with \S\ref{subsec:grenierh3}.
Define the coordinate mappings
$\bbH^2\rightarrow \Real$, for $i=1,2$, by
\startdisp
\phi_1=-\log y,\;\phi_2=x,
\finishdisp
and set
\startdisp
\phi=(\phi_1,\phi_2): \bbH^2\rightarrow\Real^3.
\finishdisp
The mapping $\phi$ is a diffeomorphism of $\bbH^2$ onto $\Real^2$,
\textit{cf.} the diffeomorphism $\phi$ of $\bbH^3$
onto $\Real^3$, defined by \eqref{eqn:phicoordsdefn}.  Let $\Gamma$
be a discrete subgroup of $\mathrm{Diffeo}(\bbH^2)$,
Let $\Gamma^{\phi_1}$ be the stabilizer of the coordinate
$\phi_1$,
$\Delta_{1,\gamma}$ the difference function, and $\sigma_2^0$
the zero section of the projection $\phi_2$.  In other words,
all have the same meanings as in Section \S\ref{subsec:grenierh3}.
In addition, the statements of the
four \textbf{A} Axioms, stated after the proof
of Lemma \ref{lem:action1cons}, remain the same, except
that we take $z\in\bbH^2$ instead of $\bbH^3$.
We now state the analogue of Theorems \ref{thm:grenierh3}
and \ref{thm:grenierh3extended} that we wish to use in this setting.
\begin{thm}\label{thm:grenierh2extended}  Let $\Gamma$ be a group of diffeomorphisms of $\bbH^2$.
Assume that the action of $\Gamma$ on $\bbH^2$ satisfies axioms \textbf{A 1}
through \textbf{A 4}.  Let $\scrG\subseteq\Real$ be a
fundamental domain for the induced action of
$\Gamma^{\phi_{[1,2]}}\backslash \Gamma^{\phi_{1}}$ on $\Real$. Suppose that
$\scrG=\overline{\intrr(\scrG)}$.  Define
\startdisp
\scrF_1=\{z\in \bbH^2\;|\;
\phi_{1}(z)\leq \phi_{1}(\gamma z),\forallindisp
\gamma\in\Gamma^{\phi_1}\}.
\finishdisp
Set
\startdisp
\scrF=\phi_{2}\inv(\scrG)\cap \scrF_1. \finishdisp
Then we have
\begin{itemize}
\item[(a)] $\scrF$ is a
fundamental domain for the action of $\Gamma^{\phi_{[1,2]}}\backslash
\Gamma$ on
$\bbH^2$ (where $\Gamma^{\phi_{[1,2]}}$ is the kernel
of the action).
\item[(b)] We have
\startdisp
\scrF_1=\overline{\intrr\scrF_1},
\finishdisp
with
\startdisp \intrr\scrF_1=\{z\in \bbH^2\;|\;
\phi_{1}(z)< \phi_{1}(\gamma x),\forallindisp\;
\gamma\in\Gamma-\Gamma^{\phi_{1}}\}, \finishdisp
and
\startdisp
\partial\scrF_1=\{z\in\scrF_1\;|\; \phi_{1}(z)=\phi_{1}(\gamma z),\;\;
\text{for some} \;\gamma\in\Gamma-\Gamma^{\phi_1}\}.
\finishdisp
\item[(c)]  We have,
further,
\startdisp
\intrr\scrF=\phi_{2}\inv(\intrr\scrG)\cap
\intrr(\scrF_1),
\finishdisp
and
\startdisp
\scrF=\overline{\intrr\scrF}.
\finishdisp
\end{itemize}
\end{thm}
As a matter of terminology, we note that for a
co-finite $\Gamma\in\mathrm{Aut}^+(\bbH^2)$ (called `Fuchsian of the first
kind' in the literature), our notion of the good Grenier
fundamental domain $\scrF$ corresponds to the \textit{Ford
fundamental domain} of the Fuchsian group $\Gamma$.
See, for example, \cite{iwaniec95}, p. 44.  However,
we use the terminology Grenier domain even in this context, in
order to stress the eventual connections with the higher-rank case.

The following results, Lemmas \ref{lem:rcasea1a2a4},
\ref{lem:rcasea3comm},  and \ref{lem:rcasea3part},
guarantee that we can apply Theorem \ref{thm:grenierh2extended}
to a wide class of subgroups of $\SL{2}{\Real}$, including
$\conj\inv(\Gamma_{\Int})$, in particular.
\begin{lem} \label{lem:rcasea1a2a4} If $\Gamma\subseteq\SL{2}{\Real}$,
acting on $\bbH^2$ from the left by fractional linear transformations,
then $\Gamma$ satisfies Axioms \textbf{A 1}, \textbf{A 2}, and \textbf{A 4}.
\end{lem}
Compare to Lemmas \ref{lem:a2} and \ref{lem:a1a2}, the corresponding
statements in the complex case.
\begin{lem} \label{lem:rcasea3comm}  For any commensurability
class of subgroups of fractional linear transformations
in $\SL{2}{\Real}$, either all the groups in the class
satisfy Axiom \textbf{A 3}, or all the groups fail to satisfy
Axiom \textbf{A 3}.
\end{lem}
\begin{lem}\label{lem:rcasea3part}  The group $\SL{2}{\Int}$,
acting on $\bbH^2$ by fractional linear transformations,
satisfies Axiom \textbf{A 3}.
\end{lem}
Lemmas \ref{lem:rcasea3comm} and \ref{lem:rcasea3part} correspond
to Corollary \ref{cor:a3} and Proposition \ref{prop:sl2zia3},
respectively, in the complex case.
As in the case of Lemmas \ref{lem:rcasea1a2a4},
the proofs are the same as those of the corresponding
statements in the complex case, except at certain
points where they are even simpler, allowing the omission
of the proofs.

Using Lemma \ref{lem:rcasea1a2a4} and
\ref{lem:rcasea3comm}, as well as Lemmas
\ref{lem:indexlemmarcase} and \ref{lem:rcasea3part} to show that
Lemma \ref{lem:rcasea3comm} applies
to this situation, we deduce the following.
\begin{cor} \label{cor:gammaz4axioms} The group $\conj\inv(\Gamma_{\Int})$ of
Definition \ref{defn:gammaz}, acting on $\bbH^2$ from the left
by fractional linear transformations, satisfies the four \textbf{A}
axioms.
\end{cor}
Because of Corollary \ref{cor:gammaz4axioms}, we can apply
Theorem \ref{thm:grenierh2extended} to obtain a general
description of the fundamental domain $\scrF_{\Real}(\scrG_{\Real})$
for the action of $\conj\inv(\Gamma_{\Int})$.  We now
proceed to calculate $\scrF_{1,\Real}$ explicitly and provide a choice
for $\scrG_{\Real}$ which will be consistent with the earlier
choice of $\scrG$.\vspace*{3mm}

\noindent \bmth\textbf{Explicit calculation of $\scrF_{1,\Real}$.}\ubmth
\hspace*{2mm}
It is very easy to see that the following Proposition holds in general.
\begin{prop} \label{prop:scrF1scontainment} Let $\Gamma\subseteq\SL{2}{\Complex}$ and define
$\Gamma_{\Int}:=\Gamma\cap\SL{2}{\Real}$.  Set
\startdisp\begin{aligned}
\scrF_1&=&&\{z\in\bbH^3\;|\; y(z)\geq y(\gamma z)\;\text{for all}\;
\gamma\in\Gamma\},\;\text{(as in \eqref{thm:grenierh3})},\;\text{and}\\
\scrF_{1,\Real}&=&&\{z\in\bbH^2\;|\; y(z)\geq y(\gamma z)\;\text{for all}\;
\gamma\in\Gamma_{\Int}\}.
\end{aligned}
\finishdisp
Then we have
\starteqn\label{eqn:scrF1scontainment}
\scrF_{1}\cap\bbH^2\subseteq\scrF_{1,\Real}.
\finisheqn
\end{prop}
\begin{myproof} {Proof} A point $z\in\bbH^2$ belongs to the left-hand
side of \eqref{eqn:scrF1scontainment} if and only if
it has maximal $y$-coordinate among points in the orbit $\Gamma z$.
So if $z\in\bbH^3$ belongs to the left-hand side of \eqref{eqn:scrF1scontainment},
then $z$ \textit{a fortiori} has maximal $y$-coordinate among
points in the orbit $\Gamma_{\Int}z$, and so $z\in\scrF_{1,\Real}$.
\end{myproof}

Our immediate aim is to show that in the particular case when
$\Gamma=\conj\inv(\SO{3}{\Int[\vecti]})$, in the situation of
\ref{prop:scrF1scontainment}, we have equality in \eqref{eqn:scrF1scontainment}.
For that, the criterion of Lemma \ref{lem:scrf1rcontainment}, below,
will suffice.  In order to set up the lemma, let $i$ be an
indexing set.  Let $\{d_i\}_{i\in I}$
be a collection complex numbers, and let $\{\rho_i\}_{i\in I}$ be
a collection of positive integers, indexed by $I$.
Let $J\subset I$ be defined by the condition
\startdisp
\text{For}\;i\in I,\; \text{we have}\;i\in J\;\text{if and only if}\;
d_i\in\Real.
\finishdisp
Also, set
\startdisp
\mathfrak{d}=\min_{i\in I-J}\left\{|\Impt (d_i)|\right\},
\finishdisp
in other words, $\mathfrak{d}>0$ is the distance to the real
line of the $d_i$ that is closest to the real line
without actually being on it.  (In
case $I=J$, we can take $\mathfrak{d}=\infty$, and the subsequent
statements will remain true without modification.)
For $\kappa>0$, we define the following condition on the pair
$\left(\{d_i\}_{i\in I},\{\rho_i\}_{i\in I}\right)$, depending
on $\kappa$.
\starteqn\label{eqn:scrF1rcondition}\text{The set $\{d_j\}_{j\in J}$
is $\kappa$-dense in $\Real$ and
$\left(\min_{i\in J}(\rho_i)\right)^2-\kappa^2
\geq\left(\max_{i\in I-J}(\rho_i)\right)^2-\mathfrak{d}^2$.}
\finisheqn
We remind the reader of our notation $\bbH^2_{\vectj}$ for the vertical
plane $x_2=0$ in $\bbH^3$, with $\bbH^2_{\vectj}=\bbH^2$.
\begin{lem} \label{lem:scrf1rcontainment} Let the pair
$\left(\{d_i\}_{i\in I},\{\rho_i\}_{i\in I}\right)$ be as above,
and assume that there exists $\kappa>0$ such that
the pair $(\{\ell_i\},\{\rho_i\})$ satisfies \eqref{eqn:scrF1rcondition}.
for the given $\kappa$.  Then we have the containment
\begin{multline}\label{eqn:scrf1rcontainment}
\left\{z\in\bbH^2_{\vectj}\;\left|\;
\left\|\begin{pmatrix}1&d_j\end{pmatrix}\begin{pmatrix}z\\
1\end{pmatrix}\right\|\geq \rho_j,\;\text{for all}\;j\in J\right.\right\}
\vspace*{0.2cm}\\
\hspace*{3cm}\bigsubseteq
\left\{z\in\bbH^2_{\vectj}\;\left|\;
\left\|\begin{pmatrix}1&d_i\end{pmatrix}\begin{pmatrix}z\\
1\end{pmatrix}\right\|\geq \rho_i,\;\text{for all}\;i\in I\right.\right\}.
\end{multline}
\end{lem}
\begin{myproof}{Proof}  Assume that $z\in\bbH^2_{\vectj}$ is contained
on the left side of \eqref{eqn:scrf1rcontainment}.  In order to show
that $z$ is contained in the right side of \eqref{eqn:scrf1rcontainment}
it will suffice to show that $z$ satisfies the condition
\starteqn\label{eqn:scrF1rcontainmentproof}
\left\|\begin{pmatrix}1&d_i\end{pmatrix}\begin{pmatrix}z\\
1\end{pmatrix}\right\|\geq \rho_i,\;\text{for all}\;i\in I-J.
\finisheqn
The assumption that $z\in\bbH^2$ belongs to the left-hand side
of \eqref{eqn:scrf1rcontainment} means that for each $j\in J$,
\startdisp
\begin{aligned}
\left\|\begin{pmatrix}1&d_j\end{pmatrix}\begin{pmatrix}z\\
1\end{pmatrix}\right\|^2&=&&||z+d_j||^2\\
&=&&y(z)^2+(x(z)+d_j)^2\\
&\geq&&\rho_j^2.
\end{aligned}
\finishdisp
Rearranging the last inequality and taking the minimum of
the $\rho_j$, we have
\startdisp
y(z)^2\geq \rho_j^2-(x(z)+d_j)^2\geq \left(\min_{j\in J}\rho_j\right)^2-
(x(z)+d_j)^2
\finishdisp
The $\kappa$-density of $\{d_j\}_{j\in J}$ in $\Real$
is equivalent to the $\kappa$-density of $\{-d_j\}_{j\in J}$.  Therefore,
by \eqref{eqn:scrF1rcondition},
we can choose $j\in J$ so that $(x(z)+d_j)^2\leq\kappa^2$.
We deduce that
\starteqn\label{eqn:scrF1rcontainmentproof2}
y(z)^2\geq\left(\min_{j\in J}\rho_j\right)^2-
\kappa^2.
\finisheqn
For any $i\in I-J$, we calculate
\startdisp
\begin{aligned}
\left\|\begin{pmatrix}1&d_i\end{pmatrix}\begin{pmatrix}z\\
1\end{pmatrix}\right\|^2&=&&||z+d_i||^2\\
&=&&y(z)^2+||\Impt(d_i)||^2+||x(z)+\Rept(d_i)||^2\\
&\geq&&\left(\min_{j\in J}\rho_j\right)^2-
\kappa^2+||\Impt(d_i)||^2+||x(z)+\Rept(d_i)||^2\;
(\text{by}\; \eqref{eqn:scrF1rcontainmentproof2})\\
&\geq &&\left(\max_{i\in I-J}(\rho_i)\right)^2-
\mathfrak{d}^2+||\Impt(d_i)||^2+
||x(z)+\Rept(d_i)||^2\; (\text{by}\;\eqref{eqn:scrF1rcondition})\\
&\geq&&\max_{i\in I-J}(\rho_i))^2,
\end{aligned}
\finishdisp
where we have used the assumption that $z\in\bbH^2$ to obtain
the second line and the definition of $\mathfrak{d}$
to obtain the final line.  The above certainly implies
\eqref{eqn:scrF1rcontainmentproof}, so by the preceding comments,
the proof is complete.
\end{myproof}

In particular, we can apply Lemma \ref{eqn:scrf1rcontainment}
with
\startdisp
(\{d_i\}_{i\in I},\{\rho_i\}_{i\in I})=(1+(1+\vecti)\Int[\vecti],\sqrt{2})
\finishdisp
For this application, the indexing is unimportant because
all the $\rho_i$ are equal.  All that we care about is that
\startdisp
\{d_i\}_{i\in I}\cap\Real=1+2\Int.
\finishdisp
In that case, we have $\mathfrak{d}=1$, and we have
\eqref{eqn:scrF1rcondition} satisfied with $\kappa=1$.
The special case of \eqref{eqn:scrf1rcontainment} obtained
from the indicated substitution is \eqref{eqn:scrf1rcontainmentsp}
in Part (a) of the following Lemma.
\begin{lem} \label{lem:scrFr1} Let $\scrF_1$ and $\scrF_{1,\Real}$
be as in Proposition
\ref{prop:scrF1scontainment}.
\begin{itemize}
\item[(a)]  We have
\begin{multline}\label{eqn:scrf1rcontainmentsp}
\left\{z\in\bbH^2_{\vectj}\;\left|\;
\left\|\begin{pmatrix}1&d\end{pmatrix}\begin{pmatrix}z\\
1\end{pmatrix}\right\|\geq 2,\;d\in1+2\Int\right.\right\}
\vspace*{0.2cm}\\
\hspace*{3cm}\bigsubseteq
\left\{z\in\bbH^2_{\vectj}\;\left|\;
\left\|\begin{pmatrix}1&d\end{pmatrix}\begin{pmatrix}z\\
1\end{pmatrix}\right\|\geq 2,\;d\in 1+(1+\vecti)\Int[\vecti]\right.\right\}.
\end{multline}
\item[(b)]  We have $\scrF_{1,\Real}$ contained in the left-hand
side of \eqref{eqn:scrf1rcontainmentsp}, while the right-hand
side of \eqref{eqn:scrf1rcontainment} equals $\scrF_1\cap\bbH^2$.
\item[(c)]  We have
\starteqn\label{eqn:scrF1Rintersectionform}
\scrF_{1,\Real}=\scrF_1\cap\bbH^2.
\finisheqn
\end{itemize}
\end{lem}
\begin{myproof}{Proof}  As noted just prior
to the Lemma, \eqref{eqn:scrf1rcontainmentsp}
is a special case of \eqref{eqn:scrf1rcontainment}.  In
order to obtain the left-hand side of \eqref{eqn:scrf1rcontainmentsp},
we use the relation
\startdisp
\Int\cap\left(1+(1+\vecti)\Int[\vecti]\right)=1+2\Int.
\finishdisp

In part (b), the equality between $\scrF_1\cap\bbH^2$
and the right-hand side of \eqref{eqn:scrf1rcontainment}
follows immediately from \eqref{eqn:scrF1gammaexplicit}.
Taken together with \eqref{eqn:scrF1scontainment},
Part (b) immediately implies Part (c).
So all that is left to verify is the containment of $\scrF_{1,\Real}$
in the left-hand side of \eqref{eqn:scrf1rcontainmentsp}.

The proof of the last containment is a routine
modification of ideas found in the proof of Lemma
\ref{lem:scrF1determination},
so we will just sketch the proof and specify the points
of departure from the proof of Lemma \ref{lem:scrF1determination}.
First, we need the analogue for $\scrF_{1,\Real}$ of the description
of $\scrF_1$ given in \eqref{eqn:scrF1altdescr}, which is
\starteqn\label{eqn:scrF1raltdescr}
\scrF_{1,\Real}=\left\{z\in\bbH^2\;\left|\;
\left\|\begin{pmatrix}0&1\end{pmatrix}\gamma
\begin{pmatrix}z\\1\end{pmatrix}\right\|^2\geq 1,\;\text{for all}\;
\gamma\in\Gamma_{\Int}\right.\right\}.
\finisheqn
The proof of \eqref{eqn:scrF1raltdescr}
is the same as that of \eqref{eqn:scrF1altdescr}.
Now, we use Proposition \ref{prop:gammaintexplicitdescrip},
the analogue for the present circumstances
of Proposition \ref{prop:inversematrixexplicitdescription},
to obtain \eqref{eqn:nonintegersetrowreal}, below.  The inclusions of
\eqref{eqn:nonintegersetrowreal}
are analogous to \eqref{eqn:nonintegersetrow} in the present
circumstances, since they state that we have
\starteqn\label{eqn:nonintegersetrowreal}
\frac{1}{\sqrt{2}}\begin{pmatrix}1&1+2\Int\end{pmatrix}\hspace*{2mm}\bigsubseteq
\hspace*{2mm}\frac{1}{\sqrt{2}}
\begin{pmatrix}0&1\end{pmatrix}
(\Xi_2)_{\Int}\begin{pmatrix}1&-1\\0&2\end{pmatrix}
\hspace*{2mm}\bigsubseteq
\hspace*{2mm}\begin{pmatrix}0&1\end{pmatrix}\Gamma_{\Int}.
\finisheqn
Then Properties \textbf{5. Scalar Function Multiple} and
\textbf{3. Inclusion Reversal} are applied to the right-hand
side of \eqref{eqn:scrF1raltdescr}, using \eqref{eqn:nonintegersetrowreal}.
In this way, one readily obtains the inclusion of
$\scrF_{1,\Real}$ in left-hand side of \eqref{eqn:scrf1rcontainmentsp},
thus completing the proof of the proposition.
\end{myproof}

\vspace*{3mm}
\noindent \bmth\textbf{Explicit Descriptions of $\scrG_{\Real}$
and $\scrF_{\Real}(\scrG_{\Real})$.}\ubmth
\begin{lem}
\begin{itemize}
\item[(a)]
We have
\startdisp
(\Gamma_{\Int})^{\phi_1}=\begin{pmatrix}1&2\Int\\
0&1\end{pmatrix}.
\finishdisp
\item[(b)]  The interval
\starteqn\label{eqn:scrGRdefn}
\scrG_{\Real}:=[0,2]
\finisheqn
is a fundamental domain for the action of $\Gamma_{\Int}^{\phi_1}$
on $\Real$ satisfying
\startdisp
\scrG_{\Real}=\overline{\intrr\scrG_{\Real}}.
\finishdisp
\item[(c)]  With $\scrG_{\Real}$ as defined in \eqref{eqn:scrGRdefn},
part (b) implies that
\starteqn\label{eqn:scrFRscrGRdescription}
\begin{aligned}
\scrF_{\Real}\left(\scrG_{\Real}\right)&=&&\{z\in\bbH^2\;
|\;0\leq x(z)\leq 2,\; y(z)^2+(x-1)^2\geq 2\}\\
&=&&\coh(\vecti,2+\vecti,\infty).
\end{aligned}
\finisheqn
\end{itemize}
\end{lem}
\begin{myproof}{Proof}
By \eqref{eqn:gammazdefn}, we have
\startdisp
\begin{aligned}
(\Gamma_{\Int})^{\phi_1}&=&&(\conj\inv(\Gamma)\cap\SL{2}{\Real})^{\phi_1}\\
&=&&\conj\inv(\Gamma)^{\phi_1}\cap\SL{2}{\Real}.\\
\end{aligned}
\finishdisp
We then apply \eqref{eqn:gammaphi1explicit} to conclude the proof
of (a).

Part (b) is immediate from part (a).

For part (c), part (b) of this lemma implies that $\scrG_{\Real}$
satisfies the hypotheses of Theorem \ref{thm:grenierh2extended}.
Using \eqref{eqn:scrF1Rintersectionform} and the definition
of $\scrF(\scrG)=\scrF_{\Real}(\scrG_{\Real})$ given in Theorem
\ref{thm:grenierh2extended}, we have
\startdisp
\begin{aligned}
\scrF_{\Real}(\scrG_{\Real})&=&&\scrF_{1,\Real}\cap\scrG_{\Real}\\
&=&&(\scrF_1\cap\bbH^2)\cap\scrG_{\Real}.
\end{aligned}
\finishdisp
The first description of $\scrF_{\Real}\left(\scrG_{\Real}\right)$
given in \eqref{eqn:scrFRscrGRdescription} then follows
from the description of $\scrF_1$ given in \eqref{eqn:scrFscrGfirstform}.
The second description in \eqref{eqn:scrFRscrGRdescription}
follows directly from the first form
because the inequalities in the first description correspond
to the half-space bounded by the three geodesics
\starteqn\label{eqn:boundinggeod}
(0,\infty)_{\bbH},\;(2,\infty)_{\bbH},\;\text{and}\;(-1,3)_{\bbH}.
\finisheqn
So the inequalities in the first description define
an (infinite) hyperbolic triangle.  We determine the corners
of the triangle by verifying that the
(finite) intersection points of the bounding
geodesics in \eqref{eqn:boundinggeod} are $\vecti,2+\vecti$.
\end{myproof}

\vspace*{3mm}
\noindent \bmth\textbf{Geometric
relation of $\scrF_{\Real}(\scrG_{\Real})$ to $\scrF(\scrG)$.}\ubmth\;
Because
\starteqn\label{eqn:scrGRunion}
\scrG_{\Real}=\left(\scrG\cup\conj(T_1)\left(R_{\frac{\pi}{2}}^2\right)\scrG
\right)\cap
\bbH^2_{\vectj},
\finisheqn
(see the proof of Proposition \ref{prop:funddomsrelation}, below)
we cannot hope that we will have the straightforward relation
\startdisp
\scrF_{\Real}\left(\scrG_{\Real}\right)=\scrF\left(\scrG\right)\cap\bbH^2_{\vectj}
\finishdisp
that we find in the classical case of  $\SL{2}{\Int[\vecti]}$ and
$\SL{2}{\Int}$.
However, we do have the next best possible relation between
the fundamental domains.

\begin{prop} \label{prop:funddomsrelation} We have the relation
\startdisp
\scrF\left(\scrG_{\Real}\right)=\left(\scrF(\scrG)\cup
\conj(T_1)\left(R_{\frac{\pi}{2}}^2\right)\scrF(\scrG)\right)\cap
\bbH^2_{\vectj}.
\finishdisp
\end{prop}

\begin{myproof}{Proof}  We begin by proving \eqref{eqn:scrGRunion}.
This is done by examining the effect of the element
$\conj(T_1)(R_{\frac{\pi}{2}}^2)$ on the endpoints $1,2$ of
$\scrG\cap\bbH^2_{\vectj}$.  We find that
\startdisp
\conj(T_1)\left(R_{\frac{\pi}{2}}^2\right)(1)=1,\quad
\text{and}\quad \conj(T_1)\left(R_{\frac{\pi}{2}}^2\right)(2)=0,
\finishdisp
so that
\startdisp
\conj(T_1)\left(R_{\frac{\pi}{2}}^2\right)\scrG_{\Real}=[0,1].
\finishdisp
Since $\scrG\cap\bbH^2_{\vectj}=[0,2]$, this completes the proof
of \eqref{eqn:scrGRunion}.

Using, successively, the definitions of $\scrF(\scrG)$ and
$\scrF_{\Real}(\scrG_{\Real})$,
given in Theorems \ref{thm:grenierh3extended} and \ref{thm:grenierh2extended},
we have
\startdisp
\begin{aligned}
\scrF_{\Real}(\scrG_{\Real})&=&&\scrF_{1,\Real}\cap\scrG_{\Real}\\
&=&&\left(\scrF_{1}\cap\bbH^2\right)\bigcap\left((\scrG\cup\conj
(T_1)(R_{\frac{\pi}{2}}^2)\scrG)\cap
\bbH^2_{\vectj}\right)
(\text{by}\;\eqref{eqn:scrF1Rintersectionform}\;\text{and}\;\eqref{eqn:scrGRunion})\\
&=&&\left(\scrF_1\cap\scrG\cap\bbH^2\right)\bigcup\left(\scrF_1\cap\conj(T_1)
(R_{\frac{\pi}{2}}^2)\scrG\cap\bbH^2\right) \\
&=&&\left(\left(\scrF_1\cap\scrG\right)\cup
\conj(T_1)(R_{\frac{\pi}{2}}^2)(\scrF_1\cap
\scrG)\right)\bigcap\bbH^2\\
&=&&\left(\scrF(\scrG)\cup\conj(T_1)\left(R_{\frac{\pi}{2}}^2\right)\scrF(\scrG)
\right)\cap\bbH^2.
\end{aligned}
\finishdisp
In the equality of the second and third line, we implicitly
use the invariance of $\scrF_1$ under the action of $\Gamma^{\phi_1}$,
which is obvious from the definition of $\scrF_1$.
\end{myproof}

\begin{wrapfigure}{r}{5cm}
\includegraphics[scale=0.5]{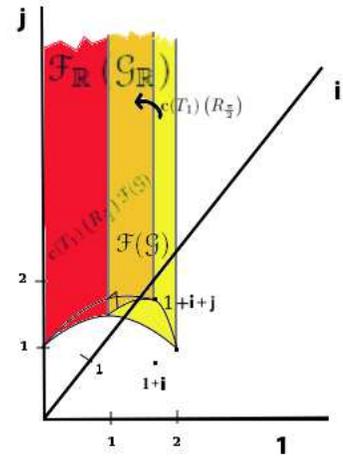}
   \caption{Relation of Fundamental domains of $\Gamma$ and $\Gamma_{\Int}$}
   \end{wrapfigure}

\begin{rem}  We also note for possible future reference that
$\scrF_{\Real}(\scrG_{\Real})$ is the \textit{normal geodesic projection}
of the union of $\scrF(\scrG)$ and one translate
$\conj(T_1)\left(R_{\frac{\pi}{2}}^2\right)\scrF(\scrG)$ of $\scrF(\scrG)$.
This relation between the fundamental domains
is connected to the one given in Proposition \ref{prop:funddomsrelation},
though neither relation implies the other, in general.
In Figure 1, we have indicated by means of a ``right-angle" symbol
at the point $1+\sqrt{2}\vectj$ that the geodesic $\bbH^1(1+\sqrt{2}\vectj,
1+\vecti+\vectj)$ is a geodesic normal to $\bbH^2_{\vectj}$.
It would take us to far afield of our main purpose to define
the concept of \textit{normal geodesic projection} precisely,
so for the moment we restrict ourselves to mentioning that this
relation between $\scrF(\scrG)$
and $\scrF_{\Real}(\scrG_{\Real})$ may be of some
use in relating spectral expansions in the complex case
to spectral expansions in the real case.
\end{rem}

\section{Volume Computations}
\label{sec.volumecomputation}
\subsection{Volumes}
The following are reasons for computing the covolume
of the lattices we are studying.  First
the volume computation is a relatively simple application of
the previous results of a group-theoretic nature
 (that is, Lemmas \ref{lem:indexlemma} and
\eqref{lem:indexlemmarcase}) to obtain quantitative
geometric results that can be compared
to previously obtained results in the literature.  Second, the volume
of the fundamental domain, or more precisely,
the reciprocal of the volume of the fundamental domain
under the Iwasawa-Haar measure, appears as a constant term in the
Laplacian-eigenfunction expansion of automorphic functions.  See, for example,
\cite{jol06}, Theorem XI.4.2, for the case of $\SL{2}{\Int[\vecti]}$
in $\SL{2}{\Complex}$.  Therefore, the result
of Corollary \ref{cor:gammacovolume} will play a direct role
in the next phase of this investigation.   Third, covolume
computations have applications in the theory of moduli space
in algebraic geometry.   For example, \cite{grithulek06}
uses a computation of the Hirzebruch-Mumford volume of an arithmetic
lattice $\Gamma$ in a real form of $G$ to compute the leading
term of growth of the space of cusp forms $S_k(\Gamma)$.
The Hirzebruch-Mumford volume, which is defined
precisely in \cite{grithulek06}, is a suitably normalized volume
of the quotient $\Gamma\backslash G/K$.
This, however, pertains to a future, projected stage of our project,
in which the geometry of the fundamental domains and the
eigenfunction expansion will be connected to the theory
of certain moduli spaces of $K3$-surfaces.
\vspace*{0.3cm}

\subsection{Ratio of covolume of $\Gamma$, resp. $\Gamma_{\Int}$,
to covolume of $\SL{2}{\Int[\vecti]}$, resp. $\SL{2}{\Int}$}
\label{subsec:ratiocovolumes}
\begin{lem} \label{lem:indexcovolumes} Let $X$ be a homogenous space
of a semi-simple Lie group $G$,
equipped with a $G$-invariant measure $\mu$.
Let $\Gamma'$, $\Gamma$ be commensurable lattices in $G$ which
act on $X$ with finite volume quotient.  Then we have
\starteqn\label{eqn:indexcovolumes}
\mu(M/\Gamma)=\mu(M/\Gamma')\frac{[\Gamma':\Gamma'\cap\Gamma]}
{[\Gamma:\Gamma'\cap\Gamma]}
\finisheqn
\end{lem}
\begin{myproof}{Proof}  The lemma follows formally from the special
case when $\Gamma'\subseteq\Gamma$, \textit{i.e.}, when
\mbox{$[\Gamma':\Gamma'\cap\Gamma]=1$}.  The reason is that
we can apply the special case of \eqref{eqn:indexcovolumes}
in the two cases
\startdisp
(\Gamma,\Gamma')=(\Gamma,\Gamma'\cap\Gamma)\;\text{resp.,}\;
=(\Gamma',\Gamma'\cap\Gamma),
\finishdisp
then divide the results
to deduce the general case of \eqref{eqn:indexcovolumes}.

So assume that $\Gamma'\subseteq\Gamma$, with
$[\Gamma:\Gamma']=r\leq\infty$.
By the assumption of commensurability there exists a finite set
of elements $\gamma_i$,
$i=1,\ldots, r$, with $r=[\Gamma:\Gamma']$, such that
\startdisp
\Gamma=\Gamma'\gamma_1\cup\hspace*{-2.5mm}\cdot\hspace{1.5mm}
\cdots\cup\hspace*{-2.5mm}\cdot\hspace{2.5mm}\Gamma'\gamma_r.
\finishdisp
Let $\scrF$ be a fundamental domain for the action
of $\Gamma$ on $X$.  Set
that
\startdisp
\scrF':=
\gamma_1\scrF\cup\cdots\cup\gamma_r\scrF,\;\text{with distinct
translates disjoint, except for their common boundary.}
\finishdisp
then $\scrF'$ is a fundamental domain for the action of $\Gamma'$ on $X$,
so that $\mu(\scrF)=\mu(\Gamma\backslash M)$.
The non-overlapping property of $\scrF$ for $\Gamma$ implies that
the $\Gamma$-translates of $\scrF$ appearing in the union defining
$\scrF'$ intersect only on their boundaries, hence in sets
of measure zero.
Thus
\startdisp
\mu(\scrF')=r\mu(\scrF)=[\Gamma:\Gamma\cap\Gamma']\mu(\scrF)
\finishdisp
Re-arranging, we obtain \eqref{eqn:indexcovolumes} in the
special case, thus completing the proof.
\end{myproof}

For the next proposition, let $\mu$ represent Iwasawa-Haar measure on $\SL{2}{\Complex}$,
or $\SL{2}{\Real}$, depending on the context.
Using Lemma \ref{lem:indexlemma}(b), resp.,
Lemma \eqref{lem:indexlemmarcase}, we deduce from
Lemma \eqref{lem:indexcovolumes} formulas giving the covolumes of
$\conj\inv(\Gamma)$ and $\conj\inv(\Gamma_{\Int})$,
in terms of the covolumes of the standard integer subgroup,
valid for any normalization of the Haar measure on $G$.
\begin{prop} \label{prop:ratiocovolumes} Let $\Gamma=\SO{3}{\Int[\vecti]}$.  Let $\conj$ the
isomorphism of $\SL{2}{\Complex}/\{\pm 1\}$ onto $\SO{3}{\Complex}$
induced by the conjugation action of $\SL{2}{\Complex}$ on its
Lie algebra.  Let $\Gamma_{\Int}\cong\SOM{2}{1}_{\Int}$
be as in Definition \ref{defn:gammaz}.  Then we have
\startdisp
\mu(\conj\inv(\Gamma)\backslash\SL{2}{\Complex})=\half
\mu(\SL{2}{\Int[\vecti]}\backslash\SL{2}{\Complex}),
\finishdisp
and
\startdisp
\mu(\conj\inv(\Gamma_{\Int})\backslash\SL{2}{\Real})=\frac{3}{2}
\mu(\SL{2}{\Int}\backslash\SL{2}{\Real}).
\finishdisp
\end{prop}

Combining the well-known Corollary \ref{cor:gammacovolume}
(resp., \ref{cor:gammazcovolume}),
below, with Proposition \ref{prop:ratiocovolumes}, we determine
the covolume of $\conj\inv(\Gamma)$ (resp., $\conj\inv(\Gamma_{\Int})$)
in $\SL{2}{\Complex}$ (resp., $\SL{2}{\Real}$).
\begin{cor} \label{cor:gammacovolumelanglandsform}  Let $\Gamma=\SO{3}{\Int[\vecti]}$.  Let $\conj$ the
isomorphism of $\SL{2}{\Complex}/\{\pm 1\}$ onto $\SO{3}{\Complex}$
induced by the conjugation action of $\SL{2}{\Complex}$ on its
Lie algebra.  Let $\Gamma_{\Int}\cong\SOM{2}{1}_{\Int}$
be as in Definition \ref{defn:gammaz}.  Then we have
\startdisp
\mu(\conj\inv(\Gamma)\backslash\SL{2}{\Complex})=\frac{1}{2}\Bflambda
\zeta_{\Rational(\vecti)}(2).
\finishdisp
and
\startdisp
\mu(\conj\inv(\Gamma_{\Int})\backslash\SL{2}{\Real})=\frac{3}{2}\Bflambda
\zeta(2).
\finishdisp
\end{cor}


\subsection{Covolumes $V_n$ and $V_{n,\Real}$ of Gaussian
and rational integer subgroups.}
\label{subsec:covolumes}
 Although the result of Corollary
\ref{cor:gammacovolume} is well-known, see \textit{e.g.}
\cite{langlandsvolumes}, the method, following Siegel's
proof of the rational-integer analogue,
is not found in any standard reference.  We include these results
as an appendix in
order to keep the treatment here of co-volumes self-contained
and in order to show that there is a completely elementary
route to the calculation of the covolumes in terms of special
values of zeta functions.
\vspace*{3mm}

\noindent\bmth\textbf{Siegel's Theorem on $V_n$.}\hspace*{0.5mm}
\ubmth
In the following calculations, we closely follow the calculations
in \S$\romtwo$.4 of \cite{posnr}.  There, one finds
an elementary calculation, following the approach of Siegel,
resulting in 4.1.3,
the covolume of $\SL{n}{\Int}$ in $\SL{n}{\Real}$.
In the following arguments,
culminating in Theorem \ref{thm:siegelvolume},
and Corollary \ref{cor:gammacovolume} we in effect
show that the arguments of \cite{posnr}, \S$\romtwo$.4,
transfer in the most direct imaginable
manner to the Gaussian integer case.
Although, for the case at hand, we only need the case $n=2$,
it is natural to treat the case of general $n$.  The treatment
of the quadratic models $\mathrm{Pos_n}(\Complex)$ and
$\mathrm{Spos}_n(\Complex)$ and the decomposition of the
$G$-invariant measures on these spaces will be useful for the extension
of the theory of the present monograph to higher rank.

\begin{defn}  The \bmth\textbf{Dedekind Zeta function $\zeta_{\Rational
(\vecti)}$ associated
to the number field $\Rational(\vecti)$}\ubmth\hspace*{.25mm} is defined
by
\starteqn\label{eqn:dedekindzetadefn}
\begin{aligned}
\zeta_{\Rational(\vecti)}(s)&=&&\sum_{I\;\text{ideal of}\;\Int[\vecti]}
N(I)^{-s}\\
&=&&\sum_{\ell\;\text{standard}}||\ell||^{-2s},
\end{aligned}
\finisheqn
where, in the first line, $N$ denotes the norm from $\Rational(\vecti)$
to $\Rational$, and in the second line, the sum is over
all standard elements of the ring $\Int[\vecti]$, in the
sense defined in \S\ref{subsec:flts}.
\end{defn}

The expression for $\zeta_{\Rational(\vecti)}$ given in the first
line of \eqref{eqn:dedekindzetadefn} is stated in terms of notions which
generalize to an arbitrary number fields $F$, and thus gives
the definition of the Dedekind zeta function associated to $F$.
In order to pass from the expression in the first line of
\eqref{eqn:dedekindzetadefn} to that of the second line,
one uses the well-known characterization of the norm,
\startdisp
N(I)=\#(\Go/I).
\finishdisp
Then one applies the one-to-one correspondence between
ideals of the principal ring $\Int[\vecti]$ and standard
integers (assigning to the ideal $I$ its unique standard generator $\ell$),
and one easily calculates
\startdisp
N(I)=\#(\Int[\vecti]/I)=\#(\Int[\vecti]/(\ell))=||\ell||^2.
\finishdisp

We tabulate some subgroups of $\SL{n}{\Complex}$ which will
allow inductive decompositions.
\starteqn\label{eqn:mirabolicdefn}
\begin{aligned}
G_n&=&&\SL{n}{\Complex}.\\
\Gamma_n&=&&\SL{n}{\Int[\vecti]}.\\
G_{n,n}&=&& \text{subgroup of $G_n$ leaving the $n^{\rm th}$ unit
 vector fixed.}\\
\Gamma_{n,n}&=&&\Gamma_n\cap G_{n,n}.
\end{aligned}
\finisheqn
\hspace*{0.3cm}

\noindent\bmth\textbf{Scholium on $G_{n,n}$.}\ubmth\hspace*{0.5mm} Subgroups
of $G$ defined by a condition of the type used in
\eqref{eqn:mirabolicdefn} to define $G_{n,n}$
are known in the literature as \textbf{mirabolic subgroups}.
Specifically, $G_{n,n}$ is known as the mirabolic subgroup
associated to the maximal parabolic preserving the the flag
$\{\Complex e_n\}$.  Directly from the definition of $G_{n,n}$
in \eqref{eqn:mirabolicdefn}, we compute that
\starteqn\label{eqn:mirabolicmatrixform}
G_{n,n}=\left\{\left.\begin{pmatrix}g'&\vectx\\0&1\end{pmatrix}\;
\right|\;\vectx^*\in\Complex^{n-1},\; g'\in\SL{n-1}{\Complex}\right\}.
\finisheqn

\hspace*{0.3cm}

To facilitate comparison with \cite{posnr},
we note here that we diverge from
their notation in keeping the subscript $n$, on $G_n$, and
its associated objects ($G_{n,n}$, $\Gamma_{n,n}$ ,etc.)
This subscript is maintained in
order to distinguish the discrete groups being considered
here (the $\SL{n}{\Int[\vecti]}$) from the
discrete group $\Gamma=\SO{3}{\Int[\vecti]}$ under consideration
in the other sections of this chapter.

There is a natural isomorphism of $G_n$-homogeneous spaces
\starteqn\label{eqn:gmodgn1}
G_{n,n}\backslash G_n\stackrel{\approx}{\longrightarrow}
(\Complex^n)^*-\{0\}
\;\;\text{given by}\; g\mapsto e_n^*g.
\finisheqn
We have a fibering
\starteqn\label{eqn:gamman1modgn1fibering}
\Int[\vecti]^{n-1}\backslash\Complex^{n-1}\rightarrow
\Gamma_{n,n}\backslash G_{n,n}\rightarrow \Gamma_{n-1}\backslash G_{n-1}
=\SL{n-1}{\Int[\vecti]}\backslash \SL{n-1}{\Complex}.
\finisheqn
arising from the coordinates $(\vectx,g')$ on $G_{n,n}$, above, in
\eqref{eqn:mirabolicmatrixform}.  From the fibration
in \eqref{eqn:gamman1modgn1fibering}, we deduce
that $\Gamma_{n,n}\backslash G_{n,n}$ has finite measure,
under the inductive assumption that $\Gamma_{n-1}\backslash G_{n-1}$
has finite measure.

Because of \eqref{eqn:gmodgn1}, we can transport
Lebesgue measure from $\Complex^{n}\approx\Real^{2n}$ to
$G_{n,n}\backslash G_n$.  We use $\vectx$ for the variable on $\Real^{2n}$,
sometimes identified with the variable in $G_{n,n}\backslash G_n$.
In an integral, we write Lebesgue measure as $\intd \vectx$.
We let $\mu_{G_{n,n}\backslash G_n}$ be the corresponding
measure on $G_{n,n}\backslash G_n$, under the isomorphism
\eqref{eqn:gmodgn1}.

Recall throughout the following discussion that a homogeneous
space of a closed unimodular linear group has an invariant
measure, unique up to constant factor.  Consider the lattice
of subgroups
\starteqn\label{eqn:groupdiamond}
\begin{diagram}
&&G_n&&\\
&\ldTo&&\rdTo& \\
G_{n,n}&&&&\Gamma_n\\
&\rdTo&&\ldTo&\\
&&\Gamma_{n,n}&&
\end{diagram}
\finisheqn
Fix a Haar measure $\intd g_n$ on $G$.  On the discrete
groups $\Gamma_n$, $\Gamma_{n,n}$, let Haar measure be the
counting measure.  Then $\intd g_n$ determines unique
measures on $\Gamma_{n}\backslash G_n$ and $\Gamma_{n,n}\backslash G_n$,
since the measure is determined locally.  We write $\intd\overline{g}_n$
for each of these induced measures from the right
hand side of the diagram in \eqref{eqn:diamondintegral}.
Passing to the left-hand
side of \eqref{eqn:diamondintegral}, having fixed $\intd g_n$ on $G_n$ and
\startdisp
\intd\mu_{G_{n,n}\backslash G_n}=\intd \vectx\;\text{on}\; G_{n,n}\backslash
G_n,
\finishdisp
\bmth\textbf{there is a unique measure $\intd G_{n,n}$ on $G_{n,n}$ such
that}
\starteqn\label{eqn:diamondintegral}
\int_{G_{n,n}\backslash G_n}\int_{\Gamma_{n,n}\backslash G_{n,n}}=
\int_{\Gamma_{n,n}\backslash G_n}=
\int_{\Gamma_{n}\backslash G_n}\int_{\Gamma_{n,n}\backslash\Gamma_n}.
\finisheqn\ubmth
The integration formula \eqref{eqn:diamondintegral} is
meant to be interpreted as two Fubini-type theorems
for the evaluation of an integral on the space $L^1(\Gamma_{n,n}\backslash
G_n)$.

\begin{lem}\label{lem:diamondintegral}  Let $f\in
L^1(\Real^{2n})\approx L^1(G_{n,n}\backslash G_n)$.  Let
$c_n=\vol(\Gamma_{n,n}\backslash G_{n,n})$, as measured
by $\intd \overline{g}_{n,n}$.  Then
\startdisp
c_n\int_{\Real^n}f(\vectx)\intd \vectx=\int_{\Gamma_n\backslash G_n}\int_{
\Gamma_{n,n}\backslash
\Gamma_n}f(\gamma g)\intd\overline{\gamma_n}\intd\overline{g}_n.
\finishdisp
\end{lem}
\begin{myproof}{Proof}  Since $L^1(G_{n,n}\backslash G_n)$
is a subspace $L^1(\Gamma_{n,n}\backslash G_n)$, we
can apply the Fubini-type integral formula \eqref{eqn:diamondintegral}
to $f$.  Because of the invariance property of $f$,
the inner integral on the left-side of \eqref{eqn:diamondintegral}
reduces to $c_n$.
\end{myproof}

Denote by \bmth$\textbf{prim}(\Int[\vecti]^n)^*$\ubmth\hspace*{0.25mm}
the set of primitive
$n$-vectors, \textit{i.e.}, integral vectors such that the GCD
of the components is $1$.  See Definition \ref{defn:GCD}
and the following properties following it for a review of the GCD in
the context of pairs Gaussian integers.  The definition and properties
have obvious extensions to $n$-tuples of Gaussian integers,
which we will use freely in what follows.
We will also make free
use of the ``Conventions regarding multiplicative structure of
$\Int[\vecti]$" in the paragraph preceding that definition.

Since $\text{prim}(\Int[\vecti]^n)^*$ is precisely the set
of vectors in $(\Int[\vecti]^n)^*$ which can be extended
to a matrix in $\Gamma_n=\SL{n}{\Int[\vecti]}$, we have
\starteqn\label{eqn:firstrowofgamman}
e_n^*\Gamma_n=e_n^*\SL{n}{\Int[\vecti]}=\text{prim}(\Int[\vecti]^n)^*.
\finisheqn
Next, we claim that we have
\starteqn\label{eqn:ntuplesprimitivedecomp}
(\Int[\vecti]^n)^*-\{0\}=\{\ell\vectv\;\text{with}\;\vectv\;\text{primitive
and}\;\ell\in\Int[\vecti]-\{0\}\;\text{standard}\},
\finisheqn
where the decomposition of $\vectw\in(\Int[\vecti]^n)^*$
into $\ell$ and $\vectv$ is understood
be unique.  In order prove the claim, for an arbitrary
element $\vectw\in(\Int[\vecti]^n)^*-\{0\}$, set $\ell$
equal to the GCD of the entries of $\vectw$.  Then
set
\startdisp
\vectv=\ell\inv\vectw.
\finishdisp
Use Properties \textbf{GCD 1} and \textbf{GCD 2} (more precisely,
the form of the properties extended to $n$-tuples) to see that $\vectv$
is primitive.
If $\ell'\vectv'$ is another decomposition of $\vectw$
in the form given in \eqref{eqn:ntuplesprimitivedecomp},
then
\startdisp
{\ell'}\inv\ell\vectv=\vectv',
\finishdisp
with both $\vectv,\,\vectv'$ primitive.  So ${\ell'}\inv\ell$
is a unit.  Since $\ell, \ell'$ are standard, the unit
in question must be $1$.  Thus $\ell=\ell'$, and we deduce the uniqueness
of the decomposition.

\begin{rem}  The decomposition \eqref{eqn:ntuplesprimitivedecomp}
corresponds
to the displayed equation immediately following (5)
in \S$\romtwo$.4 in \cite{posnr}, with the condition
that ``$k$" (the $\ell$ of our notation) is positive corresponding
to the condition that our $\ell$ is standard.  In fact,
the condition of standard-ness is an appropriate generalization of
positivity from the context
of the multiplicative theory of $\Int$ to the multiplicative theory of
$\Int[\vecti]$.
\end{rem}

Let
\startdisp
V_n=\vol(\Gamma_n\backslash G_n)=\vol(\SL{n}{\Int[\vecti]}\backslash
\SL{n}{\Complex}).
\finishdisp
If we change the Haar measure on $G$ by a constant
factor, then the volume changes by this same constant.
The volume is with respect to our fixed $\intd g_n$.  In
\eqref{eqn:Haarmeasurenormalized},
below, we shall fix a normalization of $\intd g_n$.

In the proof of Theorem \ref{thm:siegelvolume} below,
we will use the following change of variables formula.
\starteqn\label{eqn:changeofvariables}
\int_{\Complex^n}f(y\vectx)\intd \vectx=||y||^{-2n}\int_{\Complex^n}f(\vectx)\intd \vectx,\;
\text{for}\; f\in\L^1(\Complex^n),\; y\in\Complex^{\times}.
\finisheqn
The formula \eqref{eqn:changeofvariables} is derived
by computing the Jacobian factor of the multiplication endomorphism
of $\Complex^n$ given by $\vectx\mapsto y\vectx$.  Clearly, the
Jacobian of the multiplication endomorphism is independent
of $\arg(y)$, so is a function of $||y||$ alone.
The Jacobian factor is clearly mutliplicative in $y$,
hence a power of $||y||$.  Then exponent of $||y||$ in the Jacobian
is determined by noting that the  dimension of $\Complex^n$
as a real vector space is $2n$.

\begin{thm}  \textbf{After Siegel \cite{siegel45}.}\label{thm:siegelvolume}
Let $G_n=\SL{n}{\Complex}$, $\Gamma_n=\SL{n}{\Int[\vecti]}$.
  Let $\intd\vectx$ be the Lebesgue measure
on $\Complex^n\approx \Real^{2n}$.  Let $\vectw$ denote a length-$n$ vector
with entries in $\Int[\vecti]$.  Let
$f\in L^1(\Complex^n)\cong L^1(\Real^{2n})$.
Then
\startdisp
\begin{aligned}
V_n\int_{\Real^{2n}}f(\vectx)\intd \vectx&=&&\int_{\Gamma_n\backslash G_n}
\sum_{\vectw
\neq 0}f(\vectw g)\intd\overline{g}\\
&=&&\zeta_{\Rational(\vecti)}(n)\int_{\Gamma\backslash G}\sum_{\vectv\;
\text{prim}}f(\vectv g)\intd\overline{g}.
\end{aligned}
\finishdisp
Furthermore, $V_n=c_n\zeta_{\Rational(\vecti)}(n)$.  Corollary
\ref{cor:gammacovolume}, below, will determine $V_n, c_n$.
\end{thm}
\begin{myproof}{Proof}
On the right side of \eqref{eqn:diamondintegral}, we use
Lemma \ref{lem:diamondintegral} to obtain
\startdisp
\int_{\Gamma_n\backslash G_n}\sum_{\vectv\;\text{prim}}f(\vectv g)
d \overline{g}=c_n\int_{\Complex^n}f(\vectx)\intd \vectx.
\finishdisp
Replacing $f(x)$ by $f(\ell \vectx)$ with $\ell$ a \textit{standard} integer
$\ell$, and using \eqref{eqn:changeofvariables}
on the right, we find
\starteqn\label{eqn:siegelsthminter1}
\int_{\Gamma_n\backslash G_n}\sum_{\vectv\;\text{prim}}f(\ell\vectv g)
\intd\overline{g}=c_n ||\ell||^{-2n}\int_{\Complex^n}f(\vectx)\intd \vectx.
\finisheqn
Summing over all \textit{standard} $\ell\in\Int[\vecti]$, the elements
$\ell\vectv$ range over all nonzero elements $\vectw\in(\Int[\vecti]^n)^*$,
taking each value once, by \eqref{eqn:ntuplesprimitivedecomp}.
Thus, the summed form of the left-hand side of \eqref{eqn:siegelsthminter1}
gives the left-hand side of \eqref{eqn:siegelsthminter2}, below.
In order to see that the summed form of the right-hand
side of \eqref{eqn:siegelsthminter1} gives the right-hand
side of \eqref{eqn:siegelsthminter2}, we use the expression
for $\zeta_{\Rational(\vecti)}$ given in the
second line of \eqref{eqn:dedekindzetadefn}.  Thus, we obtain
\starteqn\label{eqn:siegelsthminter2}
\int_{\Gamma_n\backslash G_n}\sum_{\vectw\neq 0}f(\vectw g)\intd
\overline{g}=c_n\zeta_{\Rational(\vecti)}(n)\int_{\Complex^n}f(\vectx)\intd \vectx.
\finisheqn
Assuming that $V_n$ is finite, we shall now prove that
$V_n=c_n\zeta_{\Rational(\vecti)}(n)$. For this,
we can take a function $f$ which is continuous, $\geq 0$,
with positive integral and compact support.  We
note that for any $g\in\SL{n}{\Complex}$,
\starteqn\label{eqn:riemannsum}
\lim_{N\rightarrow\infty}\frac{1}{N^{2n}}\sum_{\vectw\neq 0}
f\left(\frac{1}{N}\vectw g\right)=\int_{\Complex^n}f(\vectx)\intd \vectx,
\finisheqn
The reason is that the left-hand side of \eqref{eqn:riemannsum},
without the limit, is the Riemann sum for the integral
on the right-hand side associated to a subdivision of the support
of $f$ into parallelotopes of side-length $\frac{1}{N}$.
The right-multiplication
by $g\in\SL{n}{\Complex}$ changes the side lengths of a parallelotope,
but preserves the volume.  Integrating \eqref{eqn:riemannsum} over
$\Gamma_n\backslash G_n$, we find
\startdisp
\begin{aligned}
V_n\int_{\Complex^n}f(\vectx)\intd \vectx&=&&\int_{\Gamma_n\backslash
G_n}\lim_{N\rightarrow\infty}\frac{1}{N^{2n}}
\sum_{\vectw\neq 0}f\left(\frac{1}{N}\vectw g\right)\intd\overline{g}\\
&=&&\lim_{N\rightarrow\infty}\frac{1}{N^{2n}}
\int_{\Gamma_n\backslash
G_n}
\sum_{\vectw\neq 0}f\left(\frac{1}{N}\vectw g\right)\intd\overline{g}\\
&=&&\lim_{N\rightarrow\infty}c_n\zeta_{\Rational(\vecti)}(n)
\frac{1}{N^{2n}}
\int_{\Complex^n}f\left(\frac{1}{N}\vectx\right)\intd \vectx&&
(\text{by}\;\eqref{eqn:siegelsthminter2})
\\
&=&&\lim_{N\rightarrow\infty}c_n\zeta_{\Rational(\vecti)}(n)
\int_{\Complex^n}f(\vectx)\intd \vectx &&(\text{by letting}\;
u=\vectx/N, \;\intd u=\intd \vectx/N^{2n}).
\end{aligned}
\finishdisp
This concludes the proof of Siegel's theorem.
\end{myproof}
\hspace*{0.3cm}

\noindent\bmth
\textbf{Decompositions of invariant measure on $\mathrm{Pos}_n(\Complex)$}
\ubmth\hspace*{0.5mm}
\begin{defn} The space \bmth$\mathrm{Pos}_n(\Complex)$\ubmth\hspace*{.25mm}
is defined to be the  Hermitian matrices of size $n$,
having positive eigenvalues.
The group $\GL{n}{\Complex}$ acts on $\mathrm{Pos}_n(\Complex)$ on the left,
according to the formula
\starteqn\label{eqn:posnraction}
[g]p\mapsto gpg^*,
\finisheqn
and $\mathrm{Pos}_n(\Complex)$ is a $\GL{n}{\Complex}$-homogeneous space.
Accordingly, it has a $\GL{n}{\Complex}$-invariant measure, which
is unique up to a constant factor.
\end{defn}

For the purposes of comparison, we introduce the following
notation for the Lebesgue measure on a Euclidean space
with a system of coordinates $Y=\left((y_{ij}^{\left(
\vecti^{\delta}\right)}\right)$
where the coordinate $y_{ij}^{(1)}$ corresponds to the real component
and $y_{ij}^{({\vecti})}$ to the imaginary component of $y_{ij}$.  We write
\startdisp
\intd\mu_{\mathrm{euc}}(Y)=
\prod\intd y_{ij}^{\left(\vecti^{\delta}\right)},\;
\text{where}\; 1\leq i\leq j\leq n,\;\text{and}\;
\delta=\begin{cases}0,1& \text{for}\;i<j\\0&\text{for}\; i=j\end{cases}\,.
\finishdisp

Here are some comments on the notation to be used below.
We shall reserve the letter $Y$ for a variable in $\mathrm{Pos}_n(\Complex)$,
and $Z$ for a variable in the space $\mathrm{Spos}_n(\Complex)$,
introduced in Definition $\ref{defn:spos}$ below.
Deviations from Lebesgue measure will be denoted
by $\intd\mu(Y)$, with $\mu$ to be specified.  If $\varphi$
is a local $C^{\infty}$ isomorphism, $J(\varphi)$ will denote
the Jacobian factor of the induced map on the measure,
so the absolute value of the determinant of the Jacobian
matrix, when expressed in terms of local coordinates.
If $g$ is a square matrix, we let $|g|$ denote its determinant
and $||g||$ the absolute value of the determinant.

\begin{prop}  A $\GL{n}{\Complex}$-bi-invariant measure
on $\mathrm{Pos}_n(\Complex)$ is given by
\starteqn\label{eqn:posncmeasuredefn}
\intd\mu_n(Y)=|Y|^{-n}\intd\mu_{\mathrm{euc}}(Y).
\finisheqn
For $g\in\GL{n}{\Complex}$, the Jacobian determinant $J(g)$
of the determinant of the transformation $[g]$ is
\startdisp
J(g)=||g||^{2n}.
\finishdisp
The invariant measure satisfies $\intd\mu_n(Y^{\inv})=\intd\mu_n(Y)$,
\textit{i.e.} it is also invariant under $Y\mapsto Y\inv$.
\end{prop}
\begin{myproof}{Proof}  We prove the second assertion first.
Note that $g\mapsto J(g)$ is multiplicative and continuous,
so that it suffices to prove the formula for a dense set of matrices $g$.
We pick the set of semisimple matrices in $\GL{n}{\Complex}$,
\textit{i.e.}, those of the form
\startdisp
gDg\inv,\;\text{with}\; D=\diag(d_1,\ldots,d_n),\;g\in\GL{n}{\Complex}.
\finishdisp
Then we readily calculate that
\startdisp
[D]Y=(d_iy_{ij}\overline{d_j}),
\finishdisp
so that, by the multiplicativity of $J$,
\startdisp
J(gDg\inv)=J(D)=\prod_{1\leq i < j\leq n}|d_id_j|^2
\prod_{1\leq i\leq n}|d_i|^2=
||D||^{2(n-1)}||D||^2=||D||^{2n}=||gDg\inv||^{2n},
\finishdisp
which, by the above comments, proves the formula for $J(g)$.  Then,
with $\intd\mu_n$ defined as in \eqref{eqn:posncmeasuredefn},
\startdisp
\begin{aligned}
\intd\mu_n([g]Y)&=&&|[g]Y|^{-n}J([g])
\intd\mu_{\rm euc}(Y)\\
&=&&||g||^{-2n}|Y|^{-n}||g||^{2n}\intd\mu_{\rm euc}(Y)\\
&=&&\intd\mu_n(Y),
\end{aligned}
\finishdisp
thus concluding the proof of left invariance.  Right
invariance follows because $J(g)=J(g^*)$.

Finally, the invariance under $Y\mapsto Y\inv$ follows
because if we let $S(Y)=Y\inv$, then for a tangent
vector $H\in\mathrm{Herm}_n$,
\startdisp
S'(Y)H=-Y\inv H Y\inv, 
\finishdisp
so $\det S'(Y)=J(Y\inv)=|Y|^{-2n}$.  Then
\startdisp
\intd\mu_n(Y\inv)=|Y|^n|Y|^{-2n}\intd\mu_{\rm euc}(Y)=
|Y|^{-n}\intd\mu_{\rm euc}(Y)=\intd\mu_n(Y),
\finishdisp
thus concluding the proof of the proposition.
\end{myproof}

We have the \textbf{first order partial Iwasawa}
decomposition of
$Y\in\mathrm{Pos}_n(\Complex)$:
\starteqn\label{eqn:firstpartialiwasawaposnc}
Y=\begin{bmatrix}I_{n-1}&\mathbf{x}\\0&1\end{bmatrix}\begin{pmatrix}
W&0\\0&v\end{pmatrix},
\finisheqn
with $v\in\Real^+,\,\mathbf{x}^*\in\Complex^{n-1},\,
W\in\mathrm{Pos}_{n-1}(\Complex)$

The decomposition \eqref{eqn:firstpartialiwasawaposnc}
gives first partial coordinates $Y=Y(W,\mathbf{x},v)$,
with the map
\startdisp
\varphi^+_{n-1,1}:\mathrm{Spos_{n-1}(\Complex)}\times{\Complex^{n-1}}^*\times
\Real^+
\rightarrow\mathrm{Spos}_n(\Complex),
\finishdisp
as in \eqref{eqn:firstpartialiwasawaposnc} as above.  Direct
multiplication in \eqref{eqn:firstpartialiwasawaposnc} yields
the explicit explicit expression
\starteqn\label{eqn:firstpartialiwasawamultposnc}
\varphi^+_{n-1,1}(W,\mathbf{x},v)=\begin{pmatrix}
v^{-\frac{1}{n-1}}W+v\mathbf{x}^*\mathbf{x}&v\mathbf{x}\\
v\mathbf{x}^*&v\end{pmatrix}.
\finisheqn
From \eqref{eqn:firstpartialiwasawamultposnc},
we see that $\varphi^+_{n-1,1}$
is bijective, because, first $v\in\Real^+$ uniquely determines
the lower right entry.  Then $\mathbf{x}\in{\Complex^{n-1}}^*$ is
uniquely determined
to give the last row (or last column), and finally $W$ is uniquely
determined by the upper-left $(n-1)\times (n-1)$ square.

We wish to compare the partial Iwasawa coordinates in
\eqref{eqn:firstpartialiwasawamultposnc} with the coordinates
on $\mathrm{Pos}_n(\Complex)$ induced by the block decomposition
\startdisp
Y=\begin{pmatrix}Y_1&\mathbf{y_2}\\ \mathbf{y_2}^*&y_3\end{pmatrix}
\finishdisp
where $Y_1$ is an $(n-1)\times(n-1)$ matrix, $\mathbf{y_2}$ is
an $(n-1)$-vector, and $y_3>0$.

\begin{prop}\label{prop:partialawasawajacobian}  The Jacobian is given by
\startdisp
J(\varphi^{+}_{n-1,1})=|V|^{2(n-1)}.
\finishdisp
For $Y=\varphi^+(W,\mathbf{x},v)$ we have the change of variable
formula
\startdisp
\intd\mu_n(Y)=|W|\inv v^{n-1}\intd\mu_{\rm euc}(X)\intd\mu_{n-1}(W)
\intd\mu_{1}(v)
\finishdisp
\end{prop}
\begin{myproof}{Proof}  We compute the Jacobian matrix and find
\startdisp
\frac{\partial(Y)}{\partial(W,\mathbf{x},v)}=
\begin{pmatrix}
I_{(n-1)^2}&*\cdots*&*\\
0&v\cdots 0&*\\
\vdots&\vdots\ddots\vdots&\vdots\\
0&0\cdots v&*\\
0&0\cdots 0&1
\end{pmatrix},
\finishdisp
with $v$ occuring $2(n-1)$ times on the diagonal.  Taking
the determinant yields the stated value.  For the change of
variable formula, we just plug in using the definitions
\startdisp
\intd\mu_{n-1}(W)=|W|^{-(n-1)}\intd\mu_{\rm euc},
\finishdisp
and similarly with $n$ and $1$, combined with the value
for the Jacobian.  The formula comes out as stated.
\end{myproof}
\hspace*{0.3cm}

\noindent\bmth
\textbf{Decompositions of invariant measure on $\mathrm{SPos}_n(\Complex)$}
\ubmth
\begin{defn} \label{defn:spos} We define
\bmth$\mathbf{Spos}_n(\Complex)$\ubmth
\hspace*{0.25mm} to be the subspace of $\mathrm{Pos}_n(\Complex)$
consisting
of unimodular matrices.  The unimodular subgroup
$\SL{n}{\Complex}$ acts on $\mathrm{Spos}_n(\Complex)$ by the
same formula \eqref{eqn:posnraction}, and $\mathrm{Spos}_n(\Complex)$
is an $\SL{n}{\Complex}$-homogeneous space.  It therefore has a
$\SL{n}{\Complex}$-invariant measure, unique up to constant
factor.  Every $Y\in\mathrm{Pos}_n(\Complex)$ can be written
uniquely in the form
\startdisp
Y=r^{1/n}Z\;\text{with}\; r>0\;\text{and}\;Z\in\mathrm{Spos}_n(\Complex).
\finishdisp
Thus, we have a product decomposition
\starteqn\label{eqn:sposproductdecomp}
\mathrm{Pos}_n(\Complex)=\Real^+\times\mathrm{Spos}_{n},
\finisheqn
in terms of the coordinates $(r,Z)$.  We denote by $\mu_n^{(1)}$
the $\SL{n}{\Complex}$-invariant measure on $\mathrm{Spos}_n(\Complex)$
such that for the product decomposition \eqref{eqn:sposproductdecomp}
we have
\starteqn\label{eqn:posnsposnmeasuredecomp}
\intd\mu_n(Y)=\frac{\intd r}{r}\intd\mu_n^{(1)}(Z).
\finisheqn
\end{defn}
We have the first order partial Iwasawa
decomposition of
$Z\in\mathrm{Spos}_n(\Complex)$:
\starteqn\label{eqn:firstpartialiwasawa}
Z=\begin{bmatrix}I_{n-1}&\mathbf{x}\\0&1\end{bmatrix}\begin{pmatrix}
v^{-1/(n-1)}W&0\\0&v\end{pmatrix},\;
\text{with $v\in\Real^+,\,\mathbf{x}^*\in\Complex^{n-1},\,
W\in\mathrm{SPos}_{n-1}(\Complex)$.}
\finisheqn
Note that, by multiplying out the right-hand side of
\eqref{eqn:firstpartialiwasawa} and comparing lower-right entries,
we obtain
\starteqn\label{eqn:lowerrightentries}
z_{nn}=v
\finisheqn
The decomposition \eqref{eqn:firstpartialiwasawa}
gives first partial coordinates $Y=Y(W,\mathbf{x},v)$,
with the map
\startdisp
\varphi^+_{n-1,1}:\Real^+\times{\Complex^{n-1}}^*\times
\mathrm{Spos_{n-1}(\Complex)}
\rightarrow\mathrm{Spos}_n(\Complex),
\finishdisp
as in \eqref{eqn:firstpartialiwasawa} as above.  Direct
multiplication in \eqref{eqn:firstpartialiwasawa} yields
the explicit expression
\starteqn\label{eqn:firstpartialiwasawamult}
\varphi^+_{n-1,1}(W,\mathbf{x},v)=\begin{pmatrix}
v^{-1/n-1}W+v\mathbf{x}^*\mathbf{x}&v\mathbf{x}\\
v\mathbf{x}^*&v\end{pmatrix}.
\finisheqn
From \eqref{eqn:firstpartialiwasawamult}, we see that $\varphi^+_{n-1,1}$
is bijective, because, first $v\in\Real^+$ uniquely determines
the lower right entry.  Then $\mathbf{x}\in{\Complex^{n-1}}^*$ is
uniquely determined
to give the last row (or last column), and finally $W$ is uniquely
determined by the upper-left $(n-1)\times (n-1)$ square.

\begin{prop} \label{prop:sposncmeasuredecomp}Associated to the first
partial Iwasawa
decomposition of \eqref{eqn:firstpartialiwasawa},
we have the measure decomposition
\startdisp
\intd\mu_n^{(1)}(Z)=v^{n}\frac{\intd v}{v}\intd
\mathbf{x}\intd\mu_{n-1}^{(1)}(W).
\finishdisp
\end{prop}
\begin{myproof}{Proof}  Write the first-partial Iwasawa decomposition
of $Y\in\mathrm{Pos}_n(\Complex)$ in the form
\startdisp
Y=\begin{bmatrix}I_{n-1}&\vectx\\&1\end{bmatrix}
\begin{pmatrix}{v'}^{-1/n-1}W&\\&v\end{pmatrix}
\finishdisp
with $W\in\mathrm{Spos}_n(\Complex)$.  Let $r=|Y|$,
so that
\startdisp
Y=r^{1/n}Z,\;\text{with}\;Z\in\mathrm{Spos}_n(\Complex).
\finishdisp
On the one hand, by \eqref{eqn:posnsposnmeasuredecomp}, we have
\starteqn\label{eqn:intdmunYfirstexp}
\intd\mu_n(Y)=\frac{\intd r}{r}\intd\mu_n^{(1)}(Z).
\finisheqn
On the other hand, we have, by \eqref{eqn:posncmeasuredefn},
\starteqn\label{eqn:sposncpropinter1}
\begin{aligned}
\intd\mu_n(Y)&=&&|Y|^{-n}\intd\mu_{\rm euc}(Y)\\
&=&&|Y|^{-n}v^{2(n-1)}\intd\mu_{\rm euc}(\vectx)\intd\mu_{\rm euc}(v)
\intd\mu_{\rm euc}({v'}^{-1/n-1}W)&&(\text{by Proposition}\;
\ref{prop:partialawasawajacobian})\\
&=&&|Y|^{-n}v^{2n-1}\frac{\intd v}{v}\intd\vectx\intd\mu_{\rm euc}
({v'}^{-1/n-1}W)\\
&=&&|Y|^{-n}v^{2n-1}|v^{-1/n-1}W|^{n-1}\intd\mu_{n-1}({v'}^{-1/n-1}W)
\frac{\intd v}{v}\intd\vectx&&(\text{by \eqref{eqn:posncmeasuredefn}}).\\
\end{aligned}
\finisheqn
Now set $r'=\frac{1}{v'}$, so that
\startdisp
{v'}^{1/n-1}W={r'}^{1/n-1}W\;\text{with}\; W\in\mathrm{Spos}_n(\Complex),
\finishdisp
and
\startdisp
\intd\mu_{n-1}({v'}^{-1/n-1}W)=\frac{\intd r'}{r'}\intd\mu_{n-1}^{(1)}(W).
\finishdisp
Thus, from \eqref{eqn:sposncpropinter1}, we have
\starteqn
\begin{aligned}
\intd\mu_n(Y)&=&&v'v^{n-1}
\left(\frac{\intd r'}{r'}\intd\mu_{n-1}^{(1)}(W)\right)
\frac{\intd v}{v}\intd\vectx\\
&=&&\left(\frac{\intd r'}{r'}\right)v'v^{n-1}\intd\mu_{n-1}^{(1)}(W)
\frac{\intd v}{v}\intd\vectx\\
&=&&\left(\frac{\intd r}{r}\right)v^{n}\intd\mu_{n-1}^{(1)}(W)
\frac{\intd v}{v}\intd\vectx,
\end{aligned}
\finisheqn
where the last line follows because $r'=r$ and $v'=v$ for
$Y\in\mathrm{Spos}_n(\Complex)$.
\end{myproof}

The measure decomposition of Proposition \ref{prop:sposncmeasuredecomp}
gives rise to a corresponding measure decomposition on the fibration.
In terms of integration over the fibers, we get the following
integral formula.
\begin{prop}\label{prop:gamman1fubini}
For a function $f$ on $\Gamma_{n,n}\backslash\mathrm{Spos}_n$,
in terms of the coordinates of Proposition \ref{prop:sposncmeasuredecomp},
we have
\startdisp
\int_{\Gamma_{n,n}\backslash\mathrm{Spos}_n(\Complex)}f(Z)\intd\mu_n^{(1)}(Z)=
\int_{\Gamma_{n-1}\backslash\mathrm{Spos}_{n-1}(\Complex)}
\int_{\Int[\vecti]^n\backslash\Complex^n}\int_{\Real^+}
f(W,\vectx,v)v^{n}\frac{\intd v}{v}\intd \vectx\intd\mu_{n-1}^{(1)}(W).
\finishdisp
\end{prop}
\hspace*{0.3cm}

\noindent\bmth
\textbf{Completion of calculation of covolume of $\SL{n}{\Int[\vecti]}$.}
\ubmth
We now fix a normalization of the Haar measure $\intd g_n$ on
$G_n=\SL{n}{\Complex}$ such that
\starteqn\label{eqn:Haarmeasurenormalized}
\int_{\mathrm{Spos}_n(\Complex)}f(Z)\intd\mu_n^{(1)}(Z)=
\int_{G_n/K_n}f(gg^*)\intd\overline{g},\;\text{for all}\;
f\in\mathrm{L}^1(\mathrm{Spos}_n(\Complex)).
\finisheqn
Because of the decomposition $G_n=U_nA_nK_n$, the
normalization fixes a normalization of the Haar measure
on $K_n$ giving $K_n$ total measure $1$.  The Haar
measure satisfying \eqref{eqn:Haarmeasurenormalized}
will be called the \textbf{symmetrically normalized
measure}.  The condition \eqref{eqn:Haarmeasurenormalized}
says that this measure is the pull-back of the measure
induced on $[\Gamma_n]\backslash\mathrm{Spos}_n(\Complex)$
by $\mu_n^{(1)}$ under the isomorphism
\startdisp
\Gamma_n\backslash G_n/K_n\rightarrow \Gamma\backslash\mathrm{Spos}_n(\Complex)=
[\Gamma_n]\backslash\mathrm{Spos}_n(\Complex).
\finishdisp
\begin{cor} \label{cor:siegelsthm} Suppose $\intd\overline{g}$ is the
symmetrically
normalized measure.  For $\varphi$ continuous (say) and in
$L^1(\Real^+)$, we have
\startdisp
\begin{aligned}
V_n\int_{\Complex^n}\varphi(\vectx^*\vectx)\intd\vectx&=&&
\int_{\Gamma_n\backslash\mathrm{Spos}_n(\Complex)}
\sum_{\vectw\neq 0}\varphi([\vectw]Z)\intd\mu_n^{(1)}(Z),\\
&=&&\zeta_{\Rational(\vecti)}(n)\int_{\Gamma\backslash\mathrm{Spos}_n(\Complex)}
\sum_{\vectv\;\text{prim}}\varphi([\vectw]Z)\intd\mu_n^{(1)}(Z),\\
\end{aligned}
\finishdisp
where the sum in the first line is over all nonzero
length-$n$ vectors with entries in $\Int[\vecti]$.
\end{cor}
\begin{myproof}{Proof}  Let $f(\vectx)=\varphi(\vectx^*\vectx)$ and apply Siegel's
formula to $f$.  Then
\startdisp
\begin{aligned}
V_n\int_{\Complex^n}f(\vectx)\intd \vectx&=&&
\int_{\Gamma_n\backslash G_n}\sum_{\vectw\neq 0}
f(\vectw g)\intd\overline{g}_n&&(\text{by Theorem}\;\ref{thm:siegelvolume})\\
&=&&\int_{\Gamma_n\backslash G_n}\sum_{\vectw\neq 0}\varphi([\vectw]g^*g)
\intd\overline{g}_n&&(\text{by definition of $\varphi$})\\
&=&&\int_{\Gamma\backslash\mathrm{Spos}_n(\Complex)}\sum_{\vectw\neq 0}
\varphi([\vectw]Z)\intd\mu_n^{(1)}(Z),
\end{aligned}
\finishdisp
by the normalization of \eqref{eqn:Haarmeasurenormalized},
thus concluding the proof of the first line of the Corollary.
The second line follows in exactly the same way, but
using the second version of Theorem \ref{thm:siegelvolume}
instead of the first in the first equality above.
\end{myproof}
\begin{prop} \label{prop:fubiniappliedtophi} For $\varphi$ on $\Real^+$
guaranteeing
convergence (for example, $\varphi\in C_c(\Real^+)$)
\startdisp
\int_{\Gamma_n\backslash\mathrm{Spos}_n}\sum_{\vectv\;\text{prim}}
\varphi([\vectv]Z)\intd\mu_n^{(1)}(Z)=V_{n-1}\int_{\Real^+}
\varphi(r)r^{n}\frac{\intd r}{r}.
\finishdisp
\end{prop}
\begin{myproof}{Proof}
We first note that the sum inside the integral,
as a function of $Z\in\mathrm{Spos}_n(\Complex)$, is
$\Gamma_n$-invariant because action by $\Gamma_n$
(on the right side of $\vectv$) simply permutes the primitive
integral vectors.  In any case, we may rewrite the left
side in the form:
\startdisp
\begin{aligned}
\text{left side}&=&&\int_{\Gamma_n\backslash\text{Spos}_n(\Complex)}
\sum_{\gamma\in\Gamma_{n,n}\backslash\Gamma_n}\varphi([e_n^*][\gamma]Z)
\intd\mu_n^{(1)}(Z)
&&(\text{by}\;\eqref{eqn:gmodgn1})\\
&=&&\int_{\Gamma_{n,n}\backslash\mathrm{Spos}_n(\Complex)}
\varphi(z_{nn})\intd\mu_n^{(1)}
(Z)&&\text{(putting $f(Z)=\varphi(z_{nn})$)}\\
&=&&\int_{\Gamma_{n-1}\backslash\mathrm{Spos}_{n-1}(\Complex)}
\int_{\Complex^n/\Int[\vecti]^n}
\int_{0}^{\infty}f(W,\vectx, v)v^n\frac{\intd v}{v}\intd
\vectx\intd\mu_{(n-1)}^{(1)}
(W)&&\text{(by Proposition \ref{prop:gamman1fubini})
and \eqref{eqn:lowerrightentries}}\\
&=&&V_{n-1}\int_0^{\infty}\varphi(v)v^n\frac{\intd v}{v},
\end{aligned}
\finishdisp
with the use of Proposition \ref{prop:gamman1fubini} in the penultimate
step.  This concludes the proof.
\end{myproof}

We shall apply the above results as in Siegel to determine
the volume of $\SL{n}{\Int[\vecti]}\backslash\SL{n}{\Complex}$,
but we need to recall some formulas from euclidean space.
We still let $\intd\vectx$ denote ordinary Lebesgue measure
on $\Complex^n\cong\Real^{2n}$.  We let $\mathbb{S}^{2n-1}$
be the unit sphere and $\mathbf{B}_{2n}$ be the unit ball.
Then we recall from calculus that
\starteqn\label{eqn:volumeof2nball}
\mu_{\rm euc}(\mathbf{B}_{2n})=\frac{\pi^n}{\Gamma(1+n)}=\frac{\pi^n}{n\Gamma(n)}.
\finisheqn
We use polar coordinates in $\Real^{2n}$, so there is a unique
decomposition
\starteqn\label{eqn:polarcoords}
\intd \vectx=r^{2n-1}\intd r\intd\mu^{(1)}_{\rm euc}(\theta),
\finisheqn
where $\intd\mu_{euc}^{(1)}$ represents a uniquely determined
measure on $\mathbf{S}^{2n-1}$ equal to $\intd\theta$ when $n=1$.
For arbitrary $n$, $\theta=(\theta_1,\ldots\theta_{2n-1})$
has $2n-1$ coordinates.  Then we find
\startdisp
\mu_{\rm euc}(\mathbf{B}_{2n})=\int_{\mathbf{B}_{2n}}\intd \vectx=
\mu^{(1)}_{\rm euc}(\mathbb{S}^{2n-1})\int_{0}^1r^{2n-1}\intd r,
\finishdisp
and therefore, using \eqref{eqn:volumeof2nball},
\starteqn\label{eqn:areaspherevolball}
\mu_{\rm euc}^{(1)}(\mathfrak{S}^{2n-1})=2n\mu_{\rm euc}(\mathbf{B}_{2n})
=\frac{2\pi^n}{n\Gamma(n)}.
\finisheqn
From \eqref{eqn:polarcoords} and \eqref{eqn:areaspherevolball}
it follows trivially that for a function $\varphi$ on $\Real^+$
one has the formula
\starteqn\label{eqn:gammafunctionintegrationformulaprelim}
\frac{\pi^{n}}{\Gamma(n)}\int_{\Real^+}\varphi(r^2)r^{2n}
2\frac{\intd r}{r}=\int_{\Real^{2n}}\varphi(\vectx^*\vectx)\intd\vectx,
\finisheqn
say for $\varphi$ continuous and in $L^1(\Real^+)$.
Now, make a change-of-coordinates in
\eqref{eqn:gammafunctionintegrationformulaprelim} from $r$ to $r'$
where
\startdisp
r^2=r'\;\text{and}\;\frac{\intd r}{r}=\half\frac{\intd r'}{r'}.
\finishdisp
Then change the name of the variable back from $r'$ to $r$
and obtain
\starteqn\label{eqn:gammafunctionintegrationformula}
\frac{\pi^{n}}{\Gamma(n)}\int_{\Real^+}\varphi(r)r^{n}
\frac{\intd r}{r}=\int_{\Real^{2n}}\varphi(\vectx^*\vectx)\intd\vectx.
\finisheqn

\begin{cor}\label{cor:gammacovolume}
Let $G_n=\SL{n}{\Complex}$ and $\Gamma_n=\SL{n}{\Int[\vecti]}$.
Let
\startdisp
\Bflambda\zeta_{\Rational(\vecti)}(s)=\pi^{-s}\Gamma(s)
\zeta_{\Rational(\vecti)}(s).
\finishdisp
Then with respect to the symmetrically normalized measure
on $G_n$, the volume $V_n$ of $\Gamma_n\backslash G_n$
is given inductively by $V_n=\Bflambda\zeta_{\Rational(\vecti)}(n)V_{n-1}$,
which yields
\startdisp
V_n=\prod_{k=2}^n\Bflambda\zeta_{\Rational(\vecti)}(k).
\finishdisp
\end{cor}
\begin{myproof}{Proof}  We start with Corollary \ref{cor:siegelsthm},
to which we apply Proposition \ref{prop:fubiniappliedtophi},
and follow up by formula \eqref{eqn:gammafunctionintegrationformula}.
The inductive relation drops out and the case $n=1$ is trivial.
\end{myproof}

Reproducing the arguments for Corollary \ref{cor:gammacovolume}
in the case of $\SL{2}{\Real}$, $\SL{n}{\Int}$, and $\mathrm{Spos}_n(\Real)$,
in place of $\SL{2}{\Complex}$, $\SL{n}{\Int[\vecti]}$
and $\mathrm{Spos}_n(\Complex)$, we obtain the statement
labelled Theorem 4.6 in Chapter II of \cite{posnr},
which we now restate for the reader's convenience.
\begin{cor}\label{cor:gammazcovolume}
Let $G_{n,\Real}=\SL{n}{\Real}$ and $\Gamma_{n,\Int}=\SL{n}{\Int}$.
Let
\startdisp
\Bflambda\zeta_{\Rational}(s)=\pi^{-s/2}\Gamma(s/2)
\zeta_{\Rational}(s).
\finishdisp
Then with respect to the symmetrically normalized measure
on $G_{n,\Real}$, the volume $V_{n,\Real}$ of $\Gamma_{n,\Int}
\backslash G_{n,\Real}$
is given inductively by $V_{n,\Real}=\Bflambda\zeta_{\Rational}(n)
V_{n-1,\Real}$,
which yields
\startdisp
V_{n,\Real}=\prod_{k=2}^n\Bflambda\zeta_{\Rational}(k).
\finishdisp
\end{cor}



\begin{singlespace}
\bibliographystyle{amsalpha}
\bibliography{memoirsbib}
\addcontentsline{toc}{chapter}{Bibliography}
\end{singlespace}
\end{document}